\documentclass[oneside,fleqn,12pt]{scrbook}               %fleqn: math. Formeln linksb�ndig setzen

\usepackage[german]{babel}

\usepackage{a4}     %Seitengr��e und Zeichensatz
\usepackage{ifthen} %\ifthenelse{ \equal{}{} }{}{}
\usepackage{version}%comment Umgebung zum auskommentieren
\usepackage{xspace} %xspace Befehl regelt Abstand nach Befehlen
\usepackage{pifont}	%zusaetzliche Zeichen ... Goosen S. 345
\usepackage{bbding} %RightHand 
\usepackage[centertags]{amsmath}    %Mathematische Zeichen und Formeln
\usepackage{amssymb}    %alle m�glichen �blichen math. Symbole
\usepackage{amstext}   %\text{bla} - Text in math. Fomeln
\usepackage{undertilde} %fuer Tildezeichen unter PI und SIGMA (LightfaceHierarchie/Analyt.Hierarchie (deskr. MngLehre))
\usepackage{amsthm}
\usepackage[all]{xy}%

\newdir{ >}{{}*!/-10pt/@{>}}    %vgl xymanual exercise 14

\newtheoremstyle{break}% name
  {9pt}%      Space above, empty = `usual value'
  {9pt}%      Space below
  {\itshape}% Body font
  {}%         Indent amount (empty = no indent, \parindent = para indent)
  {\bfseries}% Thm head font
  {.}%        Punctuation after thm head
  {\newline}% Space after thm head: \newline = linebreak  
  {}%         Thm head spec
\newcommand{\header}[1]{
\textup{\bfseries(#1)}
}
\newcommand{\point}[1]{
\medskip
{\bfseries #1:}
}
\newcommand{\Point}[1]{
\bigskip
{\bfseries #1:}
}

\newcommand{\siehe}{\HandRight\ }

\theoremstyle{plain}
\newtheorem{lemma}{Lemma}[section]%[...] ist die Z�hltiefe - weglassen bewirkt einfache Z�hlweise

\newtheorem{satz}[lemma]{Satz}%[lemma] bewirkt, da� derselbe Z�hler verwendet wird
\newtheorem{satzdefinition}[lemma]{Satz/Definition}

\newtheorem{theorem}[lemma]{Theorem}
\newtheorem{korollar}[lemma]{Korollar}
\newtheorem{bemerkung}[lemma]{Bemerkung}

\newenvironment{beweis}[1][Beweis]{\begin{proof}[\textsc{#1}]}{\end{proof}}

\theoremstyle{plain}
\newtheorem{beispiel}[lemma]{Beispiel}
\theoremstyle{definition}
\newtheorem{definition}[lemma]{Definition}%Definitionen mit Zaehlung
\newenvironment{ZITAT}
{
\makebox[1cm]{}\begin{minipage}[t]{12cm}
\begin{footnotesize}
} 
{
\end{footnotesize}
\end{minipage}
}

\newcommand{\KOM}[1][]{        %Kommentar am Rand (Kopka 1, S. 89, 90)
\ifthenelse{\equal{#1}{}}
{\marginpar{\sf \bfseries \Huge S}}
{\marginpar{\sf \bfseries #1}}        
}

\renewcommand{\emph}[1]{%<strg>+e FETT
\xspace{\bfseries{#1}}\xspace
}

\newcommand{\eemph}[1]{% SCHRÄG
\xspace{\textit{{#1}}}\xspace
}

\newcommand{\TEXT}[1]{            %gewöhnlicher Text in mathem. Formeln
    \begin{ensuremath}
    	\mbox{#1}
    \end{ensuremath}
}
\newcommand{\TT}[1]{\TEXT{#1}}

\newcommand{\gs}{         					%Gedankenstrich
\xspace--\xspace
}

\newenvironment{ITEMS}[1][equivalences]{
\newcommand{\enumeration}{eq}
\renewcommand{\enumeration}{#1}
\ifthenelse{\equal{\enumeration}{alph)}}{
                                            \enumerate}
                                           {}
\ifthenelse{\equal{\enumeration}{(alph)}}{
                                            \enumerate}
                                           {}
\ifthenelse{\equal{\enumeration}{arabic)}}{
                                            \enumerate}
                                           {}
\ifthenelse{\equal{\enumeration}{(arabic)}}{
                                            \enumerate}
                                           {}
\ifthenelse{\equal{\enumeration}{equivalences}}{
                                            \enumerate}
                                           {}
\ifthenelse{\equal{\enumeration}{enumeration}}{
                                            \enumerate}
                                           {}
\ifthenelse{\equal{\enumeration}{examples}}{
                                            \enumerate}
                                           {}
\ifthenelse{\equal{\enumeration}{definition}}{
                                            \enumerate}
                                           {}
}{
\ifthenelse{\equal{\enumeration}{alph)}}{\endenumerate}{}
\ifthenelse{\equal{\enumeration}{(alph)}}{\endenumerate}{}
\ifthenelse{\equal{\enumeration}{arabic)}}{\endenumerate}{}
\ifthenelse{\equal{\enumeration}{(arabic)}}{\endenumerate}{}
\ifthenelse{\equal{\enumeration}{equivalences}}{\endenumerate}{}
\ifthenelse{\equal{\enumeration}{enumeration}}{\endenumerate}{}
\ifthenelse{\equal{\enumeration}{examples}}{\endenumerate}{}
\ifthenelse{\equal{\enumeration}{definition}}{\endenumerate}{}

}

\newcommand{\FA}[1][]{
\begin{ensuremath}
\forall\TEXT{ } #1\TEXT{ }
\end{ensuremath}
}

\newcommand{\EX}[1][]{
\begin{ensuremath}
\exists\TEXT{ } #1\TEXT{ }
\end{ensuremath}
}

\newcommand{\UND}{
\begin{ensuremath}
\wedge
\end{ensuremath}
}

\newcommand{\ODER}{
\begin{ensuremath}
\vee
\end{ensuremath}
}

\newcommand{\NICHT}{
\begin{ensuremath}
\neg
\end{ensuremath}
}

\newcommand{\AQ}[1][]{
\begin{ensuremath}
\underset{\TT{\tiny #1}}{\Leftrightarrow}
\end{ensuremath}
}

\newcommand{\DEF}{
\begin{ensuremath}
:=
\end{ensuremath}
}

\newcommand{\EQ}[1][]{
\begin{ensuremath}
\underset{\TT{\tiny #1}}{=}
\end{ensuremath}
}

\newcommand{\NEQ}[1][]{
\begin{ensuremath}
\underset{\TT{\tiny #1}}{\neq}
\end{ensuremath}
}

\newcommand{\C}[1][]{
\begin{ensuremath}
#1^{\TT{\slshape{\scriptsize{\scriptsize{c}}}}}
\end{ensuremath}
}

\newcommand{\POW}[1][]{
\begin{ensuremath}
{\textbf P}\left( #1 \right)
\end{ensuremath}
}

\newcommand{\INT}[1][]{
\begin{ensuremath}
#1^{\circ}
\end{ensuremath}
}

\newcommand{\ABS}[1][]{
\begin{ensuremath}
\overline{#1}
\end{ensuremath}
}

\renewcommand{\o}{
\begin{ensuremath}
-
\end{ensuremath}
}

\newcommand{\sd}{
\begin{ensuremath}
\triangle
\end{ensuremath}
}

\newcommand{\LM}{
\begin{ensuremath}
\varnothing
\end{ensuremath}
}

\newcommand{\MN}[1]{
\begin{ensuremath}\left\{#1\right\}\end{ensuremath}
}

\newcommand{\MNG}[2]{
\begin{ensuremath}\{ #1\mid #2 \} \end{ensuremath}
}

\newcommand{\FL}[2][]{
\begin{ensuremath}
\left( #1 \right)_{#2}
\end{ensuremath}
}

\newcommand{\pfeil}[2]{%
\begin{ensuremath}%
\xymatrix{#1\ar@{|->}[r] &{#2}}%
\end{ensuremath}%
}%

\newcommand{\PFEIL}[3]{%
\begin{ensuremath}%
\xymatrix{#1\ar[r]^{#2} & #3}%
\end{ensuremath}%
}%

\newcommand{\IMG}{
\ensuremath{
{\textup{img}}
}
}

\newcommand{\Ers}[1][]{
\begin{ensuremath}
{\textup{\textbf{Ers}}}
\end{ensuremath}
}

\newcommand{\AC}[1][]{
\begin{ensuremath}
{\textup{\textbf{AC}}}^{#1}
\end{ensuremath}
}

\newcommand{\DC}[1][]{
\begin{ensuremath}
{\textup{\textbf{DC}}}^{#1}
\end{ensuremath}
}

\newcommand{\VL}[1][]{
\begin{ensuremath}
{\textup{\textbf{V=L}}}
\end{ensuremath}
}

\newcommand{\AD}[1][]{
\begin{ensuremath}
{\textup{\textbf{AD}}}^{#1}
\end{ensuremath}
}

\newcommand{\GCH}[1][]{
\begin{ensuremath}
{\textup{\textbf{GCH}}}
\end{ensuremath}
}

\newcommand{\ZF}[1][]{
\begin{ensuremath}
{\textup{\textbf{ZF}}}
\end{ensuremath}
}

\newcommand{\ZFC}[1][]{
\begin{ensuremath}
{\textup{\textbf{ZFC}}}
\end{ensuremath}
}

\newcommand{\BR}{
\ensuremath{{\cal N}}
}

\newcommand{\NZ}{
\ensuremath{\mathbb{N}}
}

\newcommand{\ZZ}{
\ensuremath{\mathbb{Z}}
}

\newcommand{\QZ}{
\ensuremath{\mathbb{Q}}
}

\newcommand{\IZ}{
\ensuremath{\mathbb{I}}
}

\newcommand{\RZ}{
\ensuremath{\mathbb{R}}
}

\newcommand{\ORD}{
\ensuremath{\textup{\textbf{Ord}}}
}

\newcommand{\SF}{
\ensuremath{\textup{C}}
}

\newcommand{\BIF}{
\ensuremath{\bf{\textup{R}}}
}

\newcommand{\CM}{
\ensuremath{\textup{C}}
}

\newcommand{\diam}
{
\ensuremath{\textup{diam}}
}

\newcommand{\Kugel}
{
\ensuremath{\textup{K}}
}

\newcommand{\TOP}[1]
{
\ensuremath{\underline{#1}}
}
\newcommand{\SYS}[1]
{
\ensuremath{\underline{#1}}
}
\newcommand{\SSYS}[1]
{
\ensuremath{\tilde{\underline{#1}}}
}

\newcommand{\KLEIN}[1]
{
\ensuremath{\underline{#1}}
}

\newcommand{\SALG}[1]
{
\ensuremath{\underline{#1}}
}

\newcommand{\NM}[1]
{
\ensuremath{\underline{#1}}
}

\newcommand{\BED}[1]
{
\ensuremath{\underline{#1}}
}

\newcommand{\BF}[1]
{
\ensuremath{\underline{#1}}
}

\newcommand{\erf}
{
\ensuremath{\Vdash}
}

\newcommand{\BM}[1]
{
\ensuremath{\mathcal{#1}}
}

\newcommand{\SIGMA}[3][]{
\ensuremath{
\utilde{\Sigma}^{#2}_{#3} \ifthenelse{\equal{#1}{}}{}{\left( {#1} \right)}
}
}

\newcommand{\PI}[3][]{
\ensuremath{
\utilde{\Pi}^{#2}_{#3} \ifthenelse{\equal{#1}{}}{}{\left( {#1} \right)}
}
}

\newcommand{\DELTA}[3][]{
\ensuremath{
\utilde{\Delta}^{#2}_{#3} \ifthenelse{\equal{#1}{}}{}{\left( {#1} \right)}
}
}

\newcommand{\SIGMAlf}[3][]{
\ensuremath{
\Sigma^{#2}_{#3} \ifthenelse{\equal{#1}{}}{}{\left( {#1} \right)}
}
}

\newcommand{\PIlf}[3][]{
\ensuremath{
\Pi^{#2}_{#3} \ifthenelse{\equal{#1}{}}{}{\left( {#1} \right)}
}
}

\newcommand{\DELTAlf}[3][]{
\ensuremath{
\Delta^{#2}_{#3} \ifthenelse{\equal{#1}{}}{}{\left( {#1} \right)}
}
}

\newcommand{\G}{
\ensuremath{
\textup{G}
}
}

\newcommand{\F}{
\ensuremath{
\textup{F}
}
}

\newcommand{\PROJ}[1]{
\ensuremath{
\pi \left( {#1} \right)
}
}

\newcommand{\p}{
\ensuremath{
\textup{p}
}
}

\newcommand{\lng}{%[1][]{
\ensuremath{
\textup{lng}%{(#1)}
}
}
\newcommand{\kon}[1][]{
\ensuremath{
#1{}^\smallfrown
}
}

\newcommand{\Kronecker}
{
\ensuremath{\delta}
}

\newcommand{\bs}[1][]{
\ensuremath{
\arrowvert_{#1}
}
}

\newcommand{\GQ}{
\begin{ensuremath}
\Game%gespiegeltes G - das Zeichen im Moschovakis ist ein um 120 Grad gedehtes G
\end{ensuremath}
}

\newcommand{\BMSo}[1]{
\ensuremath{
\textup{G}^{**} \ifthenelse{\equal{#1}{}}{}{\left( {#1} \right)}
}
}

\newcommand{\BMSa}[1]{
\ensuremath{
\textup{G}^{**}_p \ifthenelse{\equal{#1}{}}{}{\left( {#1} \right)}
}
}

\newcommand{\BMSb}[1]{
\ensuremath{
\tilde{\textup{G}} \ifthenelse{\equal{#1}{}}{}{\left( {#1} \right)}
}
}

\newcommand{\BMSc}[1]{
\ensuremath{
\tilde{\textup{G}}_p \ifthenelse{\equal{#1}{}}{}{\left( {#1} \right)}
}
}

\newcommand{\BMSd}[1]{
\ensuremath{
\textup{G}^{#1}
}
}

\newcommand{\BMSe}[1]{
\ensuremath{
\textup{G}_p^{#1}
}
}

\newcommand{\dd}{
\ensuremath{
d
}
}

\newcommand{\card}[1][]{
\ensuremath{
\left|{#1}\right|
}
}

\newcommand{\ID}[1]
{
\ensuremath{\frak{#1}}
}

\renewcommand{\KOM}[1][]{}

\newcommand{\clearemptydoublepage}{\newpage{\mbox{ }\pagestyle{empty}\cleardoublepage}}

\begin{document}
%\renewcommand{\KOM}[1][]{} %Randnotizen ausschalten
%==============================frontmatter==========================
\frontmatter
\thispagestyle{empty}%
\newcommand{\HRule}{\rule{\linewidth}{1mm}}%

\setlength{\parindent}{0mm}%
\setlength{\parskip}{0mm}%
    \vspace*{\stretch{1}}%
    \HRule%
    \begin{flushleft}
        \Huge Falko Weigt\\[5mm]%
        %\Huge   Spieltheoretische Verallgemeinerung mengentheoretischer Kategorisierungskonzepte%      
        %\Huge   Gr"o"sen-Konzepte und ihre Spieltheoretische Verallgemeinerung %
    		\Huge   Kategorien-Konzepte und ihre spieltheoretische Verallgemeinerung %
    		%\Huge   Kategorienkonzepte und ihre spieltheoretische Verallgemeinerung %
    		%\Huge   Kategorien-Konzepte und ihre Spieltheoretische Vereinheitlichung %
    \end{flushleft}
    \HRule
    \vspace{\stretch{3}}%

    \begin{center}%
        \Large\textsc{M\"unster 2007}
    \end{center}%
\vspace{-2cm}

\clearemptydoublepage

%======WIDMUNG======
%\begin{comment}
\newpage{
\begin{center}
%\normalsize{\scshape{F"ur meinen Hund}}
\normalsize{\scshape{Meinen Eltern}}
%\normalsize{\scshape{F"ur Mama}}
%\normalsize{\scshape{F"ur Omi}}
%\normalsize{\scshape{Papst Benedikt XVI}}
\end{center}
\pagestyle{empty}
\clearpage
}

\clearemptydoublepage
%\end{comment}

\newpage{
\chapter*{Vorwort}
%Hier kommen dann die Danksagungen hin\ldots.
%\begin{comment}
Professor Pohlers danke ich f"ur die Aufgabenstellung und seine geduldige Begleitung dieser Arbeit. Christoph Heinatsch, Christoph Duchhardt und Philipp Schlicht danke ich f"ur ihre stets offenen Ohren f"ur mathematische Fragestellungen. Meinen Eltern, Tanten und Onkeln danke ich f"ur vieles andere.
}
\clearemptydoublepage
%\end{comment}

%Inhaltsverzeichnis

\setcounter{tocdepth}{3}    %Tiefe Aufnahme der Überschriften ins Inhaltsverzeichnis

\tableofcontents

%==============================mainmatter==========================
\mainmatter

%==============================KAPITEL=============================
\chapter{Einf"uhrung}\label{Kap_Einf}
\setcounter{lemma}{0}
Zur Erfassung unserer anschaulichen Vorstellung von der \eemph{Gr"o"se} eines Objektes, genauer einer Teilmenge der reellen Zahlen oder eines anderen topologischen Raumes, gibt es in der Mathematik eine Reihe unterschiedlicher Konzepte. In Anlehnung an das von Ren\'e Baire im Jahr 1899 in \cite{Baire:1899} ver"offentlichte Gr"o"senkonzept, in dem er von Mengen \eemph{von 1.} und \eemph{von 2. Kategorie} spricht, wollen wir solche Konzepte allgemein als \eemph{Kategorien-Konzepte} bezeichnen. 

Derartige mathematischen Gr"o"senbegriffe spielen etwa in der Ma"s- und Integrationstheorie, der Mengenlehre und in der Topologie eine wesentliche Rolle. Mit ihrer Hilfe lassen sich in einem mathematisch exakt formulierten Sinne \eemph{unscharfe Begriffsbildungen} vornehmen. Beispielsweise definiert man in der Ma"s- und Integrationstheorie, f"ur einen \eemph{Ma"sraum\footnote{\emph{Ma"sraum:} Ein Tripel $(\Omega,\SYS{A},\mu)$ wobei $\Omega$ eine nicht-leere Menge, $\SYS{A}$ eine $\sigma$-Algebra auf $\Omega$ und $\mu$ ein Ma"s auf $(\Omega,\SYS{A})$ ist [\siehe Konventionen S. \pageref{Kap_Konventionen}].}} $(\Omega,\SYS{A},\mu)$ und eine Eigenschaft $E\subset\Omega$ von Elementen aus $\Omega$, da"s die Eigenschaft $E$ auf einer Menge $A\in\SYS{A}$ \eemph{$\mu$-fast "uberall} besteht, falls $E$ in $A$ gilt bis auf eine \eemph{$\mu$-Nullemenge\footnote{\emph{$\mu$-Nullmenge:} Eine Menge $A\in\SYS{A}$ mit $\mu(A)=0$ f"ur einen Ma"sraum $(\Omega,\SYS{A},\mu)$.}}: $\EX[N\in\SYS{A}](\mu(N)=0 \UND A\cap \C[N]\subset E)$. Dies ist entscheidend f"ur die Formulierung von Eindeutigkeitsaussagen, die in der Ma"s- und Integrationstheorie in der Regel nur bis auf $\mu$-Nullmengen gelten. 

In dieser Arbeit soll der Schwerpunkt auf der Untersuchung zweier Kategorien-Konzepte liegen, die sich beide mittels topologischer Begriffe bilden lassen: Die \eemph{Baire-Kategorie} und die \eemph{$\sigma$-Kategorie}. 

F"ur die spieltheoretische Charakterisierung betrachten wir beide Kategorien auf dem \eemph{Baireraum}. Das ist der Raum $\omega^\omega$ aller unendlichen Folgen nat"urlicher Zahlen mit den Mengen der Form $O_u:=\MNG{f\in\omega^\omega}{u\TT{ ist ein endl. Anfangsst"uck von } f}$ als den offenen Basismengen [\siehe Kapitel \ref{Kap_Baire_Raum}].

\bigskip Unter der \eemph{Baire-Kategorie} wollen wir das Gr"o"senkonzept von R. Baire verstehen: Als \eemph{kleine} Teilmengen eines topologischen Raumes oder speziell des Baire-Raumes $\omega^\omega$ werden die \eemph{mageren} Teilmengen ausgezeichnet. Eine Menge bezeichnet man als \eemph{mager}, wenn sie sich als abz"ahlbare Vereinigung \eemph{nirgends dichter} Teilmengen (d.h. Teilmengen, die in keiner offenen Menge dicht sind oder "aquivalent dazu: Teilmengen, deren Komplemente jeweils eine offene dichte Menge enthalten) schreiben l"a"st. 

\bigskip Die \eemph{$\sigma$-Kategorie} definieren wir direkt f"ur den Baireraum. Sie zeichnet als \eemph{kleine} Teilmengen die sogenannten \eemph{$\sigma$-beschr"ankten} Teilmengen des Baireraumes aus. Eine Teilmenge $A$ des Baireraumes wird als \eemph{$\sigma$-beschr"ankt} bezeichnet, falls es im Baireraum eine Folge $(f_i)_{i<\omega}$ gibt so, da"s 
\[
\FA[f\in A]\EX[i<\omega](f\leq f_i)
\]
wobei f"ur $f,g\in\omega^\omega$ genau dann $f\leq g$ gelte, wenn $\FA[n<\omega](f(n)\leq g(n))$ gilt (dabei sei $f(n)$ die $n$-te Stelle von $f$). Eine solche Folge von Elementen aus $\omega^\omega$ bezeichnet man auch als eine \eemph{$\sigma$-Schranke}.  

\bigskip \eemph{Der rote Faden} in dieser Arbeit ist die Behandlung der in Kapitel \ref{Kap_VereinhProblem} aufgeworfene Frage nach einer Vereinheitlichung der beiden Gr"o"senkonzepte \eemph{Baire-} und \eemph{$\sigma$-Kategorie}. Die Antwort darauf soll in einer spieltheoretisch definierten verallgemeinerten Kategorie bestehen, die die beiden erstgenannten als Spezialf"alle umfa"st. Dazu werden zun"achst in den Kapiteln \ref{Kap_BaireKategorie} und \ref{Kap_SigmaKategorie} die Baire- und die $\sigma$-Kategorie definiert und ihre grundlegenden Eigenschaften untersucht. Wie sich dabei herausstellt, lassen sich Baire- und $\sigma$-Kategorie in ganz "ahnlicher Weise spieltheoretisch charakterisieren: Dazu wird f"ur jedes der beiden Kategorien-Konzepte ein mathematisches Zwei-Personen-Spiel $\BMSo{A}$ bzw. $\BMSb{A}$ f"ur Teilmenge $A$ des Baireraumes definiert, mit dem sich ein Zusammenhang zwischen dem Kleinheitsbegriff des jeweiligen Kategorien-Konzeptes und den Gewinnstrategien in den jeweiligen Spielen  herstellen l"a"st [\siehe Theoreme \ref{SatzMagerCharSpiele} und  \ref{SatzSigmaBeschrCharSpiele}]. In beiden F"allen wird der jeweilige Kleinheitsbegriff dadurch charakterisiert, da"s eine Menge $A$ \eemph{klein} ist, genau dann, wenn Spieler II eine \eemph{Gewinnstrategie} f"ur das entsprechende Spiel mit Gewinnmenge $A$ hat. Zus"atzlich zu den Spielen $\BMSo{A}$ und $\BMSb{A}$ werden jeweils projektive Versionen $\BMSa{B}$ bzw. $\BMSc{B}$ zur Charakterisierung der jeweiligen Kategorien-Konzepte f"ur Projektionen $A=\p(B)$ f"ur eine Menge $B\subset\omega^\omega\times\omega^\omega$ bzw. $B\subset\omega^\omega\times\lambda$ f"ur eine unendliche Ordinalzahl $\lambda$ definiert [\siehe Definitionen \ref{Def_SpielMitZEugen_Baire} sowie \ref{DefBMSMitZeugen}]. Analog zu den nicht-projektiven Charakterisierungen entsprechen sich auch die Charakterisierungen f"ur Projektionen $A=\p(B)$ [\siehe Theoreme \ref{SatzMagerCharSpiele_2} sowie \ref{SatzSigmaBeschrCharSpiele_2}].

\bigskip Die \eemph{Kompaktheit} l"a"st sich f"ur abgeschlossene Teilmengen $A$ des Baireraumes wie folgt charakterisieren [\siehe Satz \ref{Satz_Heine_Borel_2}]:
\begin{equation}
\FA[s\in T_A]\EX[k<\omega]\FA[t\in\omega^{<\omega}_*](s\kon t\in T_A\Rightarrow t(0)\not > k)\label{Frml_Einl_Kompakt}
\end{equation}
wobei $T_A$ die Menge aller endlichen Anfangsst"ucke von Elementen aus $A$ ist, $\omega^{<\omega}$ die Menge aller endlichen Sequenzen nat"urlicher Zahlen und $\omega^{<\omega}_*$ gleich $\omega^{<\omega}$ ohne die leere Sequenz. Die \eemph{$\sigma$-Beschr"anktheit} einer Menge $A\subset\omega^\omega$ ist laut Satz \ref{SatzSigmaBeschr"ankt} $(ii)$ gleichbedeutend mit:
\begin{equation}
\TT{$A\subset \bigcup_{i<\omega}A_i$ f"ur eine Folge kompakter $(A_i\subset\omega^\omega)_{i<\omega}$}\label{Frml_Einl_SBeschr}.
\end{equation}
F"ur ein allgemeineres Kategorien-Konzept auf dem Baireraum nimmt man nun die beiden Charakterisierungen (\ref{Frml_Einl_Kompakt}) und (\ref{Frml_Einl_SBeschr}) f"ur Kompaktheit bzw. $\sigma$-Beschr"anktheit als Vorlagen f"ur die allgemeineren Begriffe von einer \eemph{$\BM{B}$-nirgends dichten} bzw. einer \eemph{$\BM{B}$-mageren} Teilmenge des Baireraumes (oder etwas allgemeiner des Raumes $X^\omega$ mit einer Menge $X$ mit wenigstens zwei Elementen [\siehe Definitionen \ref{DefBNirgendsDicht} und \ref{DefBMager}]. Die \eemph{$\BM{B}$-mageren} Teilmengen, da als Teilmengen abz"ahlbarer Vereinigungen \eemph{$\BM{B}$-nirgends dichter} Mengen definiert, bilden ein $\sigma$-Ideal. Die \eemph{$\BM{B}$-mageren} Teilmengen stehen also f"ur einen neuen \eemph{Kleinheitsbegriff} den wir als \eemph{verallgemeinertes Kategorien-Konzept} bezeichnen oder (in Anlehnung an die Begriffe von der \eemph{Baire-} und der \eemph{$\sigma$-Kategorie}) als die \eemph{verallgemeinerte Kategorie}. Es zeigt sich, da"s sich die Theoreme \ref{SatzMagerCharSpiele}, \ref{SatzMagerCharSpiele_2}, sowie \ref{SatzSigmaBeschrCharSpiele} und \ref{SatzSigmaBeschrCharSpiele_2} zur spieltheoretischen Charakterisierung der Kleinheitsbegriffe \eemph{mager} bzw. \eemph{$\sigma$-beschr"ankt} in entsprechenden Versionen auch f"ur die verallgemeinerte Kategorie zeigen lassen [\siehe Theoreme \ref{Theorem_Char_AllgSpiel_OhneZeugen} und \ref{Theorem_Char_AllgSpiel_MitZeugen}].

\bigskip Neben den \eemph{Kleinheitsbegriffen} charakterisieren die angesprochenen Theoreme auch Begriffe von \eemph{gro"sen} Mengen (d.h. von Komplementen kleiner Mengen) wie den \eemph{komageren Mengen} bzw. von \eemph{relativ gro"sen} Mengen (d.h. von anschaulich \glqq gro"sen\grqq\ Mengen, die aber nicht unbedingt Komplemente kleiner Mengen sind) wie etwa den \eemph{superperfekten Mengen}. 

Im Falle der spieltheoretischen Charakterisierung der Baire-Kategorie [\siehe Theoreme \ref{SatzMagerCharSpiele}, \ref{SatzMagerCharSpiele_2}] werden im Sinne der Baire-Kategorie \eemph{lokal gro"se} Mengen (d.h. \eemph{in einer offenen Menge komagere} Mengen) charakterisiert. 

Im Fall der spieltheoretischen Charakterisierung der $\sigma$-Kategorie [\siehe Theoreme \ref{SatzSigmaBeschrCharSpiele} und \ref{SatzSigmaBeschrCharSpiele_2}] werden im Sinne der $\sigma$-Kategorie \eemph{relativ gro"se} Mengen (d.h. Mengen, die eine nicht-leere \eemph{superperpfekte} Teilmenge enthalten). Als \eemph{superperfekt} bezeichnet man Teilmengen $A$ des Baireraumes, deren Baum $T_A$ (d.h. die Menge der endlichen Anfangsst"ucke von Elementen in $A$) die Eigenschaft hat, da"s es f"ur jedes Element des Baumes eine endliche Erweiterung im Baum gibt, die im n"achsten Schritt unendlich oft verzweigt. Nicht-leere superperfekte Mengen (und deren Obermengen) stehen zwar f"ur anschaulich \glqq gro"se\grqq\ Mengen, sind aber nicht unbedingt \eemph{gro"s} im Sinne der $\sigma$-Kategorie, d.h. sie sind i.a. nicht Komplement einer \eemph{$\sigma$-beschr"ankten} Menge [\siehe Beispiel \ref{Lemma_SuperPerf}]. Die \eemph{Superperfektheit} wird in Kapitel \ref{Kap_BPerf} zum Begriff der \eemph{$\BM{B}$-Perfektheit} verallgemeinert [\siehe Definition \ref{Def_BPerfTM} und Beispiel \ref{Bsp_BPerfEnthSuperperfTM}]. 

Die spieltheoretischen Charakterisierungen von Baire- und $\sigma$-Kategorie sowie der verallgemeinerten Kategorie entsprechen sich im Hinblick auf diese Gr"o"senbegriffe insofern, als Spieler I in dem jeweiligen (nicht projektiven) Spiel eine Gewinnstrategie hat genau dann, wenn auf die Gewinnmenge \eemph{lokal gro"s} bzw. \eemph{relativ gro"s} (jeweils im obigen Sinne) ist. (F"ur die entsprechenden projektiven Spiele gelten die analaogen Aussagen jeweils nur in der Hinrichtung).

\bigskip Neben der spieltheoretischen Charakterisierung von Kleinheits- und Gr"o"senbegriffen der $\sigma$- und der Baire-Kategorie in den Kapiteln \ref{KapSpielCharMager} und \ref{KapSpielCharSigmaBeschr"ankt} sowie deren Verallgemeinerung in Kapitel \ref{Kap_VerallgKategorie} leisten die Kapitel \ref{Kap_Magere_Mengen} bis \ref{Kap_Baire_Raum} eine Einordnung der Baire-Kategorie und des Baireraumes in den allgemeineren topologischen Zusammenhang: 

Um die Eigenschaft, da"s in $\RZ$ kein Intervall $[a,b]\subset\RZ$ (oder "aquivalent dazu: kein offenes Intervall $(a,b)\subset\RZ$) mit $a<b$ als Vereinigung einer Folge von Mengen von erster Kategorie dargestellt werden kann [\siehe Beispiel \ref{Bsp_mager_lauterLuecken}], f"ur topologische R"aume zu verallgemeinern, bezeichnet man einen Raum als \eemph{Baire'sch}, wenn er keine nicht-leere offene magere Menge enth"alt [\siehe Kapitel \ref{Kapitel_BairescheRaeume}].
 
Eine nicht-leere offene Teilmenge $U\subset X$ ist mager genau dann, wenn jede 
Teilmenge von $X$ \eemph{mager (und komager) in $U$} ist [\siehe Definitionen \ref{Def_Lokalisierung_Mager} und \ref{Def_KoMager}]. Die Baire'schen R"aume sind daher gerade diejenigen, in denen die \eemph{Lokalisierung} der Begriffe \eemph{mager} und \eemph{komager} sinnvoll ist.

Der \eemph{Baire'sche Kategoriensatz} identifiziert die Klassen der lokal kompakten sowie der vollst"andig metrisierbaren R"aume als Baire'sche R"aume und wird in \ref{BKS} bewiesen. 

In Satz \ref{Satz_CharBaireEig} stellen wir einen Zusammenhang her zwischen der \eemph{Determiniertheit} und der \eemph{Baire-Eigenschaft} projektiver Mengen. Anschaulich besagt die Baire-Eigenschaft "uber eine Teilmenge eines topologischen Raumes, da"s sie \glqq fast offen\grqq\gs d.h. offen bis auf eine magere Menge ist. 

In Kapitel \ref{Kap_SigmaKategorie} zeigen wir au"serdem, da"s sich einige Resultate f"ur perfekte und abz"ahlbare Teilmengen des Baireraumes analog auch f"ur superperfekte und $\sigma$-beschr"ankte Teilmengen herleiten lassen\gs vorrangig der \eemph{Satz von Cantor und Bendixson}. Die "Ubertragung \ref{CBA} des Satzes von Cantor und Bendixson f"ur $\sigma$-beschr"ankte und superperfekte Teilmengen des Baireraumes besagt, da"s sich jede abgeschlossene Teilmenge des Baireraumes eindeutig als disjunkte Vereinigung einer superperfekten und einer $\sigma$-beschr"ankten Teilmenge darstellen l"a"st. 

In Kapitel \ref{KapSpielCharSigmaBeschr"ankt} beweisen wir mittels der spieltheoretischen Charakterisierung der $\sigma$-Kategorie durch Theorem \ref{SatzSigmaBeschrCharSpiele_2} und eines Resultates von Y.N. Moschovakis sowie eines von D.A. Martin einige Definierbarkeits-Aussagen "uber $\sigma$-Schranken [\siehe S"atze \ref{SatzDefbarkeit_1} und \ref{SatzDefbarkeit_2}].  

\bigskip In den Kapiteln \ref{KapSpielCharMager}, \ref{Kap_SigmaKategorie} und \ref{Kap_VerallgKategorie} beziehen wir uns weitgehend auf S"atze und Ideen, die Alexander Kechris 1977 in seinem Artikel \eemph{\glqq On a notion of smallness for subsets of the Baire space\grqq\ } \cite{Kechris:1977} formuliert hat.

\subsection*{Inhalts"ubersicht}
\emph{Kapitel \ref{Kap_Einf}} gibt zun"achst eine Einf"uhrung in die Ideen zu dieser Arbeit. Unsere anschauliche Vorstellung von der \glqq Gr"o"se\grqq\ eines Objektes spiegelt sich im Begriff des \eemph{$\sigma$-Ideals} auf einer Menge wieder [\siehe Kapitel \ref{Kap_Kleine_Gro"se_Mengen}]. Eine solche mathematische Formulierung unserer anschaulichen Gr"o"senvorstellung nennen wir ein \eemph{Kategorien-Konzept}. In Kapitel \ref{Kap_Kategorien_Konzepte} werden drei klassische Kategorien-Konzepte\gs die Cantor- die Lebesgue- und die Baire-Kategorie\gs vorgestellt. Mittels solcher Konzepte lassen sich in der Mathematik \eemph{unscharfe Begriffe} (oder auch \eemph{Fuzzy-Begriffe}) formulieren, wie etwa der Begriff von einer \eemph{fast offenen Teilmenge} eines toplogischen Raumes [\siehe Kapitel \ref{Kap_Kategorien_Konzepte}]. Dabei weisen die Gr"o"senbegriffe der unterschiedlichen Kategorien-Konzepte strukturelle Gemeinsamkeiten auf. Viele S"atze "uber den einen Gr"o"senbegriff kann man in der entsprechenden Version auch f"ur andere Gr"o"senbegriffe zeigen. In Kapitel \ref{Kap_VereinhProblem} stellen wir daher die Frage nach der Vereinheitlichung derartiger Konzepte. Speziell geht es uns dabei um die Vereinheitlichung von Baire- und $\sigma$-Kategorie in einer verallgemeinerten spieltheoretisch konzipierten Kategorie. 

\bigskip \emph{Kapitel \ref{Kap_TopSp}\ }f"uhrt die wichtigsten spieltheoretischen Grundbegriffe ein: \eemph{topologische Spiele}, \eemph{Gewinnstrategie} und \eemph{Determiniertheit}. Die \eemph{topologischen Spiele} dienen speziell dazu, topologische Eigenschaften spieltheoretisch zu beschreiben. Dazu wird ein Zusammenhang zwischen einer topologischen Eigenschaft der Gewinnmenge des Spieles und der Determiniertheit des Spiels hergestellt. \eemph{Determiniert} nennt man ein Spiel, f"ur das einer der beiden Spieler eine Gewinnstartegie hat. Ein Spieler hat eine \eemph{Gewinnstrategie}, falls es f"ur ihn unabh"angig von den Z"ugen seines Gegners immer einen Weg gibt, das Spiel zu gewinnen. Nimmt man mit dem \eemph{Determiniertheitsaxiom} (in Kurzschreibweise $\AD$) an, da"s alle Spiele der Grundversion determiniert sind, so ist dies unvertr"aglich mit dem \eemph{Auswahlaxiom}, jedoch vertr"aglich mit dem \eemph{abz"ahlbaren Auswahlaxiom} sowie mit dem \eemph{Axiom der abh"angigen Auswahl} [\siehe Kapitel \ref{Kap_Determ_Auswahl}]. Um Inkonsistenzen zu vermeiden darf also zusammen mit dem Determiniertheitsaxiom nicht das volle Auswahlaxiom vorausgesetzt werden.

\bigskip  \emph{Kapitel \ref{Kap_BaireKategorie}\ }stellt die Baire-Kategorie vor. Die Baire-Kategorie l"a"st sich allgemein f"ur topologische R"aume definieren. Die \eemph{kleinen} Mengen der Baire-Kategorie sind die \eemph{mageren} Mengen. Sie setzen sich aus \eemph{nirgends dichten} Mengen zusammen, die sich gleichsam \glqq d"unn machen\grqq, indem ihre Komplemente offene dichte Mengen enthalten [\siehe Kapitel \ref{Kap_Magere_Mengen}]. Beispiel \ref{Bsp_RzAlsNullmengePlusMagereMenge} verdeutlicht, da"s man die Kleinheitsbegriffe verschiedener Kategorien-Konzepte keinesfalls durcheinanderw"urfeln darf\gs \eemph{klein} in dem einen Konzept bedeutet nicht unbedingt auch \eemph{klein} in einem anderen. 

Die \eemph{komageren} Mengen stehen in Baire's Kategorien-Konzept f"ur die \eemph{gro"sen} Mengen und werden daher als die Komplemente \eemph{magerer} Mengen definiert [\siehe Kapitel \ref{Kap_Komagere_Mengen}]. Daher haben die komageren Mengen die zu den mageren Mengen dualen Eigenschaften. Au"serdem wird definiert, wann eine Menge \eemph{mager} bzw. \eemph{komager in einer offenen Menge} ist. 

In einer offenen mageren Menge sind alle Teilmengen mager und komager. Sinnvollerweise betrachtet man dieses Konzept daher nur auf topologischen R"aumen, die keine offenen mageren Mengen enthalten. Derartige R"aume nennt man \eemph{Baire'sche R"aume} [\siehe Kapitel \ref{Kapitel_BairescheRaeume}]. Der Baire'sche Kategoriensatz \ref{BKS} identifiziert zwei gro"se Klassen topologischer R"aume als Baire'sch: die lokal kompakten R"aume und die vollst"andig metrisierbaren R"aume.

F"ur spieltheoretische Charakterisierungen betrachten wir in dieser Arbeit die Baire-Kategorie auf dem \eemph{Baireraum} (s.o.). In Kapitel \ref{Kap_Baire_Raum} ordnen wir den Baireraum in den allgemeineren topologischen Zusammenhang ein. Unter anderem zeigen wir, da"s der Baireraum ein \eemph{polnischer} Raum ist und somit nach dem Baire'schen Kategoriensatz \ref{BKS} $(i)$ insbesondere ein Baire'scher Raum. 

In Kapitel \ref{KapSpielCharMager} wird dann die Baire-Kategorie mittels des Spieles $\BMSo{A}$ f"ur beliebige Teilmengen $A$ des Baireraumes und mittels $\BMSa{B}$ f"ur Projektionen $A=\p(B)$ ($B\subset\omega^\omega\times\omega^\omega$) charakterisiert [\siehe Theoreme \ref{SatzMagerCharSpiele} und \ref{SatzMagerCharSpiele_2}]: Es zeigt sich, da"s eine Teilmenge $A$ des Baireraumes genau dann \eemph{klein} im Sinne der Baire-Kategorie ist, wenn Spieler II in dem Spiel $\BMSo{A}$ eine Gewinnstrategie besitzt und \eemph{gro"s} (in einer offenen Teilmenge des Baireraumes) genau dann, wenn Spieler I eine Gewinnstartegie f"ur das Spiel $\BMSo{A}$ hat [\siehe Theorem \ref{SatzMagerCharSpiele}]. F"ur Projektionen $A=\p(B)$ einer Menge $B\subset\omega^\omega\times\omega^\omega$ erh"alt man ein entsprechendes Resultat mit dem Spiel $\BMSa{B}$ in dem allerdings nur die Hinrichtungen gelten [\siehe Theorem \ref{SatzMagerCharSpiele_2}]. Andererseits bieten die Spiele $\BMSa{B}$ mit $A=\p(B)$ den entscheidenden Vorteil einer gegen"uber $A$ einfacher (mit einem $\exists$-Quantor weniger) definierten Gewinnmenge $B$. 

Kapitel \ref{Kap_BaireKategorie} abschlie"send wird als eine direkte Anwendung von Theorem \ref{SatzMagerCharSpiele_2} ein hinreichendes Kriterium f"ur die Baire-Eigenschaft projektiver Mengen angegeben [\siehe \ref{Satz_CharBaireEig}].

\bigskip  \emph{Kapitel \ref{Kap_SigmaKategorie}\ }stellt die $\sigma$-Kategorie vor. Die kompakten Teilmengen des Baireraumes werden in Kapitel \ref{Kapitel_HeineBorel} als abgeschlossen und "uberall endlich verzweigt charakterisiert [\siehe Satz \ref{Satz_Heine_Borel_2}]. Diese Charakterisierung der kompakten Mengen dient dann im sp"ateren Kapitel \ref{Kapitel_BMager} als Vorlage f"ur den allgemeineren Begriff der \eemph{$\BM{B}$-nirgends dichten} Teilmengen [\siehe Definition \ref{DefBNirgendsDicht}]. 

In Kapitel \ref{Kapitel_SigmaBeschr} wird der Kleinheits-Begriff der \eemph{$\sigma$-Beschr"anktheit} eingef"uhrt (s.o.). Die Charakterisierung der $\sigma$-beschr"ankten Mengen in Satz \ref{SatzSigmaBeschr"ankt} $(ii)$ als Teilmengen $\sigma$-kompakter Mengen dient wiederum im sp"ateren Kapitel \ref{Kapitel_BMager} als Vorlage f"ur den allgemeineren Begriff der \eemph{$\BM{B}$-mageren} Teilmengen [\siehe Definition \ref{DefBMager}]. 

Kapitel \ref{Kapitel_Superperfekt} f"uhrt mit dem Begriff der \eemph{superperfekten} und nicht-leeren Teilmenge (s.o.) eine Beschreibung f"ur \eemph{relativ gro"se} Teilmengen des Baireraumes ein. Aus dem Cantor-Bendixson-Analog in Kapitel \ref{Kapitel_CantorBendixson} folgt, da"s der Abschlu"s des Komplementes einer $\sigma$-beschr"ankten Teilmenge stets eine nicht-leere superperfekte Teilmenge enth"alt [\siehe Korollar \ref{KorollarKomplementSigmaBeschr}]. Unter der Annahme von $\AD$ enth"alt schon das Komplement einer $\sigma$-beschr"ankten Teilmenge eine nicht-leere superperfekte Teilmenge [\siehe Korollar \ref{Korollar_Komplement_Sigmabeschr}].

In Kapitel \ref{Kapitel_CantorBendixson} wird der Satz von Cantor und Bendixson f"ur h"ochstens abz"ahlbare und perfekte auf die $\sigma$-beschr"ankten und superperfekten Teilmengen des Baireraumes "ubertragen. 

In Kapitel \ref{KapSpielCharSigmaBeschr"ankt} wird dann die $\sigma$-Kategorie mittels des Spieles $\BMSb{A}$ f"ur beliebige Teilmengen $A$ des Baireraumes und mittels $\BMSc{B}$ f"ur Projektionen $A=\p(B)$ (f"ur ein $B\subset\omega^\omega\times\lambda^\omega$ und eine unendliche Ordinalzahl $\lambda$) charakterisiert [\siehe Theoreme \ref{SatzSigmaBeschrCharSpiele} und \ref{SatzSigmaBeschrCharSpiele_2}]: Analog zu den Theoremen \ref{SatzMagerCharSpiele} und \ref{SatzMagerCharSpiele_2} f"ur die Baire-Kategorie erh"alt man, da"s eine Teilmenge $A$ des Baireraumes genau dann \eemph{klein} im Sinne der $\sigma$-Kategorie ist, wenn Spieler II in dem Spiel $\BMSb{A}$ eine Gewinnstrategie besitzt und \eemph{relativ gro"s} (in dem Sinne, da"s $A$ eine nicht-leere superperfekte Teilmenge enth"alt) genau dann, wenn Spieler I eine Gewinnstartegie f"ur das Spiel $\BMSb{A}$ hat [\siehe Theorem \ref{SatzSigmaBeschrCharSpiele}]. F"ur Projektionen $A=\p(B)$ (f"ur eine Menge $B\subset\omega^\omega\times\lambda^\omega$ und eine unendliche Ordinalzahl $\lambda$) erh"alt man ein entsprechendes Resultat mit dem Spiel $\BMSc{B}$, in dem allerdings nur die Hinrichtungen gelten [\siehe Theorem \ref{SatzSigmaBeschrCharSpiele_2}]. Wie schon die Spiele $\BMSa{B}$ bieten die Spiele $\BMSc{B}$ mit $A=\p(B)$ den entscheidenden Vorteil einer gegen"uber $A$ einfacher (mit einem $\exists$-Quantor weniger) definierten Gewinnmenge $B$. 

Mittels dieser Charakterisierung (in Form von Korollar \ref{SatzSigmaBeschrCharSpiele_3}) sowie eines Resultates von D.A. Martin (\siehe\ Kapitel \ref{KapSpielCharSigmaBeschr"ankt} Fu"snote \ref{FN_GaleMartin} bzw. \cite[196ff]{KechrisSolovay:1985} oder \cite[30.10]{Kanamori:2003}) und eines weiteren Resultates von Y.N. Moschovakis (\siehe Kapitel \ref{KapSpielCharSigmaBeschr"ankt} Fu"snote \ref{FN_Moschovakis} bzw. \cite[6E.1.]{Moschovakis:1980})  lassen sich nun einige Definierbarkeitsresultate f"ur $\sigma$-Schranken (s.o.) ableiten [\siehe S"atze \ref{SatzDefbarkeit_1} und \ref{SatzDefbarkeit_2}]\footnote{\siehe \eemph{\glqq On a notion of smallness for subsets of the Baire space\grqq\ } von A. Kechris \cite[4]{Kechris:1977}.}.

\bigskip \emph{Kapitel \ref{Kap_VerallgKategorie}\ } definiert auf der Grundlage der Untersuchung von Baire- und $\sigma$-Kategorie in den vorangegangenen Kapiteln eine verallgemeinerte Kategorie und zeigt, da"s die Theoreme \ref{SatzMagerCharSpiele} und \ref{SatzMagerCharSpiele_2} sowie \ref{SatzSigmaBeschrCharSpiele} und \ref{SatzSigmaBeschrCharSpiele_2} in analoger Form auch in der verallgemeinerten Kategorie herleitbar sind.

Zun"achst definieren wir dazu in Kapitel \ref{Kap_VerallgSpiele} ein \eemph{verallgemeinertes Spiel} $\BMSd{\BM{B}}(A)$ f"ur Teilmengen $A$ des Raumes $X^\omega$. In $\BMSd{\BM{B}}(A)$ werden die vorherigen Spiele $\BMSo{A}$ und $\BMSb{A}$ dahingehend verallgemeinert, da"s Spieler II nun keine endlichen Sequenzen oder nat"urliche Zahlen mehr spielt, sondern \eemph{Bedingungen} f"ur den weiteren Spielverlauf.

In Kapitel \ref{Kapitel_BMager} dienen anschlie"send die Charakterisierungen aus Satz \ref{Satz_Heine_Borel_2} f"ur kompakte und aus Satz \ref{SatzSigmaBeschr"ankt} $(ii)$ f"ur $\sigma$-beschr"ankte Teilmengen des Baireraumes als Vorlagen f"ur die Definition der allgemeineren Begriffe \eemph{$\BM{B}$-dirgends dicht} und \eemph{$\BM{B}$-mager} [\siehe Definitionen \ref{DefBNirgendsDicht} und \ref{DefBMager}]. 

In Kapitel \ref{Kap_BPerf} wird der Begriff der \eemph{$\BM{B}$-perfekten} Teilmenge des Bairerauems eingef"uhrt. Beispiel \ref{Bsp_BPerfEnthSuperperfTM} zeigt, da"s die nicht-leeren $\BM{B}$-perfekten Teilmengen eine Verallgemeinerung von Teilmengen des Baireraumes sind, die eine nicht-leere superperfekte Teilmenge enthalten.

Kapitel \ref{KapSpielCharBMager} definiert zus"atzlich zum Spiel $\BMSd{\BM{B}}(A)$ f"ur beliebige Teilmengen $A\subset\BR$ ein entsprechendes projektives Spiel $\BMSe{\BM{B}}(C)$ f"ur Projektionen $A=\p(C)$ (f"ur eine Menge $C\subset\BR\times\lambda$ und eine beliebige unendliche Ordinalzahl $\lambda$) und zeigt, da"s sich die Theoreme aus den Kapiteln \ref{KapSpielCharMager} und \ref{KapSpielCharSigmaBeschr"ankt} zur spieltheoretischen Charakterisierung der Baire- bzw. der $\sigma$-Kategorie analog auf die verallgemeinerte Kategorie "ubertragen lassen: Analog zu den Theoremen \ref{SatzMagerCharSpiele} und \ref{SatzSigmaBeschrCharSpiele} f"ur die Baire-Kategorie bzw. die $\sigma$-Kategorie erh"alt man, da"s eine Teilmenge $A\subset X^\omega$ genau dann \eemph{klein} im Sinne der verallgemeinerten Kategorie ist, wenn Spieler II in dem Spiel $\BMSd{\BM{B}}(A)$ eine Gewinnstrategie besitzt und \eemph{relativ gro"s} (in dem Sinne, da"s $A$ eine nicht-leere $\BM{B}$-perfekte Teilmenge enth"alt) genau dann, wenn Spieler I eine Gewinnstartegie f"ur das Spiel $\BMSd{\BM{B}}(A)$ hat [\siehe Theorem \ref{Theorem_Char_AllgSpiel_OhneZeugen}]. F"ur Projektionen $A=\p(C)$ (f"ur eine Menge $C\subset\omega^\omega\times\lambda^\omega$ und eine unendliche Ordinalzahl $\lambda$) erh"alt man als Verallgemeinerung der Theoreme \ref{SatzMagerCharSpiele_2} und \ref{SatzSigmaBeschrCharSpiele_2} ein entsprechendes Resultat mit dem Spiel $\BMSe{\BM{B}}(C)$ [\siehe Theorem \ref{Theorem_Char_AllgSpiel_MitZeugen}].

%AUSKOMMENTIERT
\begin{comment}
\bigskip Es ergeben sich daher folgende Fragen: Was haben diese unterschiedlichen mathematischen Gr"o"senkonzepte gemeinsam?  sich anhand dieser Gemeinsamkeiten ein allgemeineres mathematisches Konzept aufstellen so, da"s darin auch noch die wichtigsten Theoreme? Mittels \eemph{topologischer Spiele} und speziell der sogenannten \eemph{Banach-Mazur-Spiele} werden diese Konzepte dann untersucht und Gemeinsamkeiten herausgestellt. Auf dieser Basis wird im letzten Kapitel dann eine Verallgemeinerung der Definitionen und wichtigsten S"atze erreicht. 
\end{comment}

\newpage

\section{Kleine und gro"se Mengen}\label{Kap_Kleine_Gro"se_Mengen}
Unsere anschauliche Vorstellung von der \glqq Gr"o"se\grqq\ eines Objektes\begin{comment}(in der Mathematik repr"asentiert durch eine Menge)\end{comment}, genauer einer Teilmenge der reellen Zahlen
oder eines anderen topologischen Raumes, l"a"st sich in folgenden Forderungen zusammenfassen:
\begin{ITEMS}[arabic)]
\item Die leere Menge $\LM$ ist klein,
\item Teilmengen kleiner Mengen sind klein,
\item Das Komplement kleiner Mengen ist gro"s\gs und umgekehrt: das Komplement gro"ser Mengen ist klein,
\item Endliche Vereinigungen kleiner Mengen sind klein und wegen 3) somit auch: endliche Durchschnitte 
gro"ser Mengen sind gro"s,
\item Abz"ahlbare Vereinigungen kleiner Mengen sind klein und wegen 3) somit auch: abz"ahlbare Durchschnitte 
gro"ser Mengen sind gro"s.
\end{ITEMS}

Um dies mathematisch noch pr"azieser zu formulieren, repr"asentieren wir die Begriffe \glqq klein\grqq\ und \glqq gro"s\grqq\ als Systeme $\SYS{K}$ bzw. $\SYS{G}$ von Teilmengen von $X$. Die Elemente des Systems $\SYS{K}$ nennen wir dann die \emph{kleinen} Teilmengen von $X$ und die Elemente des Systems \SYS{G} die \emph{gro"sen} Teilmengen von $X$. Da in 3) gefordert wird, da"s die gro"sen Mengen gerade die Komplemente von kleinen Mengen sein sollen, brauchen wir nur $\SYS{K}$ zu definieren so, da"s Punkt 1), 2) und 5) erf"ullt sind. Punkt 4) folgt dann bereits aus 5), da endlche Vereinigungen $A_0\cup\ldots\cup A_n$ auch als abz"ahlbare Vereinigungen $A_0\cup\ldots\cup A_n\cup\LM\cup\LM\cup\ldots$ aufgefasst werden k"onnen. Eine Pr"azisierung f"ur $\SYS{K}$ so, da"s die Punkte 1), 2) und 5) erf"ullt sind,  leistet der Begriff des $\sigma$-Ideals auf einer Menge $X$.

\begin{definition}\header{$\sigma$-Ideal}\label{Def_SigmaIdeal}\\
Sei $X$ eine Menge. Ein nicht-leeres System $\SYS{I}$ von Teilmengen von $X$ hei"st ein \emph{$\sigma$-Ideal\footnote{\emph{$\sigma$-Ideal:} \siehe etwa \cite[I 5 (11)]{Kuratowski:1976}.}} auf $X$, wenn es abgeschlossen ist gegen"uber Teilmengenbildung und abz"ahlbaren Vereinigungen, d.h. wenn gilt:
\begin{ITEMS}[arabic)]
\item $A\in\SYS{I}$ und $B\subset A\Rightarrow B\in \SYS{I}$,
\item $\SYS{S}\subset\SYS{I}$ abz"ahlbar $\Rightarrow (\bigcup_{S\in\SYS{S}} S )\in\SYS{I}$.
\end{ITEMS}
Die Elemente von $\SYS{I}$ hei"sen \emph{klein in X} oder einfach \emph{klein}. Die Komplemente kleiner Mengen in 
$X$ hei"sen \emph{gro"s in $X$} oder einfach \emph{gro"s}.
\end{definition}

Ist das Mengensystem $\SYS{I}$ nur abgeschlossen gegen"uber Teilmengen und endlichen Vereinigungen, spricht man von einem \emph{Ideal}.

\bigskip Bis auf die Abgeschlossenheit gegen"uber abz"ahlbaren Vereinigungen, verk"orpert also auch der Begriff des \eemph{Ideals} unsere Anschauung von \eemph{kleinen Mengen} (obige Punkte 1), 2) und 4) sind erf"ullt). F"ur die Auszeichnung \eemph{der kleinen Teilmengen} einer Menge $X$ werden wir nachfolgend jedoch in den meisten F"allen auf $\sigma$-Ideale zur"uckgreifen. 

\bigskip Insbesondere ist die leere Menge in jedem ($\sigma$-)Ideal auf einer Menge $X$ enthalten. Ist $\SYS{I}$ eine Ideal auf einer Menge $X$, so erh"alt man mit 
\begin{alignat*}{1}
\SYS{{I}}^\sigma:=&\MNG{A\subset X}{A=\bigcup_{i<\omega} A_i\TT{ f"ur eine Folge } (A_i\in \SYS{I})_{i<\omega}}\\
=&\MNG{A\subset X}{A\subset\bigcup_{i<\omega} A_i\TT{ f"ur eine Folge } (A_i\in \SYS{I})_{i<\omega}}
\end{alignat*}
das kleinste $\sigma$-Ideal auf $X$, das $\SYS{I}$ enth"alt.

\bigskip Vielfach ist der Ausgangspunkt f"ur einen Begriff von einer \eemph{kleinen} Menge die anschauliche Vorstellung, da"s sie sich in der gegebenen Obermenge gleichsam \glqq d"unn macht\grqq. Um dies mathematisch formulieren zu k"onnen, m"ussen die betrachteten Mengen erst mit einer geeigneten Struktur versehen werden. Eine solche Struktur liefert die Theorie der \emph{topologischen R"aume}. An dieser Stelle scheint es daher angebracht, zun"achst einige topologische Grundbegriffe zu erl"autern. Anschlie"send werden dann in Abschnitt \ref{Kap_Kategorien_Konzepte} drei klassische mathematische Gr"o"senkonzepte kurz umrissen.

\subsection*{Topologie}
Ein wichtiges Konstruktions- und Hilfsmittel aus der Analysis ist der Begriff der \eemph{Konvergenz einer Folge}\label{Def_Konv_Ansch}. Anschaulich bedeutet die Konvergenz einer Folge gegen einen Punkt, da"s in jeder \glqq Umgebung\grqq\ des Punktes \glqq fast alle\grqq\ Folgenglieder zu finden sind. 

\bigskip Um diese Anschauung mathematisch exakt fassen zu k"onnen, zeichnet man f"ur jeden Punkt $x$ einer Menge $X$ gewisse Teilmengen als \glqq Umgebungen\grqq\ dieses Punktes aus, die den folgenden sogenannten \emph{Hausdorff'schen Umgebungsaxiomen\footnote{\emph{Hausdorff'sche Umgebungsaxiome}: Benannt nach Felix Hausdorff (1868--1942). 1914 ver"offentlich in \cite[7]{Hausdorff:1914}.}} gen"ugen:
\begin{ITEMS}[arabic)]
\item $x$ geh"ort zu jeder seiner Umgebungen.
\item Jede Obermenge einer Umgebung von $x$ ist wieder eine Umgebung von $x$. $X$ ist eine Umgebung von $x$.
\item Der Durschnitt zweier Umgebungen von $x$ ist wieder eine Umgebung von $x$.
\item Jede Umgebung $U$ von $x$ enth"alt eine Umgebung $V$ von $x$ so, da"s $U$ auch eine Umgebung eines jeden Punktes von $V$ ist.
\end{ITEMS}
Auf dieser Grundlage definiert man eine Teilmenge von $X$ als \emph{offen}, falls sie Umgebung jedes ihrer Punkte ist (z.B. eine Kreisscheibe in $\RZ^2$ ohne ihren Rand). Offen sind also genau diejenigen Teilmengen, die nur aus \glqq inneren\grqq\ Punkten (und weder \glqq "au"seren Punkten\grqq\ noch \glqq Randpunkten\grqq) bestehen in dem Sinne, da"s um jeden Punkt der offenen Menge eine Umgebung gelegt werden kann, die schon ganz in der offenen Menge enthalten ist. Damit lassen sich nun folgende S"atze ableiten:

\medskip
\emph{Satz 1.} Die leere Menge und der Raum selbst sind offen.

\medskip
\emph{Satz 2.} Der Durchschnitt zweier offener Mengen ist offen.

\medskip
\emph{Satz 3.} Die Vereinigung beliebig vieler offener Mengen ist offen.

\medskip
\emph{Satz 4.} Eine Teilmenge $U$ ist genau dann eine Umgebung von $x$, wenn es eine offene Menge $O$ gibt mit $x\in O\subset U$.

\medskip
Verwendet man nun anstelle der Hausdorff'schen Umgebungsaxiome die S"atze 1 bis 3 als Axiome und Satz 4 als Definition f"ur den Umgebungs-Begriff, gelangt man zu der heute "ublichen (gleichwertigen) Definition einer Topologie:

\begin{definition}\header{Topologie}\\
Sei $X$ eine Menge. Ein System $\TOP{X}$ von Teilmengen von $X$ hei"st eine \emph{Topologie} auf $X$, wenn
gilt:
\begin{ITEMS}[arabic)]
\item $\LM\in\TOP{X}$, $X\in\TOP{X}$,
\item $O_1,O_2\in\TOP{X}\Rightarrow O_1\cap O_2\in\TOP{X}$,
\item $\SYS{S}\subset\TOP{X}\Rightarrow (\bigcup_{S\in\SYS{S}} S )\in\TOP{X}$.
\end{ITEMS}
\end{definition}
Das Paar $(X,\TOP{X})$ hei"st ein \emph{topologischer Raum} oder (falls klar ist, was gemeint ist) auch einfach \emph{Raum}. Statt $(X,\TOP{X})$ schreibt man auch nur $X$. Die Elemente von $\TOP{X}$ hei"sen \emph{offen in X} oder einfach \emph{offen}. Die Komplemente offener Mengen in $X$ hei"sen \emph{abgeschlossen in $X$} oder einfach \emph{abgeschlossen}.

\bigskip Die beiden einfachsten Beispiele einer Topoogie $\TOP{X}$ auf einer Menge $X$ sind die \emph{diskrete} und die \emph{indiskrete Topologie}. Erstere besteht aus allen Teilmengen von $X$ (also $\TOP{X}=\POW[X]$), letztere lediglich aus der leeren Menge und $X$ selber (also $\TOP{X}=\{\LM,X\})$. 

\bigskip Ein Standardverfahren, um eine Topologie auf einer Mengen $X$ zu erhalten, ist das \glqq Erzeugen\grqq\ aus einem System $\SYS{Y}$ von Teilmengen von $X$: Ist das System $\SYS{Y}$ \emph{schnittstabil} (d.h. mit $A,B\in \SYS{Y}$ ist auch $A\cap B\in \SYS{Y}$) und l"a"st sich $X$ als Vereinigung von Mengen aus $\SYS{Y}$ schreiben, so ist das System $\TOP{X}$ bestehend aus beliebigen Vereinigungen von Mengen aus $\SYS{Y}$ eine Topologie auf $X$. Man nennt $\TOP{X}$ dann die von $\SYS{Y}$ \emph{erzeugte} Topoplogie.

\bigskip Lassen sich alle offenen Mengen einer Topologie $\TOP{X}$ als Vereinigung von Mengen eines Systems $\SYS{Y}\subset\TOP{X}$ darstellen, so nennt man $\SYS{Y}$ eine \emph{Basis} der Topologie $\TOP{X}$. Die Elemente aus $\SYS{Y}$ hei"sen dann \emph{offene Basismengen}. (Dabei braucht die leere Menge nicht notwendig unter den offenen Basismengen zu sein, da man die leere Menge auch als leere Vereinigung  $\LM=\bigcup_{i\in\LM} O_i$ offener Basismengen $O_i$ erh"alt.)

\bigskip Alle Teilmengen eines topologischen Raumes $(X,\TOP{X})$ werden automatisch wieder zu topologischen R"aumen, wenn man die \glqq Spuren\grqq\ der Topologie $\TOP{X}$ auf den Teilmengen als offene Mengen auszeichnet. Unter den \glqq Spuren\grqq\ von $\TOP{X}$ auf einer Teilmenge $Y\subset X$ sind dabei die Schnitte der in $X$ offenen Mengen mit der Teilmenge $Y$ zu verstehen. Die so definierte \emph{Spurentopologie} $\TOP{Y}:=\MNG{O\cap Y}{O\subset X \TT{ offen}}$ auf $Y$ nennt man auch die von $\TOP{X}$ \emph{induzierte Topologie} oder \emph{Relativtopologie}. 

\bigskip H"aufig sind topologische R"aume aus anderen topologischen R"aumen \glqq zusammengesetzt\grqq. Auf dem kartesischen Produkt $X\times Y$ zweier R"aume $(X,\TOP{X})$ und $(Y,\TOP{Y})$ ist die Menge $\SYS{B}:=\MNG{U\times V}{U\in\TOP{X}\TT{ und }V\in\TOP{Y}}$ zwar keine Topologie, aber ein schnittstabiles System von Teilmengen, das $X\times Y$ enth"alt. Die von diesem System erzeugte Topologie $\TOP{X}\times\TOP{Y}:=\MNG{O\subset X\times Y}{O=\bigcup_i B_i\TT{ f"ur }\{B_i\}_i\subset\SYS{B}}$ nennt man die \emph{Produkttopologie}. Dies ist gerade die kleinste Topologie auf der Menge $X\times Y$ so, da"s die Pojektionen $p_X:\PFEIL{X\times Y}{}{X}$ und $p_Y:\PFEIL{X\times Y}{}{Y}$ stetig (d.h. f"ur alle $O$ offen in $X$ und $O'$ offen in $Y$ sind $p_X^{-1}(O)$ und $p_X^{-1}(O)$ wieder offen in $X\times Y$) sind. Daher nennt man die Produkttopologie auch die \emph{von den Projektionen $p_X, p_Y$ erzeugte Topologie}.

\section{Kategorien-Konzepte}\label{Kap_Kategorien_Konzepte}
\begin{ZITAT}
\flqq Dans de nombreux probl\`emes de Topologie, le r\^ole de la notation d'ensemble de I. cat\'egorie est analogue \`a celui d'ensemble de mesure nulle dans la th\'eorie de la mesure (ensembles \eemph{ n\'egligeables}).\frqq\footnotemark\ (C. Kuratowski in \cite[10 III]{Kuratowski:1958})
\end{ZITAT}\footnotetext{In vielen Problemen der Topologie spielt der Begriff einer Menge von I. Kategorie eine "ahnliche Rolle wie der einer Menge vom Ma"s Null in der Ma"stheorie (\eemph{\glq vernachl"assigbare\grq} Mengen).}\\

\bigskip Konzepte, die mathematisch exakte Begriffe unserer anschaulichen Vorstellung von der \glqq Gr"o"se\grqq\ eines Objektes liefern, sollen hier in Anlehnung an R. Baires Begriffe von der \eemph{1.} und \eemph{2. Kategorie} (\siehe \ref{Kap_Magere_Mengen}) \emph{Kategorien-Konzepte} genannt werden. Ein solches Konzept mu"s den in \ref{Kap_Kleine_Gro"se_Mengen} formulierten Anforderungen gen"ugen. In der Mathematik gibt es drei klassische Kategorien-Konzepte, die hier kurz umrissen werden sollen.

\subsection*{Cantor-Kategorie}
Das einfachste Kategorien-Konzept geht zur"uck auf Georg Cantor (1845--1918) \cite{Cantor:1883}. Eine Menge $A$ hei"st \emph{abz"ahlbar}, wenn sie gleichm"achtig zu den nat"urlichen Zahlen ist, sie hei"st \emph{h"ochstens abz"ahlbar}, falls sie endlich oder abz"ahlbar ist. Letzteres ist gleichbedeutend damit, da"s es eine injektive Abbildung $f:\PFEIL{A}{}{\omega}$ in die nat"urlichen Zahlen gibt. Eine Menge hei"st \emph{"uberabz"ahlbar}, wenn sie nicht h"ochstens abz"ahlbar ist.

\bigskip Betrachtet man die h"ochstens abz"ahlbaren als die \eemph{kleinen} Mengen und die Komplemente h"ochstens abz"ahlbarer Mengen als die \eemph{gro"sen} Mengen, so gen"ugt dieses Konzept den oben in \ref{Kap_Kleine_Gro"se_Mengen} formulierten Anforderungen. Die h"ochstens abz"ahlbaren Mengen bilden ein $\sigma$-Ideal, da h"ochstens abz"ahlbare Vereinigungen h"ochstens abz"ahlbarer Mengen wieder h"ochstens abz"ahlbar sind.

\subsection*{Lebesgue-Kategorie}
Ein weiterer klassischer Gr"o"senbegriff ist das \emph{Volumen} eines K"orpers. Die Formulierung dieses Konzepts geht zur"uck auf Henri-L\'eon Lebesgue (1875--1941) \cite{Lebesgue:1904}\footnote{\emph{Lebesguesche Integartionstheorie:} In \eemph{\glqq Le\c cons sur l'int\'egration et la recherche des fonctions primitives\grqq} (Paris 1904) ver"offentlicht Lebesgue eine Vorlesung "uber Integrationstheorie, die er im akademischen Jahr 1903--1904 am Coll\`ege de France hielt. Im Mittelpunkt steht die Frage, unter welchen Bedingungen das unbestimmte Integral eine Stammfunktion des Integranden ist. Im letzten Kapitel geht er kurz auf seinen eigenen Integralbegriff ein: Analog zum Ma"sproblem formuliert er das Integrationsproblem. Dieses f"uhrt er auf das Ma"sproblem zur"uck und gelangt mit Hilfe des Begriffs der me"sbaren Funktion zur Definition des Lebesgue-Integrals. Mit diesem Buch wurde die Lebesguesche Integartionstheorie allgemein zug"anglich.}. \KOM[LIT Leb02 (Elstrodt: S.157 2. Abschn. Lit. 2)]

\bigskip F"ur das \emph{Lebesgue-Ma"s} werden zun"achst die Volumen von einfach zu messenden Ausgangsmengen bestimmt. F"ur $a=(a_1,\ldots,a_n), b:=(b_1,\ldots,b_n)\in\RZ^n$ gelte $a\leq b$ genau dann, wenn $a_i\leq b_i$ f"ur $i=1,\ldots,n$. F"ur beliebige $a,b\in\RZ^n$ sei $(a,b]:=\MNG{z\in\RZ^n}{a<z\leq b}$. F"ur \eemph{$n$-dimensionale (links offene) Quader} $I\in \SYS{I}^n:=\MNG{(a,b]}{a,b\in\RZ^n}\subset\RZ^n$ definiert man ihr Volumen 
\[
\lambda^n(I):=\prod_{i=1}^n (b_i-a_i).
\] 
F"ur \eemph{$n$-dimensionale Figuren} $F:=\sum_{j=1}^d I_j\in \SYS{F}^n:=\MNG{\sum_{j=1}^d(a^{(j)},b^{(j)}]}{a^{(j)},b^{(j)}\in\RZ^n\TT{ f"ur } 1\leq j\leq d, d\geq 1}\subset\RZ^n$ definiert man 
\begin{equation}
\lambda^n(F):=\sum_{j=1}^d \lambda^n(I_j)\label{Frml_LebMassVonFigur}.
\end{equation}
Dabei ist $\lambda^n$ durch (\ref{Frml_LebMassVonFigur}) auf $\SYS{F}^n$ wohldefiniert (\siehe \cite[4.2]{Alsmeyer:1998}).

\bigskip In der Ma"stheorie zeigt man, da"s $\lambda^n$ eindeutig zu einem Ma"s (der Einfachheit halber ebenfalls mit $\lambda^n$ bezeichnet) auf der von $\SYS{F}^n$ erzeugten $\sigma$-Algebra $\sigma(\SYS{F}^n)=:\SYS{B}$ fortgesetzt werden kann.\footnote{Man zeigt, da"s $\lambda^n$ ein \eemph{Pr"ama"s} (d.h. ein Ma"s, das auf einem Ring [\siehe Konventionen S. \pageref{Kap_Konventionen}] definiert ist) auf dem \eemph{Ring} $\SYS{F}^n$ auf $\RZ^n$ ist. Nach dem \eemph{Fortsetzungssatz} der Ma"stheorie (\siehe \cite[3.5]{Alsmeyer:1998} bzw. \cite[5.2]{Bauer:1974}) kann $\lambda^n$ auf wenigstens eine Art zu einem Ma"s auf der von $\SYS{F}^n$ \eemph{erzeugten $\sigma$-Algebra} [\siehe Konventionen S. \pageref{Kap_Konventionen}] $\sigma(\SYS{F}^n)=:\SYS{B}$ (genannt $\sigma$-Algebra der \eemph{Borelmengen}) fortgesetzt werden. Da das Pr"ama"s $\lambda^n$ \eemph{$\sigma$-endlich} ist (d.h. es gibt eine aufsteigende Folge $\Omega_1\subset\Omega_2\subset\ldots$ von Teilmengen von $\Omega$ so, da"s $\bigcup_{i\geq 1} \Omega_i=\Omega$ und $\lambda^n(\Omega_i)<\infty$), ist diese Fortsetzung nach dem \eemph{Eindeutigkeitssatz} der Ma"stheorie (\siehe \cite[3.6]{Alsmeyer:1998}) eindeutig bestimmt.}
Die \eemph{(lebesgue-)me"sbaren Teilmengen} von $\RZ^n$ erh"alt man durch \eemph{Vervollst"andigung} der $\sigma$-Algebra $\SYS{B}$ bzgl. des Ma"ses $\lambda^n$, d.h. indem man "ubergeht zur $\sigma$-Algebra $\utilde{B}:=\MNG{A\cup N}{A\in\SYS{B}, N\subset B\TT{ f"ur ein }B\in\SYS{B}\TT{ mit }\lambda^n(B)=0}$, das Ma"s $\lambda^n$ wird auf $\utilde{B}$ fortgesetzt durch $\utilde{\lambda}^n(A\cup N):=\lambda^n(A)$ f"ur $A\in\SYS{B}$ und $N\subset B\TT{ f"ur ein }B\in\SYS{B}\TT{ mit }\lambda^n(B)=0$ (\siehe etwa \cite[II 6.3]{Elstrodt:1996}). 

\bigskip Eine me"sbare Menge $A\in\utilde{B}$ hei"st \emph{Lebesgue-Nullmenge} oder einfach \emph{Nullmenge}, wenn $\utilde{\lambda}^n(A)=0$ gilt. Der Einfachheit halber werde $\utilde{\lambda}^n$ von nun an ebenfalls mit $\lambda^n$ bezeichnet. Statt $\lambda^1$ schreiben wir auch $\lambda$.

\bigskip F"ur eine anschauliche Charakterisierung des Lebesgue-Ma"ses definiert man f"ur beliebige Teilmengen $A$ von $\RZ^n$ das \eemph{innere Ma"s} $\lambda_*(A)$ und das \eemph{"au"sere Ma"s} $\lambda^*(A)$ wie folgt:

F"ur das innere Ma"s $\lambda_*$ wird $A$ von innen her mit Borelmengen approximiert bzw. \glqq ausgesch"opft\grqq:
\[
\lambda_*(A):=\sup\MNG{\lambda^n(B)}{A\supset B\TT{ mit }B\in\SYS{B}}.
\] 
Umgekehrt wird f"ur das "au"sere Ma"s $\lambda^*$ die zu messende Menge $A$ von au"sen mit Borelmengen approximirt bzw. \glqq verpackt\grqq:
\[
\lambda^*(A):=\inf\MNG{\lambda^n(B)}{A\subset B\TT{ mit }B\in\SYS{B}}.
\] 
F"ur Teilmengen $A\subset\RZ^n$ mit $\lambda^*(A)<\infty$ gilt nun:
\[
\lambda_*(A)=\lambda^*(A)\Leftrightarrow A\in\utilde{B}
\]
und in diesem Fall ist $\lambda_*(A)=\lambda^*(A)=\lambda^n(A)$ (\siehe \cite[1.5.5]{Cohn:1982}).

\bigskip Betrachtet man die Nullmengen in $\RZ^n$ als \eemph{kleine} Mengen und die Komplemente von Nullmengen als \eemph{gro"se} Mengen, so gen"ugt dieses Konzept erneut den Anforderungen, die wir in Kapitel \ref{Kap_Kleine_Gro"se_Mengen} an ein Konzept von kleinen und gro"sen Mengen gestellt haben. Ein Beweis f"ur die $\sigma$-Idealeigenschaft der Nullmengen findet sich etwa in \cite[II,1.9]{Elstrodt:1996}.

\subsection*{Baire-Kategorie}
Im Jahr 1899 formulierte R. L. Baire sein Kategorien-Konzept \cite{Baire:1899}. Zun"achst werden \emph{nirgends dichte} Teilmengen $A$ der reellen Zahlen charakterisiert als Mengen, die in keinem Intervall dicht sind oder damit gleichbedeutend als Mengen, f"ur die jedes Intervall $I$ ein Teilintervall $J\subset I$ besitzt, das im Komplement von $A$ enthalten ist. Eine solche Menge stellt man sich am besten als Menge \glqq mit lauter L"ochern\grqq\ vor. Weitere "aquivalente Formulierungen sind: Eine Menge $A$ ist nirgends dicht genau dann, wenn das Komplement von $A$ eine offene dichte Teilmenge enth"alt und genau dann, wenn der Abschlu"s $\ABS[A]$ von $A$ keine inneren Punkte enth"alt\gs also $\INT[{\ABS[A]}]=\LM$. Mit diesen "aquivalenten Formulierungen erh"alt man eine allgemeine Definition f"ur beliebige topologische R"aume (\siehe Kapitel \ref{Kap_Magere_Mengen}).

\bigskip 
Eine abz"ahlbare Vereinigung nirgends dichter Mengen ist i.a. nicht wieder nirgends dicht: z.B. sind die rationalen Zahlen eine abz"ahlbare Vereinigung von Singleton-Mengen (also von Mengen, die nur aus einem Punkt bestehen) und solche Singleton-Mengen $\{x\}$ sind in $\RZ$ nirgends dicht (das Komplement $\RZ-\{x\}$ ist eine dichte offene Teilmenge), andererseits liegen die rationalen Zahlen aber auch dicht in den reellen Zahlen. Daher wird auf den Begriff einer nirgends dichten Menge aufbauend eine Teilmenge der reellen Zahlen als \emph{von erster Kategorie} bezeichnet, falls sie eine abz"ahlbare Vereinigung nirgends dichter Teilmengen von $\RZ$ ist. Teilmengen der reellen Zahlen, die nicht von erster Kategorie sind, hei"sen \emph{von zweiter Kategorie}. 

\bigskip Eine Menge $A$ von erster Kategorie ist nicht notwendig selber eine Menge mit \glqq lauter L"ochern\grqq\ in dem Sinne, da"s jede offene Menge das Komplement von $A$ in einer offenen Menge schneidet, kann aber durch derartige Mengen approximiert werden. Eine Menge von ertser Kategorie enth"alt anschaulich eine dichte Menge von \glqq L"ucken\grqq\gs damit ist gemeint, da"s kein Intervall $[a,b]\subset\RZ$ mit $a<b$ als Vereinigung einer Folge von Mengen von erster Kategorie dargestellt werden kann [\siehe Beispiel \ref{Bsp_mager_lauterLuecken}]. 

\bigskip
Mit den Teilmengen der reellen Zahlen von erster Kategorie als \eemph{kleine} Mengen und den Komplementen kleiner Mengen als \eemph{gro"se} Mengen beschreibt auch dieses Konzept kleine und gro"se Mengen im Sinne von Kapitel \ref{Kap_Kleine_Gro"se_Mengen}. Die Teilmengen der reellen Zahlen von erster Kategorie bilden ein $\sigma$-Ideal, da abz"ahlbare Vereinigungen abz"ahlbarer Vereinigungen wieder abz"ahlbar sind.

\subsection*{Fuzzy-Begriffe}
\begin{ZITAT}
\flqq Pour qu'une fonction born\'ee $f(x)$ soit int\'egrable, il faut et il suffit que l'ensemble de ses points de discontinuit\'e soit de mesure nulle.\frqq\footnotemark\ (H. Lebesgue in \cite[S. 45]{Lebesgue:1972:2})
\end{ZITAT}\footnotetext{Daf"ur, da"s eine beschr"ankte Funktion $f(x)$ [Riemann-]integrierbar ist, ist es notwendig und hinreichend, da"s die Menge ihrer Unstetigkeitsstellen vom Ma"s Null ist.}\\

\begin{ZITAT}
\flqq \ldots je dirai qu'une condition est remplie \eemph{presque partout} lorsqu'elle est v\'erifi\'ee en tout point, sauf aux points d'un ensemble de mesure nulle.\frqq\footnotemark\ (H. Lebesgue in \cite[S. 200]{Lebesgue:1972:2})
\end{ZITAT}\footnotetext{\ldots ich werde sagen, da"s eine Bedingung \eemph{fast "uberall} erf"ullt ist, wenn sie f"ur jeden Punkt gilt bis auf Punkte einer Menge vom Ma"s Null.}\\

\bigskip Ideale und $\sigma$-Ideale erlauben die Definition von \emph{Fuzzy-Begriffen}. Darunter sind Begriffsbildungen zu verstehen, die in einem bestimmten mathematisch exakt formulierten Sinne \glqq unscharf\grqq\ sind. Typischerweise definiert man dabei, wann ein Objekt eine vorgegebene Eigenschaft \glqq fast\grqq\ hat, indem man mittels der kleinen Objekte die gerade noch zul"assige Abweichung von dieser Eigenschaft festlegt. Hierzu sollen nun einige Beispiele gegeben werden:

\bigskip Ein klassisches Beispiel f"ur eine solche \glqq unscharfe\grqq\ Begriffsbildung ist der Konvergenzbegriff der Analysis: Eine Folge $(x_i)_{i<\omega}$ reeller Zahlen \eemph{konvergiert}\label{Def_Konv} gegen eine reelle Zahl $x$, falls f"ur jedes reelle $\epsilon>0$ \emph{fast alle} Folgenglieder einen Abstand $<\epsilon$ zu $x$ haben. Damit ist gemeint, da"s alle bis auf endlich viele Folgenglieder in einem Abstand $<\epsilon$ zu $x$ liegen. Alternativ l"a"st sich dies so formulieren: F"ur alle offenen $\epsilon$-Umgebungen $O_\epsilon:=\MNG{y\in\RZ}{\left|x-y\right|<\epsilon}\subset\RZ$ von $x$ gibt es eine endliche Menge $N\subset\RZ$ so, da"s gilt:
\[
\MNG{x_i}{i<\omega}\o N\subset O_\epsilon.
\]
Eine formal gesehen ganz "ahnliche (freilich inhaltlich v"ollig verschiedene) Definition taucht nun in der Ma"s- und Integrationstheorie auf:

\bigskip In der Ma"s- und Integrationstheorie m"ussen viele Aussagen im Hinblick auf ihre Eindeutigkeit mit einer Einschr"ankung versehen werden, bei der der Begriff der \eemph{Lebesgue-Nullmenge} eine wichtige Rolle spielt\footnote{\emph{Eindeutigkeit in der Integrationstheorie:} Dies gilt nicht nur f"ur die Lebesgue'sche Integrationstheorie, sondern allgemeiner f"ur Aussagen der Integrationstheorie in beliebigen Ma"sr"aumen (\siehe etwa \cite[10]{Alsmeyer:1998}).}:
Sei $A$ eine Lebesgue-me"sbare Teilmenge von $\RZ$ und $E\subset \RZ$ eine Eigenschaft von Elementen aus $\RZ$. Dann sagt man, die Eigenschaft $E$ \emph{gilt Lebesgue-fast "uberall} auf $A$, wenn es eine Lebesgue-Nullmenge $N$ gibt so, da"s gilt:
\[
A\o N\subset E.
\]
Dies erlaubt es nun, Begriffe und Aussagen in der Integrationstheorie zu formulieren, die im obigen Sinne eindeutig bis auf Lebesgue-Nullmengen sind. Beispielsweise wird die Gleichheit von Funktionen $f,g$ auf $\RZ$ ausgedr"uckt durch:
\[
f=g\TT{ Lebesgue-fast "uberall}
\]
oder die Endlichkeit der Werte einer Funktion $f$ auf $\RZ$ durch:
\[
|f|<\infty\TT{ Lebesgue-fast "uberall}.
\]
Dies ist notwendig, da es f"ur die Aussagen in der Integrationstheorie i.a. auf Nullmengen nicht ankommt.

\bigskip In der deskriptiven Mengenlehre sagt man, eine Teilmenge $A$ eines Baire'schen Raumes $X$ habe die \emph{Baire-Eigenschaft}, falls sie \emph{fast offen} ist, d.h.:
\[
A\sd O\TT{ ist mager}. 
\]
f"ur eine offene Teilmenge $O$ von $X$.

\bigskip Allgemeiner l"a"st sich dies f"ur eine Menge $X$ mit einem $\sigma$-Ideal $\SYS{I}\subset \POW[X]$ und einer Eigenschaft $\SYS{E}\subset\POW[X]$ von Teilmengen von $X$ definieren: Eine Teilmenge $A$ von $X$ habe \emph{$\SYS{I}$-fast} die Eigenschaft $\SYS{E}$, falls es eine Menge $B\in \SYS{E}$ (d.h. $B$ hat die Eigenschaft $\SYS{E}$) gibt so, da"s gilt: $A\sd B\in\SYS{I}\TT{ (d.h. $A\sd B$ ist $\SYS{I}$-klein})$.
F"ur eine Eigenschaft $E\subset X$ von Elementen von $X$ l"a"st sich allgemein formulieren: Eine Teilmenge $A$ von $X$ hat \emph{$\SYS{I}$-fast "uberall} die Eigenschaft $E$, falls es eine Menge $N\in\SYS{I}$ gibt so, da"s gilt: $A\o N\subset E$.

\bigskip Unscharfe Begriffe spielen in der Mathematik eine wichtige Rolle, wenn es um die pr"azise Formulierung von Aussagen geht, f"ur die es \glqq auf kleine Abweichungen nicht ankommt\grqq. 

\begin{comment}
\bigskip Neben der Mathematik spielen derartige Konzepte auch in Bereichen wie der Informatik und dem Maschinenbau eine gro"se Rolle. Speziell in der Me"s- und Regeltechnik, wenn es darum geht, einer Maschine beizubringen, wann sich zwei Situationen oder "au"sere Reize \glqq nicht wesentlich voneinander unterscheiden\grqq, so da"s die Maschine darauf in gleicher oder "ahnlicher Weise reagieren kann.\KOM[VERS\\DATA MINIG?]
\end{comment}

\bigskip Nach diesen Beispielen von Kategorien-Konzepten wollen wir nun die Frage nach der Vereinheitlichung derartiger Konzepte stellen. Dabei beschr"anken wir uns auf Kategorien-Konzepte f"ur den Baireraum.

\section{Vereinheitlichungsproblem}\label{Kap_VereinhProblem}
\begin{ZITAT}
\flqq It is natural to ask whether these notions of smallness are related.\frqq\ (J.C. Oxtoby in \cite[S. 4]{Oxtoby:1980})\footnotemark
\end{ZITAT}\footnotetext{Oxtoby bezieht sich dabei auf das Verh"altnis zwischen Nullmengen und Mengen von erster Kategorie im Raum der reellen Zahlen. Wir wollen diese Frage nun "ubertragen auf Mengen von erster Kategorie und $\sigma$-beschr"ankte Mengen im Baireraum.}\\

\bigskip In den Kapiteln \ref{Kap_BaireKategorie} und \ref{Kap_SigmaKategorie} werden das Kategorien-Konzept von Baire und die $\sigma$-Kategorie vorgestellt. Es zeigt sich dass die Gr"ossenbegriffe viele strukturelle Gemeinsamkeiten aufweisen. Viele S"atze "uber einen Gr"ossenbegriff kann man in der entsprechenden Version auch f"ur andere Gr"ossenbegriffe zeigen [\siehe die Kapitel \ref{KapSpielCharMager} und \ref{KapSpielCharSigmaBeschr"ankt}]. Daher stellt sich nun folgende Frage: 

\subsection*{Frage}
Wie lassen sich Baire- und $\sigma$-Kategorie in einem allgemeineren Kategorien-Konzept vereinheitlichen so, da"s sich darin die wichtigsten Ergebnisse zu den beiden Konzepten reproduzieren lassen?
\begin{comment}
Gibt es ein allgemeineres Kategorien-Konzept, das Baire- und $\sigma$-Kategorie (und m"oglichst noch andere) als Spezialf"alle enth"alt und in dem sich die wichtigsten Ergebnisse zu den beiden Konzepten reproduzieren lassen?
\end{comment}

\subsection*{Vorgehensweise} 
F"ur die Erstellung eines allgemeineren Kategorien-Konzeptes, wird hier zun"achst nur die Frage nach einer vereinheitlichten Fassung der Baire- und der $\sigma$-Kategorie gestellt. In den Kapiteln \ref{Kap_BaireKategorie} und \ref{Kap_SigmaKategorie} werden dann die Baire- und die $\sigma$-Kategorie spieltheoretisch charakterisiert. Auf den Gemeinsamkeiten dieser Charakterisierungen aufbauend wird anschlie"send in Kapitel \ref{Kap_VerallgKategorie} ein verallgemeinertes Kategorien-Konzept vorgestellt, das die beiden erstgenannten vereinheitlicht und auch noch weitere Spezialf"alle umfasst. In dem allgemeineren Konzept erh"alt man Theoreme analog zu denen der Kapitel \ref{Kap_BaireKategorie} und \ref{Kap_SigmaKategorie}. 

\bigskip Die spieltheoretischen Charakterisierungen werden dabei mit Hilfe spezieller topologischer Spiele vorgenommen\gs den sogenannten \eemph{Banach-Mazur-Spielen}. 
\begin{comment}
Auf diese soll nun im n"achsten Kapitel eingegangen werden.
\end{comment}

\begin{comment}
\subsection*{Untersuchung mittels Spielen}\KOM[GEF"ALLT MIR NOCH NICHT SO]
In der Spieltheorie gibt es eine gro"se Vielzahl von Spielen f"ur die verschiedensten zu untersuchenden Gegenst"ande. Der "Ubersicht halber teilt man diese nach ihren gemeinsamen Eigenschaften in Klassen ein. 
\end{comment}

\begin{comment}
\bigskip Der Grundtypus, auf den sich alle hier verwendeten mathematischen Spiele zur"uckf"uren lassen, sieht wie folgt aus: Eine Menge $A\subset \RZ$ wird als \emph{Gewinnmenge} eines Spieles betrachtet. Spieler I und II spielen abwechselnd nat"urliche Zahlen. Das Spielergebnis wird bestimmt, indem die gespielten nat"urlichen Zahlen in der Reihenfolge, in der sie gespielt wurden, als Dezimaldarstellung einer reellen Zahl $x$ in $\RZ$ aufgefasst werden. Liegt $x$ in $A$, so gewinnt Spieler I, andernfalls gewinnt Spieler II. 
\end{comment}

\bigskip Derartige Spiele erlauben dann eine spieltheoretische Charakterisierung topologischer Eigenschaften $E_1$ und $E_2$ der Menge $A$, indem ein Zusammenhang zu den den Gewinnstrategien der Spieler f"ur ein Spiel mit Gewinnmenge $A$ hergestellt wird\gs etwa derart, da"s gilt: \glqq Spieler I hat eine Gewinnstrategie f"ur das Spiel mit Gewinnmenge $A$ genau dann, wenn die Mange $A$ die topologische Eigenschaft $E_1$ besitzt. Spieler II hat eine Gewinnstrategie f"ur das Spiel mit Gewinnmenge $A$ genau dann, wenn die Menge $A$ die topologische Eigenschaft $E_2$ besitzt.\grqq

\bigskip Im nachfolgenden Kapitel werden nun derartige topologischen Spiele vorgestellt und einige wichtige Eigenschaften besprochen.

%==============================KAPITEL=============================
\chapter{Topologische Spiele}\label{Kap_TopSp}
\setcounter{lemma}{0}

Die Spieltheorie ist eine eigenst"andige mathematische Disziplin und dient als wissenschaftliche Grundlage f"ur zahlreiche empirisch arbeitende Wissenschaften. In den Wirtschafts- und Sozialwissenschaften werden mit spieltheoretischen Methoden wirtschaftliche, gesellschaftliche und milit"arische Vorg"ange modelliert, in der Biologie nutzt man Spieltheorie zur Analyse von Ph"anomenen wie der Ausbreitung von Genen oder Spezies. In der Informatik und im Maschinenbau kommen spieltheoretische Methoden etwa bei der Verifikation \eemph{reaktiver Systeme\footnote{\emph{Reaktive Systeme:} Computer-Systeme, die in st"andiger Interaktion mit ihrer Umgebung stehen, indem sie Reize aus ihrer Umgebung empfangen und verarbeiten und als Reaktion darauf steuernd und regelnd auf diese Umgebung einwirken. Sie werden eingesetzt bei unterschiedlichsten Anwendungen wie "Uberwachung und Steuerung von medizinischen Ger"aten, Fahrzeugen, Produktionsanlagen und Kraftwerken. Da viele dieser Anwendungen sicherheitskritisch sind, ist man sehr daran interessiert, schon im Vorfeld bei der Entwicklung\gs etwa durch Methoden der spieltheoretischen Verifikation\gs Fehlerquellen auszuschlie"sen.}} zum Einsatz. 
%(\siehe etwa \cite{Manna:1992} oder \cite{Manna:1995}

\begin{comment}
Vor allem im \eemph{Operations Research\footnote{\emph{Operations Research:} Der Begriff kommt urspr"unglich aus dem Milit"arwesen aus der Zeit des zweiten Weltkrieges. Eine Fragestellung damals war z.B. die optimale Anzahl von Schiffen und Begleitschutz f"ur Schiffskonvois. Heute "uberwiegen ingenieurs- und wirtschaftsinformatische Fragestellungen. Bekannte klassische Probleme aus dem Operations Research sind z.B. das K"onigsberger Br"uckenproblem oder das Problem eines Handelsreisenden.}} 
\end{comment}

\bigskip Die mathematische Spieltheorie betrachtet eine Vielzahl von Spielen unterschiedlichster Typen. Ein allgemeiner Spielbegriff findet sich etwa in \cite[7]{vNeumann:1944}.

\bigskip Die \eemph{topologischen Spiele} dienen im allgemeinen dazu, topologische Eigenschaften spieltheoretisch zu beschreiben. Dazu wird ein Zusammenhang zwischen einer topologischen Eigenschaft der Gewinnmenge und den Gewinnstrategien des Spiels hergestellt. 

\bigskip In einem \eemph{topologischen Spiel} w"ahlen die Spieler Objekte, die in Zusammnheng stehen mit der Toplogie eines Raumes\gs dies k"onnen etwa Punkte, offene oder abgeschlossene Teilmengen, abgeschlossene H"ullen oder dergleichen sein (\siehe etwa \cite{Telarsky:1987}). Zu dieser Art Spiel geh"oren auch die sp"ater noch zu definierenden Spiele $\BMSo{A}$ und $\BMSa{B}$ aus Kapitel \ref{Kap_BaireKategorie}, $\BMSb{A}$ und $\BMSc{B}$ aus Kapitel \ref{Kap_SigmaKategorie} und $\BMSd{\BM{B}}(A)$ und $\BMSe{\BM{B}}(C)$  aus Kapitel \ref{Kap_VerallgKategorie}. Der gemeinsame Grundtyp dieser Spiele sind die \eemph{unendlichen Zwei-Personen-Spiele mit perfekter Information\footnote{\emph{Unendliche Zwei-Personen-Spiele mit perfekter Information:} Die Spiele hei"sen \glqq unendlich\grqq, da unendlich viele Spielz"uge gemacht werden. Sie hei"sen \glqq mit perfekter Information\grqq, da die Spieler jederzeit "uber alle gemachten Spielz"uge im Bilde sind. Dies w"are z.B. nicht der Fall, wenn die Spieler einige ihrer Z"uge verdeckt oder gleichzeitig abg"aben.}} aus Kapitel \ref{Kap_UnendlZweiPersSpiel}.

\bigskip In Kapitel \ref{Kap_Determ_Auswahl} wird das \eemph{Determiniertheitsaxiom} eingef"uhrt. Dieses besagt, da"s gewisse mathematische Spiele immer determinieren (d.h. es gibt stets eine Gewinnstrategie f"ur einen der beiden Spieler). Das Determiniertheitsaxiom ist nicht vertr"aglich mit dem so elementar erscheinenden \eemph{Auswahlaxiom} [\siehe Lemma \ref{Lemma_Inkonsistenz_AD_AC}]. Das Auswahlaxiom besagt, da"s es f"ur jedes System nicht-leerer Mengen m"oglich ist, aus jeder Menge dieses Systems ein Element auszuw"ahlen. Allerdings ist das Determiniertheitsaxiom vertr"aglich mit dem schw"acheren \eemph{abz"ahlbaren Auswahlaxiom} [\siehe Lemma \ref{Lemma_Konsistenz_AD_AC_omega}]. Dieses garantiert die Auswahl von Elementen f"ur abz"ahlbare Mengensysteme. Um Inkonsistenzen zu vermeiden, mu"s man also stets darauf achten, das Determiniertheitsaxiom nicht zusammen mit dem vollen Auswahlaxiom zu verwenden. 

\bigskip Bevor jedoch auf diese Zusammenh"ange in Kapitel \ref{Kap_Determ_Auswahl} n"aher eingegangen wird, sollen hier zun"achst in den Kapiteln \ref{Kap_UnendlZweiPersSpiel} und \ref{Kap_Banach_Mazur_Spiele} die nachfolgend immer wieder ben"otigten grundlegenden spieltheoretischen Begriffe eingef"uhrt werden.

\section{Unendliche Zwei-Personen-Spiele mit perfekter Information}\label{Kap_UnendlZweiPersSpiel}

Die Grundvariante eines unendlichen Zwei-Personenspieles mit perfekter Information sieht wie folgt aus:

\begin{definition}\header{Das Spiel $\textup{G}(A)$}\label{Def_SpielGrundvariante}\\
Sei $X$ eine beliebige Menge. Zu einer beliebigen Teilmenge $A\subset X^\omega$ definiert man das \emph{Spiel $\textup{G}(A)$} wie folgt:
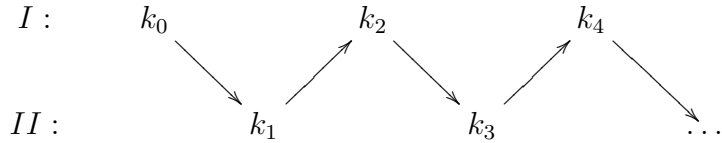
\[
\xymatrix{
I: &k_0 \ar[dr] & 					 	 &k_2 \ar[dr]  &						&k_4 \ar[dr]  &\\
II:&		 				&k_1 \ar[ur]	 &						 &k_3 \ar[ur]	&							&\ldots
}
\]
Spieler I und Spieler II spielen abwechselnd\gs Spieler I spielt ein $k_0\in X$, Spieler II spielt eine $k_1\in X$, Spieler I spielt ein $k_2\in X$, Spieler II spielt ein $k_3\in X$ u.s.w.. Sei nun $f=(k_0,k_1,\ldots)\in X^\omega$\gs man nennt $f$ auch die \emph{Spielfolge} oder \emph{Partie}. 
Man sagt, da"s \emph{Spieler I gewinnt}, falls gilt:
\[
f\in A.
\]
Ansonsten sagt man, da"s \emph{Spieler II gewinnt}. $A$ nennt man dabei auch die \emph{Gewinnmenge}. Hat die Gewinnmenge eine Eigenschaft $\SYS{E}\subset\POW[X]$ (etwa $A$ offen, abgeschlossen, $\SIGMA{1}{n}{}(X)$ f"ur einen polnischen Raum $X$\ldots), so sagt man auch: \emph{das Spiel $\textup{G}(A)$ hat die Eigenschaft $\SYS{E}$} oder bezeichnet $\textup{G}(A)$ als ein \emph{$\SYS{E}$-Spiel}.

\bigskip Ein solches Spiel bezeichnet man auch als \emph{unendliches Zwei-Personen-Spiel mit perfekter Information}, da es nicht nach endlich vielen Spielz"ugen endet und jeder der beiden Spieler zu jedem Zeitpunkt die komplette Information dar"uber hat, wie der bisherige Spielverlauf aussieht\gs insbesondere wie der andere Spieler gespielt hat.
\end{definition}

\begin{comment}
\bigskip F"ur $X=\omega\times\ldots\times\omega$ kann man nun auch Spiele $\textup{G}(A)$ f"ur $A\subset X^\omega=(\omega\times\ldots\times\omega)^\omega=\omega^\omega\times\ldots\times\omega^\omega$\KOM\KOM[F"UR\\DEFINIERBARK\\KAPITEL] betrachten, indem man die gespielten $k_i\in \omega\times\ldots\times\omega$ in nat"urliche Zahlen kodiert.
\end{comment}

\begin{definition}\header{(Gewinn-)Strategie f"ur I}\label{Def_Gewinnstrategie_I}\\
Eine Strategie f"ur Spieler I ist eine Funktion
\begin{alignat*}{1}
&\sigma:\PFEIL{ \MNG{s\in X^{<\omega}}{\lng(s) \TT{ ist gerade}} }{}{X}.
\end{alignat*}
Statt $\sigma(s)$ schreiben wir auch kurz $\sigma s$.
Man sagt, Spieler I \emph{spielt nach der Strategie $\sigma$}, falls gilt:
\begin{alignat*}{1}
&k_0=\sigma()\\
&k_2=\sigma(k_0,k_1)\\
&k_4=\sigma(k_0,k_1,k_2,k_3)\\
&\ldots
\end{alignat*}
Die $k_0,k_2,k_4\ldots$ nennt man dann auch \emph{$\sigma$-gespielt}. Falls Spieler I nach der Strategie $\sigma$ spielt, ist die Spielfolge bereits festgelegt durch $\sigma$ und $p:=(k_1,k_3,k_5,\ldots)$ und wird mit $\sigma * p$ bezeichnet.

\bigskip $\sigma$ ist eine \emph{Gewinnstrategie f"ur I}, falls 
\begin{alignat*}{1}
&\FA[p\in\omega^{<\omega}](\sigma * p\in A).
\end{alignat*}
\end{definition}

Analog definiert man (Gewinn-)Strategien f"ur Spieler II:

\begin{definition}\header{(Gewinn-)Strategie f"ur II}\label{Def_Gewinnstrategie_II}\\
Eine Strategie f"ur Spieler II ist eine Funktion
\begin{alignat*}{1}
&\tau:\PFEIL{ \MNG{s\in X^{<\omega}}{\lng(s) \TT{ ist ungerade}} }{}{X}.
\end{alignat*}
Statt $\tau(s)$ schreiben wir auch kurz $\tau s$. Man sagt, Spieler II \emph{spielt nach der Strategie $\tau$}, falls gilt:
\begin{alignat*}{1}
&k_1=\tau(k_0)\\
&k_3=\tau(k_0,k_1,k_2)\\
&k_5=\tau(k_0,k_1,k_2,k_3,k_4)\\
&\ldots
\end{alignat*}
Die $k_1,k_3,k_5\ldots$ nennt man dann auch \emph{$\tau$-gespielt}. Falls Spieler II nach der Strategie $\tau$ spielt, ist die Spielfolge bereits festgelegt durch $\tau$ und $q:=(k_0,k_2,k_4,\ldots)$ und wird mit $q * \tau$ bezeichnet.

\bigskip $\tau$ ist eine \emph{Gewinnstrategie f"ur II}, falls 
\begin{alignat*}{1}
&\FA[q\in\omega^{<\omega}](q * \tau\in \C[A]).
\end{alignat*}
\end{definition}

\section{Banach-Mazur-Spiele}\label{Kap_Banach_Mazur_Spiele}

Die Banach-Mazur-Spiele sind spezielle unendlichen Zwei-Personenspiele mit perfekter Information. Da die Spieler in ihnen topologische Objekte (z.B. offene Basismengen) spielen, z"ahlen sie zu den topologischen Spielen. Die Grundvariante der Banach-Mazur-Spiele sieht wie folgt aus
\footnote{\emph{Banach-Mazur-Spiel:} Der polnische Mathematiker Stanis\l av Mazur (1905--1981) entwickelte 1935 folgendes Spiel: Zu einer beliebigen vorgegebenen Teilmenge $A$ des Einheitsintervalls $[0,1]\subset\RZ$ spielen Spieler I und Spieler II abwechselnd Intervalle $I_n$ und $J_n\subset [0,1]$ wobei $I_n\supset J_n\supset I_{n+1}$. Man sagt, da"s \emph{Spieler I gewinnt}, falls gilt:
\begin{alignat*}{1}
\bigcap_{n<\omega}I_n\cap A\neq\LM &\TT{ oder gleichbedeutend: }\bigcap_{n<\omega}J_n\cap A\neq\LM
\end{alignat*}
Ansonsten sagt man, da"s \emph{Spieler II gewinnt}. Die unter \ref{Def_Banach_Mazur_Spiel} gegebene Definition ist eine Verallgemeinerung dieses Spiels f"ur beliebige topologische R"aume von John C. Oxtoby (1910--1991).}:

\begin{definition}\header{Das Spiel $\textup{BM}(X,A)$}\label{Def_Banach_Mazur_Spiel}\\
Sei $X$ ein nicht-leerer topologischer Raum. Zu einer beliebigen Teilmenge $A\subset X$ definiert man das \emph{Banach-Mazur-Spiel $\textup{BM}(X,A)$} wie folgt:
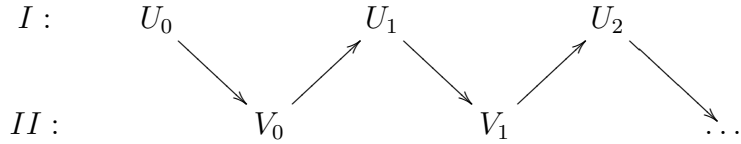
\[
\xymatrix{
I: &U_0 \ar[dr] & 					 	 &U_1 \ar[dr]  &						&U_2 \ar[dr]  &\\
II:&		 				&V_0 \ar[ur]	 &						 &V_1 \ar[ur]	&							&\ldots
}
\]
Spieler I und Spieler II spielen abwechselnd offene nicht-leere Mengen $U_n$ und $V_n\subset X$ wobei $U_n\supset V_n\supset U_{n+1}$. 
Man sagt, da"s \emph{Spieler I gewinnt}, falls gilt:
\begin{alignat*}{1}
\bigcap_{n<\omega}U_n\cap A\neq\LM &\TT{ oder gleichbedeutend: }\bigcap_{n<\omega}V_n\cap A\neq\LM
\end{alignat*}
Ansonsten sagt man, da"s \emph{Spieler II gewinnt}.

\bigskip Falls klar ist, in welchem Raum $X$ gespielt wird, bezeichnet man das Spiel auch einfach mit \emph{$\textup{BM}(A)$}.
\end{definition}

Die topologischen Spiele $\BMSo{A}$ und $\BMSa{B}$ aus Kapitel \ref{Kap_BaireKategorie} sowie $\BMSb{A}$ und $\BMSc{B}$ aus Kapitel \ref{Kap_SigmaKategorie} werden weiter unten f"ur Teilmengen $A$ des Baire-Raumes bzw. Teilmengen $B$ von ${\BR\times\BR}$ (im Falle von $\BMSa{B}$) und Teilmengen $B$ von ${\BR\times\lambda^\omega}$ f"ur eine unendliche Ordinalzahl $\lambda$ (im Falle von $\BMSa{B}$) mit $A=\p(B)$ definiert. Sie lassen sich auch als Varianten dieses Grundtypes formulieren. Beispielsweise entsprechen in dem Spiel $\BMSb{A}$ [\siehe Definietion \ref{DefBMSOhneZeugen}] die von I gespielten endlichen Sequenzen $u_n\in\omega^{<\omega}$ offenen Basismengen $U_n:=O_{u_0\kon\ldots\kon u_n}\subset\omega^\omega$ und die von II gespielten nat"urlichen Zahlen $k_n$ entsprechen den offenen Mengen $V_n:=\bigcup_{u(0)> k_n}O_{u_0\kon\ldots\kon u_{n-1}\kon u}\subset\omega^\omega$.

\section{Determiniertheit und Auswahl}\label{Kap_Determ_Auswahl}

Die Determiniertheit von Spielen in Abh"angigkeit von ihrer Gewinnmenge ist ein entscheidendes Thema in der Spieltheorie. Dabei beziehen wir uns mit dem Begriff \eemph{Determiniertheit} stets auf die in Kapitel \ref{Kap_UnendlZweiPersSpiel} defnierten Spiele $\textup{G}(A)$ f"ur eine Teilmenge $A\subset \omega^\omega$. Reden wir "uber die Determiniertheit von anderen Spielen, so gehen wir davon aus, da"s diese auf den Grundtyp $\textup{G}(A)$ zur"uckgef"uhrt werden k"onnen (bei einem Spiel, in dem endliche Sequenzen nat"urlicher Zahlen gespielt werden: etwa durch die Kodierung dieser Sequenzen in nat"urliche Zahlen).

\begin{definition}\header{$\textup{G}(A)$ determiniert}\label{Def_determiniert_G_A}\\
Sei $A\subset X^\omega$. Das Spiel $\textup{G}(A)$ hei"st \emph{determiniert}, falls einer der beiden Spieler eine Gewinnstrategie hat.
\end{definition}

\begin{definition}\header{Determiniertheitsaxiom $\AD$}\\
Das \emph{Determiniertheitsaxiom} $\AD$ besagt, da"s f"ur jedes $A\subset \omega^\omega$ das zugeh"orige Spiel $\textup{G}(A)$ determiniert ist\footnote{\emph{Determiniertheitsaxiom:} 1962 von Jan Mycielski und Hugo Steinhaus in ihrem Aufsatz \eemph{\glqq A mathematical axiom contradicting the axiom of choice\grqq} (\siehe \cite{Mycielski:1962}) formuliert als \glqq Axiom of Determinateness\grqq: \glqq Every set of reals is determined\grqq\ (\siehe auch \cite{Steinhaus:1965}). Unter der Annahme von $\AD$ l"a"st sich zeigen, da"s jede Menge reeller Zahlen \eemph{lebesgue-me"sbar} ist und die \eemph{Baire-Eigenschaft} sowie die \eemph{Perfekte-Mengen-Eigenschaft} besitzt (\siehe etwa \cite[33.]{Jech:2003}). Diese Resultate aus dem Jahr 1964 waren urspr"unglich der Hauptanreiz f"ur eine Besch"aftigung mit $\AD$. Nachdem Robert M. Solovay gezeigt hatte, da"s $\aleph_1$ unter der Annahme von $\AD$ eine me"sbare Kardinalzahl ist, richtete sich dann die Aufmerksamkeit auf das Verh"altnis zwischen $\AD$ und \eemph{gro"sen Kardinalzahlen} (\siehe etwa \cite[28.]{Kanamori:2003}).}.
\end{definition}

Das Auswahlaxiom wird manchmal auch als \eemph{volles Auswahlaxiom} bezeichnet, um den Unterschied zu schw"acheren Formen der Auswahl deutlich zu machen. Das volle Auswahlaxiom ist wie folgt definiert:

\begin{definition}\header{Auswahlsaxiom $\AC$}\\
Das \emph{Auswahlsaxiom} $\AC$ besagt, da"s f"ur jede Familie $S$ von Mengen mit $\LM\not\in S$ eine Auswahlfunktion 
\begin{alignat*}{1}
&f:\PFEIL{S}{}{\bigcup S}
\end{alignat*}
mit $\FA[X\in S](f(X)\in X)$ existiert\footnote{\emph{Auswahlaxiom:} 1904 durch Ernst Zermelo (1871--1953) in \cite{Zermelo:1904} formuliert, um zu zeigen, da"s jede Menge wohlgeordnet werden kann.}.
\end{definition}

Das Determiniertheitsaxiom $\AD$ und das volle Auswahlaxiom $\AC$ sind nicht miteinander vertr"aglich:

\begin{lemma}\header{Inkonsistenz von $\AD$ mit $\AC$}\label{Lemma_Inkonsistenz_AD_AC}\\
Wenn $\AC$ gilt, dann gibt es ein $A\subset \omega^\omega$ so, da"s das zugeh"orige Spiel $\textup{G}(A)$ nicht determiniert ist.
\end{lemma}
\begin{beweis}\KOM[Jech 628]
Seien $\MNG{\sigma_\alpha}{\alpha<2^{\aleph_0}}$ und $\MNG{\tau_\alpha}{\alpha<2^{\aleph_0}}$ die Mengen der Strategien f"ur Spieler I bzw. II (in Spielen $\textup{G}(A)$ mit $A\subset\BR$ beliebig). Dann konstruiert man die Mengen
\begin{alignat*}{1}
&X:=\MNG{x_\alpha}{\alpha<2^{\aleph_0}}\\
&Y:=\MNG{y_\alpha}{\alpha<2^{\aleph_0}}
\end{alignat*}
wie folgt:

F"ur gegebene $\MNG{x_\xi}{\xi<\alpha}$ und $\MNG{y_\xi}{\xi<\alpha}$ w"ahle man mit Hilfe von $\AC$ ein $y_\alpha:=\sigma_\alpha*b\not\in\MNG{x_\xi}{\xi<\alpha}$ f"ur ein $b\in\BR$ (solch ein $y_\alpha$ existiert, da die Kardinalit"at von $\MNG{\sigma_\alpha*b}{b\in\BR}$ gleich $2^{\aleph_0}>\alpha$ ist). Danach w"ahlt man ein $x_\alpha:=a*\tau_\alpha\not\in\MNG{y_\xi}{\xi\leq\alpha}$ (dies existiert aus demselben Grund, wie das $y_\alpha$). 

\point{1) Es gilt $X\cap Y=\LM$}
Angenommen es gibt ein $z\in X\cap Y$\gs etwa
\begin{alignat*}{1}
z	&=y_\alpha=\sigma_\alpha*b\not\in\MNG{x_\xi}{\xi<\alpha}\TEXT{ und}\\
z	&=x_{\alpha'}=a*\tau_{\alpha'}\not\in\MNG{y_\xi}{\xi\leq{\alpha'}}.
\end{alignat*}
Wegen $z=x_{\alpha'}\not\in\MNG{x_\xi}{\xi<\alpha}$ gilt dann $\alpha'\geq \alpha$ und wegen $z=y_\alpha\not\in\MNG{y_\xi}{\xi\leq{\alpha'}}$ gilt $\alpha'<\alpha$\gs Widerspruch.

\point{2) Es gilt $\FA[\alpha<2^{\aleph_0}]\EX[(a,b)\in\BR\times\BR](a*\tau_\alpha\in X\UND\sigma_\alpha*b\not\in X)$}
Nach der Konstruktion von $X$ und $Y$ gibt es f"ur jedes $\alpha<2^{\aleph_0}$ ein $a*\tau_\alpha\in X$ und ein $\sigma_\alpha*b\in Y$. Da $X\cap Y=\LM$ mu"s dann gelten $\sigma_\alpha*b\not\in X$.

\bigskip2) bedeutet, da"s weder Spieler I noch II eine Gewinnstrategie f"ur das Spiel $\textup{G}(X)$ haben. Das Spiel $\textup{G}(X)$ ist also nicht determiniert.
\end{beweis}

Eine schw"achere Form der Auswahl ist die abz"ahlbare Auswahl:

\begin{definition}\header{Axiom der abz"ahlbaren Auswahl ${\AC}^\omega$}\\
Das \emph{Axiom der abz"ahlbaren Auswahl} ${\AC}^\omega$ besagt, da"s f"ur jede abz"ahlbare Familie $S$ von Mengen mit $\LM\not\in S$ eine Auswahlfunktion 
\begin{alignat*}{1}
&f:\PFEIL{S}{}{\bigcup S}
\end{alignat*}
mit $\FA[X\in S](f(X)\in X)$ existiert.

Sei $M$ eine Menge. Das \emph{Axiom der abz"ahlbaren Auswahl auf $M$} ${\AC}^\omega_M$ besagt, da"s f"ur jede abz"ahlbare Familie $S$ von Teilmengen in $M$ mit $\LM\not\in S$ eine Auswahlfunktion 
\begin{alignat*}{1}
&f:\PFEIL{S}{}{\bigcup S}
\end{alignat*}
mit $\FA[X\in S](f(X)\in X)$ existiert.
\end{definition}

Das Axiom der abz"ahlbaren Auswahl ${\AC}^\omega_{\BR}$ ist vertr"aglich mit dem Determiniertheitsaxiom $\AD$\gs es folgt sogar daraus:

\begin{lemma}\header{Konsistenz von $\AD$ mit ${\AC}^\omega_{\BR}$}\label{Lemma_Konsistenz_AD_AC_omega}\\
Aus $\AD$ folgt ${\AC}^\omega_{\BR}$.
\end{lemma}
\begin{beweis}\KOM[Jech 628]
Sei $S=\MNG{X_n}{n<\omega}$ eine Familie nicht-leerer Teilmengen in $\BR$. 

\point{Es gilt $\EX[f:\PFEIL{S}{}{\bigcup S}](f(X_n)\in X_n)$}

Um dies zu zeigen, wird zun"achst das Spiel $\textup{G}_S$ wie folgt definiert: 
\begin{alignat*}{1}
&\TT{Spieler I spielt }(a_0,a_1,a_2,\ldots)\in\BR,\\
&\TT{Spieler II spielt }(b_0,b_1,b_2,\ldots)=:b\in\BR.
\end{alignat*}
Spieler II gewinne genau dann, wenn $b\in X_{a_0}$ gilt. (Mit $A:=\MNG{(a_0,b_0,a_1,b_1,\ldots)}{(b_0,b_1,b_2,\ldots\in X_{a_0})}$ l"a"st sich $\textup{G}_S$ als Spiel $\textup{G}(A)$ der Grundvariante aus Definition \ref{Def_SpielGrundvariante} auffassen und es ist klar, da"s aus $\AD$ folgt, da"s $\textup{G}_S$ determiniert ist.)

Spieler I kann f"ur $\textup{G}_S$ keine Gewinnstrategie haben, da Spieler II zu $a_0$ nur ein $b\in X_{a_0}\neq\LM$ zu spielen braucht, um zu gewinnen (Spieler I hat nach seinem ersten Zug $a_0$ keinerlei Einflu"s mehr auf den Ausgang des Spieles). Wegen $\AD$ hat somit Spieler II eine Gewinnstrategie $\tau$. Mit dieser l"a"st sich dann eine Auswahlfunktion 
\begin{alignat*}{1}
&f(X_n):=\left\langle (n,0,0,\ldots)*\tau \right\rangle_{\TEXT{\tiny{II}}}
\end{alignat*}
f"ur $S$ definieren, wobei $\left\langle (a_0,a_1,a_2,\ldots)*\tau \right\rangle_{\TEXT{\tiny{II}}}:=b$ mit $b\in X_{a_0}$ die Z"uge von Spieler II in $\textup{G}_S$ seien. Es gilt $f(X_n)\in X_n$ da $\tau$ eine Gewinnstrategie von II f"ur $\textup{G}_S$ ist.
\end{beweis}

Um nicht in Widerspr"uche zu geraten, ist also darauf zu achten, $\AD$ und $\AC$ niemals gleichzeitig vorauszusetzen. $\AC[\omega]$ ist jedoch f"ur die Untersuchung von Spielen des Typs $\textup{G}(A)$ f"ur eine Teilmenge $A\subset \omega^\omega$ in der Regel schon ausreichend. Eine andere schw"achere Form der Auswahl ist das Axiom der \eemph{abh"angigen Auswahl} $\DC$. Die abh"angige Auswahl folgt aus $\AC[\omega]$ und ist ebenfalls vertr"aglich mit $\AD$ (\siehe vgl. Kapitel \ref{Kap_SigmaKategorie} Fu"snote \ref{FN_GaleMartin}). An den Stellen, an denen wir $\AD$ annehmen, soll dies jedesmal ausdr"ucklich erw"ahnt werden.

\clearemptydoublepage

%==============================KAPITEL=============================
%\chapter{Magere und komagere Mengen}\KOM\KOM[einiges\\erg"anzen\\Diss. 1.Kap.]
\chapter{Baire-Kategorie}\label{Kap_BaireKategorie}
\setcounter{lemma}{0}

In diesem Kapitel soll nun die Baire-Kategorie vorgestellt werden. Die Baire-Kategorie l"a"st sich allgemein f"ur topologische R"aume definieren. Die \eemph{kleinen} Mengen der Baire-Kategorie sind die \eemph{mageren} Mengen. Sie setzen sich aus \eemph{nirgends dichten} Mengen zusammen, die sich gleichsam \glqq d"unn machen\grqq, indem ihre Komplemente offene dichte Mengen enthalten [\siehe Kapitel \ref{Kap_Magere_Mengen}]. Die nirgends dichten Mengen sind anschaulich gesehen Mengen mit \glqq lauter L"ochern\grqq, da jede offene Menge das Komplement von $A$ in einer offenen Menge schneidet. Eine magere Menge $A$ ist in dem Sinne nicht notwendig selber eine Menge mit \glqq lauter L"ochern\grqq, kann aber durch derartige Mengen approximiert werden. 
Offene Intervalle $(a,b)\neq\LM$ in $\RZ$ sind nie mager [\siehe Beispiel \ref{Bsp_mager_lauterLuecken}]. Diese Eigenschaft motiviert die Definition \eemph{Baire'scher R"aume} [\siehe Definition \ref{Kapitel_BairescheRaeume}]: Ein Raum hei"st \eemph{Baire'sch}, falls er keine nicht-leeren mageren offenen Mengen enth"alt. 
\begin{comment}
Eine magere Menge in $\RZ$ oder allgemeiner in einem \eemph{Baire'schen Raum} [\siehe Definition \ref{Kapitel_BairescheRaeume}] ist anschaulich eine Menge mit einer dichten Menge von \glqq L"ucken\grqq. Diese Anschauung bezieht sich darauf, da"s sich kein Intervall $I$ in $\RZ$ als Vereinigung $M$ einer Folge von Mengen von erster Kategorie darstellen l"a"st [\siehe Beispiel \ref{Bsp_mager_lauterLuecken}]. Da in $\RZ$ der Durchschnitt abz"ahlbar vieler offener dichter Mengen nicht leer ist, folgt daraus: Falls $I\supset M$ dann $\FA[U\neq\LM\TT{ offen }][U\subset I\Rightarrow\EX[x\in U](x\not\in M)]$.\KOM[Oxt S. 4]%Magere Mengen sind durch diese Eigenschaft nat"urlich nicht vollst"andig charakterisiert.
\end{comment}

Beispiel \ref{Bsp_RzAlsNullmengePlusMagereMenge}\footnote{Bis auf Beispiel \ref{Bsp_RzAlsNullmengePlusMagereMenge} (aus \cite[11.8]{Jech:2003}), \ref{Bsp_RatZahlenUndIrratZahlen}, \ref{BspMagerErstTrennungsAx_etc} und \ref{Bsp2KatMetrRaumEndl_etc} stammen die Beispiele in Kapitel \ref{Kap_Magere_Mengen} aus \cite[II 6, VI 2]{Koehnen:1988}} verdeutlicht, da"s man die Kleinheitsbegriffe verschiedener Kategorien-Konzepte keinesfalls durcheinanderw"urfeln darf\gs \eemph{klein} in dem einen Konzept bedeutet nicht unbedingt auch \eemph{klein} in einem anderen. 

Die \eemph{komageren} Mengen stehen in Baire's Kategorien-Konzept f"ur die \eemph{gro"sen} Mengen und werden daher als die Komplemente \eemph{magerer} Mengen definiert [\siehe Kapitel \ref{Kap_Komagere_Mengen}]. Daher haben die komageren Mengen die zu den mageren Mengen dualen Eigenschaften. Au"serdem wird definiert, wann eine Menge \eemph{mager} bzw. \eemph{komager in einer offenen Menge} ist.

Eine magere Menge in $\RZ$ enth"alt anschaulich eine dichte Menge von \glqq L"ucken\grqq\gs damit ist gemeint, da"s keine Intervall $[a,b]\subset\RZ$ (oder "aquivalent dazu: kein offenes Intervall $(a,b)\subset\RZ$) mit $a<b$ als Vereinigung einer Folge von Mengen von erster Kategorie dargestellt werden kann [\siehe Beispiel \ref{Bsp_mager_lauterLuecken}]. Um diese Eigenschaft f"ur topologische R"aume zu verallgemeinern, bezeichnet man einen Raum als \eemph{Baire'sch}, wenn er keine nicht-leere offene magere Menge enth"alt [\siehe Kapitel \ref{Kapitel_BairescheRaeume}].

Eine nicht-leere offene Teilmenge $U\subset X$ ist mager genau dann, wenn jede 
Teilmenge von $X$ mager (und komager) in $U$ ist. Die Baire'schen R"aume sind daher gerade diejenigen, in denen die \eemph{Lokalisierung} der Begriffe \eemph{mager} und \eemph{komager} sinnvoll ist.

 Der Baire'sche Kategoriensatz \ref{BKS}\footnote{Der Bairesche Kategoriensatz \ref{BKS} findet sich in dieser Form etwa in \cite[II.2]{Haworth:1977}, Teil $(ii)$ findet sich auch in \cite[7.3.12]{Preuss:1970}.} identifiziert zwei gro"se Klassen topologischer R"aume als Baire'sch: die lokal kompakten R"aume und die vollst"andig metrisierbaren R"aume.

F"ur spieltheoretische Charakterisierungen betrachten wir in dieser Arbeit die Baire-Kategorie auf dem \eemph{Baireraum}. In Kapitel \ref{Kap_Baire_Raum} ordnen wir den Baireraum in den allgemeineren topologischen Zusammenhang\footnote{Bei Begriffen wie \eemph{quasikompakt}, \eemph{kompakt} und \eemph{lokal (quasi-) kompakt} halten wir uns an die Definitionen von \cite[7.1, 7.3]{Preuss:1970}. \siehe etwa \cite[4]{Preuss:1970} f"ur die \eemph{Trennungsaxiome} und den Begriff \eemph{regul"ar}.} ein. Unter anderem zeigen wir, da"s der Baireraum ein \eemph{polnischer} Raum ist und somit nach dem Baire'schen Kategoriensatz \ref{BKS} $(i)$ insbesondere ein Baire'scher Raum. 

In Kapitel \ref{KapSpielCharMager} wird dann die Baire-Kategorie mittels des Spieles $\BMSo{A}$ f"ur beliebige Teilmengen $A$ des Baireraumes und mittels $\BMSa{B}$ f"ur Projektionen $A=\p(B)$ ($B\subset\omega^\omega\times\omega^\omega$) charakterisiert [\siehe Theoreme \ref{SatzMagerCharSpiele} und \ref{SatzMagerCharSpiele_2}]\footnote{\siehe etwa \cite[III 33]{Jech:2003} f"ur Theorem \ref{SatzMagerCharSpiele}.  \siehe \cite[S. 194]{Kechris:1977} f"ur Theorem \ref{SatzMagerCharSpiele_2}.} . Es zeigt sich, da"s eine Teilmenge $A$ des Baireraumes genau dann \eemph{klein} im Sinne der Baire-Kategorie ist, wenn Spieler II in dem Spiel $\BMSo{A}$ eine Gewinnstrategie besitzt und \eemph{lokal gro"s} (\eemph{gro"s} in einer offenen Teilmenge des Baireraumes) genau dann, wenn Spieler I eine Gewinnstartegie f"ur das Spiel $\BMSo{A}$ hat [\siehe Theorem \ref{SatzMagerCharSpiele}]. F"ur Projektionen $A=\p(B)$ einer Menge $B\subset\omega^\omega\times\omega^\omega$ erh"alt man ein entsprechendes Resultat mit dem Spiel $\BMSa{B}$ [\siehe Theorem \ref{SatzMagerCharSpiele_2}].

Kapitel \ref{Kap_BaireKategorie} abschlie"send wird als eine direkte Anwendung von Theorem \ref{SatzMagerCharSpiele_2} ein hinreichendes Kriterium f"ur die Baire-Eigenschaft projektiver Mengen angegeben [\siehe Satz \ref{Satz_CharBaireEig}]\footnote{\siehe \cite[S. 193]{Kechris:1977} f"ur Satz \ref{Satz_CharBaireEig}.}.

\section{Magere Mengen}\label{Kap_Magere_Mengen}
Die \eemph{kleinen} Mengen der Baire-Kategorie sind die \eemph{mageren} Mengen. Sie setzen sich aus \eemph{nirgends dichten} Mengen zusammen, die sich gleichsam \glqq d"unn machen\grqq. Im Gegensatz dazu stehen die \eemph{"uberall dichten} Mengen, die sich gleichsam \glqq breit machen\grqq, indem sie jede offene Menge schneiden:

\begin{definition}\header{"uberall dicht}\\
Eine Menge $A \subset X$ hei"st \emph{"uberall dicht} oder einfach \emph{dicht} genau dann, wenn
sie mit jeder offenen nicht-leeren Menge $U\subset X$ einen nicht-leeren Schnitt hat. 
\end{definition}

\begin{definition}\header{dicht in einer offenen Menge $U$}\label{DefDichtInOffMng}\\
Eine Menge $A \subset X$ hei"st \emph{dicht in $U$} f"ur eine offene Menge $U\subset X$ genau dann, wenn sie in $U$ dicht bez"uglich der Relativtopologie ist.
\end{definition}

Beispielsweise liegen die rationalen Zahlen dicht in den reellen Zahlen: Jede reelle Zahl l"a"st sich als unendliche Dezimalzahl darstellen. Diese l"a"st sich als unendliche Reihe im Sinne der Analysis auffassen. Die Partialsummen dieser Reihe sind rationale Zahlen. Also l"a"st sich jede reelle Zahl beliebig genau durch rationale Zahlen approximieren.

\begin{comment}
\begin{beispiel}
Die Menge der auf einem abgeschlossenen Intervall $[a,b]\subset\RZ$ definierten st"uckweise linearen Funktionen
liegt dicht in der Menge $\SF$$[a,b]$ der auf $[a,b]$ definierten stetigen Funktionen.
\end{beispiel}
\begin{beweis}
K"ohnen 136\KOM[(FEHLT)]
\end{beweis}

\begin{beispiel}
Die Menge $M = \MNG{x\in\QZ}{x=\frac{m}{2^n}, m\in\ZZ, n\in\NZ}$ liegt dicht in $\RZ$.
\end{beispiel}
\begin{beweis}
\KOM[K"ohnen 133]\KOM[(FEHLT)]
\end{beweis}

\begin{beispiel}
Im euklidischen metrischen Raum $\RZ^2$ liegt die Menge aller Elemente in $\RZ^2$ mit rationalen Koordinaten dicht.
\end{beispiel}
\begin{beweis}
\KOM[K"ohnen 134]\KOM[(FEHLT)]
\end{beweis}
\end{comment}

\bigskip Im Gegensatz dazu machen sich die nirgends dichten Mengen gleichsam \glqq d"unn\grqq\ dadurch, da"s ihr Komplement eine offene dichte Menge enth"alt:

\begin{definition}\header{nirgends dicht}\label{Def_Nirg_Dicht}\\
Eine Menge $A\subset X$ hei"st \emph{nirgends dicht} genau dann, wenn $\C[A]$ eine offene dichte Menge $G\subset X$ enth"alt.
\end{definition}

Anschaulich erh"alt man eine \eemph{nirgends dichte} Teilmenge $\tilde{Y}$ eines Raumes $X$, indem man aus einer Teilmenge $Y$ von $X$ iterativ aus jeder $X$-offenen Menge $U\subset Y$ eine offene Teilmenge $V$ aus $U$ entfernt. Eine nirgends dichte Menge stellt man sich daher am besten als eine Menge \glqq mit lauter L"ochern\grqq\ vor. 

\bigskip Aus Definition \ref{Def_Nirg_Dicht} folgt direkt, da"s Teilmengen nirgends dichter Mengen ebenfalls nirgends dicht sind. Die leere Menge $\LM$ ist nirgends dicht, da $\C[\LM] = X$ offen und dicht ist. Das Innere $\INT[A]$ einer nirgends dichten Menge $A$ ist leer, da $A$ keine offene Menge enthalten kann (eine offene Menge in $A$ k"onnte ja nicht die dichte Menge $G \subset \C[A]$ schneiden). Insbesondere ist $A$ also nicht offen. Au"serdem kann ein nirgends dichtes $A$ nicht dicht sein in einer Menge $B \subset X$, die die dichte Menge $G \subset \C[A]$ schneidet, da $A$ sonst $G \cap B\subset \C[A]$ schneiden m"usste\gs insbesondere kann $A$ somit nicht dicht sein in $X$ oder einer offenen Teilmenge $U$ von $X$. Eine abgeschlossene Menge $A\subset X$, die das Komplement einer dichten Menge ist, mu"s nirgends dicht sein, da ihr Komplement offen und dicht ist. Es gilt also: Ist $\INT[{\ABS[A]}]=\LM$, so ist $A$ (als Teilmenge der nirgends dichten Menge $\ABS[A]$) nirgends dicht. Insbesondere ist $A$ nirgends dicht, falls abgeschlossen und $\INT[A]=\LM$.

\begin{satz}\header{Ideal der nirgends dichten Mengen}\\
Die nirgends dichten Teilmengen eines topologischen Raumes bilden ein Ideal. 
\end{satz}
\begin{beweis}
Es bleibt zu zeigen, da"s die Vereinigung endlich vieler nirgends dichter Mengen nirgends dicht ist.
Sei also $\bigcup_{i=1}^{n}{A_i}$ die Vereinigung endlich vieler nirgends dichter Teilmengen $A_i$ eines 
topologischen Raumes $X$. Dann liegt im Komplement dieser Vereinigung der endliche Schnitt
$\bigcap_{i=1}^{n}{G_i}$ offener dichter Teilmengen $G_i\subset X$. Endliche Durchschnitte offener 
dichter Teilmengen sind aber stets selbst wieder offen und dicht.
\end{beweis}

\begin{satz}\header{nirgends dicht}\label{satz_nirgdicht}\\
F"ur eine Menge $A\subset X$ sind "aquivalent:
\begin{ITEMS}
\item $A$ nirgends dicht.
\item $\C[{\ABS[A]}]$ dicht.
\item $A$ ist in keiner offenen nicht-leeren Menge $U\subset X$ dicht.
\item F"ur jede offene nicht-leere Menge $U\subset X$ gibt es eine offene nicht-leere Menge $V\subset X$ so, da"s $V\subset U$ und $V\cap A = \LM $.
\item $\INT[{\ABS[A]}]=\LM$.
\end{ITEMS}
\end{satz}

\begin{beweis}\KOM[Diss. 6]
$(i)\Rightarrow(ii):$ Sei $A$ nirgends dicht. Dann ist $\INT[{\ABS[A]}]=\LM$. $\ABS[A]$ enth"alt also keine offene nicht-leere Menge. Folglich gilt f"ur jede offene nicht-leere Menge $U\subset X$: $U\cap \C[{\ABS[A]}]\neq\LM$. Das hei"st, $\C[{\ABS[A]}]$ ist dicht.
\begin{comment}
$(i)\Rightarrow(ii):$ Sei $A$ nirgends dicht\gs etwa $G\subset \C[A]$ mir $G$ offen und dicht in $X$. Angenommen $A$ ist dicht in einem $U\subset X$ offen. Da $G$ dicht ist, ist dann $U\cap G\neq\LM$. Da $G$ und $U$ offen sind, ist $U\cap G$ offen. Und da $A$ dicht in $U$ und $U\cap G$ offen ist mu"s $A\cap U\cap G\neq\LM$ gelten. Widerspruch, da $G\subset\C[A]$.

\medskip Die Annahme war also falsch und es gilt: $A$ ist in keiner offenen Menge $U\subset X$ dicht.
\end{comment}

$(ii)\Rightarrow(iii):$ Sei $\C[{\ABS[A]}]$ dicht. Angenommen $A$ ist dicht in einer offenen nicht-leeren Menge $U\subset X$. Da $\C[{\ABS[A]}]$ dicht ist, ist dann $U\cap \C[{\ABS[A]}]\neq\LM$. Da $\C[{\ABS[A]}]$ und $U$ offen sind, ist $U\cap \C[{\ABS[A]}]$ offen. Und da $A$ dicht in $U$ und $U\cap \C[{\ABS[A]}]$ offen ist, mu"s $A\cap U\cap \C[{\ABS[A]}]\neq\LM$ gelten. Widerspruch, da $\C[{\ABS[A]}]\subset\C[A]$.

\medskip Die Annahme war also falsch und es gilt: $A$ ist in keiner offenen nicht-leeren Menge $U\subset X$ dicht.

\bigskip $(iii)\Rightarrow(iv)$ Sei $U$ eine offene nicht-leere Teilmenge von $X$ und gelte $(iii)$. Angenommen alle offenen nicht-leeren Mengen $V\subset U$ haben einen nicht-leeren Schnitt mit $A$. Das hei"st $A$ liegt dicht in der offenen nicht-leeren Menge $U$\gs im Widerspruch zu $(iii)$. 

\medskip Die Annahme war also falsch, d.h. es gibt eine offene nicht-leere Menge $V\subset U$ mit $V\cap A=\LM$.

\bigskip $(iv)\Rightarrow(i)$ Gelte $(iv)$. Angenommen $A$ ist nicht nirgends dicht. Dann gilt $\INT[{\ABS[A]}]\neq\LM$. Da $\INT[{\ABS[A]}]\subset A$ gilt dann 
\begin{alignat*}{1}
&\FA[{V\subset {\INT[{\ABS[A]}]} \TT{ offen } \neq\LM}](V\cap A\neq\LM)
\end{alignat*}
und da $\INT[{\ABS[A]}]$ offen ist, widerspricht dies $(iv)$.

\medskip Die Annahme war also falsch und $A$ ist somit nirgends dicht.

\bigskip $(v)\Leftrightarrow(ii)$
In ${\ABS[A]}$ ist keine offene Menge enthalten genau dann, wenn sich jede offene Menge mit $\C[{\ABS[A]}]$ schneidet.
\end{beweis}

\begin{comment}
\begin{beispiel}
Die Cantor-Menge $\CM$ liegt nirgends dicht in $\RZ$. 
\end{beispiel}
\begin{beweis}K"ohnen 145\\Elstrodt 70f\KOM[(FEHLT)]

\end{beweis}
\end{comment}

\begin{beispiel}\label{Bsp_GraphStetNirgDicht}
Der Graph $\textup{G}_f:=\MNG{(x,y)\in\RZ^2}{f(x)=y}$ einer stetigen, reellwertigen Funktion $f$ ist eine im euklidischen Raum\footnote{\emph{euklidischer Raum:} \siehe Konventionen S. \pageref{Def_Eukl_Raum}.} $(\RZ^2,\textup{d}^2)$ nirgends dichte Menge.
\end{beispiel}
\begin{beweis}\KOM[K"ohnen 146]
Nach Satz \ref{satz_nirgdicht} gen"ugt es zu zeigen, da"s $\ABS[{\textup{G}_f}]$ keine offene Menge enth"ahlt.

\medskip Um zu zeigen, da"s $\textup{G}_f$ gleich $\ABS[{\textup{G}_f}]$, gen"ugt es zu zeigen, da"s $\textup{G}_f\supset\ABS[{\textup{G}_f}]$. Sei $a=(x,y)\in\ABS[{\textup{G}_f}]$ beliebig. Da $\ABS[{\textup{G}_f}]$ abgeschlossen ist gibt es dann eine Folge $(a_i)_i=(x_i,y_i)_i$ in ${\textup{G}_f}$, die gegen $a$ konvergiert und es gilt $\lim_i(x_i)=x$ und  $\lim_i(y_i)=y$. Da $f$ stetig ist, gilt dann
\[
f(x)=\lim_i f(x_i)=\lim_i y_i=y
\]
und somit mu"s $a=(x,y)$ in $\textup{G}_f$ liegen.

\medskip Nun gen"ugt es zu zeigen, da"s $\textup{G}_f$ keine offene $\epsilon$-Kugel $\Kugel(a,\epsilon):=\MNG{b\in\RZ}{\textup{d}^2(a,b)<\epsilon}$ f"ur $\epsilon>0$ enth"alt. Dies ist aber klar, da der Graph $\textup{G}_f$ von $f$ mit einer $\epsilon$-Kugel $\Kugel(a=(x,y),\epsilon)$ auch die Menge 
\[
L:=\MNG{(x',y')\in\RZ^2}{x'=x \UND y-\epsilon<y'<y+\epsilon}\subset \Kugel(a,\epsilon)
\]
enthalten m"u"ste (was nicht geht, da $\epsilon>0$ und $f$ eine Funktion ist).

\medskip Insgesamt ist damit gezeigt, da"s $\INT[{\ABS[{\textup{G}_f}]}]=\LM$. Somit ist $\textup{G}_f$ nach Satz \ref{satz_nirgdicht} nirgends dicht.
\end{beweis}

Eine abz"ahlbare Vereinigung nirgends dichter Mengen kann dicht in einer offenen Menge sein und somit nicht wieder nirgends dicht (Beispiel: Die Menge der rationalen Zahlen $\QZ=\bigcup_{q\in\QZ}{\{q\}}$ ist eine abz"ahlbare Vereinigung nirgends dichter Teilmengen $\{q\}\subset\RZ$ aber dicht in $\RZ$. Nirgends dichte Mengen bilden also i.a. kein $\sigma$-Ideal\gs um ein solches zu erhalten definiert man:

\begin{definition}\header{mager}\\
Eine Menge $A \subset X$ hei"st \emph{mager in $X$} oder einfach \emph{mager} 
genau dann, wenn $A$ als abz"ahlbare Vereinigung $\bigcup_{i \geq 0}{A_i}$ von nirgends 
dichten Mengen $A_i$ geschrieben werden kann. Magere Mengen bezeichnet man auch als
 \emph{von erster Kategorie} und nicht magere Mengen als \emph{von zweiter Kategorie}
\footnote{\emph{mager:} Bei R. L. Baire hei"st es in \cite{Baire:1899} \emph{ensemble de premi\`ere cat\'egorie} bzw. \emph{ensemble de deuxi\`eme cat\'egorie}, korrekt "ubersetzt also \emph{Mengen von erster Kategorie} bzw. \emph{Mengen von zweiter Kategorie}. "Ublich ist aber auch die Kurzform \emph{1. Kategorie} bzw. \emph{2. Kategorie}. Die franz"osische Mathematikergruppe Bourbaki hat vorgeschlagen, die sprachlich etwas schwerf"alligen und unanschaulichen Ausdrucksweisen \emph{von 1. Kategorie} und \emph{von 2. Kategorie} fallenzulassen, und stattdessen eine Menge von 1. Kategorie \emph{mager} zu nennen. Der Vorschlag setzte sich aber nicht recht durch, und so sind bis heute beide Bezeichnungsweisen gebr"auchlich.}.
\end{definition}

Der Begriff \eemph{mager} l"a"st sich wie folgt lokalisieren:

\begin{definition}\header{mager in einer offenen Menge $U$}\label{Def_Lokalisierung_Mager}\\
Sei $U$ eine offene Menge in $X$, dann hei"st eine Menge $A \subset X$ \emph{mager in $U$} 
genau dann, wenn $A \cap U$ mager ist. Da $U$ offen ist, ist dies "aquivalent dazu, da"s $A \cap U$
mager bez"uglich der Relativtopologie von $U$ ist.
\end{definition}

Nirgends dichte Mengen sind also mager\gs insbesondere die leere Menge. Abz"ahlbare Vereinigungen 
magerer Mengen sind ebenfalls mager, da eine abz"ahlbare Vereinigung abz"ahlbarer Vereinigungen 
von Mengen wieder abz"albar ist. Teilmengen magerer Mengen sind mager\gs sei etwa $B \subset A$ und $A$ 
mager, dann ist $A = \bigcup_{i\geq 0}{A_i}$ mit $A_i$ nirgends dicht und somit 
$B = \bigcup_{i\geq 0}{B\cap A_i}$ mit $B\cap A_i \subset A_i$ nirgends dicht. Mithin gilt:

\begin{satz}\header{$\sigma$-Ideal der mageren Mengen}
Die mageren Teilmengen eines topologischen Raumes bilden ein $\sigma$-Ideal.
\end{satz}

Die mageren Mengen erf"ullen also die nat"urlichen Eigenschaften kleiner Mengen.

\bigskip Obwohl \eemph{mager} wie \eemph{Lebesgue-Nullmenge} jeweils f"ur \glqq vernachl"assigbar klein\grqq\ stehen, l"a"st sich $\RZ$ als Vereinigung einer mageren Menge und einer Nullmenge darstellen. Die Kleinheitsbegriffe sind also wohl zu trennen\gs \glqq klein\grqq\ in dem einen Sinne bedeutet nicht automatisch auch \glqq klein\grqq\ in einem anderen Sinne:

\begin{beispiel}\label{Bsp_RzAlsNullmengePlusMagereMenge}
Es gibt eine Lebesgue-Nullmenge $N\subset \RZ$ und eine magere Menge $M\subset \RZ$ so, da"s $\RZ=N\cup M$ gilt.
\end{beispiel}
\begin{beweis}
Sei $q_1,q_2,\ldots$ eine Aufz"ahlung der rationalen Zahlen. F"ur jedes $n\geq 1$ und $k\geq 1$ sei $I_{n,k}$ das offene Intervall mit Mittelpunkt $q_n$ von der L"ange $\frac{1}{k 2^n}$. Sei $D_k:=\bigcup_{n\geq 1} I_{n,k}$ und $A=\bigcap_{k\geq 1} D_k$. Da die $I_{n,k}$ offen sind ist dann jedes $D_k$ offen. Au"serdem ist jedes $D_k$ auch dicht in $\RZ$, weil $\QZ\subset D_k$ gilt und $\QZ$ dicht in $\RZ$ ist. Somit ist $\C[A]$ mager. $A$ ist eine Nullmenge, da $\lambda(D_k)\leq\frac{1}{k}$ und somit $\lambda(A)=0$ gilt (mit Lebesgue-Ma"s $\lambda$). 
\end{beweis}

Um eine Anschauung zu bekommen, wie \glqq allt"agliche\grqq\ magere Mengen aussehen k"onnen und auf welche Weise sie sich \glqq d"unn machen\grqq, sollen abschlie"send noch einige weitere Beispiele betrachtet werden.

\begin{beispiel}\label{Bsp_AbzbReelleTM_Mager}
Jede h"ochstens abz"ahlbare Teilmenge von $\RZ$ ist mager in $\RZ$.
\end{beispiel}
\begin{beweis}
Dies ergibt sich sofort daraus, da"s einelementige Teilmengen von $\RZ$ nirgends dicht in $\RZ$ sind. Eine einelementige Teilmenge $\{x\}$ von $\RZ$ ist stets nirgends dicht in $\RZ$, weil ihr Komplement $\C[\{x\}]$ offen und dicht in $\RZ$ ist.
\begin{comment}
Die Singulum-Teilmengen von $\RZ$ sind nirgends dicht: Sei $\{x\}$ eine beliebige Singulum-Teilmenge von $\RZ$. Da $\RZ$ metrisierbar ist, gibt es f"ur jede reele Zahl $y\neq x$ eine offene $\epsilon$-Kugel $\Kugel(y,\epsilon)$, die $y$ aber nicht $x$ enth"alt\gs $\C[\{x\}]$ ist also offen. Da $\RZ$ keine isolierten Punkte (Punkte, die eine Umgebung haben, in der keine weiteren Punkte aus $\RZ$ enthalten sind) hat, kann $\{x\}$ nicht offen sein\gs somit schneidet jede offene Teilmenge von $\RZ$ das Komplement von $\{x\}$, d.h. $\C[\{x\}]$ ist dicht. Insgesamt ist $\C[\{x\}]$ also offen und dicht und $\{x\}$ somit nirgends dicht.

\medskip Eine abz"ahlbare Teilmenge von $\RZ$ l"a"st sich dann als abz"ahlbare Vereinigung nirgends dichter Mengen darstellen:
\[
A=\bigcup_{x\in A}\{x\}
\]
und ist somit mager.

\medskip Eine endliche Teilmenge von $\RZ$ ist Teilmenge einer abz"ahlbaren Teilmenge von $\RZ$ und somit ebenfalls mager.

\medskip Also sind h"ochstens abz"albare Teilmengen von $\RZ$ mager.
\end{comment}
\end{beweis}

\begin{beispiel}\label{Bsp_RatZahlenUndIrratZahlen}
Die Menge der rationalen Zahlen $\QZ$ ist mager in $\RZ$ und die Menge der irrationalen Zahlen $\IZ$ ist nicht mager in $\RZ$.
\end{beispiel}
\begin{beweis}
Da $\QZ\subset\RZ$ abz"ahlbar ist, ist $\QZ$ nach Beispiel \ref{Bsp_AbzbReelleTM_Mager} mager. Da $\RZ=\QZ\cup\IZ$ gilt und $\RZ$ nicht mager ist, ist $\IZ$ nicht mager. (In Satz \ref{BKS} $(i)$ zeigen wir, da"s in vollst"andig metrisierbaren R"aumen der Durchschnitt von offenen dichten Teilmengen wieder dicht ist. Da $\RZ$ vollst"andig metrisierbar ist folgt daraus, da"s $\RZ$ nicht mager sein kann, da sonst: $\RZ=\bigcup_{i<\omega}\underset{\TT{\tiny{nirg. dicht}}}{A_i}\Rightarrow \LM=\C[\RZ]\supset\bigcap_{i<\omega}\underset{\substack{\TT{\tiny{off. dicht}}}}{G_i}\NEQ[\ref{BKS} $(i)$]\LM$\gs Widerspruch). 
\end{beweis}

\bigskip In Beispiel \ref{Bsp_AbzbReelleTM_Mager} wurde lediglich benutzt, da"s $\RZ$ keine isolierten Punkte besitzt und das \eemph{erste Trennungsaxiom} (je zwei Punkte besitzen Umgebungen, die den jeweils anderen Punkt nicht enthalten [\siehe sp"ateres Kapitel \ref{Kapitel_BairescheRaeume} Formel (\ref{Formel_T_1})] erf"ullt. Mit der Argumentation aus Beispiel \ref{Bsp_AbzbReelleTM_Mager} erh"alt man somit auch:

\begin{beispiel}\label{BspMagerErstTrennungsAx_etc}
Sei $X$ ein topologischer Raum, der das erste Trennungsaxiom erf"ullt und keine isolierten Punkte besitzt. Dann ist jede h"ochstens abz"ahlbare Teilmenge von $X$ mager.  
\end{beispiel}

\begin{beispiel}
Die Vereinigungsmenge $G$ der Graphen aller Polynome $f$ mit rationalen Koeffizienten ist mager im euklidischen
Raum $\RZ^2$.
\end{beispiel}
\begin{beweis}\KOM[K"ohnen 356]
Jedes der Polynome $f$ ist eine auf $\RZ$ stetige Funktion. Nach Beispiel \ref{Bsp_GraphStetNirgDicht} sind somit die Graphen $\textup{G}_f$ der Polynome nirgends dicht in $\RZ^2$. Da die Menge $L:=\bigcup_{i<\omega} \QZ[X_1,\ldots,X_i]$ der Polynome mit rationalen Koeffizienten abz"ahlbar ist, ist dann $G=\bigcup_{f\in L} \textup{G}_f$ als abz"ahlbare Vereinigung nirgends dichter Mengen eine magere Menge.
\end{beweis}

\begin{comment}
\begin{beispiel}
Die Menge $\BIF [a,b]$ der "uber einem abgeschlossenen reellen Intervall $[a,b]\subset \RZ$ beschr"ankten 
Riemann-integrierbaren Funktionen ist mager in sich selbst.
\end{beispiel}
\begin{beweis}\KOM[(FEHLT)]
\end{beweis}
\end{comment}

\begin{beispiel}
Sei $A\subset\RZ$ dicht und die Funktion $f:\PFEIL{\RZ}{}{\RZ}$ in $A$ stetig\footnote{\emph{stetig in $A\subset\RZ$:} Eine Funktion $f:\PFEIL{\RZ}{}{\RZ}$ hei"st \emph{stetig in $A\subset\RZ$}, falls $f$ stetig in allen Punkten $x\in A$ ist. F"ur ein $x\in\RZ$ hei"st $f$ \emph{stetig im Punkt $x$}, falls $\FA[\epsilon>0]\EX[\delta>0]\FA[y\in\RZ](\dd(x,y)<\delta\Rightarrow\dd(f(x),f(y))<\epsilon)$.}. Dann ist die Menge $M$ der
Unstetigkeitsstellen von $f$ mager.
\end{beispiel}
\begin{beweis}
Sei $M$ die Menge der Unstetigkeitsstellen von $f$. Dann l"a"st sich $M$ schreiben als $M=\bigcup_{n\geq 1} M_n$,
wobei eine reelle Zahl $x$ genau dann zu $M_n$ geh"ore, wenn es eine Folge $(x_k)_{k\geq 1}$ gibt, die in $\RZ$ gegen
$x$ konvergiert und f"ur die $|f(x_k)-f(x)|>\frac{1}{n}$ f"ur alle $k\geq1$.

Die $M_n$ sind nirgends dicht (f"ur alle $n\in\NZ$): Angenommen f"ur ein $n_0\geq 1$ ist $M_{n_0}$ in $\RZ$ 
\eemph{nicht} nirgends dicht. Wegen $\ABS[A]=\RZ$ gibt es einen Stetigkeitspunkt $y_0\in A$ von $f$, der 
Ber"uhrungspunkt von $M_{n_0}$ (d.h. in $M_{n_0}$ gibt es eine gegen $y_0$ konvergente Folge) ist. Da $f$ in $y_0$ 
stetig ist, gibt es zu $\frac{1}{2n_0}$ ein $\delta>0$ mit
\[
\FA[y\in\RZ](|y-y_0|<\delta\Rightarrow|f(y)-f(y_0)|<\frac{1}{2n_0}).
\]
Da $y_0$ Ber"uhrungspunkt von $M_{n_0}$ ist, gibt es ein $x_0\in M_{n_0}$ so, da"s $|x_0-y_0|<\delta$. Sei nun $(x_k)$ die oben genannte Folge, die in $M_{n_0}$ gegen $x_0$ konvergiert (existiert laut Definition von $M_{n_0}$).
Dann gilt f"ur hinreichend gro"ses $k\geq1$:
\begin{alignat*}{2}
|f(x_k)-f(x_0)|		&\leq	\underbrace{|f(x_k)-f(y_0)|}_{<\frac{1}{2n_0}}	
										+\underbrace{|f(y_0)-f(x_0)|}_{<\frac{1}{2n_0}}\\
									&<\frac{1}{n_0}
\end{alignat*}
Insgesamt gilt dann f"ur hinreichend gro"ses $k\geq1$:
\[
|f(x_k)-f(x_0)|>\frac{1}{n_0} \TT{ und } |f(x_k)-f(x_0)|<\frac{1}{n_0}.
\]
Widerspruch\gs also gibt es kein $n_0\in\NZ$, f"ur das $M_{n_0}$ nicht nirgends dicht w"are.\\
Somit ist dann $M=\bigcup_{n\geq 1} M_n$ mit nirgends dichten $M_n$\gs also ist $M$ mager.
\end{beweis}

\begin{beispiel}\label{Bsp2KatMetrRaumEndl_etc}
In einem metrischen Raum $(X,\dd)$ mit nur endlich vielen Elementen ist jede nicht-leere Teilmenge $A\subset X$ und insbesondere $X$ selber von 2. Kategorie.
\end{beispiel}
\begin{beweis}
Dies ergibt sich sofort daraus, da"s in einem metrischen Raum $X$ mit nur endlich vielen Punkten jede Singulum-Menge $\{x\}\subset X$ offen ist.
\begin{comment}
Da $X$ eine endliche Menge ist, enth"alt jede nicht-leere Teilmenge $A$ von $X$ eine offene $\epsilon$-Kugel. Dann ist $A$ nicht nirgends dicht. Somit ist klar, da"s keine nicht-leere Teilmenge in $X$ als abz"ahlbare Vereinigung nirgends dichter Teilmengen von $X$ dargestellt werden kann. Also sind alle Teilmengen von $X$ und insbesondere $X$ selbst von 2. Kategorie.
\end{comment} 
\end{beweis}

Allgemeiner gilt (nach der selben Argumentation), da"s jede Teilmenge eines topologischen Raumes $X$, die nur aus isolierten Punkten (d.h. aus Punkten in $X$, gegen die keine Folge in $X$ konvergiert [\siehe Kapitel \ref{Kap_SigmaKategorie}] besteht, von 2. Kategorie ist.

\begin{beispiel}\label{Bsp_mager_lauterLuecken}
Sei $I$ ein Intervall (also $I$ gleich $(a,b),(a,b],[a,b)$ oder $[a,b]$ mit $a\lneq b$ in $\RZ$) des euklidischen Raumes $\RZ$. Dann ist $I$ von 2. Kategorie in $\RZ$.
\end{beispiel}
\begin{beweis}
Der Baire'sche Kategoriensatz [\siehe Satz \ref{BKS}$(i)$] zeigt, da"s in vollst"andigen metrischen R"aumen (wie $\RZ$) offene Mengen niemals mager sind. Ist $I$ also offen, so kann $I$ nicht mager in $\RZ$ sein. Ist $I$ abgeschlossen oder halboffen, so enth"alt es ein offenes Intervall und kann daher ebenfalls nicht mager sein. 
\end{beweis}

\section{Komagere Mengen}\label{Kap_Komagere_Mengen}
Die \eemph{komageren} Mengen stehen in Baire's Kategorien-Konzept f"ur die \glqq gro"sen\grqq\ Mengen und werden daher als die Komplemente \glqq kleiner\grqq\ Mengen definiert. Analog zur Lokalisierung des Begriffes einer \eemph{mageren} Mengen wird auch definiert, wann eine Menge \eemph{komager in einer offenen Menge} ist:

\begin{definition}\header{komager in $X$, komager in $U\subset X$}\label{Def_KoMager}\\
Eine Menge $A \subset X$ hei"st \emph{komager in $X$} oder einfach \emph{komager} 
genau dann, wenn $\C[A]$ mager ist. Eine Menge $A \subset X$ hei"st \emph{komager in $U$} f"ur eine 
offene Menge $U\subset X$ genau dann, wenn $U\o A$ mager ist \gs das ist gleichbedeutend damit, da"s 
$\C[A]$ mager in $U$ ist.
\end{definition}

\begin{lemma}\header{komager in $U\subset X$}\label{Lemma_Komager}\\
Sei $U \subset X$ offen. Dann ist $A \subset X$ komager in $U$ genau dann, wenn es eine Folge 
$(G_i)_{i<\omega}$ offener dichter Mengen in $X$ gibt so, da"s
\begin{equation*} U \cap A \supset \bigcap_{i<\omega}{U \cap G_i}\end{equation*}
gilt.
\end{lemma}
\begin{beweis}
Seien $U\subset X$ offen und $A\subset X$. Dann gilt
\begin{alignat*}{2}
A \TT{ komager in U }
							&\AQ U\cap{\C[A]}		&					&\TT{mager}\\
						 	&\AQ U\cap{\C[A]} 	&= 				&\bigcup_{i<\omega}{U\cap A_i}\TT{ mit $A_i$ nirgends dicht}\\
							&\AQ U\cap A 			&= 				&\bigcap_{i<\omega}{U\cap {\C[A_i]}}\\
							&\AQ U\cap A 			&\supset 	&\bigcap_{i<\omega}{U\cap G_i}\TT{ mit $G_i\subset \C[A_i]$ offen und dicht.}
\end{alignat*}
\end{beweis}

Insbesondere gilt, da"s $A\subset X$ genau dann komager ist, wenn $A$ einen abz"ahlbaren Schnitt 
$\bigcap_{i<\omega}{G_i}$ offener dichter Teilmengen $G_i$ von $X$ enth"alt.

\bigskip In Kapitel \ref{Kap_Magere_Mengen} wurde bereits gezeigt, da"s magere Mengen \eemph{klein} im Sinne von Kapitel \ref{Kap_Kleine_Gro"se_Mengen} sind (d.h. sie bilden ein $\sigma$-Ideal). Als Komplemente magerer Mengen verhalten sich die komageren Mengen somit genau so, wie wir es in Kapitel \ref{Kap_Kleine_Gro"se_Mengen} von \eemph{ gro"sen} Mengen gefordert haben: Sie sind abgeschlossen gegen"uber Obermengen und abz"ahlbaren Schnitten.

\section{Baire'sche R"aume}\label{Kapitel_BairescheRaeume}

F"ur die spieltheoretischen Untersuchungen werden wir in dieser Arbeit fast ausschlie"slich den Baire-Raum verwenden. Dieser wird in Abschnitt \ref{Kap_Baire_Raum} eingef"uhrt und einige seiner grundlegenden Eigenschaften umrissen. Da die Baire-Kategorie jedoch nicht blo"s f"ur den Baire-Raum, sondern f"ur allgemeine topologische R"aume definiert ist, erscheint eine Einordnung in den allgemeineren Zusammenhang der Topologie angebracht.

\bigskip In diesem Abschnitt werden die Bair'schen R"aume allgemein definiert und einige ihrer grundlegende topologischen Eigenschften vorgestellt. Der Baire'sche Kategoriensatz \ref{BKS} identifiziert zwei gro"se Klassen von Topologischen R"aumen als Baire'sch: die lokal kompakten und die vollst"andig metrisierbaren R"aume. 

\bigskip Offene (nicht-leere) Intervalle in $\RZ$ und somit beliebige nicht-leere Teilmengen von $\RZ$ sind niemals mager [\siehe Beispiel \ref{Bsp_mager_lauterLuecken}]. Um diese Eigenschaft f"ur topologische R"aume zu verallgemeinern definiert man:

\begin{definition}\header{Baire'scher Raum}\label{Def_BairescherRaum}\\
Ein topologischer Raum $X$ hei"st \emph{Baire'sch} genau dann, wenn er
keine nicht-leere offene magere Menge enth"alt.
\end{definition}

Eine nicht-leere offene Teilmenge $U\subset X$ ist mager genau dann, wenn jede 
Teilmenge von $X$ mager (und komager) in $U$ ist. Die Baire'schen R"aume sind daher gerade diejenigen, in denen die \eemph{Lokalisierung} der Begriffe \eemph{mager} und \eemph{komager} sinnvoll ist.

\begin{lemma}\header{Baire'scher Raum}\\
F"ur einen topologischen Raum $X$ ist "aquivalent
\begin{ITEMS}[equivalences]
\item $X$ ist Baire'sch.
\item F"ur jede Folge $(D_i)_{i<\omega}$ offener dichter Teilmengen von $X$ ist der Schnitt 
$\bigcap_{i\geq 0}{D_i}$ dicht in $X$.
\item Jede komagere Menge in $X$ ist dicht.
\item F"ur jede Folge $(A_i)_{i<\omega}$ abgeschlossener Teilmengen von $X$ mit $\INT[A_i]=\LM$ ist das Innere der Vereinigung $\INT[(\bigcup_{i<\omega}{A_i})]$ ebenfalls leer.
\end{ITEMS}
\end{lemma}
\begin{beweis}
$(i)\Rightarrow(ii):$ Sei $X$ Baire'sch und sei $(D_i)_{i<\omega}$ eine Folge offener nicht-leerer dichter Teilmengen in $X$. Angenommen $\bigcap_{i<\omega}D_i$ ist nicht dicht in $X$. Sei etwa $U\subset X$ offen nicht-leer und schnittfremd mit $\bigcap_{i<\omega}D_i$. Dann gilt
\begin{alignat*}{1}
&U=U\o \bigcap_{i<\omega}D_i=\bigcup_{i<\omega}(U\o D_i)
\end{alignat*}
und die $U\o D_i$ sind nirgends dicht (laut Definition \ref{Def_Nirg_Dicht}, denn das Komplement enth"alt die offene in $X$ dichte Teilmenge $D_i$). Somit ist $U$ mager und offen\gs im Widerspruch dazu, da"s $X$ als Baire'sch vorausgestzt war.

\smallskip Also ist obige Annahme falsch, d.h. $\bigcap_{i<\omega}D_i$ ist dicht in $X$.

\bigskip$(ii)\Rightarrow(iii):$ Sei $A\subset X$ komager und es gelte $(ii)$. Angenommen es gibt eine offene nicht-leere Menge $U\subset X$ mit $U\cap A=\LM$. Da $A$ komager ist, gilt laut Lemma \ref{Lemma_Komager}: $A\supset\bigcap_{i<\omega}G_i$ f"ur einige $G_i\subset X$ offen und dicht. Wegen $U\cap A=\LM$ gilt dann auch 
\begin{alignat*}{1}
&U\cap\bigcap_{i<\omega}G_i=\LM.
\end{alignat*}
Da die $G_i$ offen und dicht sind, ist jedoch $\bigcap_{i<\omega}G_i$ nach $(ii)$ dicht\gs es gilt also 
\begin{alignat*}{1}
&U\cap\bigcap_{i<\omega}G_i\neq\LM.
\end{alignat*}
Widerspruch.

\smallskip Somit war obige Annahme falsch\gs es gibt also keine offene nicht-leere Menge $U\subset X$, die $A$ nicht schneidet. Das hei"st $A$ ist dicht.

\bigskip$(iii)\Rightarrow(iv):$ Sei $(A_i)_{i<\omega}$ eine Folge abgeschlossener Teilmengen von $X$ mit $\INT[A_i]=\LM$ und es gelte $(iii)$. Angenommen $\INT[(\bigcup_{i<\omega}{A_i})]\neq\LM$. Sei etwa $U\subset \bigcup_{i<\omega}{A_i}$ offen und nicht-leer. Dann gilt f"ur das Komplement 
\begin{alignat*}{1}
&\C[(\bigcup_{i<\omega}{A_i})]\TT{ nicht dicht}.
\end{alignat*}
Andererseits ist aber $\bigcup_{i<\omega}{A_i}$ mager (da $A_i$ abgeschlossen mit $\INT[A_i]=\LM$ also $\C[A_i]$ offen und dicht also $A_i$ nirgends dicht) und somit gilt nach $(iii)$: 
\begin{alignat*}{1}
&\C[(\bigcup_{i<\omega}{A_i})]\TT{ dicht}.
\end{alignat*}
Widerspruch.

\smallskip Obige Annahme war also falsch und es gilt $\INT[(\bigcup_{i<\omega}{A_i})]=\LM$.

\bigskip$(iv)\Rightarrow(i):$
Sei $U\subset X$ offen nicht-leer und es gelte $(iv)$. Angenommen $U$ ist mager\gs etwa $U=\bigcup_{n<\omega}A_i$ mit $A_i$ nirgends dicht. Dann ist das Innere des Abschlusses $\INT[{\ABS[A_i]}]$ leer f"ur alle $A_i$. Nach $(iv)$ gilt dann
\begin{alignat*}{1}
&\INT[{(\bigcup_{i<\omega}{A_i})}]=\LM.
\end{alignat*}
Da $U=\bigcup_{i<\omega}{A_i}$ offen und nicht-leer ist gilt aber auch
\begin{alignat*}{1}
&\INT[{(\bigcup_{i<\omega}{A_i})}]\neq\LM.
\end{alignat*}
Widerspruch.

\smallskip Obige Annahme war also falsch, es gibt also keine offene magere Menge $U\subset X$\gs das hei"st $X$ ist Baire'sch.
\end{beweis}

Endliche topologische R"aume sind immer Baire'sch, da sie nur endlich viele offene Mengen besitzen und Schnitte endlich
vieler offener dichter Mengen immer dicht sind. Wichtige Beispiele f"ur Baire'sche R"aume sind lokal kompakte Hausdorff-R"aume und vollst"andig metrisierbare R"aume\gs insbesondere die polnischen R"aume. Zun"achst ist es jedoch notwendig einige weitere topologische Grundbegriffe zu definieren:

\begin{definition}\header{Hausdorff'scher Raum oder $T_2$-Raum}\label{Def_T2}\\
Ein topologischer Raum $(X,\TOP{X})$ hei"st \emph{Hausdorff'sch} oder \emph{$T_2$-Raum}, wenn je zwei verschiedene Punkte von $X$ disjunkte Umgebungen besitzen ("aquivalent: disjunkte offene Umgebungen besitzen), d.h.:
\begin{alignat}{1}
&\FA[x\neq y]\EX[U,V\in\TOP{X}](x\in U\UND y\in V\UND U\cap V=\LM)\label{Formel_Hausdorff_Eigenschaft}.
\end{alignat}
Die Formel (\ref{Formel_Hausdorff_Eigenschaft}) nennt man auch die \emph{Hausdorff-Eigenschaft} oder das \emph{zweite Trennungsaxiom} (auch kurz \emph{$T_2$} genannt).
\end{definition}

\bigskip F"ur den Beweis des Baire'schen Kategoriensatzes in diesem Kapitel spielen f"ur uns neben dem zweiten Trennungsaxiom auch noch das sogenannte \emph{erste} und das \emph{dritte Trennungsaxiom} eine gewisse Rolle. Sei wie oben $(X,\TOP{X})$ ein topologischer Raum. Nach dem \emph{ersten Trennungsaxiom} (auch kurz \emph{$T_1$} genannt) haben je zwei unterschiedliche Punkte offene Umgebungen, die den jeweils anderen Punkt nicht enthalten, d.h.:
\begin{alignat}{1}
&\FA[x\neq y]\EX[U,V\in\TOP{X}](x\in U\UND y\in V\UND x\not\in V\UND y\not\in U)\label{Formel_T_1}.
\end{alignat}
Das \emph{dritte Trennungsaxiom} (auch kurz \emph{$T_3$} genannt) besagt, da"s jede abgeschlossene Menge und jeder Punkt, der nicht in dieser Menge liegt, schnittfremde offene Umgebungen besitzen, d.h.:
\begin{alignat}{1}
&\FA[A\TT{ abg.}]\FA[x\not\in A]\EX[U,V\in\TOP{X}](A\subset U\UND x\in V\UND U\cap V=\LM)\label{Formel_T_3}.
\end{alignat}
Einen Raum, der das erste und dritte Trennungsaxiom erf"ullt, nennt man \emph{regul"ar}.

\subsection*{Lokal kompakte R"aume}
Der Begriff der \emph{Kompaktheit} ist in der Topologie von zentraler Bedeutung. Dabei werden gewisse Endlichkeitseigenschaften von "Uberdeckungen gefordert. Auf diese Weise erh"alt man die \emph{kompakten} bzw. die \emph{lokal kompakten} R"aume. Zugrunde liegt diesen Begriffen die f"ur die reelle Analysis wichtige topologische Eigenschaft, da"s man abgeschlossene und beschr"ankte Teilmengen des $\RZ^n$ durch "Uberdeckungen charakterisieren kann (Heine-Borel'scher "Uberdeckungssatz\footnote{\emph{Heine-Borel'scher "Uberdeckungssatz:} In der deutschen Literatur nach \'Emile Borel (1871--1956) und Eduard Heine (1821--1881) benannter "Uberdeckungssatz, nach dem die kompakten Teilmengen von $\RZ^n$ gerade die beschr"ankt und abgeschlossenen sind. Allerdings formulierte und bewies Borel diesen Satz lediglich f"ur abz"ahlbare "Uberdeckungen und weist in einer Notiz "uber seine wissenschaftlichen Arbeiten selbst darauf hin, da"s H. Lebesgue den Beweis f"ur beliebige "Uberdeckungen erbracht hat. Daher wird dieser Satz in der franz"osischen Literatur zurecht als Satz von Borel und Lebesgue bezeichnet. (\siehe vgl. \cite[IV, 2,2]{Bourbaki:1998:1})}).\KOM[LIT Elstrodt S.39]

\bigskip Eine Familie $(A_i)_{i\in I}$ von Teilmengen einer Menge $X$ hei"st eine \emph{"Uberdeckung} von $X$, wenn $X\subset\bigcup_{i\in I}A_i$ gilt. Der Bequemlichkeit halber nennt man manchmal auch die Darstellung $X\subset\bigcup_{i\in I}A_i$ eine \emph{"Uberdeckung} von $X$. Ist $I$ endlich (abz"ahlbar), so hei"st die "Uberdeckung $(A_i)_{i\in I}$ \emph{endlich (abz"ahlbar)}. Ist $(X,\TOP{X})$ ein topologischer Raum, $(A_i)_{i\in I}$ eine "Uberdeckung von $X$ und sind alle $A_i$ offen (abgeschlossen), so hei"st die "Uberdeckung $(A_i)_{i\in I}$ \emph{offen (abgeschlossen)}. Ist $(A_i)_{i\in I}$ eine "Uberdeckung einer Menge $X$, $J\subset I$ und $(A_j)_{j\in J}$ ebenfalls eine "Uberdeckung von $X$, so hei"st $(A_j)_{j\in J}$ eine \emph{Teil"uberdeckung} von $(A_i)_{i\in I}$.

\bigskip Ein topologischer Raum hei"st \emph{quasikompakt}, falls jede offene "Uberdeckung des Raumes eine endliche Teil"uberdeckung besitzt. Ein quasikompakter $T_2$-Raum hei"st \emph{kompakt\footnote{\emph{(quasi-)kompakt:} Viele Autoren verwenden den Begriff \emph{kompakt} anstelle von \emph{quasikompakt}. Wir halten uns hier an die Notation von \cite{Preuss:1970} und fordern f"ur \emph{kompakt} zus"atzlich noch die Hausdorff-Eigenschaft.}}. Eine Teilmenge eines topologischen Raumes hei"st \emph{(quasi-) kompakt}, wenn sie als Unterraum (quasi-)kompakt ist.

\bigskip In einem quasikompakten Raum $(X,\TOP{X})$ l"a"st sich wie folgt von lokalen auf globale Eigenschaften schlie"sen: Sei $\SYS{E}\subset\TOP{X}$ eine Eigenschaft von offenen Teilmengen von $X$ so, da"s f"ur beliebige offene Teilmengen $U,V$ von $X$ gilt:
\[
U,V\in\SYS{E}\Rightarrow U\cup V\in\SYS{E}.
\]
Hat dann $X$ die Eigenschaft $\SYS{E}$ \emph{lokal} (d.h. f"ur jedes $x\in X$ gibt es eine offene Umgebung $U_x$ von $x$ mit der Eigenschaft $\SYS{E}$), dann hat schon der ganze Raum $X$ die Eigenschaft $\SYS{E}$. Diese Schlu"sweise ist deshalb korrekt, weil die Vereinigung $\bigcup_{x\in X}U_x\supset X$ aller $U_x$ mit jeweils der Eigenschaft $\SYS{E}$ eine "Uberdeckung von $X$ ist. Da $X$ quasikompakt ist, gibt es dazu eine endliche Teil"uberdeckung $\bigcup_{j=1}^n U_{x_j}\supset X$ die die Eigenschaft $\SYS{E}$ auf $X$ "ubertr"agt.

\bigskip Quasikompakte R"aume sind alternativ durch folgende Eigenschaft charakterisiert:

\begin{lemma}\header{quasikompakt}\label{Lemma_Char_Quasikompakt}\\
F"ur einen topologischen Raum $(X,\TOP{X})$ ist "aquivalent:
\begin{ITEMS}
\item $(X,\TOP{X})$ ist quasikompakt.
\item In jeder Familie $(A_i)_{i\in I}$ abgeschlossener Teilmengen von $X$ mit $\bigcap_{i\in I} A_i=\LM$ gibt es endlich viele Glieder $A_{i_1},\ldots,A_{i_n}$ mit $\bigcap_{k=1}^n A_k=\LM$.
\end{ITEMS}
\end{lemma}
\begin{beweis}
Die Aussagen $(i)$ und $(ii)$ sind zueinander dual: Die Behauptung, da"s der Durchschnitt abgeschlossener Teilmengen $\bigcap_i A_i$ leer ist, ist gleichbedeutend damit, da"s $\bigcup_i \C[A_i]$ gleich $X$ ist, d.h.: da"s die Komplemente der abgschlossenen Mengen eine offene "Uberdeckung von X sind.
\begin{comment}
$(i)\Rightarrow(ii):$ 
Sei $(X,\TOP{X})$ quasikompakt und $(A_i)_i$ eine beliebige Familie von abgeschlossenen Teilmengen in $X$ mit $\bigcap_{i} A_{i}=\LM$. 

\smallskip\emph{Es ist zu zeigen:} Es gibt endlich viele $A_{i_1},\ldots,A_{i_n}$ mit $\bigcap_{j=1}^n A_{i_j}=\LM$.

\smallskip Es gilt:
\begin{alignat*}{3}
						\bigcap_{i} A_{i}=\LM
\Leftrightarrow	&\bigcup_{i} \C[A_{i}]=X
\underset{\tiny{X \TT{ quasikomp.}}}{\Rightarrow}	&\bigcup_{j=1}^n \C[A_{i_j}]=X
\Rightarrow	&\bigcap_{j=1}^n A_{i_j}=\LM
\end{alignat*}
und somit folgt die Behauptung.

\bigskip $(i)\Leftarrow(ii):$ 
Gelte in $X$, da"s f"ur jede Familie $(A_i)_i$ abgeschlossener Teilmengen in $X$ mit $\bigcap_{i} A_{i}=\LM$ es schon endlich viele $A_{i_1},\ldots,A_{i_n}$ gibt mit $\bigcap_{j=1}^n A_{i_j}=\LM$. Sei $(O_i)_i$ eine beliebige offene "Uberdeckung von $X$.

\smallskip\emph{Es ist zu zeigen:} Zu $(O_i)_i$ gibt es eine endliche Teil"uberdeckung von $X$.

\smallskip Es gilt:\KOM[FEHLT]
\begin{alignat*}{3}
									\bigcup_{i} O_{i}=X
\Rightarrow	&\bigcap_{i} \C[O_{i}] = \LM
\underset{\tiny{\TT{Vorauss.}}}{\Rightarrow}	&\bigcap_{j=1}^n \C[O_{i_j}]=\LM
\Rightarrow	&\bigcup_{j=1}^n O_{i_j}=X
\end{alignat*}
und somit folgt die Behauptung.
\end{comment}
\end{beweis}

F"ur den Beweis des Baire'schen Kategoriensatzes \ref{BKS} sollen nun schnell noch einige weitere Lemmata bewiesen werden.

\begin{lemma}\header{hinr. Krit. f"ur abgeschlossen}\label{Lemma_Abgeschlossen}\\
Kompakte Teilmengen von $T_2$-R"aumen sind abgeschlossen.
\end{lemma}
\begin{beweis}
Sei $X$ ein $T_2$-Raum und $K\subset X$ kompakt. 

\point{Es gen"ugt zu zeigen} $\C[K]$ ist offen bzw.: f"ur jedes $y\in\C[K]$ liegt eine offene Umgebung $V_y$ von $y$ in $\C[K]$.

\bigskip Sei $y\in \C[K]$ beliebig. Da $X$ ein $T_2$-Raum ist, gibt es zu jedem $x\in K$ schnittfremde offene Umgebungen $O_x$ von $x$ und $V_y^x$ von $y$. Es gilt dann $K\subset\bigcup_{x\in K} O_{x}$ und
\[
\bigcup_{x\in K} O_{x}\cap \bigcap_{x\in K} V_{y}^x = \LM.
\]

Da $K$ kompakt ist, gibt es dazu schon eine endliche Teil"uberdeckung $\bigcup_{i=1}^n O_{x_i}\supset K$ und es gilt: 
\[
\bigcup_{i=1}^n O_{x_i}\cap \bigcap_{i=1}^n V_{y}^{x_i} = \LM.
\]
Da $\bigcap_{i=1}^n V_{y}^{x_i}$ eine offene Umgebung von $y$ ist und ganz in $\C[K]$ liegt, folgt dann die Behauptung.
\end{beweis}

Anders als in $\RZ^n$ reicht in kompakten R"aumen die Abgeschlossenheit zur Charakterisierung kompakter Teilmengen aus:

\begin{lemma}\header{hinr. Krit. f"ur (quasi-)kompakt}\label{Lemma_Quasikompakt}\\
Abgeschlossene Teilmengen von (quasi-)kompakten R"aumen sind (quasi-)kompakt.
\end{lemma}
\begin{beweis}
Sei $X$ ein quasikompakter Raum und $A\subset X$ abgeschlossen. Sei $\bigcup_i U_i\supset A$ eine beliebige offene "Uberdeckung von $A$. Dann ist
\[
\bigcup_i U_i\cup\C[A]=X 
\]
eine offene "Uberdeckung von $X$. Da $X$ quasikompakt ist existiert dazu eine endliche offene Teil"uberdeckung
\[
\bigcup_{i=1}^n U_i\cup\C[A]=X 
\]
von $X$. Dann mu"s $\bigcup_{i=1}^n U_i\supset A$ eine endliche Teil"uberdeckung von $A$ zu $(U_i)_i$ sein. Da die offene "Uberdeckung $(U_i)_i$ beliebig gew"ahlt war, ist $A$ somit quasikompakt.

\smallskip Da sich die Hausdorff-Eigenschaft auf Teilmengen "ubertr"agt, gilt die Aussage analog auch f"ur abgeschlossene Teilmengen kompakter R"aume.
\end{beweis}

\begin{lemma}\header{regul"ar}\label{Lemma_Regul"ar}\\
Ein kompakter Raum ist immer auch regul"ar.
\end{lemma}
\begin{beweis}
Sei $X$ ein kompakter Raum.

\point{$X$ erf"ullt $T_1$} Ein kompakter (d.h. quasikompakter und Hausdorff'scher) Raum erf"ullt immer das erste Trennungsaxiom, weil $T_2\Rightarrow T_1$ gilt [\siehe Formeln (\ref{Formel_Hausdorff_Eigenschaft}) und (\ref{Formel_T_1})].

\point{$X$ erf"ullt $T_3$} Sei $A\subset X$ eine beliebige abgeschlossene Teilmenge von $X$ und sei $x\in X$ ein beliebiger Punkt in $X$, der nicht in $A$ liegt. Dann gibt es nach $T_2$ f"ur jeden Punkt $y\in A$ schnittfremde Umgebungen $U_y$ von $y$ und $U_x^y$ von $x$. Die Vereinigung $\bigcup_{y\in A} U_y$ ist eine offene "Uberdeckung von $A$ und 
\begin{alignat*}{1}
&\bigcup_{y\in A}U_y\cup\C[A]=X
\intertext{ist eine offene "Uberdeckung von $X$. Da $X$ quasikompakt ist gibt es zu dieser "Uberdeckung schon eine endliche Teil"uberdeckung von $X$:}
&\bigcup_{j=1}^n U_{y_j}\cup\C[A]=X.
\end{alignat*}
Mit $U_A:=\bigcup_{j=1}^n U_{y_j}$ und $U_x:=\bigcap_{j=1}^n U_x^{y_j}$ erh"alt man so zwei schnittfremde offene Umgebungen von $A$ und $x$. D.h. $X$ erf"ullt $T_3$.
\end{beweis}

Einen Raum, der lokal aussieht wie ein kompakter Raum, nennt man \emph{lokal kompakt}: 

\begin{definition}\header{lokal (quasi-)kompakt, relativ kompakt}\label{Def_Lokal_Quasikompakt}\\
Ein topologischer Raum $(X,\TOP{X})$ hei"st \emph{lokal quasikompakt}, wenn jeder Punkt $x\in X$ eine quasikompakte Umgebung (d.h. eine quasikompakte Teilmenge, die Umgebung des Punktes ist) besitzt.

\smallskip Ein lokal quasikompakter $T_2$-Raum hei"st \emph{lokal kompakt} 

\smallskip Eine Teilmenge $A$ eines topologischen Raumes $(X,\TOP{X})$ hei"st \emph{relativ kompakt}, wenn $\ABS[A]$ kompakt ist.
\end{definition}

In einem lokal kompakten Raum enth"alt jede Umgebung eines Punktes schon eine kompakte Umgebung dieses Punktes:

\begin{lemma}\header{lokal kompakt}\label{Lemma_Lokal_Kompakt}\\
In einem lokal kompakten Raum $(X,\TOP{X})$ gilt:
\[
\FA[x\in X]\FA[U\TT{ Umgeb. von $x$}]\EX[\TT{komp. Umgeb. $K$ von $x$}](K\subset U).
\]
\end{lemma}
\begin{beweis}
Sei $(X,\TOP{X})$ ein lokal kompakter Raum und sei $x\in X$ beliebig mit einer beliebigen offenen Umgebung $U$ von $x$. Zu $x$ gibt es dann eine kompakte Umgebung $K$ von $x$. Sei $V:=\INT[(U\cap K)]$. $V$ ist also eine in $K$ enthaltene $X$-offene Umgebung von $x$ ($x\in V$, da $K$ komp. Umgebung von $x$\gs etwa $x\in G\subset K$ mit $G$ offen, somit folgt $x\in G\cap U\subset\INT[(U\cap K)]$). $K$ ist kompakt also nach Lemma \ref{Lemma_Regul"ar} auch regul"ar\gs erf"ullt also insbesondere $T_3$. Nun ist $K\o V$ abgeschlossen in $K$ und $x$ liegt nicht in $K\o V$. Demnach gibt es nach $T_3$ eine in $K$ offene Umgebung $O_{K\o V}$ von der $K$-abgeschlossenen Menge $K\o V$ und eine $K$-offene Umgebung $H\subset V$ vom Punkt $x$ so, da"s $O_{K\o V}\cap H=\LM$. Sei $O_{K\o V}$ etwa gleich $O\cap K$ f"ur $O\subset X$ offen. Dann ist $W:=\C[O]\cap K$ eine $K$-abgeschlossene Umgebung von $x$ in $K$ ($W$ ist in $K$ eine Umgebung von $x$, da $H\subset W$ $K$-offen). Die Menge $W$ ist eine kompakte Umgebung von $x$ in $X$:

\point{$W$ ist kompakt} Da $W$ eine abgeschlossene Teilmenge des kompakten Raumes $K$ ist, ist $W$ nach Lemma \ref{Lemma_Quasikompakt} quasikompakt. Die Hausdorff-Eigenschaft "ubertr"agt sich von $K$ auf $W$. Also ist $W$ insgesamt kompakt. 

\point{$W$ ist eine Umgebung von $x$ in $X$} Da $W$ eine Umgebung von $x$ in $K$ ist, gibt es eine in $X$ offene Menge $G\subset X$ mit $x\in G\cap K\subset W \subset V\subset K$. Dann ist $G\cap K\cap V\EQ[\tiny{$G\cap V\subset K$}] G\cap V$ und wegen $G\cap K\subset W$ in $W$ enthalten und eine in $X$ offene Umgebung von $x$. Es gilt also: $W$ ist ein Umgebung von $x$ in $X$.
\end{beweis}

Nun fehlen nur noch einige Grundbegriffe "uber metrische R"aume, und der Baire'schen Kategoriensatz kann bewiesen werden.

\subsection*{Vollst"andig metrisierbare R"aume}
In der Mathematik wird die anschauliche Vorstellung vom \glqq Abstand\grqq\ zweier Punkte durch den Begriff der \emph{Metrik} beschrieben. F"ur eine Menge $X$ ist eine Abbildung $d:\PFEIL{X\times X}{}{\RZ}$ in die reellen Zahlen eine Metrik, falls gilt:
\begin{ITEMS}[arabic)]
\item $\dd(x,y)=0\Leftrightarrow x=y$
\item $\dd(x,y) = \dd(y,x)$
\item $\dd(x,z) \leq \dd(x,y) + \dd(y,z)$
\end{ITEMS}
f"ur alle $x,y,z\in X$. Das Paar $(X,\dd)$ nennt man einen \emph{metrischen Raum\footnote{\emph{metrischer Raum:} Von Maurice Fr\'echet 1906 in \cite{Frechet:1906} eingef"uhrt (ebenda eine Axiomatisierung des Begriffes der \eemph{Konvergenz}).}}.  F"ur je zwei Elemente $x,y\in X$ hei"st $\dd(x,y)$ der \emph{Abstand} von $x$ und $y$. 

\bigskip Beispielsweise ist $(\RZ^n,\textup{d}^n)$ f"ur $n<\omega$ mit 
\[
\textup{d}^n((x,\ldots,x_n),(y,\ldots,y_n)):=\sqrt{(x_1-y_1)^2+\ldots+(x_n-y_n)^2}
\]
ein metrischer Raum. $(\RZ^n,\textup{d}^n)$ hei"st der $n$-dimensionale \emph{euklidische Raum} und die Metrik $\textup{d}^n$ die zugeh"orige \emph{euklidische Metrik}. 

\bigskip Die Begriffe Metrik und Topologie stehen miteinander in einem engen Zusammenhang, denn jeder metrische Raum ist gleichzeitig auch ein topologischer Raum. Die offenen $\epsilon$-Kugeln $O_\epsilon(x):=\MNG{y\in X}{\dd(x,y)<\epsilon}$ f"ur Elemente $x\in X$ und reelle $\epsilon>0$ bilden die Basis einer Topologie auf $X$. Man nennt diese Topologie auch die \emph{von $\dd$ induzierte Topologie}. 

\bigskip Umgekehrt l"a"st sich nicht selten f"ur eine gegebene Topologie $\TOP{X}$ auf einer Menge $X$ eine Metrik $\dd$ auf $X$ angeben, die die Topologie $\TOP{X}$ induziert. In diesem Fall nennt man den topololgischen Raum $(X,\TOP{X})$ \emph{metrisierbar}. Ist ein Raum metrisierbar, lassen sich also alle topologischen Aussagen wahlweise auch in der Sprache der Metrik formulieren.

\bigskip Ein metrischer Raum $(X,\dd)$ hei"st \emph{vollst"andig}, wenn in ihm jede Cauchyfolge konvergiert. Dabei bezeichnet man eine Folge von Punkten $(x_i)_{i\in I}$ in $X$ als \emph{Cauchyfolge\footnote{\emph{Cauchyfolge:} Benannt nach Augustin Cauchy (1789-1857).}}, wenn 
\[
\FA[\epsilon>0]\EX[n_0<\omega]\FA[n,m>n_0](\dd(x_n,x_m)<\epsilon).
\] 
Ein topololgischer Raum $(X,\TOP{X})$ hei"st \emph{vollst"andig metrisierbar}, wenn es eine Metrik $\dd$ auf $X$ gibt so, da"s $\dd$ die Topologie $\TOP{X}$ induziert und $(X,\dd)$ vollst"andig ist. Ein topologischer Raum hei"st  \emph{separabel}, wenn er eine abz"ahlbare dichte Teilmenge enth"alt. Einen separabelen und vollst"andig metrisierbaren Raum nennt man auch einen \emph{polnischen Raum\footnote{\emph{polnischer Raum:} "Aquivalent: Vollst"andig metrisierbar und mit einer abz"ahlbaren Basis (da f"ur metrische R"aume $X$ gilt: $X$ separabel $\Rightarrow$ $X$ hat abz"ahlbare Basis).}}.

\bigskip Anschlie"send soll der Baire'schen Kategoriensatzes f"ur die vollst"andig metrisierbaren und die lokal kompakten R"aume bewiesen werden. Der Beweis f"ur die vollst"andig metrisierbaren R"aume beruht auf dem sogenannten \glqq Schachtelungsprinzip\grqq\ (f"ur den lokal kompakten Fall ben"otigen wir an entsprechender Stelle die Charakterisieung quasikompakter R"aume durch Lemma \ref{Lemma_Char_Quasikompakt}). Das Schachtelungsprinzip besagt, da"s in vollst"andig metrisierbaren R"aumen eine absteigende Folge abgeschlossener Teilmengen, deren Durchmesser gegen null geht, genau ein Element enth"alt\gs insbesondere also nicht leer ist.

\bigskip F"ur Teilmengen $Y$ eines metrischen Raumes $(X,\dd)$ definiert man den \emph{Durchmesser} von $Y$ als
\[
\diam(Y):=\sup\MNG{\dd(x,y)}{x,y\in Y}.
\] 

\begin{lemma}\header{Schachtelungsprinzip}\label{Lemma_Schachtelungsprinzip}\\
Sei $(X,\TOP{X})$ ein vollst"andig metrisierbarer Raum und $A_0\supset A_1\supset A_2\supset\ldots$ eine absteigende Folge abgeschlossener Teilmengen von $X$, f"ur die $\diam(A_i)$ gegen $0$ geht. Dann ist $\bigcap_{i<\omega} A_i=\{x\}$ und insbesondere:
\[
\bigcap_{i<\omega} A_i\neq\LM.
\]
\end{lemma}
\begin{beweis}
Sei $(x_i)_i$ eine Folge mit $(x_i\in A_i)_{i<\omega}$ (eine solche Folge erh"alt man bereits mittels $\AC[\omega]_X$). Dann ist $(x_i)_{i<\omega}$ eine Cauchy-Folge, da gilt:
\[
\dd(x_m,x_n)\leq\diam(A_N)\TT{ f"ur $n,m\geq N$}.
\] 
Da $X$ vollst"andig ist, konvergiert $(x_i)_{i<\omega}$ gegen ein $x\in X$. Da die $A_i$ abgeschlossen sind ist auch $\bigcap_{i<\omega} A_i$ abgeschlossen und $x$ mu"s in $\bigcap_{i<\omega} A_i$ liegen. Also gilt $\bigcap_{i<\omega} A_i\neq\LM$. Da nun die Durchmesser der $A_i$ gegen null gehen, kann kein weiteres Element $y\in X$ im Schnitt der $A_i$ liegen. Also $\bigcap_{i<\omega} A_i=\{x\}$
\end{beweis}

\subsection*{Baire'scher Kategoriensatz}
Die vollst"andig metrisierbaren und die lokal kompakten R"aume sind zwei sehr unterschiedliche Klassen topologischer R"aume. Dennoch zeigt der folgende Satz, da"s sie beide zur Klasse der Baire'schen R"aume geh"oren. Dies ist auch ein Grund daf"ur, da"s man sich in der Topologie f"ur Baire'sche R"aume an sich interessiert (\siehe etwa \cite[IX, 5,3]{Bourbaki:1998:2}\KOM[LIT Dis.: 13 p. 193]).

\begin{satz}\header{Baire'scher Kategoriensatz}\label{BKS}
\mbox{}
\begin{ITEMS}
\item Jeder vollst"andig metrisierbare Raum ist Baire'sch.
\item Jeder lokal kompakte Raum ist Baire'sch.
\end{ITEMS}
\end{satz}
\begin{beweis}
\KOM[Kechris:41,42 oder Diss.:19 oder Preu"s]
$(i)$ Sei $X$ ein vollst"andig metrisierbarer Raum. Sei $(D_n)_n$ eine beliebige Folge offener dichter
Teilmengen in $X$.

\point{Es ist zu zeigen} $\FA[O\subset X\TT{ offen, nicht-leer}](O\cap \bigcap_n D_n\neq\LM)$.

\medskip Sei $U_1\subset X$ offen und nicht-leer. Dann ist $U_1\cap D_1\neq\LM$ (da $D_1$ dicht ist). Sei etwa $x_1\in U_1\cap D_1$, d.h. $U_1\cap D_1$ ist eine offene Umgebung von $x_1$. Da $X$ metrisch ist, gibt es dann eine abgeschlossene Kugel $U_2:=\ABS[{O_{\epsilon_2}(x_1)}]\subset U_1\cap D_1$. Dann ist $\INT[U_2]\cap D_2\neq\LM$ (da $D_2$ dicht ist). Sei etwa $x_2\in\INT[U_2]\cap D_2$. Da $X$ metrisch ist, gibt es dann erneut eine abgeschlossene Kugel $U_3:=\ABS[{O_{\epsilon_3}(x_2)}]\subset \INT[U_2]\cap D_2$. 

\medskip Induktiv fortfahrend erh"alt man so eine absteigende Folge $U_2\supset U_3\supset\ldots$ nicht-leerer abgeschlossener Kugeln in $X$ mit Radien $(0<\epsilon_i<\epsilon_{i-1}/2)_{i>2}$, d.h. $\diam(U_i)$ geht gegen $0$, f"ur die gilt:
\begin{alignat}{1}
&U_n\subset \INT[U_{n-1}]\cap D_{n-1}\label{Formel_BCT_3}
\end{alignat}
f"ur $n\geq 2$. Da wir die $U_n$ f"ur $n\geq2$ so gew"ahlt haben, da"s sie abgeschlossen sind und auch die "ubrigen Voraussetzungen des Schachtelungsprinzips \ref{Lemma_Schachtelungsprinzip} erf"ullen, gilt dann $\bigcap_{n\geq 2} U_n\neq\LM$ und wegen $U_1\supset U_2\supset\ldots$ auch: 
\begin{alignat}{1}
&\bigcap_{n\geq 1} U_n\neq\LM.\label{Formel_BCT_4}
\end{alignat}

\medskip Wegen (\ref{Formel_BCT_3}) gilt $U_1\supset U_2\supset \INT[U_2]\supset U_3\supset\INT[U_3]\supset\ldots$ $(*)$ und somit $\bigcap_{n\geq 1}U_n\EQ[$(*)$]\bigcap_{n\geq 1}\INT[U_n]\underset{\TT{\tiny{(\ref{Formel_BCT_3})}}}{\subset} U_1\cap \bigcap_n D_n$. Wegen (\ref{Formel_BCT_4}) gilt dann 
\begin{alignat*}{1}
&U_1\cap \bigcap_n D_n\neq\LM.
\end{alignat*}
Es folgt also die Behauptung.

\bigskip $(ii)$ Sei $X$ ein lokal kompakter Raum. Sei $(D_n)_n$ eine beliebige Folge offener dichter
Teilmengen in $X$.

\point{Es ist zu zeigen} $\FA[O\subset X\TT{ offen, nicht-leer}](O\cap \bigcap_n D_n\neq\LM)$.

\medskip Sei $O_1\subset X$ offen und nicht-leer. Dann ist $O_1\cap D_1\neq\LM$ (da $D_1$ dicht ist). Sei etwa $x\in O_1\cap D_1$, d.h. $O_1\cap D_1$ ist eine offene Umgebung von $x$. Dann gibt es eine kompakte Umgebung $K_1$ von $x$ (d.h. eine kompakte Menge $K_1$, die eine offene Umgebung $U$ von $x$ enh"alt) mit $K_1\subset O_1\cap D_1$ (da $X$ lokal kompakt ist und Lemma \ref{Lemma_Lokal_Kompakt}). Sei
\begin{alignat*}{1}
&O_2:=\INT[K_1].\\
\intertext{Dann gilt:}
&\ABS[O_2]\EQ[Def.]\ABS[{\INT[K_1]}]\subset\ABS[K_1]\EQ[$K_1$ abg.]K_1\subset O_1\cap D_1
\end{alignat*}
denn $K_1$ ist als kompakte Teilmenge eines $T_2$-Raumes abgeschlossen [\siehe Lemma \ref{Lemma_Abgeschlossen}]. Dabei ist $\ABS[O_2]$ als abgeschlossene Teilmenge des kompakten Raumes $K_1$ kompakt [\siehe Lemma \ref{Lemma_Quasikompakt}], d.h. $O_2$ ist relativ kompakt. 

\medskip Wendet man die gleiche "Uberlegung auf $O_2$ und $D_2$ an, erh"alt man eine relativ kompakte offene  Umgebung $O_3$ von $x$ mit $\ABS[O_3]\subset O_2\cap D_2$.

\smallskip Induktiv fortfahrend ergibt sich so eine absteigende Folge $(O_n)_n$ relativ kompakter offener (nicht-leerer) Umgebungen $O_n$ von $x$ mit
\begin{alignat}{1}
&\ABS[O_n]\subset O_{n-1}\cap D_{n-1}\label{Formel_BCT_1}
\end{alignat}
f"ur $n\geq 2$. Da f"ur $n\geq 2$ alle $\ABS[O_n]$ als abgeschlossene Teilmengen des kompakten Raumes $\ABS[O_2]$ aufgefa"st werden k"onnen und je endlich viele Glieder der Folge $(O_n)_n$ einen nicht leeren Durchschnitt haben (die Folge ist absteigend), gilt nach Lemma \ref{Lemma_Char_Quasikompakt}, da"s $\bigcap_{n\geq 2} \ABS[O_n]\neq\LM$ und wegen $\ABS[O_1]\supset \ABS[O_2]\supset\ldots$ auch: 
\begin{alignat}{1}
&\bigcap_{n\geq 1} \ABS[O_n]\neq\LM.\label{Formel_BCT_2}
\end{alignat}

\medskip Wegen (\ref{Formel_BCT_1}) gilt $\ABS[O_1]\supset O_1\supset\ABS[O_2]\supset O_2\supset\ldots$ $(*)$ und somit $\bigcap_{n\geq 1}\ABS[O_n]\EQ[$(*)$]\bigcap_{n\geq 1}O_n\underset{\TT{\tiny{(\ref{Formel_BCT_1})}}}{\subset} O_1\cap \bigcap_n D_n$. Wegen (\ref{Formel_BCT_2}) gilt dann 
\begin{alignat*}{1}
&O_1\cap \bigcap_n D_n\neq\LM.
\end{alignat*}
Es folgt also die Behauptung.
\end{beweis}

Der Baire-Raum $\omega^\omega$ (mit der Topologie, die durch die diskreten Topologien auf den $\omega$ induziert wird)
ist ein Beispiel f"ur einen nicht lokal kompakten Baire'schen Raum [\siehe Kapitel \ref{Kap_Baire_Raum}]. Auch sind Baire'sche R"aume nicht unbedingt 
Hausdorff'sch (Beispiel: der topologische Raum $X$ mit mindestens zwei Punkten und der indiskreten Topologie 
$\MN{X,\LM}$ ist Baire'sch jedoch nicht Hausdorff'sch).
\KOM[Bsp. f"ur\\Baire'sch\\und nicht\\polnisch?]

\section{Die Baire-Eigenschaft}\label{Kap_BaireEig}
Die anschauliche Vorstellung einer \glqq fast offenen\grqq\ Menge wird durch den Begriff der \eemph{Baire-Eigenschaft} beschrieben:

\begin{definition}\header{Baire-Eigenschaft}\label{Def_BaireEigenschaft}\\
Sei $X$ ein Baire'scher Raum. Eine Menge $A\subset X$ habe die  \emph{Baire-Eigenschaft}, falls es eine offene Menge
$U\subset X$ gibt so, da"s $A\sd U$ mager ist. Dazu "aquivalent ist, da"s $A = U\sd P$ f"ur eine offene Menge $U$ 
und eine
magere Menge $P$ in $X$ ist (setze $P = A\sd U$).
\end{definition}

Jede magere Menge hat die Baire-Eigenschaft, da sie in diesem Sinne fast leer und die
leere Menge offen ist. Offene Mengen haben die Baire-Eigenschaft, da die leere Menge nirgends dicht 
und somit mager ist.

\bigskip Ein System $\SYS{A}$ von Teilmengen einer Menge $X$ nennt man eine \emph{$\sigma$-Algebra}, falls $\SYS{A}$ die Menge $X$ enth"alt und abgeschlossen ist gegen"uber Komplementbildung und abz"ahlbaren Vereinigungen, d.h. wenn gilt:
\begin{ITEMS}[arabic)]
\item $X\in\SYS{A}$,
\item $A\in\SYS{A}\Rightarrow\C[A]\in\SYS{A}$,
\item $\SYS{S}\subset\SYS{A}$ abz"ahlbar $\Rightarrow (\bigcup_{S\in\SYS{S}} S )\in\SYS{A}$.
\end{ITEMS}

Damit gilt nun folgendes Lemma:

\begin{lemma}
Die Teilmengen eines topologischen Raumes $X$, die die Baire-Eigenschaft besitzen, bilden eine $\sigma$-Algebra.
\end{lemma}
\begin{beweis}
\point{Es ist zu zeigen}$X$ selber hat die Baire-Eigenschaft und die Teilmengen von $X$, die die Baire-Eigenschaft
besitzen, sind abgeschlossen gegen"uber Komplementbildung und abz"ahlbaren Vereinigungen.

\medskip $X$ ist offen und hat somit die Baire-Eigenschaft.

\medskip Da beliebige (insbesondere abz"ahlbare) Vereinigungen offener Mengen offen sind, "ubertr"agt sich die Abgeschlossenheit magerer Mengen gegen"uber abz"ahlbaren Vereinigungen direkt auf Mengen mit der Baire-Eigenschaft. 

\medskip Es bleibt die Abschlu"seigenschaft gegen"uber der Komplementbildung zu zeigen. Habe $A\subset X$ die Baire-Eigenschaft\gs etwa $A\sd U$ mager f"ur ein $U\subset X$ offen. Dann ist auch $A\sd \ABS[U]$ mager:

Zun"achst folgt aus $U$ offen, da"s $\ABS[U]\o U$ nirgends dicht ist, da $\INT[{(\ABS[U]\o U)}]=\LM$ gilt und $\C[{(\ABS[U]\o U)}]$ offen aus $U$ offen folgt. Sei etwa $A\sd \ABS[U]=(A\sd U) \cup B$ mit $B\subset \ABS[U]\o U$ nirgends dicht. Dann ist $A\sd \ABS[U]$ mager. \begin{comment}Damit sieht man wegen $A\sd \ABS[U] = (A\o \ABS[U]) \cup  (\ABS[U]\o A)$ und $A\o \ABS[U] \subset A\o U$ mager sowie $\ABS[U]\o A = (U\o A) \cup {\left[ (\ABS[U]\o U) \o A \right]}$ mager (wegen $U\o A$ mager und $\left[ (\ABS[U]\o U) \o A \right] \subset (\ABS[U]\o U)$ mager), da"s $A\sd \ABS[U]$ mager ist.\end{comment} 

Somit hat auch $\C[A]$ die Baire-Eigenschaft, da $(\C[A])\sd (\C[{\ABS[U]}]) = \C[({\C[ {(A\sd \ABS[U])} ]})] = A\sd \ABS[U]$ mager ist und $\C[ {\ABS[U]} ]$ offen.
\end{beweis}

In einem Baire'schen Raum bilden die Mengen mit der Baire-Eigenschaft die kleinste $\sigma$-Algebra, die die
Borelmengen und die mageren Mengen enth"alt. Die Baire-Eigenschaft l"a"st sich f"ur den Baireraum auch gut spieltheoretisch charakterisieren: in Satz \ref{Satz_CharBaireEig} geben wir ein hinreichendes Kriterium f"ur die Baire-Eigenschaft projektiver Mengen an.

\section{Der Baire-Raum}\label{Kap_Baire_Raum}\KOM[1. RxR \\nicht hom"om. \\zu R\\ABER\\$\BR\times\BR$ \\h. zu \\$\BR$\\2. BorelIso-Zeugs]

Der Baireraum $\omega^\omega$ und der Cantorraum $2^\omega$ werden in der Mengenlehre anstelle der reellen Zahlen untersucht\footnote{\emph{Baire- und Cantorraum:} Benannt nach Ren\'e Baire (1874--1932) bzw. Georg Cantor (1845--1918).}. Die Zuordnung einer reellen Zahl $x\in[0,1]$ zur Folge $b\in n^\omega$ ihrer Nachkommastellen einer g-adischen Darstellung ist i.a. nicht eindeutig. So bezeichnen bei Dezimaldarstellung die Folgen $0,0999\ldots$ und $0,1000\ldots$ dieselbe reelle Zahl. Diese Uneindeutigkeiten treten weder im Baire- noch im Cantor-Raum auf. Diese Folgenr"aume haben auf der anderen Seite ganz "ahnliche Eigenschaften wie die reellen Zahlen: Wie man eine reelle Zahl durch Angabe von immer mehr Nachkommastellen immer genauer bestimmt, so werden auch Elemente $f$ der Folgenr"aume durch Angabe von immer l"angeren Anfangsst"ucken immer besser approximiert. Dabei werden im Gegensatz zur p-adischen Darstellung einer reellen Zahlen $x\in\RZ$ zwei Folgen $(f(0),f(1),\ldots)$ und $(g(0),g(1),\ldots)$ aus Baire- oder Cantor-Raum wirklich nur dann miteinander identifiziert, wenn sie Glied f"ur Glied "ubereinstimmen, d.h. $f(n)=g(n)$ f"ur alle $n<\omega$.

\bigskip Zwei Folgen $f$ und $g$ aus dem Baire- oder Cantor-Raum sind intuitiv "ahnlich, wenn sie in einem langen Anfangsst"uck "ubereinstimmen. Auf diese Weise erh"alt man einen Begriff von \emph{$f$ liegt nahe bei $g$} f"ur zwei Folgen aus dem Baire- oder Cantor-Raum, wie man ihn auch von den reellen Zahlen her kennt. Die mathematische Pr"azisierung dieser Vorstellung liefert dann R"aume, die den reellen Zahlen sehr "ahneln, und die zudem f"ur die Untersuchungen in der Mengenlehre besser geeignet sind.

\bigskip Dar"uberhinaus "ubertragen sich aber viele Resultate, die man im Baire-Raum erh"alt nicht nur auf die reellen Zahlen, sondern zus"atzlich auch auf polnische R"aume. In der deskriptiven Mengenlehre zeigt man, da"s je zwei "uberabz"ahlbare polnische R"aume (insbesondere der Baireraum \siehe Satz \ref{BR->polnisch} und jeder andere polnische Raum) zueinander \eemph{borel-isomorph\footnote{\emph{borel-isomorph:} Zwei topologische R"aume hei"sen \eemph{borel-isomorph}, falls es zwischen ihnen einen \eemph{Borel-Isomorphismus} gibt. Eine Abbildung $f:\PFEIL{X}{}{Y}$ zwischen zwei topologischen R"aumen hei"st \eemph{Borel-Isomorphismus}, falls $f$ bijektiv ist und sowohl $f$ als auch die Umkehrabbildung $f^{-1}$ \eemph{borel-me"sbar} sind. Eine Abbildung $f:\PFEIL{X}{}{Y}$ zwischen zwei topologischen R"aumen hei"st \eemph{borel-me"sbar}, falls f"ur eine Borelmenge $A\subset Y$ ihr Urbild $f^{-1}A$ stets wieder eine Borel-Menge in $X$ ist (\siehe etwa \cite[11]{Kechris:1995}).}} sind (\siehe etwa \cite{Kechris:1995}). 

\bigskip Ein Modell $V$ von $\ZF$ enth"alt mit $\omega$ stets nat"urliche Zahlen im Sinne der Peano'schen Axiome und somit auch den Baire-Raum $\omega^\omega$. Dieser ist hom"oomorph zum Raum der Irrationalzahlen aufgefa"st als Teilraum der reellen Zahlen (\siehe etwa \cite[4.6.2]{Alexandroff:1994}).

\bigskip Sei $\omega^{<\omega}$ die Menge der endlichen Sequenzen nat"urlicher Zahlen:
\[
\omega^{<\omega} \DEF \bigcup_{n< \omega}{\omega^n}
\]
und $\omega_*^{<\omega}$ die Menge der nicht-leeren endlichen Sequenzen nat"urlicher Zahlen. Sei $\lng(u)$ die L"ange der endlichen Sequenz $u\in \omega^{<\omega}$:
\[
\lng(u):=m\TT{ f"ur }u=(u(0),\ldots,u(m-1)).
\]
F"ur $u, v\in \omega^{<\omega}$ bedeute $u\prec v$, da"s $u$ ein Anfangsst"uck von $v$ ist:
\[
u\prec v:\Leftrightarrow u=(v(0),\ldots,v(\lng(u)-1))
\]
wobei $\lng(u)\leq\lng(v)$ sei\gs d.h. $u$ ist ein echtes Anfangsst"uck von $v$ oder $u$ ist gleich $v$.\footnote{A. Kechris definiert in \cite[S. 192]{Kechris:1977} davon abweichend $u\prec v:\Leftrightarrow \lng(u)\gneq\lng(v)\UND\FA[i\leq\lng(v)](u(i)=v(i))$\gs d.h. $v$ ist ein echtes Anfangsst"uck von $u$.}\label{DefPrec}
Ist $f\in\omega^\omega$ und $u\in\omega^{<\omega}$ ein endliches Anfangsst"uck der unendlichen Sequenz $f$, schreibt man ebenfalls $u\prec f$. 

\bigskip F"ur eine unendliche Sequenz $f\in\omega^\omega$ bezeichne $f\bs[n]$ das endliche Anfangsst"uck von $f$ der L"ange $n$:
\[
f\bs[n]:=(f(0),\ldots,f(n-1)).
\]
\begin{definition}\header{Baireraum}\label{BR}\\
Der  \emph{Baire-Raum} ist der Raum $\BR := \omega^\omega$ aller unendlichen Folgen nat"urlicher Zahlen mit der Topologie, die durch die offenen Basismengen $O_s \DEF \MNG{f\in{\BR}}{s\prec f}$ mit $s\in\omega^{<\omega}$ definiert ist. Die Mengen $O_s$ nennt man auch  \emph{elementare offene Teilmengen}. Da der Baireraum $\BR=O_{()}$ selber eine elementare offene Menge ist, und das System der elementaren offenen Teilmengen des Baireraumes schnittstabil ist, bilden die elementaren offenen Teilmengen tats"achlich die Basis einer Topologie.
\end{definition}

Das ist gerade die Produkttopologie der diskreten Topologien auf den $\omega$. 

\bigskip Nachfolgend werden nun einige elementare topologische Eigenschaften des Baireraumes dargestellt. Au"serdem wird der Baireraum als polnischer Raum der Klasse der Baire'schen R"aume zugeordnet [\siehe S"atze \ref{BR->polnisch} und \ref{Satz_Baireraum_Bairesch}].

\subsection*{Topologische Eigenschaften}

Singleton-Mengen des Baireraumes sind nirgends dicht: Das Komplement einer Singleton-Menge $\{f\}\subset\BR$ l"a"st sich schreiben als $\C[\{f\}]=\bigcup_{u\not\prec f}O_u$. Also ist $\C[\{f\}]$ offen. Da jede offene Basismenge $O_v$ des Baireraumes ein $g\neq f$ enth"alt, gilt $g\in O_u$ f"ur ein $u\not\prec f$ und somit $g\in\C[\{f\}]$. Also ist $\C[\{f\}]$ nicht nur offen sondern auch dicht. Die Singleton-Menge $\{f\}$ ist also nirgends dicht.

\bigskip Der Baireraum selber ist eine offene Basismenge, da er sich schreiben l"a"st als ${\BR}=O_{()}$. Hingegen l"a"st sich die leere Menge nicht als offene Basismenge auffassen: Angenommen die leere Menge ist eine offene Basismenge. Gelte etwa $\LM=O_u$ f"ur ein $u\in\omega^\omega$. Dann folgt $\neg\EX[f\in\omega^\omega](f\succ u)$\gs Widerspruch.

\bigskip F"ur die Arbeit im Baireraum spielen die offenen Basismengen eine besondere Rolle. Sie sind nicht nur offen sondern immer auch abgeschlossen und haben einige weitere besondere Eigenschaften:

\begin{lemma}\label{Lemma_Baireraum_OffeneBasismengen}
Je zwei offene Basismengen in $\BR$ liegen entweder ineinander oder haben einen leeren Durschnitt - genauer: 
Seien $O_s$ und $O_t$ offene Basismengen in $\BR$ mit $O_s\cap O_t \neq \LM$, dann gilt
\begin{alignat*}{1}
&O_s\subset O_t \TT{ oder } O_t\subset O_s.
\end{alignat*}
\end{lemma}
\begin{beweis}
Seien $O_s$ und $O_t$ beliebige offene Basismengen des Baireraumes und gelte $O_s\cap O_t\neq\LM$. Sei etwa $f\in O_s\cap O_t$. Dann gilt 
\begin{alignat*}{1}
&s\prec f \TT{ und }t\prec f\\ 
\Rightarrow &s\prec t \TT{ oder }t\prec s\\
\Rightarrow &O_s\supset O_t \TT{ oder }O_t\supset O_s 
\end{alignat*}
und somit folgt die Behauptung.
\end{beweis}

F"ur zwei beliebige offene Basismengen gilt stets
\[
O_s\subset O_t\Leftrightarrow s\succ t.
\] 
Liegen zwei offene Basismengen $O_s$ und $O_t$ nicht ineinander, schreibt man auch
\[
O_s\perp O_t\TT{ oder }s\perp t.
\]

\begin{lemma}\label{BR:offeneBasismengen->abgeschlossen}
Die offenen Basismengen des Baireraumes sind auch abgeschlossen.
\end{lemma}
\begin{beweis}
Das Komplement einer offene Basismenge $O_u\subset\BR$ ist offen, da es sich schriebenl"a"st als
$\bigcup_{s\perp u}O_s$. Daher ist $O_u$ abgeschlossen.
\end{beweis}

\begin{lemma}\label{BR:offeneBasismengen->nichtkompakt}
Die nicht-leeren offenen Teilmengen des Baireraumes sind nicht kompakt.
\end{lemma}
\begin{beweis}
Eine offene Basismenge $O_s$ des Baireraumes ist nicht kompakt, da es zur offenen "Uberdeckung $(O_{s\kon (m)})_{m<\omega}$ keine endliche Teil"uberdeckung gibt. Somit sind auch beliebige offene Teilmengen $U$ des Baireraumes nicht kompakt, da es f"ur eine beliebige offene Basismenge $O_s\subset U$ zur offenen "Uberdeckung $(O_{s\kon (m)})_{m<\omega}\cup U\o O_s$ keine endliche Teil"uberdeckung von $U$ gibt.
\end{beweis}

Insbesondere ist wegen ${\BR}=O_{()}$ der Baireraum nicht kompakt. Au"serdem sieht man, da"s kompakte Mengen $A\subset\BR$ nirgends dicht sind, da $\FA[O_u](O_u\not\subset A)$ gelten mu"s ($O_u\subset A$ abgeschlossen und $A$ kompakt w"urde sonst nach Lemma \ref{Lemma_Quasikompakt} folgen: $O_u$ kompakt\gs im Widerspruch zu Lemma \ref{BR:offeneBasismengen->nichtkompakt}). 

\bigskip Eine gegen"uber der Kompaktheit schw"achere Anforderung stellt man an einen Lindel"of-Raum. Ein topologischer Raum $X$ hei"st \emph{Lindel"of-Raum} oder  \emph{Lindel"of'sch} genau dann, wenn es zu jeder offenen "Uberdeckung $\bigcup_{i\in I}{U_i}$ von $X$ eine abz"ahlbare Teil"uberdeckung $\bigcup_{i<\omega}{U_i}$ gibt. Beispiele sind R"aume mit abz"ahlbarer Basis:

\begin{satz}\header{Lindel"of}\label{Lindel"of}\\
Sei $X$ ein topologischer Raum mit abz"ahlbarer Basis und $U\subset X$. Dann ist $U$ mit der Relativtopologie ein
Lindel"of-Raum. Insbesondere ist jeder topologische Raum mit abz"ahlbarer Basis ein Lindel"of-Raum.
\end{satz}
\begin{beweis}
\KOM[(Preu"s\\S.40)]
Sei das System $\SYS{B}$ eine abz"ahlbare Basis der Topologie von $X$. Sei weiter $(O_i)_{i\in I}$ eine beliebige offene "Uberdeckung von $X$. 

\point{Es ist zu zeigen} Es gibt eine abz"ahlbare Teil"uberdeckung von $(O_i)_{i\in I}$.

\medskip Jedes $O_i$ l"a"st sich darstellen als Vereinigung von Mengen aus $\SYS{B}$. Bestehe $\SYS{B}^*$ aus allen offenen Basismengen, die zu Darstellung von den $O_i$ ben"otigt werden. Dann ist $\SYS{B}^*$ abz"ahlbar, da $\SYS{B}^*\subset \SYS{B}$. Zu jedem $B^*\in\SYS{B}^*$ l"a"st sich nun ein $i_{B^*}\in I$ so ausw"ahlen, da"s $B^*\subset O_{i_{B^*}}$. Dann ist $(O_j)_{j\in J}$ mit $J:=\MNG{i_{B^*}}{B^*\in\SYS{B}^*}\subset I$ eine abz"ahlbare Teil"uberdeckung von $(O_i)_{i\in I}$.
\end{beweis}

\begin{lemma}\label{BR->2.Axiom}
Der Baire-Raum hat eine abz"ahlbare Basis.
\end{lemma}
\begin{beweis}
Da $\omega^{<\omega}$ als abz"ahlbare Vereinigung abz"ahlbarer Mengen $\omega^n$ mit $n\in \omega$ abz"ahlbar ist, ist auch die Menge $\MNG{O_u}{u\in\omega^{<\omega}}$ der offenen Basismengen des Baire-Raumes abz"ahlbar.
\end{beweis}

Daraus ergibt sich nun:

\begin{korollar}\label{BR->Lindel"of'sch}
Der Baire-Raum ist ein Lindel"of-Raum.
\end{korollar}
\begin{beweis}
Nach Lemma \ref{BR->2.Axiom} hat der Baire-Raum eine abz"ahlbare Basis und ist somit wegen Satz \ref{Lindel"of} Lindel"off'sch.
\end{beweis}

\begin{lemma}\label{2.Axiom->separabel}
Jeder topologische Raum mit einer abz"ahlbaren Basis ist separabel.
\end{lemma}
\begin{beweis}
\KOM[(Preu"s\\S.39)]
Sei $X$ ein topologischer Raum mit einer abz"ahlbaren Basis. W"ahlt man aus jeder offenen Basismenge $B$ ein $x_B$
(hierzu verwendet man $AC^{\omega}$), so ist $D=\MNG{x_B}{B \TT{ offene Basismenge}}$ abz"ahlbar und jede nicht leere 
offene Menge schneidet sich mit $D$.
\end{beweis}

\begin{satz}\label{BR->polnisch}
Der Baire-Raum ist polnisch.
\end{satz}
\begin{beweis}
\point{$\BR$ ist separabel}Nach Lemma \ref{BR->2.Axiom} und Satz \ref{2.Axiom->separabel} mu"s $\BR$ separabel sein (dies sieht man auch direkt, da die Teilmenge $\{u\kon (0,0,0\ldots)\in{\BR}\mid u\in\omega^\omega\}$ dicht in $\BR$ und abz"ahlbar ist).

\point{$\BR$ ist metrisierbar}
Falls die Abbildung $d:\PFEIL{\omega^\omega\times\omega^\omega}{}{\RZ}$ mit
\begin{equation*}
\dd(f,g):=\sum_{n=0}^\infty\frac{1}{2^{n+1}}\Kronecker(f(n),g(n))\\
\Kronecker(i,j):=
\begin{cases}
0 &\TT{falls }i=j\\
1 &\TT{falls }i\neq j
\end{cases}
\end{equation*}
wohldefiniert ist (d.h. die Reihe $\dd(f,g)$ konvergiert f"ur beliebige $f,g\in\omega^\omega$), dann definiert sie eine Metrik auf $\omega^\omega$. F"ur beliebige $f,g,h\in\omega^\omega$ gilt n"amlich:
\begin{ITEMS}[arabic)]
\item $\dd(f,f)=0$, \\da $\Kronecker(f(n),f(n))=0$ gilt,
\item $\dd(f,g)=\dd(g,f)$, \\da $\Kronecker(f(n),g(n))=\Kronecker(g(n),f(n))$ gilt,
\item $\dd(f,h)\leq \dd(f,g)+\dd(g,h)$, \\da $\Kronecker(f(n),h(n))\leq\Kronecker(f(n),g(n))+\Kronecker(g(n),h(n))$ gilt. 
\end{ITEMS}
Dabei ist bei Punkt $3)$ zu beachten, da"s $\sum_{n=0}^\infty\frac{1}{2^{n+1}}(\Kronecker(f(n),g(n))+\Kronecker(g(n),h(n)))$ absolut konvergiert und daher beliebig umgeordnet werden kann.

\point{{\mdseries $\dd$} ist wohldefiniert} Die geometrische Reihe $\sum_{n=0}^\infty \frac{1}{2^{n+1}}$ hat nur nicht-negative Glieder und konvergiert (gegen $1$). Wegen $0\leq\Kronecker(f(n),g(n))\leq 1$ gilt:
\[
\left|\frac{1}{2^{n+1}}\Kronecker(f(n),g(n))\right|\leq\frac{1}{2^{n+1}}.
\]
Daher bezeichnet man die Reihe $\sum_{n=0}^{\infty} \frac{1}{2^{n+1}}$ auch als eine \emph{Majorante} von $\sum_{n=0}^{\infty} \frac{1}{2^{n+1}}\Kronecker(f(n),g(n))$. Nach dem \emph{Majorantenkriterium} aus der Analysis konvergiert dann $\sum_{n=0}^{\infty} \frac{1}{2^{n+1}}\Kronecker(f(n),g(n))$ absolut und somit auch im gew"ohnlichen Sinne.

\point{Die Metrik {\mdseries $\dd$} erzeugt die Topologie des Baireraumes}
Die Topologie des Baireraumes wird erzeugt von den offenen Basismengen $O_u=\MNG{f\in{\BR}}{u\prec f}$ mit $u\in\omega^{<\omega}$. Es gen"ugt also zu zeigen, da"s sich jede offene Basismenge $O_u$ eine Darstellung als eine Vereinigung offener $\epsilon$-Kugeln $\Kugel(f,\epsilon_f)$ hat:
\[
O_u=\bigcup_{f}\Kugel(f,\epsilon_f).
\]
Sei $O_u$ eine beliebige offene Basismenge und $u=(u(0),\ldots,u(n))$ d.h. $\lng(u)=n+1$. W"ahlen wir
\[
\epsilon\leq1-\sum_0^n\frac{1}{2^{i+1}},
\]
so gilt:
\[
O_u=K:=\bigcup_{f\in O_u}\Kugel(f,\epsilon).
\]
\point{Es gilt $K\supset O_u$} Weil $\FA[f\in O_u](f\in\Kugel(f,\epsilon)\subset K)$.

\point{Es gilt $K\subset O_u$} Sei $g\in K$ beliebig. 
\begin{alignat*}{1}
\Rightarrow 	&\EX[f\in O_u](g\in\Kugel(f,\epsilon))\\
\Rightarrow 	&\EX[f\in O_u](\dd(f,g)<\epsilon\leq 1-\sum_0^n\frac{1}{2^{i+1}})\\
\Rightarrow	  &u\prec g\\
\Rightarrow	  &g\in O_u.
\end{alignat*}

\begin{comment}%Falsch...
W"ahlt man nun ein beliebiges $f$ aus $O_u$ als \glqq Mittelpunkt\grqq\ so gilt f"ur beliebiges $g\in \BR$:
\begin{alignat*}{1}
&g\in O_u\\
\Leftrightarrow &u\prec g\\
\Leftrightarrow &\dd(f,g)<1-\sum_{i=0}^n \frac{1}{2^{i+1}}=:\epsilon,
\end{alignat*}
wegen
\begin{alignat*}{1}
\dd(f,g)&=\sum_{i=0}^\infty \frac{1}{2^{i+1}}\Kronecker(f(i),g(i))\\
&=\underbrace{\sum_{i=0}^\infty \frac{1}{2^{i+1}}}_{=1}-\sum_{f(i)=g(i)} \frac{1}{2^{i+1}}\TT{ und }f(i)=g(i)\TT{ f"ur }i=0,\ldots,n.
\end{alignat*}
Somit gilt:
\begin{alignat}{1}
O_u=\MNG{g\in{\BR}}{\dd(f,g)<\epsilon}=:\Kugel(f,\epsilon).\label{Frml_O_u=K_epsilon}
\end{alignat}
Die $O_u$ entsprechen also den offenen (und da die $O_u$ auch abgeschlossen sind ebenso den abgeschlossenen) \emph{$\epsilon$-Kugeln} $\Kugel(f,\epsilon)$ des Baireraumes mit der Metrik $\dd$. 
\end{comment}
Somit erzeugt $\dd$ die Topologie des Baireraumes.

\point{Der metrische Raum $({\BR},\dd)$ ist vollst"andig}
Es ist zu zeigen, da"s jede Cauchy-Folge in $({\BR},\dd)$ bez"uglich $\dd$ konvergiert. Sei $(f_i)_i<\omega$ eine beliebige Cauchy-Folge in $({\BR},\dd)$\gs das hei"st:
\[
\FA[\epsilon>0]\EX[N<\omega]\FA[i,j\geq N](\dd(f_i,f_j)<\epsilon).
\]
W"ahlt man nun:
\begin{alignat*}{1}
&\epsilon\leq 1-\sum_{i=0}^n\frac{1}{2^{i+1}},
\end{alignat*}
dann gilt 
\begin{alignat*}{1}
&\EX[N<\omega]\FA[i,j\geq N]((f_i(0),\ldots,f_i(n))=(f_j(0),\ldots,f_j(n)))
\end{alignat*}
F"ur immer kleiner werdendes $\epsilon$ erhalten wir auf diese Weise immer l"anger werdende Anfangsst"ucke, die fast allen (allen bis auf endlich vielen) Folgengliedern $f_i$ gemeinsam sind. 

\medskip Diese immer l"anger werdenden jeweils fast allen $f_i$ gemeinsamen endlichen Anfangsst"ucke definieren ein $f\in\BR$, gegen das die Cauchy-Folge konvergiert.
\begin{comment}
: Zur Definition von $f(n)$ w"ahle man $\epsilon:=1-\sum_{i=0}^{n-1}\frac{1}{2^{i+1}}$. Dann gilt (da $(f_i)$ eine Cauchy-Folge ist):
\[
\EX[N_\epsilon<\omega]\FA[i,j\geq N_\epsilon](\dd(f_i,f_j)<1-\sum_{i=0}^{n-1}\frac{1}{2^{i+1}}).
\]
Dann gilt auch
\[
\dd(f_i,f_j)\leq 1-\sum_{i=0}^n\frac{1}{2^{i+1}}
\]
und das ist genau dann der Fall, wenn $(f_i(0),\ldots,f_i(n))=(f_j(0),\ldots,f_j(n))$. Nun definiert man 
\[
f(n):=k<\omega$ mit $\FA[i\geq N_\epsilon](f_i(n)=k).
\]
Dann ist klar, da"s die Cauchy-Folge gegen $f$ konvergiert: F"ur beliebiges $\epsilon>0$ existiert ein $N\in\omega$ so, da"s f"ur alle $i,j>N$ gilt $\dd(f_i,f_j)<\epsilon$. Dies bedeutet gerade, da"s $(f_i(0),\ldots,f_i(n))=(f_j(0),\ldots,f_j(n))=(f(0),\ldots,f(n))$, f"ur alle $n<\omega$ mit $\epsilon\leq 1-\sum_{i=0}^n\frac{1}{2^{i+1}}$. Daher gilt
\[
\FA[\epsilon>0]\EX[N<\omega]\FA[i\geq N](\dd(f_i,f)<\epsilon).
\] 
Das hei"st, da"s die Folge $(f_i)_{i<\omega}$ gegen $f$ konvergiert.
\end{comment}
\end{beweis}

Als metrischer Raum ist der Baire-Raum somit auch Hausdorff'sch: F"ur zwei beliebige Punkte $f,g$ des Baire-Raumes sind $\Kugel(f,\epsilon_f)$ und $\Kugel(g,\epsilon_g)$ mit $\epsilon_f,\epsilon_g<\frac{1}{2}\dd(f,g)$ schnittfremde offene Umgebungen von $f$ bzw. $g$ (angenommen $h\in\Kugel(f,\epsilon_f)\cap\Kugel(g,\epsilon_g)$, dann $\dd(f,h)+\dd(h,g)\leq \epsilon_f+\epsilon_g < \dd(f,g)\leq \dd(f,h)+\dd(h,g)$\gs Widerspruch).

\bigskip Der Baire-Raum ist weder endlich (klar) noch lokal kompakt: Angenommen der Baireraum ist lokal kompakt. Das hei"st: f"ur jeden Punkt des Baireraumes gibt es eine kompakte Teilmenge, die eine offene Basisumgebung dieses Punktes enth"alt [\siehe Definition \ref{Def_Lokal_Quasikompakt}]. Diese offene Basisumgebung ist nach Lemma \ref{BR:offeneBasismengen->abgeschlossen} auch abgeschlossen und somit nach Lemma \ref{Lemma_Quasikompakt} als abgeschlossene Teilmenge einer kompakten Menge selber kompakt. Offene Teilmengen des Baireraumes sind aber niemals kompakt [\siehe Lemma \ref{BR:offeneBasismengen->nichtkompakt}].

\bigskip Wie f"ur alle polnischen R"aume gilt nach dem Bairesch'en Kategorienatz \ref{BKS} f"ur den Baire-Raum:

\begin{satz}\label{Satz_Baireraum_Bairesch}
Der Baire-Raum ist Baire'sch.
\end{satz}
\begin{beweis}
Der Baire-Raum ist nach Satz \ref{BR->polnisch} insbesondere vollst"andig metrisierbar und somit nach Satz \ref{BKS} Baire'sch.
\end{beweis}

\subsection*{B"aume}

Ein wichtiges kombinatorisches Hilfsmittel zur Untersuchung von Folgenr"aumen $X^\omega$ wie etwa des Baireraumes $\omega^\omega$ oder des Cantorraumes $2^\omega$ ist das Konzept von B"aumen. In der deskriptiven Mengenlehre wie auch in der Graphtheorie und etwa in der theoretischen Informatik spielen Konzepte von B"aumen eine grundlegende Rolle. 

\bigskip An dieser Stelle sollen nun B"aume $T\subset X^\omega$ allgemein auf beliebigen nicht-leeren Mengen $X$ definiert werden. In den nachfolgenden Kapiteln gebrauchen wir dieses Baum-Konzept zun"achst nur f"ur den Fall $X=\omega$. Im letzten Kapitel, in dem allgemeinere Folgenr"aume $X^\omega$ vorkommen, werden wir dann von den B"aumen auf beliebigen nicht-leeren Mengen $X$ gebrauch machen. 

\bigskip Sei nun $X^{<\omega}$ die Menge der endlichen Sequenzen und $X_*^{<\omega}$ die Menge der nicht-leeren endlichen Sequenzen von Elementen in $X$. F"ur $u,v\in X^{<\omega}$ und $f\in X^\omega$ seien die Ausdr"ucke $u\prec v$, $\lng(u)$, $u\prec f$ und $f\bs[n]$ analog zu den Definitionen f"ur den Baireraum [\siehe S. \pageref{DefPrec}].\label{DefEndlSeqEtc}

\bigskip F"ur eine beliebige nicht-leere Menge $X$ ist die \emph{Standardtopologie} f"ur den Folgenraum $X^\omega$ analog zu der Topologie des Baireraumes definiert. Die offenen Basismengen sind die Teilmengen der Gestalt
\[
O_u:=\MNG{f\in X^\omega}{u\prec f}.
\]
Analog zum Baireraum gilt f"ur zwei beliebige offene Basismengen stets $O_s\subset O_t\Leftrightarrow s\succ t$. Liegen zwei offene Basismengen $O_s$ und $O_t$ nicht ineinander, schreibt man auch $O_s\perp O_t\TT{ oder }s\perp t$. 

\bigskip Als einen \emph{Baum} auf einer nicht-leeren Menge $X$ bezeichnet man nun eine Teilmenge $T\subset X^{<\omega}$, die abgeschlossen ist bez"uglich Anfangsst"ucken:
\[
u\in T\UND v\prec u\Rightarrow v\in T.
\] 
Die Elemente von $T$ hei"sen auch \emph{Knoten}. Ein \emph{unendlicher Ast} oder \emph{unendlicher Pfad} in $T$ ist eine unendliche Sequenz $f\in X^\omega$ so, da"s $f\bs[n]\in T$ f"ur alle $n<\omega$. 

\bigskip Die zu einem Baum $T$ geh"orige Teilmenge $[T]$ von $X^\omega$ oder auch der \emph{K"orper} von $T$ ist die Menge aller unendlichen Pfade in $T$:
\[
[T]:=\MNG{f\in X^\omega}{\FA[n<\omega](f\bs[n]\in T)}.
\]
Einen Knoten $u\in T$, der keine echte Erweiterung $v\succ u$, $v\neq u$ in $T$ besitzt, nennt man einen \emph{endlichen Pfad} oder auch \emph{Blatt} in $T$. Ein Baum $T$ hei"st \emph{blattlos} oder \emph{beschnitten}, falls er keine endlichen Pfade enth"alt:
\[
\FA[u\in T]\EX[v\in T](u\prec v\UND u\neq v).
\]
Nicht beschnittene B"aume nennt man auch B"aume \emph{mit Bl"attern}.

\bigskip F"ur eine Teilmenge $A$ des Folgenraumes $X^\omega$ sei $T_A\subset\omega^{<\omega}$ die Menge aller endlichen Anfangsst"ucke von Elementen in $A$: 
\begin{equation*}
T_A\DEF\MNG{u\in X^{<\omega}}{\EX[f\in A]\FA{n<\lng(u)} (f(n)=u(n)) }.
\end{equation*}
$T_A$ hei"st der \emph{Baum der Menge $A$}.

\bigskip Durch die Abbildung  $[\ \_ \ ]: \PFEIL{\{\TT{B"aume auf }X\}}{}{\POW[X^{\omega}]}$ geht f"ur nicht beschnittene B"aume die Information der endlichen Pfade verloren. F"ur nicht beschnittene B"aume $T$ ist $T_{[T]}\subset T$ daher eine echte Inklusion.

\bigskip Umgekehrt gilt f"ur die Abbildung $T_{\_}: \PFEIL{\POW[X^\omega]}{}{\{\TT{B"aume auf }X\}}$ und f"ur nicht abgeschlossenen Teilmengen $A\subset X^\omega$, da"s $[T_A]\supset A$ eine echte Inklusion ist. Da $A$ n"amlich nicht abgeschlossen ist, gibt es eine Folge $(g_n)_{n<\omega}$ in $A$, die gegen ein $f\not\in A$ konvergiert. In $[T_A]$ ist nun aber auch $f$ enthalten, da die endlichen Anfangsst"ucke von $f$ gleichzeitig endliche Anfangsst"ucke der $g_n$ sind, die in $T_A$ liegen.

\bigskip Beschr"ankt man nun die Abbildung $[\ \_ \ ]$ auf die Menge der beschnittenen B"aume auf $X$ und die Abbildung $T_{\_}$ auf die abgeschlossenen Teilmengen von $X^\omega$, so ergibt sich eine Bijektion. Die abgeschlossenen Teilmengen des Folgenraumes entsprechen also genau den blattlosen B"aumen:

\begin{lemma}
Die Abbildungen 
\[
\xymatrix{
{\{\TT{blattlose B"aume}\subset X^{<\omega}\}} \ar@<1ex>[r]^(0.55){[\ \_ \ ]} &{\{\TT{abg. Teilm.}\subset X^\omega\}} \ar@<1ex>[l]^(0.45){T_{\_}}
}
\]
sind zueinander invers.
\end{lemma}
\begin{beweis}
Es ist zu zeigen, da"s $T=T_{[T]}$ f"ur beliebige blattlose B"aume $T\subset X^{<\omega}$ und $A=[T_A]$ f"ur beliebige abgeschlossene Teilmengen $A$ von $X^\omega$.

\point{F"ur $T\subset X^{<\omega}$ blattloser Baum gilt $T\supset T_{[T]}$}
Jedes $u\in X^{<\omega}$, das ein endliches Anfangsst"uck eines unendlichen Pfades in $T$ ist, liegt in $T$: Die unendlichen Pfade $f$ in $T$ sind ja gerade dadurch definiert, da"s ihre endlichen Anfangsst"ucke $f\bs[n]$ f"ur alle $n<\omega$ in $T$ liegen.

\point{F"ur $T\subset X^{<\omega}$ blattloser Baum gilt $T\subset T_{[T]}$}
Da $T$ keine endlichen Pfade hat, ist jedes $u\in T$ ein endliches Anfangsst"uck eines unendlichen Pfades $f\in [T]$. Dann gilt aber $u\in T_{[T]}$.

\point{F"ur $A\subset X^\omega$ abgeschlossen gilt $A\supset[T_A]$}
F"ur ein $f\in [T_A]$ gilt $f\bs[n]\in T_A$ f"ur alle $n<\omega$. Also kann $f$ beliebig gut durch unendliche Sequenzen $g_n\in A$ mit $f\bs[n]\prec g_n$ approximiert werden. Anders ausgedr"uckt: es gibt immer eine Folge $(g_n)_n$ in $A$, die gegen $f$ konvergiert. Da $A$ abgeschlossen ist, mu"s dann $f$ in $A$ liegen. 

\point{F"ur $A\subset X^\omega$ abgeschlossen gilt $A\subset[T_A]$}
F"ur ein $f$ aus $A$ enth"alt $T_A$ alle endlichen Anfangsst"ucke $f\bs[n]$ mit $n<\omega$ von $f$. Also mu"s $[T_A]$ auch wieder $f$ enthalten. 
\end{beweis}

Demnach ist eine Menge $A\subset X^\omega$ genau dann abgeschlossen, wenn die Relation $f\in A$ schon durch die endlichen Anfangsst"ucke von $f$ charakterisiert ist, d.h. genau dann, wenn
\[
\FA[i<\omega](f\bs[i]\in T_A)\Rightarrow f\in A.
\] 
Das ist genau dann der Fall, wenn $A=[T]$ f"ur einen Baum $T$ auf $X$ gilt.  

\bigskip Eine aufsteigende Folge endlicher Sequenzen $s_0\prec s_1\prec s_2\prec\ldots$ in $X^{<\omega}$ definiert stets eine Cauchyfolge $f_0\succ s_0, f_1\succ s_1, f_2\succ s_2,\ldots$ in $X^\omega$ die nach Satz \ref{BR->polnisch} gegen ein $f\in X^\omega$ konvergiert. In diesem Fall sagen wir dann ebenfalls von der aufsteigenden Folge $(s_n)_{n<\omega}$, da"s sie gegen $f$ \emph{konvergiert}.

\section{Charakterisierung durch Banach-Mazur Spiele}\label{KapSpielCharMager}

Die Baire-Kategorie soll nun spieltheoretisch charakterisiert werden. Theorem \ref{SatzMagerCharSpiele} charakterisiert die Baire-Kategorie mittels der Banach-Mazur-Spiele $\BMSo{A}$ f"ur beliebige Teilmengen $A$ des Baireraumes. In Theorem \ref{SatzMagerCharSpiele_2} wird dann die Baire-Kategorie mittels der Banach-Mazur-Spiele $\BMSa{B}$ speziell f"ur sogenannte projektive Teilmengen $A=\p(B)$ des Baireraumes charakterisiert\footnote{\emph{Banach-Mazur-Spiele $\BMSo{A}$ und $\BMSa{B}$:} Die $**$-Notation f"ur diese Varianten des Banach-Mazur-Spieles kommt von J. Mycielski (\siehe \cite{Mycielski:1964:1}). Die Notation $_p$ kennzeichnet \emph{projektive} Spiele\gs auch genannt Spiele \emph{mit Zeugen}.}. 

\bigskip An dieser Stelle erscheint es daher angebracht den von nun an immer wieder benutzten Begriff der projektiven Menge im Zusammenhang mit den Boldface-Punktklassen einzuf"uhren. Bei dieser Gelegenheit werden auch gleich die analytischen Mengen und die Lightface-Punktklassen definiert, die sp"ater f"ur die Definierbarkeitsresultate in Kapitel \ref{KapSpielCharSigmaBeschr"ankt} gebraucht werden.

\subsection*{Boldface-Punktklassen}
Mit den Boldface-Punktklassen wird eine Systematik f"ur gewisse Teilmengen polnischer R"aume eingef"uhrt. Die Boldface-Punktklassen teilen sich ein in die Borel-Hierarchie und die projektive Hierarchie (auch Lusin-Hierarchie). \begin{comment}
Die Lightface-Punktklassen (auch Kleene-Punktklassen genannt) f"uhren eine Systematik f"ur gewisse Teilmengen von R"aumen $\BR^n$ mit $n<\omega$ ein. Die Lightface-Punktklassen teilen sich ein in die arithmetische und die analytische Hierarchie.
\end{comment}

\bigskip Die Borel-Hierarchie beginnt mit den offenen Teilmengen und deren Komplementen und wird mittels abz"ahlbarer Vereinigung und Komplementbildung induktiv fortgef"uhrt.  

\bigskip F"ur jedes $n<\omega$ sind die Klassen $\SIGMA{0}{n}$, $\PI{0}{n}$ und $\DELTA{0}{n}$ von Teilmengen polnischer R"aume wie folgt definiert\footnote{\emph{$\F_\sigma$ und $\G_\delta$:} Entsprechend der alten Notation der Borel-Hierarchie schreibt man auch $\G_\delta$ anstatt $\PI{0}{2}$ und $\F_\sigma$ anstatt $\SIGMA{0}{2}$ (\siehe etwa \cite[VII.4]{Kuratowski:1976}). Eine $\G_\delta$-Menge ist also ein abz"ahlbarer Schnitt offener Mengen und eine $\F_\sigma$-Menge eine abz"ahlbare Vereinigung abgeschlossener Mengen.\label{Def_GDeltaFSigma}}:
\begin{alignat*}{1}
\SIGMA{0}{1}{}&=\TT{ Klasse aller offenen Mengen},\\
\PI{1}{0}{}&=\TT{ Klasse aller abgeschlossenen Mengen},\\
\SIGMA{0}{n+1}{}&=\TT{ Klasse der $A=\bigcup_{i<\omega} A_i$ f"ur $\PI{0}{n}{}$-Mengen $A_i$,}\\
\PI{0}{n}{}&=\TT{ Klasse der Komplemente von $\SIGMA{0}{n}{}$-Mengen in $X$},\\
						&=\TT{ Klasse der $A=\bigcap_{i<\omega} A_i$ f"ur $\SIGMA{0}{n}{}$-Mengen},\\
\DELTA{0}{n}{}&=\SIGMA{0}{n}{}\cap\PI{0}{n}{}.
\end{alignat*}
Die Elemente aus $\SIGMA{0}{n}{}$ oder $\PI{0}{n}{}$ hei"sen \emph{Borel-Mengen}. Es ergibt sich eine Hierarchie $\DELTA{0}{n}{}\subset\SIGMA{0}{n}{}\subset\DELTA{0}{n+1}{}$ sowie $\DELTA{0}{n}{}\subset\PI{0}{n}{}\subset\DELTA{0}{n+1}{}$ mit echten Inklusionen (\siehe etwa \cite[11, S.140f]{Jech:2003}). Diese bezeichnet man als die \emph{Borel-Hierarchie}. Die einzelnen Klassen $\SIGMA{0}{n}{}$, $\PI{0}{n}{}$ und $\DELTA{0}{n}{}$ bezeichnet man als \emph{Borel-Punktklassen}. F"ur die Borel-Punktklassen eines polnischen Raumes $X$ schreibt man entsprechend $\SIGMA{0}{n}{}(X)$, $\PI{0}{n}{}(X)$ und $\DELTA{0}{n}{}(X)$.

\bigskip Die Klasse
\[
\TT{\bfseries{\textup{B}}}:=\bigcup_{n<\omega}\SIGMA{0}{n}{}=\bigcup_{n<\omega}\PI{0}{n}{}
\]
ist abgeschlossen gegen"uber abz"ahlbaren Vereinigungen, abz"ahlbaren Schnitten und Komplementbildung. \siehe etwa \cite[1D]{Moschovakis:1980}, \cite[3.12, S. 146f]{Kanamori:2003} oder \cite[11, S. 140f]{Jech:2003} zur Definition der Borel-Hierarchie.

\bigskip Die projektive Hierarchie beginnt mit den analytischen Teilmengen und deren Komplementen und wird mittels Projektionen und Komplementbildung induktiv fortgef"uhrt. Eine Teilmenge $A$ eines polnischen Raumes $X$ hei"st \emph{analytisch}, falls $A$ die Projektion einer abgeschlossenen Teilmenge $B$ von $X\times\BR$ ist: 
\[
A=\p(B):=\MNG{f\in X}{\EX[{g\in{\BR}}]((f,g)\in B)}
\]
f"ur eine abgeschlossene Teilmenge $B\subset X\times\BR$.
Das ist "aquivalent dazu, da"s $A$ das stetige Bild einer Borel-Menge eines polnischen Raumes ist (\siehe etwa \cite[S. 142]{Jech:2003}). Also sind die Borel-Mengen und insbesondere die offenen und die abgeschlossenen Teilmengen eines polnischen Raumes analytisch.

\bigskip F"ur jedes $n<\omega$ sind die Klassen $\SIGMA{1}{n}$, $\PI{1}{n}$ und $\DELTA{1}{n}$ von Teilmengen polnischer R"aume wie folgt definiert:
\begin{alignat*}{1}
\SIGMA{1}{1}{}&=\TT{ Klasse aller analytischen Mengen},\\
\PI{1}{1}{}&=\TT{ Klasse der Komplemente analytischer Mengen},\\
\SIGMA{1}{n+1}{}&=\TT{ Klasse der Projektionen von $\PI{1}{n}{}$-Mengen in $X\times\BR$},\\
\PI{1}{n}{}&=\TT{ Klasse der Komplemente von $\SIGMA{1}{n}{}$-Mengen in $X$},\\
\DELTA{1}{n}{}&=\SIGMA{1}{n}{}\cap\PI{1}{n}{}.
\end{alignat*}
Die Elemente aus $\SIGMA{1}{n}{}$ oder $\PI{1}{n}{}$ hei"sen \emph{projektive Mengen} (oder auch \emph{Lusin-Mengen}). Es ergibt sich eine Hierarchie $\DELTA{1}{n}{}\subset\SIGMA{1}{n}{}\subset\DELTA{1}{n+1}{}$ sowie $\DELTA{1}{n}{}\subset\PI{1}{n}{}\subset\DELTA{1}{n+1}{}$ mit echten Inklusionen (\siehe etwa \cite[11, S.145]{Jech:2003}). Diese bezeichnet man als die \emph{projektive Hierarchie} (oder auch \emph{Lusin-Hierarchie\footnote{\emph{Lusin-Hierarchie}: Von Nikolai Lusin (1883--1950) und Wac\l aw Sierpi\'nski (1882--1969).}}). Die einzelnen Klassen $\SIGMA{1}{n}{}$, $\PI{1}{n}{}$ und $\DELTA{0}{1}{}$ bezeichnet man als \emph{projektive Punktklassen} (oder auch \emph{Lusin-Punktklassen}). F"ur die projektiven Punktklassen eines polnischen Raumes $X$ schreibt man entsprechend $\SIGMA{1}{n}{}(X)$, $\PI{1}{n}{}(X)$ und $\DELTA{1}{n}{}(X)$.
 
\bigskip Eine projektive Menge $A$ aus $\SIGMA[\BR]{1}{n+1}$ ist also die Projektion $\p(B)$ einer Menge $B$ aus $\PI[\BR\times\BR]{1}{n}$:
\begin{alignat*}{1}
&A=\p(B)=\MNG{f\in{\BR}}{\EX[{g\in{\BR}}]((f,g)\in B)}.
\end{alignat*} 

\siehe etwa \cite[1E]{Moschovakis:1980}, \cite[3.12, S. 148f]{Kanamori:2003} oder \cite[11, S. 144f]{Jech:2003} zur Definition der projektiven Hierarchie.

\bigskip Gelegentlich (etwa f"ur Satz \ref{Satz_CharBaireEig} oder die Definierbarkeitsresultate in Kapitel \ref{KapSpielCharSigmaBeschr"ankt}) ben"otigen wir noch das folgende Lemma "uber die Abschlu"seigenschaften der projektiven Punktklassen.

\begin{lemma}\header{Abschlu"seigenschaften}\label{Lemma_Sigma_SchnittStabil}
\begin{ITEMS}
\item Die $\SIGMA{1}{n}{}$ sind abgeschlossen gegen stetige Urbilder, abz"ahlbare Vereinigungen (d.h. $\exists$-Quantifizierung "uber $\omega$), abz"ahlbare Schnitte (d.h. $\forall$-Quantifizierung "uber $\omega$) und stetige Bilder (insbesondere gegen Projektionen, d.h. $\exists$-Quantifizierung "uber polnischen R"aumen).
\item Die $\PI{1}{n}{}$ sind abgeschlossen gegen stetige Urbilder, abz"ahlbare Vereinigungen (d.h. $\exists$-Quantifizierung "uber $\omega$), abz"ahlbare Schnitte (d.h. $\forall$-Quantifizierung "uber $\omega$) und Urbilder von Projektionen (d.h. $\forall$-Quantifizierung "uber polnischen R"aumen).
\begin{comment}
\item Die $\SIGMAlf{1}{n}{}$ und die $\PIlf{1}{n}{}$ sind abgeschlossen gegen abz"ahlbare Vereinigungen (d.h. $\exists$-Quantifizierung "uber $\omega$), abz"ahlbare Schnitte (d.h. $\forall$-Quantifizierung "uber $\omega$).
\end{comment}
\end{ITEMS}
\end{lemma}
\siehe etwa \cite[37.1]{Kechris:1995}.% f"ur $(i), (ii)$ und \cite[3E.2]{Moschovakis:1980} f"ur $(iii)$.

\subsection*{Banach-Mazur-Spiel ohne Zeugen $\BMSo{A}$}
Um eine spieltheoretische Charakterisierung der Baire-Kategorie zu erhalten, definiert man f"ur beliebiges $A\subset\BR$ ein Spiel $\BMSo{A}$, in dem Spieler I genau dann eine Gewinnstrategie hat, wenn $A$ komager in einer nicht-leeren offenen Menge ist, und II eine Gewinnstrategie hat genau dann, wenn $A$ mager ist.

\begin{definition}\header{Banach-Mazur-Spiel ohne Zeugen $\BMSo{A}$}\\
Zu einer beliebigen Menge $A\subset\BR$ definiert man das \emph{Banach-Mazur-Spiel $\BMSo{A}$} wie folgt:
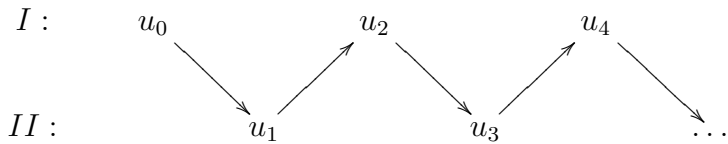
\[
\xymatrix{
I: &u_0 \ar[dr] & 					 	 &u_2 \ar[dr]  &						&u_4 \ar[dr]  &\\
II:&		 				 			&u_1 \ar[ur]	 &						       &u_3 \ar[ur]	&										&\ldots
}
\]
Spieler I und Spieler II spielen abwechselnd\gs Spieler I spielt eine nicht-leere endliche Sequenz nat"urlicher Zahlen $u_0\in\omega_*^{<\omega}$, Spieler II spielt eine nicht-leere endliche Sequenz $u_1\in\omega_*^{<\omega}$, Spieler I spielt eine nicht-leere endliche Sequenz $u_2\in\omega_*^{<\omega}$, Spieler II spielt eine nicht-leere endliche Sequenz $u_3\in\omega_*^{<\omega}$ u.s.w.. Sei nun $f:=u_0\kon u_1\kon u_2\kon\ldots\in\BR$. Man sagt, da"s \emph{Spieler I gewinnt}, falls gilt:
\[
f\in A.
\]
Ansonsten sagt man, da"s \emph{Spieler II gewinnt}. $A$ nennt man dabei auch die \emph{Gewinnmenge}.
\end{definition}

Die Begriffe \emph{(Gewinn-)Strategie f"ur I (bzw. f"ur II)} und \emph{determiniert} sind f"ur die Spiele $\BMSo{A}$ analog zu den Definitionen \ref{Def_Gewinnstrategie_I} (bzw. \ref{Def_Gewinnstrategie_II}) und \ref{Def_determiniert_G_A}  definiert.

\bigskip Mit einer Kodierung der in $\BMSo{A}$ gespielten Sequenzen in die nat"urlichen Zahlen l"a"st sich dieses Spiel auch als ein Spezialfall der Grundvariante [\siehe Definition \ref{Def_SpielGrundvariante}] auffassen und es gilt:

\begin{lemma}\header{Determiniertheit}\label{Lemma_Determiniertheit_G**}\\
Falls $\AD$ gilt, so sind auch die obigen Spiele $\BMSo{A}$ determiniert.
\end{lemma}
\begin{beweis}\KOM[KURZ]
Kodiert man die in $\BMSo{A}$ gespielten endlichen Sequenzen nat"urlicher Zahlen in nat"urliche Zahlen, so gibt es eine Menge $C\subset \omega^\omega$ so, da"s Spieler I (bzw. II) eine Gewinnstrategie in $\BMSo{A}$ hat genau dann, wenn Spieler I (bzw. II) eine Gewinnstrategie in dem Spiel $\textup{G}(C)$ hat.
\end{beweis}

In dem Spiel $\BMSo{A}$ spielen Spieler I und II also abwechselnd endliche Sequenzen $s_i$\gs Spieler I strebt dabei an, da"s $f:=s_0\kon s_1\kon\ldots$ in $A$ liegt und Spieler II, da"s $f$ in $\C[A]$ liegt. Dadurch l"a"st sich nun ein Zusammenhang herstellen zwischen den Gewinnstrategien der Spieler I und II in dem Spiel $\BMSo{A}$ und den topologischen Eigenschaften der Menge $A$ bzgl. der Baire-Kategorie:

\begin{theorem}\header{Charakterisierung durch Spiele ohne Zeugen}\label{SatzMagerCharSpiele}\\
Sei $A\subset\BR$. Dann gilt:
\begin{ITEMS}
\item $I$ hat Gewinnstrategie in $\BMSo{B}$ $\Leftrightarrow$ $A$ ist komager in einer offenen nicht-leeren Menge $G$.
\item $II$ hat Gewinnstrategie in $\BMSo{B}$ $\Leftrightarrow$ $A$ ist mager.
\end{ITEMS}
\end{theorem}

\begin{beweis}$(i)\Leftrightarrow:$ \KOM[Jech 631] Spieler I hat eine Gewinnstrategie f"ur $\BMSa{B}$ genau
dann, wenn es eine endliche Sequenz $s_0\in\omega_*^{<\omega}$ gibt so, da"s Spieler II eine Gewinnstrategie f"ur folgendes Spiel hat: I spielt $s_1\succ s_0$ nicht-leer, II spielt $s_2\succ s_1$ nicht-leer, I spielt $s_3\succ s_2$ nicht-leer etc. und I gewinnt, wenn $f:=s_0\kon s_1\kon s_2\kon\ldots\in O_{s_0}\o A$. Nach Theorem \ref{SatzMagerCharSpiele} $(ii)$ hat II eine Gewinnstrategie in diesem Spiel genau dann, wenn $O_{s_0}\o A$ mager ist\gs also $A$ komager in der offenen Menge $O_{s_0}$ ist.

\bigskip$(ii)\Leftarrow:$ Sei $A$ mager. Dann gilt nach Lemma \ref{Lemma_Komager}:
\[
\C[A]\supset\bigcap_{i<\omega} G_i
\] 
mit offnen dichten Mengen $G_i\subset{\BR}$. Der Schnitt endlich vieler offener dichter Mengen ist stets wieder offen und dicht $(*)$. Dann kann II folgenderma"sen spielen:
\begin{alignat*}{1}
						&\TT{I spielt beliebiges }s_0\\
\underset{}{\Rightarrow}	&O_{s_0}\cap G_0\neq\LM\TT{ offen (wegen $G_0$ offen und dicht)}\\
\Rightarrow	&\EX[O_{s_0\kon s_1}](O_{s_0\kon s_1}\subset O_{s_0}\cap G_0)\\
						&\TT{II spielt }s_1\\
						&\TT{I spielt beliebiges }s_2\\
\underset{(*)}{\Rightarrow}	&O_{s_0\kon\ldots\kon s_2}\cap G_0\cap G_1\neq\LM\TT{ offen}\\
\Rightarrow	&\EX[O_{s_0\kon\ldots\kon s_2\kon s_3}](O_{s_0\kon\ldots\kon s_3}\subset O_{s_0\kon\ldots\kon s_2}\cap G_0\cap G_1)\\
						&\TT{II spielt }s_3\\
						&\ldots\\
						&\TT{I spielt beliebiges }s_{2n}\\
\underset{(*)}{\Rightarrow}	&O_{s_0\kon\ldots\kon s_{2n}}\cap G_0\cap\ldots\cap G_n\neq\LM\TT{ offen}\\
\Rightarrow	&\EX[O_{s_0\kon\ldots\kon s_{2n}\kon s_{2n+1}}](O_{s_0\kon\ldots\kon s_{2n}\kon s_{2n+1}}\subset O_{s_0\kon\ldots\kon s_{2n}}\cap G_0\cap\ldots\cap G_n)\\
						&\TT{II spielt }s_{2n+1}\\
						&\ldots			
\end{alignat*}
Das Spielen von $s_1,s_3,\ldots$ in dieser Weise bedeutet f"ur II eine Gewinnstrategie, da
\[
s_0\kon s_1\kon s_2\kon\ldots\in \bigcap_{i<\omega} G_i\subset \C[A].
\]

\bigskip$(ii)\Rightarrow:$
Habe II eine Gewinnstrategie $\tau$. 

\point{Gute Sequenz} 
Eine Sequenz $u:=(s_0,s_1,\ldots,s_{2n},s_{2n+1})$ mit $s_i\in \omega^{<\omega}_*$ f"ur $i=0,\ldots, 2n+1$ hei"se \emph{gut}, falls f"ur alle $i=0,\ldots ,n$ die $s_{2i+1}$ gem"a"s der Gewinnstrategie $\tau$ gespielt wurden. Die leere Sequenz sei definitionsgem"a"s gut. 

\medskip Sei nun $f\in A$ beliebig. Dann mu"s es eine Sequenz $u=(s_0,s_1,\ldots,s_{2n},s_{2n+1})$ (eventuell $u=()$) geben mit:
\begin{alignat*}{1}
&u\TT{ gut } \UND\\
&s_0\kon\ldots\kon s_{2n+1}\prec f\ \UND\\
&\neg\EX[(s_{2n+2},s_{2n+3})](u\kon(s_{2n+2},s_{2n+3})\TT{ gut }\UND s_0\kon\ldots\kon s_{2n+3}\prec f)
\end{alignat*}
Ansonsten k"onnte Spieler I gewinnen, obwohl Spieler II mit seiner Gewinnstrategie $\tau$ spielt (Widerspruch). Dann gilt $f\in K_u$ mit
\begin{alignat}{2}
K_u &&:= \{f'\in \omega^\omega\mid &s_0\kon\ldots\kon s_{2n+1}\prec f'\UND\notag\\
			  && &\neg\EX[(s_{2n+2},s_{2n+3})]\notag\\
			  && &(u\kon  (s_{2n+2},s_{2n+3})\TT{ gut}\UND s_0\kon\ldots\kon s_{2n+3}\prec f') \}\label{Frml_G**_Ku_1}\\	
				&&= \{f'\in \omega^\omega\mid &s_0\kon\ldots\kon s_{2n+1}\prec f'\UND\notag\\
			  && &\neg\EX[(s_{2n+2},s_{2n+3})]\notag\\
			  &&&(s_{2n+3} \TT{ $\tau$-gesp.}\UND s_0\kon\ldots\kon s_{2n+3}\prec f') \}\label{Frml_G**_Ku_2}
\end{alignat}
wobei $s_{2n+3}$ \glqq$\tau$-gesp.\grqq\ hei"st $s_{2n+3}=\tau(u\kon (s_{2n+2}))$. Damit gilt dann
\begin{alignat*}{1}
&A\subset\bigcup_{u\TT{\scriptsize{ gut}}} K_u.
\end{alignat*}
Die Vereinigung ist abz"ahlbar (weil $\omega^{<\omega}$ abz"ahlbar ist). Um zu zeigen, da"s $A$ mager ist, gen"ugt es nun zu zeigen, da"s die $K_u$ nirgends dicht sind (denn dann ist $A$ als Teilmenge der mageren Menge $\bigcup_u K_u$ mager). Das ist "aquivalent dazu, da"s die $K_u$ abgeschlossen sind und $\INT[{K_u}]=\LM$ gilt (nach Satz \ref{satz_nirgdicht} $(v)$).

\point{Die $K_u$ sind abgeschlossen}
Sei $K_u$ wie oben konstruiert mit $u=(s_0,s_1,\ldots,s_{2n},s_{2n+1})$. Dann l"a"st sich $\C[{K_u}]$ nach (\ref{Frml_G**_Ku_2}) darstellen als
\begin{alignat*}{1}
&\C[K_u]=\C[{O_{(s_0,\ldots,s_{2n+1})}}]\cup\bigcup_{\substack{(s_{2n+2},s_{2n+3})\in\BR\times\BR\\s_{2n+3}\ \tau\TT{\scriptsize{-gespielt}}}}O_{s_0\kon\ldots\kon s_{2n+1}\kon s_{2n+2}\kon s_{2n+3}}.
\end{alignat*}
Also ist $\C[K_u]$ offen und somit $K_u$ abgeschlossen.

\point{Das Innere der $K_u$ ist leer}
Angenommen eine offene Basismenge liegt in $K_u$:
\begin{alignat*}{1}
&O_{s_0\kon\ldots\kon s_{2n+1}\kon s_{2n+2}}\subset K_u.
\end{alignat*}
Dann gilt f"ur $s_{2n+3}:=\tau(u\kon (s_{2n+2}))$, da"s $s_0\kon\ldots\kon s_{2n+3}\prec f'$ f"ur ein $f'\in O_{s_0\kon\ldots\kon s_{2n+1}\kon s_{2n+2}}\subset K_u$ ist. Es gibt also ein $f'\in K_u$ mit:
\begin{alignat*}{1}
&\EX[(s_{2n+2},s_{2n+3})](s_{2n+3}=\tau(u\kon (s_{2n+2})\UND s_0\kon\ldots\kon s_{2n+3}\prec f')
\end{alignat*}
im Widerspruch zur Definition von $K_u$.

\medskip Also liegt keine offene Menge in $K_u$, d.h. $\INT[K_u]=\LM$.

\medskip Damit ist die Hinrichtung von $(ii)$ bewiesen.
\end{beweis}

\subsection*{Banach-Mazur-Spiel mit Zeugen $\BMSa{B}$}\label{Kap_Baire_SpielMitZeugen}
Manchmal ist es n"utzlich, eine gegen"uber der zu charakterisierenden Menge $A$ vereinfachte Gewinnmenge $B$ betrachten zu k"onnen [\siehe etwa Satz \ref{Satz_CharBaireEig}].
Um eine spieltheoretische Charakterisierung der Baire-Kategorie speziell f"ur Teilmengen $A$ des Bairerauems mit $A=\p(B)$ f"ur eine Teilmenge $B\subset\BR\times\BR$ zu erhalten, definiert man daher ein Spiel $\BMSa{B}$ auf $B$ mit der Eigenschaft, da"s $A$ komager in einer offenen nicht-leere Menge ist, falls Spieler I eine Gewinnstrategie hat, und da"s $A$ mager ist, falls II eine Gewinnstrategie hat. 
\begin{comment}
Um eine spieltheoretische Charakterisierung der Baire-Kategorie speziell f"ur Teilmengen $A$ des Bairerauems mit $A=\p(B)$ f"ur eine Teilmenge $B\subset\BR\times\BR$ zu erhalten, definiert man ein Spiel $\BMSa{B}$ auf $B$ mit der Eigenschaft, da"s $A$ komager in einer offenen nicht-leere Menge ist, falls Spieler I eine Gewinnstrategie hat, und da"s $A$ mager ist, falls II eine Gewinnstrategie hat. 
\end{comment}

\begin{definition}\header{Banach-Mazur-Spiel mit Zeugen $\BMSa{B}$}\label{Def_SpielMitZEugen_Baire}\\
Sei $A\subset\BR$ und $B\subset\BR\times\BR$ und es gelte:
\begin{alignat*}{1}
&A=\p(B).%:=\MNG{f\in{\BR}}{\EX[{g\in{\BR}}]((f,g)\in B)}.
\end{alignat*}
Dann definiert man das \emph{Banach-Mazur-Spiel $\BMSa{B}$} wie folgt:
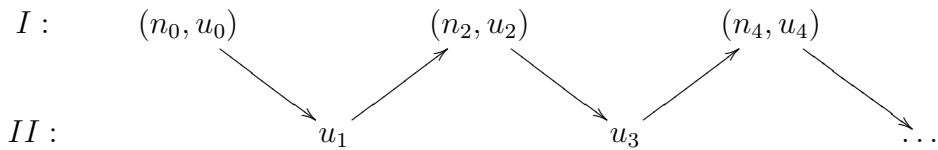
\[
\xymatrix{
I: &(n_0,u_0) \ar[dr] & 					 	 &(n_2,u_2) \ar[dr]  &						&(n_4,u_4) \ar[dr]  &\\
II:&		 				 			&u_1 \ar[ur]	 &						       &u_3 \ar[ur]	&										&\ldots
}
\]
Spieler I und Spieler II spielen abwechselnd\gs Spieler I spielt eine nat"urliche Zahl $n_0\in\omega$ und eine nicht-leere endliche Sequenz nat"urlicher Zahlen $u_0\in\omega_*^{<\omega}$, Spieler II spielt eine nicht-leere endliche Sequenz $u_1\in\omega_*^{<\omega}$, Spieler I spielt eine nat"urliche Zahl $n_2\in\omega$ und eine nicht-leere endliche Sequenz $u_2\in\omega_*^{<\omega}$, Spieler II spielt eine nicht-leere endliche Sequenz $u_3\in\omega_*^{<\omega}$ u.s.w.. Sei nun $g=(n_0,n_2,\ldots)\in\BR$ und $f=u_0\kon u_1\kon u_2\kon\ldots\in\BR$. 
Man sagt, da"s \emph{Spieler I gewinnt}, falls gilt:
\[
(f,g)\in B.
\]
Ansonsten sagt man, da"s \emph{Spieler II gewinnt}. $B$ nennt man dabei auch die \emph{Gewinnmenge}.

Da Spieler I nicht nur endliche Sequenzen sondern auch Elemente aus $\omega$ spielt, die zu einem \eemph{Zeugen} $g\in\omega^\omega$  zusammengesetzt werden\gs wobei Spieler I $(f,g)\in B$ anstrebt\gs geh"ort das Spiel $\BMSa{B}$ zu den Banach-Mazur-Spielen \emph{mit Zeugen}.
\end{definition}

Die Begriffe \emph{(Gewinn-)Strategie f"ur I (bzw. f"ur II)} und \emph{determiniert} sind f"ur die Spiele $\BMSa{B}$ analog zu den Definitionen \ref{Def_Gewinnstrategie_I} (bzw. \ref{Def_Gewinnstrategie_II}) und \ref{Def_determiniert_G_A}  definiert.

\bigskip Mit einer Kodierung der in $\BMSa{B}$ gespielten Sequenzen in die nat"urlichen Zahlen l"a"st sich dieses Spiel auch als ein Spezialfall der Grundvariante [\siehe Definition \ref{Def_SpielGrundvariante}] auffassen und es gilt:

\begin{lemma}\header{Determiniertheit}\label{Lemma_Determiniertheit_G**p}\\
Falls $\AD$ gilt, so sind auch die obigen Spiele $\BMSa{B}$ determiniert.
\end{lemma}
\begin{beweis}
Analog zu Lemma \ref{Lemma_Determiniertheit_G**}.
\end{beweis}

F"ur die Mengen $A$ und $B$ aus obiger Definition gilt dann:

\begin{theorem}\header{Charakterisierung durch Spiele mit Zeugen}\label{SatzMagerCharSpiele_2}\\
Sei $A\subset\BR$ und und $A=\p(B)$ f"ur ein $B\subset\BR\times\BR$.  Dann gilt:
\begin{ITEMS}
\item $I$ hat Gewinnstrategie in $\BMSa{B}$ $\Rightarrow$ $A$ komager in einer offenen nicht-leeren Menge $G$.
\item $II$ hat Gewinnstrategie in $\BMSa{B}$ $\Rightarrow$ $A$ mager.
\end{ITEMS}
\end{theorem}

\begin{beweis}$(i)$ Habe Spieler I eine Gewinnstrategie $\sigma$ f"ur das Spiel $\BMSa{B}$. Wegen $A=\p(B)$ hat Spieler I dann eine Gewinnstrategie $\tilde{\sigma}$ f"ur das Spiel $\BMSo{A}$. Mit Theorem \ref{SatzMagerCharSpiele} $(i)$ folgt dann, da"s $A$ komager in einer offenen nicht-leeren Menge $G$ ist.

\bigskip$(ii)$
Habe II eine Gewinnstrategie $\tau$. 

\point{Gute Sequenz} 
Eine Sequenz $u:=(k_0,s_0,s_1,\ldots,k_{2n},s_{2n},s_{2n+1})$ mit $s_0\in \omega_*^{<\omega}$, $s_i\in \omega^{<\omega}_*$ f"ur $i=1,\ldots, 2n+1$ und $k_{2i}\in\omega$ f"ur $i=0,\ldots,n$ hei"se \emph{gut}, falls f"ur alle $i=0\ldots n$ die $s_{2i+1}$ gem"a"s der Gewinnstrategie $\tau$ gespielt wurden. Die leere Sequenz sei definitionsgem"a"s gut. 

\medskip Sei nun $(f,g)\in B$ beliebig. Dann mu"s es eine Sequenz $u=(\underbrace{k_0,s_0}_I,\underbrace{s_1}_{II},\ldots,\underbrace{k_{2n},s_{2n}}_I,\underbrace{s_{2n+1}}_{II})$ (eventuell $u=()$) geben mit:
\begin{multline*}
u\TT{ gut } \UND\\
(k_0,\ldots,k_{2n})=(g(0),\ldots,g(n))\prec g \UND s_0\kon\ldots\kon s_{2n+1}\prec f\ \UND\\
\neg\EX[(k_{2n+2},s_{2n+2},s_{2n+3})](u\kon(k_{sn+2},s_{2n+2},s_{2n+3})\TT{ gut }\\\UND (k_0,\ldots,k_{2n},k_{2n+2})\prec g\UND s_0\kon\ldots\kon s_{sn+3}\prec f)
\end{multline*}
Ansonsten k"onnte Spieler I gewinnen, obwohl Spieler II mit seiner Gewinnstrategie $\tau$ spielt (Widerspruch). Sei $m:=k_{2n+2}:=g(n+1)$, dann gilt $f\in M_{(u,m)}$ mit
\begin{alignat}{2}
M_{(u,m)} &&:= \{f'\in \omega^\omega\mid &s_0\kon\ldots\kon s_{2n+1}\prec f'\UND\notag\\
			  && &\neg\EX[(s_{2n+2},s_{2n+3})]\notag\\
			  && &(u\kon  (m,s_{2n+2},s_{2n+3})\TT{ gut}\UND s_0\kon\ldots\kon s_{2n+3}\prec f')\}\label{Frml_Gp**_Mum_1}\\	
				&&= \{f'\in \omega^\omega\mid &s_0\kon\ldots\kon s_{2n+1}\prec f'\UND\notag\\
			  && &\neg\EX[(s_{2n+2},s_{2n+3})]\notag\\
			  &&&(s_{2n+3} \TT{ $\tau$-gesp.}\UND s_0\kon\ldots\kon s_{2n+3}\prec f')\}\label{Frml_Gp**_Mum_2}
\end{alignat}
wobei $s_{2n+3}$ \glqq$\tau$-gesp.\grqq\ hei"st $s_{2n+3}=\tau(u\kon  (m, s_{2n+2}))$. Damit gilt dann
\begin{alignat*}{1}
&A\subset\bigcup_{(u,m)}M_{(u,m)}
\end{alignat*}
wobei "uber alle $(u,m)\in\omega^{<\omega}\times\omega$ mit $u$ gut und $m$ wie oben vereinigt wird. Die Vereinigung ist abz"ahlbar (weil $\omega^{<\omega}\times\omega$ abz"ahlbar ist). 

\medskip Da wir die Mengen $M_{(u,m)}$ ganz "ahnlich wie die Menge $K_u$ im Beweis der Hinrichtung von $(ii)$ von Theorem \ref{SatzMagerCharSpiele} definiert haben, k"onnen wir von nun an analog dazu fortfahren.

\medskip Um zu zeigen, da"s $A$ mager ist, gen"ugt es (wie im Beweis von Theorem \ref{SatzMagerCharSpiele}) zu zeigen, da"s die $M_{(u,m)}$ nirgends dicht sind (denn dann ist $A$ als Teilmenge der mageren Menge $\bigcup_{(u,m)}M_{(u,m)}$ mager). Das ist "aquivalent dazu, da"s die $M_{(u,m)}$ abgeschlossen sind und $\INT[{M_{(u,m)}}]=\LM$ gilt (nach Satz \ref{satz_nirgdicht} $(v)$).

\point{Die $M_{(u,m)}$ sind abgeschlossen}
Sei $M_{(u,m)}$ wie oben konstruiert mit $u=(k_0,s_0,s_1,\ldots,k_{2n},s_{2n},s_{2n+1})$ und $m=k_{2n+2}\in\omega$. Dann l"a"st sich $\C[M_{(u,m)}]$ nach (\ref{Frml_Gp**_Mum_2}) darstellen als
\begin{alignat*}{1}
&\C[M_{(u,m)}]=\C[O_{(s_0,\ldots,s_{2n+1})}]\cup\bigcup_{\substack{(s_{2n+2},s_{2n+3})\in\BR\times\BR\\s_{2n+3}\ \tau\TT{\scriptsize{-gespielt}}}}O_{s_0\kon\ldots\kon s_{2n+1}\kon s_{2n+2}\kon s_{2n+3}}.
\end{alignat*}
Also ist $\C[M_{(u,m)}]$ offen und somit $M_{(u,m)}$ abgeschlossen.

\point{Das Innere der $M_{(u,m)}$ ist leer}
Angenommen eine offene Basismenge liegt in $M_{(u,m)}$:
\begin{alignat*}{1}
&O_{s_0\kon\ldots\kon s_{2n+1}\kon s_{2n+2}}\subset M_{(u,m)}.
\end{alignat*}
Dann gilt f"ur $s_{2n+3}:=\tau(u\kon (m,s_{2n+2}))$, da"s $s_0\kon\ldots\kon s_{2n+3}\prec f'$ f"ur ein $f'\in O_{s_0\kon\ldots\kon s_{2n+1}\kon s_{2n+2}}\subset M_{(u,m)}$ ist. Es gibt also ein $f'\in M_{(u,m)}$ mit:
\begin{alignat*}{1}
&\EX[(s_{2n+2},s_{2n+3})](s_{2n+3}=\tau(u\kon (m,s_{2n+2})\UND s_0\kon\ldots\kon s_{2n+3}\prec f')
\end{alignat*}
im Widerspruch zur Definition von $M_{(u,m)}$.

\medskip Also liegt keine offene Menge in $M_{(u,m)}$, d.h. $\INT[M_{(u,m)}]=\LM$.

\medskip Damit ist $(ii)$ bewiesen.
\end{beweis}

\subsection*{Determiniertheit und Baire-Eigenschaft}
Nun lassen sich auch Begriffe wie etwa die Baire-Eigenschaft [\siehe Definition \ref{Def_BaireEigenschaft}], die mittels magerer Mengen definiert wurden, gut spieltheoretisch charakterisieren: Der folgende Satz gibt ein hinreichendes Kriterium f"ur die Baire-Eigenschaft projektiver Mengen an.

\begin{satz}\header{Determiniertheit und Baire-Eigenschaft}\label{Satz_CharBaireEig}\\
Wenn jede $\SIGMA{1}{2n}$-Menge in $\BR$ determiniert ist, dann hat jede $\SIGMA{1}{2n+1}$-Menge
in $\BR$ die Baire-Eigenschaft.
\end{satz}

\begin{beweis} 
Die Idee f"ur den Beweis besteht darin, die die Baire-Eigenschaft definierenden Begriffe mager und komager f"ur $\SIGMA{1}{2n+1}$-Mengen in $\BR$ auf die Determiniertheit von $\SIGMA{1}{2n}$-Mengen in $\BR$ zur"uckzuf"uhren. Dies geschieht mittels des Spieles $\BMSa{B}$. Falls nicht anders erw"ahnt seien hier alle Mengen in $\BR$.

\medskip Sei $A\in\SIGMA[\BR]{1}{2n+1}$ beliebig, $B\in\PI[\BR\times\BR]{1}{2n}$ mit $\PROJ{B}=A$ und es sei jede
$\SIGMA{1}{2n}$-Menge determiniert (und somit auch jede $\PI{1}{2n}$-Menge). F"ur das Spiel $\BMSa{B}$ hat dann 
entweder Spieler I oder Spieler II eine Gewinnstrategie. Nach obigem Lemma ist $A$ dann entweder mager und hat somit
auch die Baire-Eigenschaft, oder $A$ ist komager auf einer offenen nicht-leeren Menge $G$. F"ur den letzteren Fall
bleibt also zu zeigen: $A$ hat die Baire-Eigenschaft, d.h.: $\EX{U}\TT{offen }(A\o U) \cup (U\o A)$ mager.

\medskip Sei $\tilde{G} = \bigcup\MNG{O_s}{A \TT{ komager auf }O_s}$, dann ist $(A\o \tilde{G}) \cup (\tilde{G}\o A)$ mager,
da $(A\o \tilde{G})$ und $(\tilde{G}\o A)$ mager sind:\\
\point{$\tilde{G}\o A$ mager}
\begin{alignat*}{1}
											\tilde{G}\o A		  &=\tilde{G}\cap \C[A]\\
											  								&=(\bigcup_{abzb.}{O_s}{}) \cap \C[A]\\			
											  								&=\bigcup_{abzb.}{{( O_s\cap \C[A] )}}\\
											  								&=\bigcup_{abzb.}{{(\underbrace{ O_s\o A }_{mager})}},
\end{alignat*}
wobei die Ausdr"ucke $O_s\o A$ mager sind (da $A$ komager in $N_s$). Die Vereinigung ist endlich, da $s\in\omega^{<\omega}$ und $\omega^{<\omega}$ abz"ahlbar ist.

\point{$A\o\tilde{G}$ mager} 
Zun"achst mu"s mit $A\in \SIGMA{1}{2n+1}$ auch $(A\o\tilde{G}) \in \SIGMA{1}{2n+1}$ gelten, weil $\tilde{G}$ offen ist ($\tilde{G}$ offen $\Rightarrow$ $\tilde{G}\in \SIGMA{0}{1}$ und $\C[\tilde{G}]\in \PI{0}{1}\subset\SIGMA{1}{2n+1}$, da nun $A\in\SIGMA{1}{2n+1}$ und $\SIGMA{1}{n+1}$ nach Lemma \ref{Lemma_Sigma_SchnittStabil} schnittstabil ist, folgt $A\cap{\C[\tilde{G}]}\in\SIGMA{1}{2n+1}$).
											
\medskip Auf $A\o\tilde{G}$ treffen nun wieder die Voraussetzungen des obigen Lemmas zu\gs allerdings kann $A\o\tilde{G}$ nicht komager auf einer offenen nicht-leeren Menge $U$ sein, weil sonst $U\o(A\o \tilde{G})$ und somit auch $U\o A$ mager w"are\gs laut Definition von $\tilde{G}$ g"alte dann $U\subset\tilde{G}$. Damit lie"se sich folgern $U\o(A\o \tilde{G})$ mager $\Rightarrow U\o(A\o U)(= U)$ mager. $U$ ist aber offen und nicht-leer und kann somit nicht mager sein.

\medskip Also gilt nach obigem Lemma $A\o\tilde{G}$ mager.
\end{beweis}

\chapter{$\sigma$-Kategorie}\label{Kap_SigmaKategorie}
\setcounter{lemma}{0}
In diesem Kapitel sollen unter anderem einige Ergebnisse f"ur das Begriffpaar \eemph{h"ochstens abz"ahlbar} und \eemph{perfekt} analog auf die Begriffe \eemph{$\sigma$-beschr"ankt} und \eemph{superperfekt} "ubertragen werden.

\bigskip Der Begriff eines perfekten Raumes ist ein sehr anschaulicher: Ein Raum hei"st \emph{perfekt}, wenn er keine \eemph{isolierten} Punkte enth"alt. Dabei bezeichnet man einen Punkt eines Raumes als \emph{isoliert}, falls es eine Umgebung des Punktes gibt, die au"ser dem Punkt selber keinen weiteren Punkt des Raumes enth"alt. Ein solcher Punkt ist dann auch im anschaulichen Sinne \glqq isoliert\grqq\ von den "ubrigen Punkten des Raumes. "Aquivalente Formulierungen f"ur einen Punkt $x$ eines topologischen Raumes $X$ sind:
\begin{ITEMS}
\item $x$ ist isoliert.
\item $x$ ist kein H"aufungspunkt von $X$.
\item $\{x\}$ ist offen in $X$.
\item Es gibt keine Folge $(y_n)_n$ in $X$, die gegen $x$ konvergiert (\siehe S. \pageref{Def_Konv_Ansch}, \pageref{Def_Konv}).
\end{ITEMS}
Dabei bezeichnet man einen Punkt $x$ als einen \emph{H"aufungspunkt} von $X$, falls jede Umgebung von $x$ einen Punkt $y\neq x$ aus $X$ enth"alt. Dazu "aquivalent ist, da"s $x$ mit Punkten $y\neq x$ aus $X$ approximiert werden kann. 

\bigskip Eine Teilmenge eines Raumes hei"st \emph{perfekt}\label{Def_perfekteTM}, falls sie ein perfekter Raum und zus"atzlich abgeschlossen ist.

\bigskip Zun"achst werden in Kapitel \ref{Kapitel_HeineBorel} die kompakten Teilmengen des Baireraumes charakterisiert als abgeschlossen und "uberall endlich verzweigt [\siehe Satz \ref{Satz_Heine_Borel_2}]. Diese Charakterisierung der kompakten Mengen dient dann im sp"ateren Kapitel \ref{Kapitel_BMager} als Vorlage f"ur den allgemeineren Begriff der $\BM{B}$-nirgends dichten Teilmengen [\siehe Definition \ref{DefBNirgendsDicht}]. 

\bigskip In Kapitel \ref{Kapitel_SigmaBeschr} wird der Kleinheits-Begriff der \eemph{$\sigma$-Beschr"anktheit} eingef"uhrt. Laut Satz \ref{Satz_SigmaBeschr_SigmaIdeal} bilden die $\sigma$-beschr"ankten Mengen ein $\sigma$-Ideal\gs verhalten sich also wirklich wie \eemph{kleine} Mengen. Bis auf die leere Menge gibt es im Baireraum keine offene Menge, die $\sigma$-beschr"ankt ist [\siehe Satz \ref{Satz_SigmaBeschr_Topologie}]. Das ist eine weitere Gemeinsamkeit der $\sigma$-beschr"ankten mit den mageren Teilmengen des Baireraumes. Die Charakterisierung der $\sigma$-beschr"ankten Mengen in Satz \ref{SatzSigmaBeschr"ankt} $(ii)$ als Teilmengen $\sigma$-kompakter Mengen dient dann im sp"ateren Kapitel \ref{Kapitel_BMager} als Vorlage f"ur den allgemeineren Begriff der $\BM{B}$-mageren Teilmengen [\siehe Definition \ref{DefBMager}]. 

\bigskip Kapitel \ref{Kapitel_Superperfekt} f"uhrt mit dem Begriff der \eemph{superperfekten} (und nicht-leeren) Teilmenge eine Beschreibung f"ur \eemph{relativ gro"se} Teilmengen des Baireraumes ein. Jedoch sind nicht alle superperfekten Teilmengen des Baireraumes \eemph{gro"s} im Sinne der $\sigma$-Kategorie (d.h. Komplemente $\sigma$-beschr"ankter Mengen) [\siehe Beispiel \ref{Lemma_SuperPerf}]. Aus dem Cantor-Bendixson-Analog in Kapitel \ref{Kapitel_CantorBendixson} folgt jedoch, da"s der Abschlusses des Komplementes einer $\sigma$-beschr"ankten Teilmenge stets eine nicht-leere superperfekte Teilmenge enth"alt [\siehe Korollar \ref{KorollarKomplementSigmaBeschr}]. Unter der Annahme von $\AD$ enth"alt schon das Komplement einer $\sigma$-beschr"ankten Teilmenge eine nicht-leere superperfekte Teilmenge [\siehe Korollar \ref{Korollar_Komplement_Sigmabeschr}].

\bigskip In Kapitel \ref{Kapitel_CantorBendixson} wird der Satz von Cantor und Bendixson f"ur h"ochstens abz"ahlbare und perfekte auf die $\sigma$-beschr"ankten und superperfekten Teilmengen des Baireraumes "ubertragen.\footnote{\siehe \cite[2B]{Kechris:1977}.} 

\bigskip Abschlie"send wird in Kapitel \ref{KapSpielCharSigmaBeschr"ankt} das Begriffspaar \eemph{$\sigma$-beschr"ankt} und \eemph{Obermenge einer nicht-leeren superperfekten Menge} dann auf zwei verschiedene Arten spieltheoretisch charakterisiert [\siehe Theoreme \ref{SatzSigmaBeschrCharSpiele} und \ref{SatzSigmaBeschrCharSpiele_2}]\footnote{\siehe \cite[3.3, 3.1]{Kechris:1977} f"ur Theorem \ref{SatzSigmaBeschrCharSpiele} sowie die Aussage des Theorems \ref{SatzSigmaBeschrCharSpiele_2}.}. 

Mittels dieser Charakterisierung (in Form von Korollar \ref{SatzSigmaBeschrCharSpiele_3}) sowie eines Resultates von D.A. Martin (\siehe Kapitel \ref{KapSpielCharSigmaBeschr"ankt} Fu"snote \ref{FN_GaleMartin} bzw. \cite[196ff]{KechrisSolovay:1985} oder \cite[30.10]{Kanamori:2003}) und eines weiteren Resultates von Y.N. Moschovakis (\siehe Kapitel \ref{KapSpielCharSigmaBeschr"ankt} Fu"snote \ref{FN_Moschovakis} bzw. \cite[6E.1.]{Moschovakis:1980}) werden dann noch einige Definierbarkeitsresultate f"ur $\sigma$-Schranken abgeleitet [\siehe S"atze \ref{SatzDefbarkeit_1} und \ref{SatzDefbarkeit_2}]\footnote{\siehe \cite[4]{Kechris:1977}.}.

\bigskip Die Definitionen sowie teilweise auch die Aussagen der Lemmata und S"atze (bis auf das Cantor-Bendixson-Analog \ref{CBA} und die Definierbarkeits-Resultate \ref{SatzDefbarkeit_1} und \ref{SatzDefbarkeit_2} ohne Beweise) finden sich in \eemph{\glqq On a notion of smallness for subsets of the Baire space\grqq\ } von A. Kechris \cite[2--4]{Kechris:1977}. Theorem \ref{SatzSigmaBeschrCharSpiele} sowie die Aussage von Theorem \ref{SatzSigmaBeschrCharSpiele_2} finden sich in \cite[3.1, 3.3]{Kechris:1977}.

\section{Heine-Borel-Analog}\label{Kapitel_HeineBorel}

Im euklidischen Raum $\RZ^n$ ($n\in \NZ$ gr"o"ser $0$) sind die \eemph{kompakten} Teilmengen vollst"andig charakterisiert als \eemph{beschr"ankt und abgeschlossene} Teilmengen. Eine Teilmenge $A$ eines metrischen Raumes $(X,\dd)$ hei"st \emph{beschr"ankt}, falls
\begin{equation}\label{DefBeschr"ankt}
\diam(A)<\infty\\
\TT{wobei }\diam(A)=\sup\MNG{\dd(x,y)\in \RZ}{x,y\in A}.
\end{equation}
Mit dieser Definition ist dann in einem metrischen Raum $(X,\dd)$ z.B. jede abgeschlossene Kugel
$\ABS[\Kugel(x,r)]:=\MNG{y\in X}{\dd(x,y)\leq r}$ beschr"ankt und abgeschlossen (wobei $x\in X$ und $r\in \RZ$ gr"o"ser 
$0$). 

\bigskip Im Baireraum hingegen lassen sich die kompakten Teilmengen nicht als beschr"ankt und abgeschlossen charakterisiern.\label{BR:KompNichtBeschrUAbg} Metrisiert man den Baire'schen Raum $\BR$ mittels der Metrik $d:\PFEIL{\BR}{}{\RZ}$ mit\KOM[(Kechris S. 3)]
\begin{equation*}
\dd(f,g):=\sum_{n=0}^{\infty}\frac{1}{2^{n+1}}\Kronecker(f(n),g(n))\\
\Kronecker(i,j):=
\begin{cases}
0 &\TT{falls }i=j\\
1 &\TT{falls }i\neq j\\
\end{cases}
\end{equation*}
so sind die offenen Basismengen $O_u\subset\BR$ beschr"ankt, da
\[
\FA[f,g\in\omega^\omega](\dd(f,g)\leq 1)
\]
gilt und abgeschlossen [\siehe Lemma \ref{BR:offeneBasismengen->abgeschlossen}]. Nun sind im Baireraum aber offene Mengen nicht kompakt [\siehe Lemma \ref{BR:offeneBasismengen->nichtkompakt}].

Insbesondere sind also die beschr"ankten und abgeschlossenen $O_u$ im Baireraum nicht kompakt. F"ur 
Teilmengen $A$ des Baire-Raumes kann \eemph{kompakt} somit nicht gleichbedeutend mit \eemph{beschr"ankt} (im
Sinne von (\ref{DefBeschr"ankt})) und \eemph{abgeschlossen} sein. Die kompakten Teilmengen des Baire-Raumes $\BR$ k"onnen jedoch mit Hilfe einer partiellen Ordnung als $\leq$-beschr"ankt und abgeschlossen charakterisiert werden. 

\bigskip Eine Teilmenge $A$ des Baireraumes hei"se \emph{$\leq$-beschr"ankt}, falls es ein $f\in\BR$ gibt mit $g\leq f$ f"ur alle $g\in A$\gs dabei sei
\begin{equation*}
g\leq f \TT{ genau dann, wenn f"ur jedes } n< \omega \TT{ gilt, da"s } g(n)\leq f(n).
\end{equation*}
Ein solches $f$ hei"st eine \emph{$\leq$-Schranke f"ur $A$}. Nun erh"alt man folgende Charakterisierung kompakter Teilmengen des Baireraumes:

\begin{satz}\header{Heine-Borel-Analog}\label{Satz_Heine_Borel_1}\\
Eine Menge $A\subset\BR$ ist kompakt genau dann, wenn sie $\leq$-beschr"ankt und abgeschlossen ist.
\end{satz}
\begin{beweis}
\point{$\Rightarrow$} Sei $A\subset\BR$ kompakt. 

\point{$A$ ist abgeschlossen} Da der Baireraum metrisierbar ist [\siehe Satz \ref{BR->polnisch}], ist der Baireraum ein Hausdorff-Raum. In einem Hausdorff-Raum sind kompakte Mengen jedoch immer abgeschlossen [\siehe Lemma \ref{Lemma_Abgeschlossen}]. 

\point{$A$ ist $\leq$-beschr"ankt} Da $A$ kompakt ist, gibt es f"ur die offene "Uberdeckung von $A$ mit allen offenen Basismengen $O_u$ des Baireraumes eine endliche Teil"uberdeckung. Daraus folgt, da"s sich der Baum $T_A$ von $A$ auf "uberall nur endlich oft verzweigt\gs genauer: f"ur jedes $u\in T_A$ ist die Menge $\MNG{m\in \omega}{u\kon (m)\in T_A}$ endlich. Sei $T_A^{n}$ die Menge der endlichen Sequenzen in $T_A$ der L"ange $n\in\omega$. Dann ist $f\in \BR$ mit $f(n)=\max\MNG{m\in\omega}{u\kon (m)\in T_A^{n}\TT{ f"ur ein $u\in T_A$}}$ eine Schranke f"ur $A$.

\point{$\Leftarrow$}Sei umgekehrt $A\subset\BR$ $\leq$-beschr"ankt und abgeschlossen. 

\point{$A$ ist kompakt} Angenommen es gibt eine offene "Uberdeckung von $A$ zu der es keine endliche Teil"uberdeckung gibt. Liege diese offene "Uberdeckung ohne Einschr"ankung in Form offener Basismengen $\FL[O_{u_i}]{i\in I}$ vor (gibt es keine endliche Teil"uberdeckung der beliebigen offenen "Uberdeckung so auch keine der sie erzeugenden offenen Basismengen). Da $A$ $\leq$-beschr"ankt ist, verzweigt sich der Baum $T_A$ auf jeder Stufe nur endlich oft (s.o.). Dann l"a"st sich aus den endlich vielen Sequenzen einer L"ange $n<\omega$ ein  $v_1\in T_A$ heraussuchen so, da"s $O_{v_1}\cap A$ keine endliche Teil"uberdeckung hat (ansonsten g"abe es eine endliche Teil"uberdeckung f"ur ganz $A$). Aus dem so definierten Teilbaum $T_A\cap T_{O_{v_1}}$ l"a"st sich dann  analog eine endliche Sequenz $v_2\succ v_1$ heraussuchen so, da"s $A\cap O_{v_2}$ ebenfalls  keine endliche Teil"uberdeckung besitz und so fort. (Dabei sind die $O_{v_i}$ i.a. nicht unbedingt aus der obigen "Uberdecklung $\FL[O_{u_i}]{i\in I}$.) Sei $f\in\BR$ mit $f\bs[n]\DEF v_n$ f"ur $n\in\omega$ und seien $g_1\in O_{v_1}\cap A, g_2\in O_{v_2}\cap A,\ldots$. Dann konvergiert die Folge $\FL[g_n]{n< \omega}$ gegen $f$. Da $g_n\in A$ f"ur $n< \omega$ und $A$ abgeschlossen ist, mu"s also $f$ in $A$ liegen. Somit liegt $f$ aber schon in einer der offenen Basismengen $O_{u_f}$aus der "Uberdeckung $\FL[O_{u_i}]{i\in I}$ und es gibt ein $N\in\omega$ so, da"s $u_f\prec v_n$ f"ur alle $n\geq N$. F"ur $n>N$ haben die $O_{v_n}\cap A$ dann aber die endliche Teil"uberdeckung $O_{u_f}$\gs Widerspruch.
\end{beweis}

Eine Teilmenge $A$ des Baireraumes ist $\leq$-beschr"ankt genau dann, wenn ihr Baum "uberall (d.h. f"ur jedes $s\in T_A$) endlich verzweigt ist. Eine andere Formulierung des obigen Satzes lautet daher:\KOM[WICHTIG]

\begin{satz}\header{Heine-Borel-Analog 2}\label{Satz_Heine_Borel_2}\\
Eine Menge $A\subset\BR$ ist kompakt genau dann, wenn sie abgeschlossen und "uberall endlich verzweigt ist\gs das hei"st: $A$ ist abgeschlossen und es gilt:
\begin{alignat}{1}
&\FA[s\in T_A]\EX[k<\omega]\FA[t\in\omega^{<\omega}_*](s\kon t\in T_A\Rightarrow t(0)\not > k)\label{Omega_Nirgends_Dicht}
\end{alignat}
wobei $\omega^{<\omega}_*$ gleich $\omega^{<\omega}$ ohne die leere Sequenz sei.\footnote{F"ur die Definition von $\omega^{<\omega}$ und $\omega_*^{<\omega}$ \siehe S. \pageref{DefEndlSeqEtc} oder Anhang S. \pageref{DefMengenR"aumeEtc}.}
\end{satz}

Die Formulierung (\ref{Omega_Nirgends_Dicht}) wird in Kapitel \ref{Kapitel_BMager} als Vorlage f"ur eine  verallgemeinerte Definition der Kompaktheit dienen.

\section{$\sigma$-beschr"ankte Mengen}\label{Kapitel_SigmaBeschr}
Die kleinen Mengen der $\sigma$-Kategorie sind die $\sigma$-beschr"ankten Teilmengen des Baireraumes. Sie werden wie folgt definiert:

\begin{definition}\header{$\sigma$-beschr"ankt}\\
Eine Teilmenge $A\subset \BR$ hei"st \emph{$\sigma$-beschr"ankt} genau dann, wenn es
eine Folge $\FL[f_i]{i< \omega}$ in $\BR$ gibt so, da"s es f"ur jedes $f\in A$ ein $i< \omega$
gibt mit $f\leq f_i$.
\end{definition}
Teilmengen einer $\sigma$-beschr"ankten Menge $A\subset\BR$ sind $\sigma$-beschr"ankt, da eine $\sigma$-Schranke f"ur
$A$ auch eine $\sigma$-Schranke f"ur jede Teilmenge von $A$ ist.
Eine abz"ahlbare Vereinigung $A=\bigcup_{i<\omega}A_i$ $\sigma$-beschr"ankter Mengen $A_i\subset\BR$ ist 
$\sigma$-beschr"ankt, da die abz"ahlbare Vereinigung $(f_{ij})_{i,j<\omega}$ der $\sigma$-Schranken
$(f_{ij})_{j<\omega}$ der $A_i$ abz"ahlbar und somit eine $\sigma$-Schranke f"ur $A$ ist. Es gilt also:

\begin{satz}\label{Satz_SigmaBeschr_SigmaIdeal}
Die $\sigma$-beschr"ankten Teilmengen des Baireraumes bilden ein $\sigma$-Ideal.
\end{satz}

In der Definition eines Baire'schen Raumes wird gefordert, da"s es keine nicht-leeren mageren offenen Teilmengen gibt. In Theorem \ref{BKS} und Satz \ref{Satz_Baireraum_Bairesch} wurde dann gezeigt, da"s der Baireraum ein Baire'scher Raum ist. Hierzu analog gilt der folgende Satz:

\begin{satz}\header{offen, nicht-leer $\Rightarrow$ nicht $\sigma$-beschr"ankt}\label{Satz_SigmaBeschr_Topologie}\\
Im Baireraum gibt es keine nicht-leeren offenen $\sigma$-beschr"ankten Teilmengen.
\end{satz}
\begin{beweis}
Mit einem Diagonalargument sieht man zun"achst, da"s offene Basismengen des Baire-Raumes niemals $\sigma$-beschr"ankt sind:

\medskip Sei $O_u$ eine offene Basismenge des Baire-Raumes (mit $u\in\omega^{<\omega}$). Sei $(f_i)_{i<\omega}$ eine beliebige Folge von Elementen des Baireraumes (eine potentielle $\sigma$-Schranke von $O_u$). Dann liegt $g$ mit
\begin{alignat*}{1}
g		&:=u\kon (n,f_{\lng(u)+1}(\lng(u)+1)+1,f_{\lng(u)+2}(\lng(u)+2)+1,\ldots)\\
n		&:=\max(f_0(\lng(u))+1,\ldots,f_{\lng(u)}(\lng(u))+1)
\end{alignat*}
in $O_u$. Falls nun $i\leq\lng(u)$ ist, so gilt:
\begin{alignat*}{1}
&f_i(\lng(u))<\max(f_0(\lng(u))+1,\ldots,f_{\lng(u)}(\lng(u))+1)=g(\lng(u))\\
\intertext{und falls $i>\lng(u)$, so gilt:}
&f_i(i)<f_{i}(i)+1=g(i).
\end{alignat*}
Insgesamt gilt damit $\FA[i<\omega](g\not\leq f_i)$. Das bedeutet, da"s $(f_i)_i$ keine $\sigma$-Schranke von $O_u$ und somit $O_u$ (da $(f_i)_i$ beliebig gew"ahlt war) nicht $\sigma$-beschr"ankt sein kann.

\medskip Dann gibt es aber keine nicht-leere offene Teilmenge des Baireraumes, die $\sigma$-beschr"ankt ist. Denn Teilmengen $\sigma$-beschr"ankter Mengen sind ja wieder $\sigma$-beschr"ankt und die nicht-leeren offenen Teilmengen des Baireraumes enthalten offene Basismengen (diese sind aber wie eben gesehen nicht $\sigma$-beschr"ankt).
\end{beweis}

Die Begriffe \eemph{$\sigma$-beschr"ankt} und \eemph{mager} beschreiben beide kleine Mengen (im Sinne von Kapitel \ref{Kap_Kleine_Gro"se_Mengen}). Dar"uber hinaus teilen sie nach Satz \ref{Satz_SigmaBeschr_Topologie} die Eigenschaft, da"s ihre $\sigma$-Ideale mit der Topologie des Baireraumes nur die leere Menge gemeinsam haben.

\bigskip Zur besseren Handhabe wird nun der Begriff der $\sigma$-Beschr"anktheit durch einige "aquivalente Formulierungen charakterisiert. Anschlie"send sollen einige Beispiele den Begriff veranschaulichen. 

\bigskip Topologisch lassen sich $\sigma$-beschr"ankte Mengen als Teilmengen $\sigma$-kompakter Mengen beschreiben. Eine \emph{$\sigma$-kompakte} Teilmenge eines topologischen Raumes ist eine abz"ahlbare Vereinigung kompakter Teilmengen.

\bigskip Mittels einer partiellen Ordnung $\leq_{\TEXT{\tiny{fin}}}$ lassen sich die $\sigma$-beschr"ankten Teilmengen des Baireraumes auch als die $\leq_{\TEXT{\tiny{fin}}}$-beschr"ankten Teilmengen beschreiben. Dabei gelte f"ur Elemente $f, g$ des Baireraumes $f\leq_{\TEXT{\tiny{fin}}} g$ genau dann, wenn f"ur alle bis auf endlich viele $n<\omega$ gilt $f(n)\leq g(n)$\gs genauer:
\begin{equation*}
f\leq_{\TEXT{\tiny{fin}}} g \TT{ genau dann, wenn } \EX[n_0<\omega]\FA[n\geq n_0] (f(n)\leq g(n)).
\end{equation*}
Eine Menge $A\subset\BR$ hei"st \emph{$\leq_{\TEXT{\tiny{fin}}}$-beschr"ankt}, falls es ein $f\in\BR$ gibt mit $g\leq_{\TEXT{\tiny{fin}}} f$ 
f"ur alle $g\in A$. Ein solches $f$ hei"st eine \emph{$\leq_{\TEXT{\tiny{fin}}}$-Schranke f"ur $A$}.

\bigskip Damit l"a"st sich nun die $\sigma$-Beschr"anktheit wie folgt charakterisieren:

\begin{satz}\header{$\sigma$-beschr"ankt}\label{SatzSigmaBeschr"ankt}\\
F"ur eine Teilmenge $A$ des Baireraumes ist "aquivalent:
\begin{ITEMS}[equivalences]
\item $A$ ist $\sigma$-beschr"ankt.
\item $A\subset B$ f"ur eine $\sigma$-kompakte Menge $B\subset\BR$.
\item $A$ ist $\leq_{\TEXT{\tiny{fin}}}$-beschr"ankt.
\end{ITEMS}
\end{satz}
\begin{beweis}
$(i)\Rightarrow(ii)$: Sei $A$ etwa durch $(f_i)_{i<\omega}$ $\sigma$-beschr"ankt. Sei $A_i:=\MNG{g\in A}{g\leq f_i}$. Dann sind die $A_i$ $\leq$-beschr"ankt und wegen
\[
\C[A_i]=\bigcup_{\EX[n<\lng(u)](u(n)>f_i(n))}O_u
\]
auch abgeschlossen\gs nach Satz \ref{Satz_Heine_Borel_1} also kompakt. Somit ist $A\subset \bigcup_{i<\omega} A_i$ $\sigma$-kompakt.
\begin{comment}%auch richtig
Sei $A\subset\BR$ $\sigma$-beschr"ankt durch die Folge $\FL[f_i]{i<\omega}$ in $\BR$. Dann l"a"st sich $A$ schreiben als $A=\bigcup_{i<\omega}{B_i\cap A}$ mit $B_i\DEF \MNG{{g\in\BR}}{g\leq f_i}$. Da $B_i\cap A$ beschr"ankt ist (durch $f_i$) ist $\ABS[B_i\cap A]$ $\leq$-beschr"ankt und abgeschlossen also kompakt [\siehe Lemma \ref{Satz_Heine_Borel_1}].
\end{comment}

\begin{comment}
$(ii)\Rightarrow(i)$: Sei umgekehrt $A\subset B$ mit $B$ $\sigma$-kompakt\gs etwa $B=\bigcup_{i<\omega}{B_i}$ mit 
$B_i$ kompakt. Dann sind die $B_i$ (nach obigem Lemma) beschr"ankt und somit $A$ $\sigma$-beschr"ankt. 
\end{comment}

\bigskip $(ii)\Rightarrow(iii)$: Sei $A\subset B$ f"ur eine $\sigma$-kompakte Menge $B\subset\BR$ und etwa 
$B=\bigcup_{i<\omega}B_i$ mit kompakten $B_i\subset\BR$ ($i<\omega$). Dann sind die einzelnen $B_i$ jeweils
$\leq$-beschr"ankt (Heine-Borel-Analog)\gs etwa durch $f_i\in\BR$. Definiert man $f\in\BR$ durch:
$f(i):=\max(f_0(i),\ldots,f_i(i))$ f"ur $i<\omega$, dann gilt $\FA[g\in B](g\leq_{\TEXT{\tiny{fin}}} f)$.

\bigskip $(iii)\Rightarrow(i)$: Sei $A$ $\leq_{\TEXT{\tiny{fin}}}$-beschr"ankt durch ein $f\in\BR$, d.h. es gilt
\begin{alignat*}{1}
\FA[g\in A]\EX[n_0<\omega]\FA[n\geq n_0]{(g(n)\leq f(n))}.			
\end{alignat*}
Mit $A_m:=\MNG{g\in A}{\FA[n\geq m](g(n)\leq f(n))}$ die Menge der $g\in A$, die ab $m<\omega$
durch $f$ majorisiert werden, gilt dann
\begin{alignat*}{1}
&A=\bigcup_{m<\omega}A_m.
\end{alignat*}
Dabei kann das endliche Anfangsst"uck $u:=(g(0),\ldots,g(m-1))$ eines $g\in A_m$ Werte in $\omega^m$\gs also nur 
abz"ahlbar viele\gs annehmen. Dann gilt
\begin{alignat*}{1}
A&=\bigcup_{m<\omega}A_m\\
 &=\bigcup_{m<\omega}\bigcup_{u<\omega^m}(A_m\cap O_u).
\end{alignat*}
Diese Vereinigung ist abz"ahlbar (als abz"ahlbare Vereinigung abz"ahlbar Vereinigungen) und jedes der $A_m\cap O_u$ ist beschr"ankt durch $f^u:=u\kon (f(m),f(m+1),\ldots)$. Insgesamt ist $A$ also $\sigma$-beschr"ankt.
\end{beweis}

Punkt $(ii)$ des Satzes \ref{SatzSigmaBeschr"ankt} wird in Kapitel \ref{Kapitel_BMager} als Vorlage einer verallgemeinerten Definition kleiner Mengen dienen\gs den \emph{$\BM{B}$-nirgends dichten} Mengen f"ur eine vorgegebene Bedingungsmenge $\BM{B}$ [\siehe Kapitel \ref{Kapitel_BMager}].\KOM[WICHTIG]

\bigskip Die $\sigma$-Beschr"anktheit charakterisiert laut Satz \ref{SatzSigmaBeschr"ankt} $(ii)$ genau die Teilmengen des Baireraumes, die Teilmenge einer $\sigma$-kompakten Menge sind. Ein solcher Begriff kann nur in solchen R"aumen sinnvoll sein, die nicht selber $\sigma$-kompakt sind (denn in einem $\sigma$-kompakten Raum sind alle Teilmengen ebenfalls $\sigma$-kompakt). Metrische R"aume $X$, in denen die abgeschlossenen Kugeln $\ABS[\Kugel(x,\epsilon)]$ (mit $x\in X$ und $\epsilon>0$ in $\RZ$) kompakt sind, sind stets $\sigma$-kompakt, da dann gilt $X=\bigcup_{r\geq 1}\ABS[\Kugel(x,r)] $. Dies trifft z.B. auf den Raum der reellen Zahlen zu. Der Umstand, da"s die abgeschlossenen $\epsilon$-Kugeln $\Kugel(f,\epsilon)=\ABS[\Kugel(f,\epsilon)]$ des Baireraumes nicht kompakt sind [\siehe S. \pageref{BR:KompNichtBeschrUAbg}], ist also notwendig daf"ur, da"s der Begriff $\sigma$-beschr"ankt im Baireraum "uberhaupt sinnvoll ist. In einem $\sigma$-kompakten Raum wie den reellen Zahlen w"urde ein Kleinheitsbegriff, der wie die $\sigma$-Beschr"anktheit in Satz \ref{SatzSigmaBeschr"ankt} $(ii)$ charakterisiert ist, eben keinen rechten Sinn ergeben.

\bigskip Um den Begriff der $\sigma$-Kompaktheit noch etwas zu veranschaulichen, sollen nun noch einige Beispiele gegeben werden.

\begin{beispiel}
Der euklidische Raum $\RZ$ ist $\sigma$-kompakt aber nicht kompakt.
\end{beispiel}
\begin{beweis}
Die Menge der reelen Zahlen l"a"st sich schreiben als $\RZ=\bigcup_{n<\omega} [-n,n]$, wobei die $[-n,n]$ f"ur 
$n\in\NZ$ abgeschlossene und beschr"ankte also kompakte Intervalle in $\RZ$ sind.\\
Da $\RZ=\bigcup_{n<\omega}(-n,n)$ eine offene "Uberdeckung von $\RZ$ ist, zu der es keine endliche Teil"uberdeckung
gibt, ist $\RZ$ nicht kompakt.
\end{beweis}

\begin{beispiel}
Jede abz"ahlbare Teilmenge $A\subset\BR$ ist $\sigma$-beschr"ankt.
\end{beispiel}
\begin{beweis}
$A$ bildet eine $\sigma$-Schranke f"ur sich selbst.
\end{beweis}

\begin{beispiel}
Jede Teilmenge $A\subset\BR$, die eine offene Menge $U\subset\BR$ enth"alt, ist nicht $\sigma$-beschr"ankt.
\end{beispiel}
\begin{beweis}
W"are $A$ $\sigma$-beschr"ankt so auch die in $A$ enthaltene offene Menge $U$. Das kann nach Satz \ref{Satz_SigmaBeschr_Topologie} nicht sein.
\end{beweis}

\begin{beispiel}
Lokal quasikompakte R"aume mit abz"ahlbarer Basis sind $\sigma$-kompakt.
\end{beispiel}
\begin{beweis}
Sei der Raum $X$ lokal quasikompakt und habe eine abz"ahlbare Basis. F"ur jedes $x\in X$ gibt es dann [laut Definition \ref{Def_Lokal_Quasikompakt}, S. \pageref{Def_Lokal_Quasikompakt}] eine kompakte Umgebung $K_x$ von $x$, die eine offene Umgebung $U$ von $x$ enh"alt. $U$ enth"alt eine offene Basismenge $O_x$ von $x$. Da $X$ eine abz"ahlbare Basis besitzt, l"a"st sich $X$ dann so schreiben: 
\[
X=\bigcup_{x\in X}O_x
\]
und wegen $O_x\subset K_x$ gilt dann
\[
X=\bigcup_{x\in X}K_x.
\]
Da die offenen Basismengen $O_x$ abz"ahlbar sind, ist diese Vereinigung abz"ahlbar und somit $X$ $\sigma$-kompakt.
\end{beweis}

In den Baire'schen R"aumen, die eine abz"ahlbare Basis besitzen (wie der Baireraum) und zudem noch lokal quasikompakt sind, macht also ebenfalls ein wie in Satz \ref{SatzSigmaBeschr"ankt} $(ii)$ charakterisierter Kleinheitsbegriff keinen Sinn. In Kapitel \ref{Kap_Baire_Raum} hatten wir allerdings bereits gesehen, da"s der Baireraum nicht lokal quasikompakt sein kann.

\section{Superperfekte Mengen}\label{Kapitel_Superperfekt}

Die \eemph{gro"sen} Mengen der $\sigma$-Kategorie sind die Komplemente $\sigma$-beschr"ankter Teilmengen des Baireraumes. Der Begriff der \eemph{superperfekten} (und nicht-leeren) Teilmenge des Bairerauems kennzeichnet \eemph{relativ gro"se} Teilmengen des Baireraumes, die jedoch nicht unbedingt \eemph{gro"s} im Sinne der $\sigma$-Kategorie sein m"ussen, d.h. sie m"ussen  in $\BR$ nicht Komplemente von $\sigma$-beschr"ankten Mengen sein [\siehe Beispiel \ref{Lemma_SuperPerf}].

\begin{definition}\header{superperfekt}\\
Sei $X$ eine beliebige Menge und $T$ ein Baum "uber $X$. Ein $v\in T$ hei"st \emph{$\sigma$-Knoten in $T$}, falls $v$ unendlich viele direkte Erweiterungen in $T$ hat\gs genauer: $\MNG{{m\in X}}{v\kon m \in T}$ ist unendlich.

\medskip Ein Baum $T$ "uber $X$ hei"st \emph{superperfekt} genau dann, wenn es f"ur jedes $u\in T$ eine
Erweiterung $v\succ u$ in $T$ gibt, die ein $\sigma$-Knoten in $T$ ist. Eine solche Erweiterung $v\in T$ 
hei"st eine \emph{$\sigma$-Erweiterung von $u$ in $T$}. Eine Teilmenge $A\subset X^\omega$ hei"st \emph{superperfekt} 
genau dann, wenn $A=[T]$ f"ur einen Baum $T$ auf $X$ und der Baum $T$ superperfekt ist.
\end{definition}

Die leere Menge $\LM$ ist superperfekt, da $T_{\LM}=\MNG{u\in X^{<\omega}}{\EX[f\in \LM](u\prec f)}=\LM$. Jede offene Basismenge $O_u$ des Baireraumes ist superperfekt. 

\begin{beispiel}\header{superperfekt $\not\Rightarrow$ Komplement einer $\sigma$-beschr. Menge}\label{Lemma_SuperPerf}\\
Die offenen Basismengen $O_u$ des Baireraumes sind superperfekt aber nicht Komplement $\sigma$-beschr"ankter Mengen. 
\end{beispiel}
\begin{beweis}
Die offenen Basismengen $O_u\subset\BR$ sind superperfekt, jedoch ist $\C[O_u]$ nicht $\sigma$-beschr"ankt (beispielsweise ist $O_v$ mit $v\perp u$ Teilmenge von $\C[O_u]$ und laut Satz \ref{Satz_SigmaBeschr_Topologie} nicht $\sigma$-beschr"ankt\gs also ist auch $\C[O_u]$ nicht $\sigma$-beschr"ankt).
\end{beweis}

Superperfekte Mengen lassen sich in folgender Weise charakterisieren:

\begin{satz}
Eine Menge $A\subset\BR$ ist superperfekt genau dann, wenn $A$ abgeschlossen ist und f"ur jedes $f\in A$ und jede offene Menge $U\subset \BR$,
die $f$ enth"alt, $U\cap A$ nicht in einer kompakten Menge $B\subset \BR$ enthalten ist.
\end{satz}
\begin{beweis}
$\Rightarrow:$ Sei $A\in \BR$ superperfekt, d.h. $A$ abgeschlossen und $T_A$ superperfekt.
Angenommen es existieren ein $f\in A$ und eine offene Basismenge $O_u\subset \BR$, die $f$ enth"alt, d.h. $u\prec f$ 
(insbesondere ist $u$ dann aus $T_A$), und $(O_u\cap A)\subset B$ f"ur eine kompakte Menge $B\subset \BR$. 
Da $B$ kompakt ist, ist $B$ nach obigem Lemma beschr"ankt. Also ist auch $(O_u\cap A)\subset B$ beschr"ankt\gs etwa
durch $g\in \BR$. Dann kann aber keine Erweiterung $v$ von $u$ in $T_A$ unendlich viele direkte Erweiterungen haben, 
da $g$ sonst keine Schranke f"ur $(O_u\cap A)$ w"are. Das steht aber im Widerspruch dazu, da"s $T_A$ superperfekt ist.
Statt der offenen Basismenge $O_u$ kann man auch von einer beliebigen offenen Menge $U$, die $f$ enth"alt, ausgehen
\gs $f$ mu"s dann ja schon in einer offenen Basismenge $O_u\subset U$ enthalten sein.

\medskip $\Leftarrow:$ Sei umgekehrt $A$ abgeschlossen und gelte, da"s f"ur jedes $f\in A$ und jede offene Menge $U\subset\BR$, 
die $f$ enth"alt, $U\cap A$ nicht in einer kompakten Menge $B\subset \BR$ liegt. Dann mu"s $T_A$ superperfekt sein:\\
Angenommen $T_A$ ist nicht superperfekt, d.h. es gibt ein $u\in T_A$, dessen s"amtliche Erweiterungen $v\succ u$ 
in $T_A$ jeweils wieder nur endlich viele direkte Erweiterungen $v\kon m\in T_A$ haben. Zun"achst ist dann
(wegen $u\in T_A$) $u$ ein Anfangsst"uck eines $f\in A$\gs es gibt also ein $f\in (O_u\cap A)$. Weiter ist $(O_u\cap
 A)$ abgeschlossen (da $O_u$ und $A$ abgeschlossen) und beschr"ankt:\\ \KOM[Argument oben\\ schon mal\\ benutzt\gs evtl. \\eigenes Lemma]
$T_{(O_u\cap A)}$ kann auf jeder Ebene nur endlich oft verzweigen (da alle Erweiterungen $v$ von $u$ laut Annahme nur
endlich viele direkte Erweiterungen in $T_A$ haben). Das hei"st die Mengen $T_{(O_u\cap A)}^{n}$ (alle endlichen
Sequenzen in $T_{(O_u\cap A)}$ der L"ange $n<\omega$) und damit auch die Mengen 
$\MNG{m<\omega}{ u\kon m\in T_{(O_u\cap A)}^{n}}$ (mit $n<\omega$) sind jeweils endlich. Somit l"a"st sich mit
$g(n):=\max\MNG{m\in\omega}{ u\kon m\in T_{(O_u\cap A)}^{n}}$ ($n<\omega$) eine Schranke f"ur $(O_u\cap A)$
definieren.\\
%$M_i := \MNG{m<\omega}{ v\succ u \TT{ in } T_A \TT{ und } \lng(v)=i \TT{ und } v\kon m \in T_A}$ ($i<\omega$) sind %jeweils endlich 
%(da alle Erweiterungen $v$ von $u$ laut Annahme nur endlich viele direkte Erweiterungen in $T_A$ haben).
%Somit l"a"st sich durch $g := (max(M_i))_{i<\omega}$ eine Schranke f"ur $(O_u\cap A)$ definieren.\\
Nach obigem Satz \ref{Satz_Heine_Borel_1} ist $(O_u\cap A)$ dann kompakt (da beschr"ankt und abgeschlossen). Das steht jedoch im 
Widerspruch zur Vorraussetzung, da"s f"ur eine offene Menge $U\subset \BR$, die $f$ enth"alt, $U\cap A$ nicht in einer
kompakten Menge enthalten ist.
\end{beweis}

\begin{comment}
Folgende Definition soll f"ur jedes $u$ in einem superperfekten Baum $T_A$ eine eindeutige $\sigma$-Erweiterung $v\succ u$ mit minimaler L"ange heraussuchen:

\begin{definition}\header{$A$-minimale $\sigma$-Erweiterung}\\
Sei $A\subset\BR$. 
Eine Erweiterung $v\in T_A$ von $u\in T_A$ habe \emph{$A$-minimale Breite}, wenn f"ur alle $v'\in T_A$ mit $v'\succ u$ und $\lng(v')=\lng(v)$ gilt, da"s $v(n_0)<v'(n_0)$ f"ur $n_0:=\min\MNG{n<\lng(v)}{v(n)\neq v'(n)}$. Falls f"ur ein $u\in T_A$ eine $\sigma$-Erweiterung in $T_A$ existiert, so hei"st diejenige mit minimaler L"ange und falls mehrere mit minimaler L"ange exisiteren zus"atzlich mit $A$-minimaler Breite die \emph{$A$-minimale $\sigma$-Erweiterung von $u$}.
\end{definition}

Damit l"a"st sich nun folgender Satz zeigen:
\end{comment}

\begin{satz}\header{Zu $\BR$ hom"oomorphe Teilmengen superperfekter Mengen}\\
Jede nicht-leere superperfekte Menge $A\subset\BR$ enth"alt eine zu $\BR$ hom"oomorphe superperfekte Teilmenge.
\end{satz}
\begin{beweis}
Sei $A$ eine superperfekte Teilmenge des Baireraumes. Es ist zu zeigen, da"s $A$ eine superperfekte (insbesondere abgeschlossene) Teilmenge $B$ enth"alt, die zu $\BR$ hom"oomorph ist.

\medskip Dazu definieren wir zun"achst einen Teilbaum $T^\omega$ von $T_A$ wie folgt:

\medskip Sei $u_0\in T_A$ ein beliebiger $\sigma$-Knoten. Sei $T^0:=\MNG{s\in\omega^{<\omega}}{s\prec u_0}\subset T_A$. Dann enth"alt $T^0$ als einziges Blatt $u_0$, und dieses ist ein $\sigma$-Knoten in $T_A$. Sei der Baum $T^n\subset T_A$ bereits definiert und seien alle Bl"atter in $T^n$ $\sigma$-Knoten in $T_A$, dann erh"alt man den Baum $T^{n+1}\supset T^n$ wie folgt: Wir definieren f"ur jedes Baltt $u\in T^n$ eine gegen Anfangsst"ucke abgeschlossene Menge $A_u\subset T_A$ und setzen dann
\[
T^{n+1}:=(\bigcup_{u\TT{\scriptsize{ Blatt in }} T^n}A_u)\cup T^n.
\]
Dabei erh"alt man die Menge $A_u$ f"ur ein Blatt $u$ in $T^n$ wie folgt:

\medskip Sei $u$ ein Blatt in $T^n$. Da alle Bl"atter in $T^n$ $\sigma$-Knoten in $T_A$ sind, hat $u$ abz"ahlbar viele direkte Erweiterungen $u\kon (i)$ ($i<\omega$) in $T_A$. Da $T_A$ supeperfekt ist, gilt:
\[
\FA[i<\omega]\EX[\sigma(i)\in T_A](\sigma(i)\TT{ $\sigma$-Knoten in $T_A$}\UND u\kon(i)\prec \sigma(i)).
\]
Man w"ahle f"ur jede direkte Erweiterung $u\kon (i)$ von $u$ in $T_A$ \eemph{genau einen} $\sigma$-Knoten $\sigma(i)\in T_A$ mit $u\kon(i)\prec \sigma(i)$ aus. Dann definieren wir $A_u$ als die Menge dieser $\sigma$-Erweiterungen von $u$ in $T_A$ zusammen mit ihren Anfangsst"ucken: 
\[
A_u:=\MNG{s\in \omega^{<\omega}}{\EX[i<\omega](s\prec\sigma(i))}.
\]
Insbesondere enth"alt $T^{n+1}$ so wieder nur Bl"atter, die $\sigma$-Knoten in $T_A$ sind.
 
\medskip Sei nun 
\[
T^\omega:=\bigcup_{n<\omega} T^n.
\]

Mit dieser Konstruktion sind die Bl"atter der Teilb"aume $T^n\subset T^\omega$ (also $\sigma$-Knoten in $T_A$) auch $\sigma$-Knoten in $T^\omega$ und bilden die einzigen Verzweigungen in $T^\omega$.

\medskip Sei 
\[
\sigma_{(0)},\sigma_{(1)},\sigma_{(2)},\ldots\in T^1
\]
eine Abz"ahlung der f"ur $T^1$ ausgew"ahlten $\sigma$-Erweiterungen von $u_0\in T^0$,
\[
\sigma_{(0,0)},\sigma_{(0,1)},\sigma_{(0,2)},\ldots\in T^2
\]
eine Abz"ahlung der f"ur $T^2$ ausgew"ahlten $\sigma$-Erweiterungen von $\sigma_{(0)}\in T^1$,
\[
\sigma_{(1,0)},\sigma_{(1,1)},\sigma_{(1,2)},\ldots\in T^2
\]
eine Abz"ahlung der f"ur $T^2$ ausgew"ahlten $\sigma$-Erweiterungen von $\sigma_{(1)}\in T^1$ 

u.s.w..

\medskip Bildet man nun $f\in\BR$ auf dasjenige $g\in [T^\omega]$ ab, gegen das die Folge $\sigma_{(f(0))},\sigma_{(f(0),f(1))},\sigma_{(f(0),\ldots,f(2))},\ldots$ von $\sigma$-Erweiterungen von $u$ konvergiert ($g$ liegt in $[T^\omega]$, da $[T^\omega]$ abgeschlossen ist), so ergibt sich auf diese Weise eine Abbildung:
\begin{alignat*}{2}
	\varphi: 	&\BR	&\PFEIL{}{}{}	&[T^\omega]\\
					   &f 	&\pfeil{}{}		&\lim_{n<\omega}\sigma_{(f(0),\ldots,f(n))}.
\end{alignat*}
Da wir im Baum $T^\omega$ keinerlei Verzweigungen zwischen den in $T^\omega$ liegenden $\sigma$-Knoten aufgenommen haben gilt: $\FA[n<\omega](\sigma_{(f(0),\ldots,f(n))}\prec\lim_{i<\omega}\sigma_{(f(0),\ldots,f(i))})$ und jedes $g\in [T^\omega]$ definiert eine eindeutige Folge $(\sigma_{(f(0),\ldots,f(n))}\prec g)_{n<\omega}$. Damit erhalten wir die zu $\varphi$ inverse Abbildung:
\begin{alignat*}{2}
	\varphi^{-1}: 	&[T^\omega]					  &\PFEIL{}{}{}	&\BR\\
					    		&g								  	&\pfeil{}{}		&f
\end{alignat*}
und $\varphi$ ist somit bijektiv. Ferner gilt:

\point{$\varphi$ ist stetig}
Sei $u_0$ der $\sigma$-Knoten in $T^{\omega}$ mit $\lng(u_0)$ minimal und $O_v\cap[T^\omega]$ beliebige offene Basismenge in $[T^\omega]$ (mit $v\in\omega^{<\omega}$). Dann mu"s laut Konstruktion von $[T^\omega]$ gelten, da"s $v\prec u_0$ oder $v\succ u_0\UND v\neq u_0$.

\medskip 1) Falls $v\prec u$, so ist $O_v\cap[T^\omega]=[T^\omega]$ und somit $\varphi^{-1}(O_v\cap[T^\omega])=\BR$ offen.

\medskip 2) Falls $v\succ u\UND v\neq u$, dann ist $\varphi^{-1}(O_v\cap[T^\omega])=\bigcup_{\sigma_{(i_1,\ldots,i_n)}\succ v}O_{(i_1,\ldots,i_n)}$ offen.

\medskip Da die Urbilder aller offenen Basismengen in $[T^\omega]$ offen in $\BR$ sind, gilt gleiches f"ur alle offenen Mengen in $[T^\omega]$, $\varphi$ ist also stetig.

\point{$\varphi^{-1}$ ist stetig}
Sei $O_v\subset\BR$ beliebige offene Basismenge. Dann ist $\varphi(O_v)=O_{\sigma_v}\cap [T^\omega]$ und somit offen in $[T^\omega]$:

\point{Es gilt: $\varphi(O_v)\subset O_{\sigma_v}\cap [T^\omega]$}\\ $v\prec f\Rightarrow \sigma_v\prec g$ mit $(\sigma_{(f_0,\ldots,f_n)})_n$ konvergiert gegen $g$.

\point{Es gilt: $\varphi(O_v)\supset O_{\sigma_v}\cap [T^\omega]$}\\ Sei $\sigma_v\prec g$ und $\varphi(f)=g$ f"ur ein $f\in\BR$, d.h. $\sigma_{(f(0))}\prec \sigma_{(f(0),f(1))}\prec \sigma_{(f(0),\ldots,f(2))}\prec \ldots$ konvergiert gegen $g$ f"ur $f\in \BR$. Dann mu"s $v$ ein Anfangsst"uck von $f$ sein und somit auch $f\in\varphi(O_v)$gelten.

\medskip Also sind unter $\varphi$ alle Bilder offener Mengen in $\BR$ offen. $\varphi^{-1}$ ist also stetig.

\medskip Insgesamt ist damit bis hierhin gezeigt, da"s $\BR$ hom"omorph ist zur Teilmenge $[T^\omega]$ der nicht-leeren superperfekten Menge $A\subset\BR$. 

\medskip Die Menge $[T^\omega]$ ist abgeschlossen und nach Konstruktion von $T^\omega$ auch superperfekt.
\end{beweis}

\section{Cantor-Bendixson-Analog}\label{Kapitel_CantorBendixson}
Der Satz von Cantor und Bendixson, der sich auf abz"ahlbare und perfekte Mengen bezieht, l"a"st sich auf
$\sigma$-beschr"ankte und superperfekte Mengen "ubertragen. 
Dazu soll im Beweis des folgenden Satzes eine superperfekte Teilmenge aus einer abgeschlossenen Teilmenge $A$ des Baireraumes konstruiert werden. Dies geschieht durch \glqq schrittweises\grqq\ Aussortieren von Elementen aus dem Baum von $A$. 
%Eben wegen dieses schrittweisen Vorgehens finden dabei die Ordinalzahlen Verwendung.
Die Klasse aller Ordinalzahlen sei mit $\ORD$ bezeichnet.

\begin{satz}\header{Cantor-Bendixson-Analog\footnote{\siehe etwa \cite[4.6]{Jech:2003} f"ur Cantor-Bendixson.}}\label{CBA}\\
Sei $A\subset \BR$ abgeschlossen. Dann kann $A$ eindeutig geschrieben werden als $A=P\cup C$ mit $P$ superperfekt,
$C$ $\sigma$-beschr"ankt und $P\cap C=\LM$.
\end{satz}
\begin{beweis}
Die Idee ist, durch iteratives \glqq Aussortieren\grqq\ mittels einer Funktion $\partial$ aus dem Baum $T$ eine superperfekte Teilmenge von $A$ zu erhalten. Da man aber im vorhinein nicht wei"s, wieviele Aussortierungsschritte $\partial$ daf"ur n"otig sind, definiert man $\partial^\alpha T$ f"ur alle Ordinalzahlen $\alpha$\gs das hei"st anschaulich: f"ur jede \glqq Anzahl\grqq\ Aussortierungsschritte, die man "uberhaupt machen kann: 

\medskip F"ur einen Baum $T$ definiert man die \emph{Ableitung} $\partial T$ von $T$ als die Menge derjenigen Elemente von $T$,
die in $T$ eine $\sigma$-Erweiterung haben:%Erweiterung haben, die wiederum unendlich viele direkte Erweiterungen \eemph{in $T$} hat:
\begin{alignat*}{1}
							\partial T 	&:=\MNG{u\in T}
							{
								\EX[v\in T]( v \TT{ ist $\sigma$-Erweiterung von $u$ in $T$})
							}\\
												&\EQ[Def.] \MNG{u\in T}
							{
								\EX[v\in T]( v\succ u \TT{ und } \MNG{m<\omega}{v\kon (m) \in T}\TT{ unendlich})
							}.			
\end{alignat*}
Mittels transfiniter Rekursion (\siehe etwa \cite[2.15]{Jech:2003}) definiert man:
\begin{alignat*}{1}
							&\partial^0 T = T\\
							&\partial^{\alpha +1}T = \partial\partial^{\alpha}T \TT{ } (\alpha\in\ORD)\\
							&\partial^{\lambda}T = \bigcap_{\alpha<\lambda} \partial^\alpha T \TT{ } (\lambda \TT{ Limes-Ordinalzahl}).
\end{alignat*}
Aus der Definition von $\partial T$ und $\partial^\alpha T$ f"ur beliebige $\alpha\in\ORD$ folgt direkt, da"s f"ur einen Baum $T$ auch $\partial^\alpha T$ stets wieder ein Baum ist. 

\medskip Sei nun $A$ abgeschlossen, dann ist also $A=[T]$ f"ur einen Baum $T$. Bei jedem Inklusions-Schritt \glqq$\supset$\grqq\ in der Kette 
\begin{alignat}{1}
&\partial^0 T\supset \partial^1 T\supset\ldots\supset\partial^{\lambda}T\supset\partial^{\lambda+1}T\supset\ldots\label{Frml_Kette}
\end{alignat}
werden entweder Elemente aussortiert oder die Kette ist an der Stelle konstant, d.h. $\partial^\alpha T = \partial^{\alpha+1}T$\gs dann bleibt sie ab dieser Stelle auch konstant (nach obiger induktiver Definition). Sei $\alpha_0\in \ORD$ minimal mit $\partial^{\alpha_0}T=\partial^{\alpha}T$ f"ur alle Ordinalzahlen $\alpha >\alpha_0$. Dieses $\alpha_0$ mu"s abz"ahlbar sein, da $T$ abz"ahlbar ist und man aus einer abz"ahlbaren Menge nur abz"ahlbar oft Elemente aussortieren kann. Setzt man nun $P:=[\partial^{\alpha_0}T]$ und $C:=A\o P$, so ist $P\cap C=\LM$ und es gilt: $P$ ist superperfekt, $C$ ist $\sigma$-beschr"ankt und $P$ sowie $C$ sind eindeutig durch diese Eigenschaften bestimmt:

\point{$P$ ist superperfekt}Angenommen $P$ ist nicht superperfekt, dann ist der Baum $\partial^{\alpha_0}T$ nicht superperfekt und es gibt ein $u\in \partial^{\alpha_0}T$ zu dem es keine $\sigma$-Erweiterung $v\succ u$ in $\partial^{\alpha_0}T$ gibt. Laut Definition von $\partial$ wird ein solches $u$ aber im n"achsten Schritt $\partial^{\alpha_0}T\supset\partial^{\alpha_0+1}T$ aussortiert, d.h. $\partial^{\alpha_0}T\neq\partial^{\alpha_0+1}T$ im Widerspruch zur Wahl von $\alpha_0$.

\point{$C$ ist $\sigma$-beschr"ankt}
In jedem Schritt der absteigenden Kette (\ref{Frml_Kette}) werden gegen"uber den Obermengen in der Kette h"ochstens abz"ahlbar viele Elemente aussortiert\gs n"amlich solche $s\in\partial^\alpha T$, die keine $\sigma$-Erweiterung in $\partial^\alpha T$ haben. F"ur solche $s$ ist dann $O_s\cap [\partial^\alpha T]$ endlich verzweigt.\KOM[FRAGE:\\SCHRITTE\\EINZELN\\DURCHRECHN]

\medskip Die Menge $C$ l"a"st sich nun wie folgt darstellen:
\begin{alignat}{1}
C&:=A\o P\notag\\
 &=A\cap \C[P]\notag\\
 &=A\cap\bigcup_{\alpha<\alpha_0}\bigcup_{\substack{s\TT{ \scriptsize{ohne}}\\\sigma\TT{\scriptsize{-Erweiterung }}\\\TT{\scriptsize{in }}\partial^\alpha T}}(\underbrace{O_s\cap[\partial^\alpha T]}_{=:A_{(\alpha,s)}})\notag%\label{Frml_P}
\end{alignat}
wobei die Vereinigungen beide abz"ahlbar sind und die $A_{(\alpha,s)}:=O_s\cap[\partial^\alpha T]$ "uberall endlich verzweigt (s.o.) und abgeschlossen sind. Nach Satz \ref{Satz_Heine_Borel_2} sind die $A_{(\alpha,s)}$ damit kompakt. Die Menge $C$ ist also in einer abz"ahlbaren Vereinigung kompakter Mengen enthalten:
\[
C\subset\bigcup_{(\alpha,s)}A_{(\alpha,s)}
\]
und somit nach Satz \ref{SatzSigmaBeschr"ankt} $\sigma$-beschr"ankt.

\point{$P$ und $C$ eindeutig}
Seien $A=P\cup C=P'\cup C'$ zwei Darstellungen von $A$ mit $P, P'$ superperfekt, $C, C'$ $\sigma$-beschr"ankt und jeweils $P$ und $C$ sowie $P'$ und $C'$ schnittfremd. Angenommen $P\subsetneq P'$. Etwa $f\in P'\o P$. Dann mu"s $T_{P'}$ damit schon einen superperfekten Teilbaum $T_{sp}\subset T_{P'}$ enthalten mit $T_{sp}\not\subset T_P$. Demnach mu"s $T_{sp}\subset T_C$ gelten. Das ist ein Widerspruch, da eine $\sigma$-beschr"ankte Menge keine superperfekte Teilmenge enthalten kann. Demnach mu"s $P=P'$ und somit auch $C=C'$ gelten.
\end{beweis}

\begin{korollar}\label{KorollarKomplementSigmaBeschr}
Sei $A$ das Komplement einer $\sigma$-beschr"ankten Menge $B\subset\BR$, dann enth"alt $\ABS[A]$ eine nicht-leere
superperfekte Teilmenge. 
\end{korollar}
\begin{beweis}
Sei $A=\C[B]$ mit $B$ $\sigma$-beschr"ankt. Laut \ref{CBA} ist $\ABS[A]=P\cup C$ mit einer superperfekten Menge 
$P$ und einer $\sigma$-beschr"ankten Menge $C$. Dann kann $P$ nicht leer sein, weil:

Angenommen $P=\LM$. Dann ist $\ABS[A]=C$ $\sigma$-beschr"ankt und somit $\BR$$=C \cup B$ $\sigma$-beschr"ankt\gs
Widerspruch (der Baire-Raum ist nicht $\sigma$-beschr"ankt). Also mu"s $P\neq\LM$ gelten.

Das hei"st $\ABS[A]$ enth"alt eine superperfekte Teilmenge.
\end{beweis}

Unter der Annahme von $\AD$ gilt sogar, da"s $A$ als Komplement einer $\sigma$-beschr"ankten Menge eine nicht-leere superperfekte Teilmenge enth"alt [\siehe unten Korollar \ref{Korollar_Komplement_Sigmabeschr}]. 

Der Satz \ref{CBA} liefert insbesondere, da"s abgeschlossene Mengen entweder \glqq gro"s\grqq\ oder \glqq klein\grqq\ im Sinne der $\sigma$-Kategorie sind:

\begin{korollar}\label{KorollarSigmaBeschrUndSuperperf}
Sei $A\subset\BR$ abgeschlossen und nicht-leer. Dann ist $A$ entweder $\sigma$-beschr"ankt oder enth"alt eine nicht-leere superperfekte Teilmenge.
\end{korollar}
\begin{beweis}
Sei $A\subset\BR$ abgeschlossen und nicht-leer. Dann kann $A$ nach Satz \ref{CBA} geschrieben werden als $A=P\cup C$, wobei $P$ superperfekt und $C$ $\sigma$-beschr"ankt ist, und es gilt:

\point{Fall 1} $P=\LM$. Dann ist $A=C$ $\sigma$-beschr"ankt.

\point{Fall 2} $P\neq\LM$. Dann enth"alt $A$ eine nicht-leere superperfekte Teilmenge (n"amlich $P$).
\end{beweis}

\section{Charakterisierung durch Banach-Mazur-Spiele}\label{KapSpielCharSigmaBeschr"ankt}
In Kapitel \ref{KapSpielCharMager} haben wir die Baire-Kategorie mittels der $\BMSo{}$- und $\BMSa{}$-Spiele beschrieben. In "ahnlicher Weise soll nun die $\sigma$-Kategorie zun"achst f"ur beliebige Teilmengen $A$ des Baireraumes mittels der Spiele ohne Zeugen $\BMSb{A}$ charakterisiert werden [\siehe Theorem \ref{SatzSigmaBeschrCharSpiele}]. Um vereinfachte Gewinnmengen zu erhalten definieren wir zus"atzlich die Spiele mit Zeugen $\BMSc{B}$ und erhalten so eine Charakterisierung der $\sigma$-Kategorie f"ur Projektionen $A=\PROJ{B}\subset{\BR}$ [\siehe Theorem \ref{DefBMSMitZeugen}].

\subsection*{Banach-Mazur-Spiel ohne Zeugen $\BMSb{A}$}
Um eine spieltheoretische Charakterisierung der $\sigma$-Kategorie f"ur beliebige Teilmengen $A$ des Baireraumes zu erhalten, definiert man ein Spiel $\BMSb{A}$ mit der Eigenschaft, da"s Spieler I genau dann eine Gewinnstrategie hat, wenn $A$ eine nicht-leere superperfekte Teilmenge enth"alt und II eine Gewinnstrategie hat genau dann, wenn $A$ $\sigma$-beschr"ankt ist.

\begin{definition}\header{Banach-Mazur-Spiel ohne Zeugen $\BMSb{A}$}\label{DefBMSOhneZeugen}\\
Zu einer beliebigen Menge $A\subset\BR$ definiert man das \emph{Banach-Mazur-Spiel $\BMSb{A}$} wie folgt:
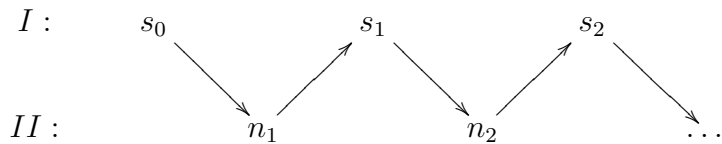
\[
\xymatrix{
I: &s_0 \ar[dr]  & 					 	 &s_1 \ar[dr]  &						&s_2 \ar[dr]\\
II:&		 				 &n_1 \ar[ur]	 &						 &n_2 \ar[ur]	&						&\ldots
}
\]
Spieler I und Spieler II spielen abwechselnd\gs Spieler I spielt eine nicht-leere endliche Sequenz $s_0\in\omega_*^{<\omega}$, Spieler II spielt eine nat"urliche Zahl $n_1<\omega$, Spieler II spielt eine nicht-leere endliche Sequenz $s_1\in\omega_*^{<\omega}$, Spieler I spielt eine nat"urliche Zahl $n_2<\omega$ u.s.w.. Man sagt, da"s \emph{Spieler I gewinnt}, falls gilt:
\begin{ITEMS}[arabic)]
\item $s_0\kon s_1\kon\ldots\in A$,
\item $s_i(0)>k_i$ f"ur $i\geq 1$. 
\end{ITEMS}
Ansonsten sagt man, da"s \emph{Spieler II gewinnt}. $A$ nennt man dabei auch die \emph{Gewinnmenge}.
\end{definition}

Die Begriffe \emph{(Gewinn-)Strategie f"ur I (bzw. f"ur II)} und \emph{determiniert} sind f"ur die Spiele $\BMSb{A}$ analog zu den Definitionen \ref{Def_Gewinnstrategie_I} (bzw. \ref{Def_Gewinnstrategie_II}) und \ref{Def_determiniert_G_A}  definiert.

\bigskip Mit einer Kodierung der in $\BMSb{A}$ gespielten Sequenzen in die nat"urlichen Zahlen l"a"st sich dieses Spiel auch als ein Spezialfall der Grundvariante [\siehe Definition \ref{Def_SpielGrundvariante}] auffassen und analog zu Lemma \ref{Lemma_Determiniertheit_G**p} gilt:

\begin{lemma}\header{Determiniertheit}\label{Lemma_Determiniertheit_G_tilde}\\
Falls $\AD$ gilt, so sind auch die obigen Spiele $\BMSb{A}$ determiniert.
\end{lemma}
\begin{beweis}\KOM[KURZ]
Analog zu Lemma \ref{Lemma_Determiniertheit_G**p}.
\end{beweis}

In dem Spiel $\BMSb{A}$ darf Spieler II also nach jedem Zug von Spieler I endlich viele direkte Nachfolger der bisher 
von I gebildeten endlichen Sequenz $s_0\kon s_1\kon\ldots s_n$ f"ur den weiteren Spielverlauf ausschlie"sen. Dies
wirkt sich genau in gew"unschter Weise auf die Gewinnstrategien der Spieler I und II in Abh"angigkeit von $A$ aus:

\begin{theorem}\header{Charakterisierung durch Spiele ohne Zeugen}\label{SatzSigmaBeschrCharSpiele}\\
Sei $A\subset\BR$. Dann gilt:
\begin{ITEMS}
\item Spieler I hat eine Gewinnstrategie in $\BMSb{A}$ $\Leftrightarrow$ $A$ enth"alt eine nicht-leere superperfekte Teilmenge $P$.
\item Spieler II hat eine Gewinnstrategie in $\BMSb{A}$ $\Leftrightarrow$ $A$ ist $\sigma$-beschr"ankt.
\end{ITEMS}
\end{theorem}
\begin{beweis}
$(i)\Leftarrow$: Enthalte $A$ eine nicht-leere superperfekte Teilmenge $P$. Spieler I spiele als erstes einen beliebigen
$\sigma$-Knoten $s_0\in T_P$ ($T_P$ enth"alt $\sigma$-Knoten, da $P$ superperfekt und nicht-leer ist). Beim n"achsten
Zug spiele I eine $\sigma$-Erweiterung $s_1\succ s_0$ in $T_P$ mit $s_0>n_1$ (das ist m"oglich, da I zuletzt die
$\sigma$-Erweiterung $s_0$ gespielt hat) u.s.w.. Dann liegt $s_0\kon s_1\kon\ldots$ in $A$ (weil $A$ als superperfekte
Menge laut Definition abgeschlossen ist) und $s_i(0)>n_i$ f"ur alle $i<\omega$. Spieler I hat somit eine
 Gewinnstrategie.

\bigskip$(i)\Rightarrow$: Habe Spieler I eine Gewinnstartegie $\sigma$ f"ur $\BMSb{A}$. Sei $T_\sigma:=\MNG{s\in T_A}{s\prec s_0\kon\ldots\kon s_n\TT{ f"ur mit $\sigma$ gespielte }s_0,\ldots,s_n}$. Dann mu"s die Menge $P:=[T_\sigma]$ superperfekt sein: Angenommen $P$ ist nicht superperfekt\gs das hei"st $T_\sigma$ enth"alt ein $s$, zu dem es keine $\sigma$-Erweiterung $s\kon s'$ in $T_\sigma$ gibt. F"ur beliebiges $s\kon s'\in T_\sigma$ kann Spieler II dann
\[
\max\MNG{m<\omega}{s\kon s'\kon(m)\in T_\sigma}
\]
spielen und gewinnt. Das steht jedoch im Widerspruch dazu, da"s $s$ in $T_\sigma$ und $\sigma$ eine Gewinnstrategie von Spieler $I$ ist.

%%%%%%%%%%%%%%
%%%%%%%%%%%%%%
\begin{comment}
(ii)$\Leftarrow$: (mit $\AD$) Sei $A$ $\sigma$-beschr"ankt. Dann enth"alt $A$ keine superperfekte Teilmenge [\siehe Korollar 
\ref{KorollarSigmaBeschrUndSuperperf}]. Somit hat Spieler I laut (i) keine Gewinnstrategie f"ur $\BMSb{A}$. 
\KOM[$\AD$\\vorrauss.?]Mit $\AD$ folgt somit, da"s Spieler II eine
Gewinnstrategie f"ur $\BMSb{A}$ besitzt.
%%%%%%%%%%%%%%
%%%%%%%%%%%%%%
\end{comment}

\bigskip$(ii)\Leftarrow$: Sei $A$ $\sigma$-beschr"ankt. Nach Satz \ref{SatzSigmaBeschr"ankt} gilt dann:
\begin{alignat*}{1}
\EX[f_0\in\BR]\FA[f\in A]\EX[n_0<\omega]\FA[n\geq n_0](f(n)\leq f_0(n)).
\end{alignat*}
Dann gilt f"ur beliebiges Ergebnis $f:=s_0\kon s_1\kon s_2\ldots$ in $A$:
\[
\EX[n_0<\omega]\FA[n\geq n_0](s_n(0)\leq f_0(\sum_{i=0}^{n-1}\lng(s_i))).
\]
Das Spielen von 
\begin{alignat*}{1}
&f_0(\lng(s_0)),\\
&f_0(\lng(s_0)+\lng(s_1)),\\
&f_0(\lng(s_0)+\lng(s_1)+\lng(s_2)),\\
&\ldots
\end{alignat*}
bedeutet f"ur Spieler II also eine Gewinnstrategie.

\bigskip$(ii)\Rightarrow$:  Habe Spieler II eine Gewinnstrategie $\tau$.

\point{Gute Sequenz}Eine Sequenz $u=(s_0,k_1,s_1,k_2,\ldots,s_n,k_{n+1})$ mit $s_i\in\omega_*^{<\omega}$ f"ur $i=0,\ldots,n$ sowie $k_i\in\omega$ f"ur $i=1,\ldots,n+1$ hei"se \emph{gut}, falls gilt:
\begin{alignat*}{1}
&s_i(0)>k_i\TT{ f"ur }i=1,\ldots,n\quad\UND\\
&k_j\TT{ mit $\tau$ gespielt f"ur }j=1,\ldots,n+1.
\end{alignat*}
Die leere Sequenz sei definitionsgem"a"s gut.

\point{Sequenz gut f"ur $f\in A$}$u$ hei"se \emph{gut f"ur} $f\in A$, falls gilt:
\begin{alignat*}{1}
&\EX[s_{n+1}](s_0\kon\ldots\kon s_n\kon s_{n+1}\prec f\UND s_{n+1}(0)>k_{n+1})
\end{alignat*}
Insbesondere ist die leere Sequenz gut f"ur jedes $f\in A$.

\medskip Sei nun $f\in A$ beliebig. Dann mu"s es eine Sequenz $u=(s_0,k_1,s_1,k_2,\ldots,s_n,k_{n+1})$ geben (eventuell $u=()$) mit $u$ ist gut und $u$ ist gut f"ur $f$ und:
\begin{alignat*}{1}
&\neg\EX[(s_{n+1},k_{n+2})\in\BR\times\omega](u\kon(s_{n+1},k_{n+2})\TT{ gut}\UND\TT{ gut f"ur } f)
\end{alignat*}
\begin{comment}
so, da"s:
\begin{multline*}
u\TT{ gut } \UND s_0\kon,\ldots,\kon s_n=(f(0),\ldots,f(m-1))\prec f \UND f(m)>k_{n+1} \UND\\
\neg\EX[(s_{n+1},k_{n+2})\in\BR\times\omega](u\kon(s_{n+1},k_{n+2})\TT{ gut}\UND\\ s_0\kon,\ldots,\kon s_n,s_{n+1}=(f(0),\ldots,f(m'-1))\prec f\UND f(m')>k_{n+2})
\end{multline*}
\end{comment}
Ansonsten k"onnte Spieler I gewinnen, obwohl Spieler II mit seiner Gewinnstrategie $\tau$ spielt (Widerspruch). Dann gilt $f\in K_u$ mit
\begin{alignat}{2}
&K_u&:=\{f'\in \omega^\omega\mid &u\TT{ gut f"ur } f'\UND\notag\\
		&& &\neg\EX[(s_{n+1},k_{n+2})](u\kon (s_{n+1},k_{n+2}) \TT{ gut }\UND\TT{ gut f"ur } f')\}\notag\\
		%&& &u\kon (s_{n+1},k_{n+2}) \TT{ gut f"ur } f')\}\notag\\
&		&=\{f'\in \omega^\omega\mid &\EX[s_{n+1}](s_0\kon\ldots\kon s_n\kon s_{n+1}\prec f'\UND s_{n+1}(0)>k_{n+1})\UND\notag\\%label{Ku_01}\\
		&& &\neg\EX[(s_{n+1},k_{n+2})]((s_{n+1}(0)>k_{n+1}\UND k_{n+2}\TT{ $\tau$-gesp.})\UND\notag\\
		&& &\EX[s_{n+2}](s_0\kon\ldots\kon s_n\kon s_{n+1}\kon s_{n+2}\prec f'\UND s_{n+2}(0)> k_{n+2}))\}\label{Ku_02}
\end{alignat}
\begin{comment}
\begin{alignat*}{1}
K_u:=&\{f'\in A\mid s_0\kon\ldots\kon s_n=(f'(0),\ldots,f'(m-1))\prec f'\UND f'(m)>k_{n+1}\UND\\
		&\neg\EX[(s_{n+1},k_{n+2})\in \BR\times\omega](u\kon (s_{n+1},k_{n+2}) \TT{ gut }\UND\\
		&s_0\kon\ldots\kon s_n\kon s_{n+1}=(f'(0),\ldots,f'(l-1))\prec f'\UND f'(l)>k_{n+2}\}
\end{alignat*}
\end{comment}
wobei \glqq$k_{n+2}\TT{ $\tau$-gesp.}$\grqq\ hei"st: $k_{n+2}=\tau(u\kon (s_{n+1}))$.
Damit gilt dann
\begin{alignat*}{1}
&A\subset\bigcup_{u\TT{\scriptsize{ gut}}}K_u.
\end{alignat*}
Da $u\in\omega^{\omega}$ und $\omega^{<\omega}$ abz"ahlbar ist, ist die Vereinigung abz"ahlbar. 

\medskip Um zu zeigen, da"s $A$ $\sigma$-beschr"ankt ist, gen"ugt es nun zu zeigen, da"s die $K_u$ $\sigma$-kompakt sind, denn dann ist die abz"ahlbare Vereinigung $\bigcup_{u}K_u$ $\sigma$-kompakter Mengen $\sigma$-kompakt und $A$ als Teilmenge einer $\sigma$-kompakten Menge $\sigma$-beschr"ankt.

\point{Die $K_u$ sind $\sigma$-kompakt}
Sei $K_u$ wie oben konstruiert mit $u=(s_0,k_1,s_1,k_2,\ldots,s_n,k_{n+1})$. Dann l"a"st sich $K_u$ schreiben als
\[
K_u=\bigcup_{\widetilde{s_{n+1}}(0)>k_{n+1}}\underbrace{O_{s_0\kon\ldots\kon s_n\kon \widetilde{s_{n+1}}}\cap \C[{(\bigcup_{\substack{s_{n+1}(0)>k_{n+1}\\\EX[k_{n+2}](k_{n+2}\TT{\scriptsize{ $\tau$-gesp.}}\UND s_{n+2}(0)>k_{n+2})}} O_{s_0\kon\ldots\kon s_n\kon s_{n+1}\kon s_{n+2}})}]}_{=:L_{\widetilde{s_{n+1}}}}.
\]
Nun gen"ugt es zu zeigen, da"s die $L_{\widetilde{s_{n+1}}}$ kompakt sind. Dies ist nach Satz \ref{Satz_Heine_Borel_2} genau dann der Fall, wenn die $L_{\widetilde{s_{n+1}}}$ abgeschlossen und "uberall endlich verzweigt sind.

\point{Die $L_{\widetilde{s_{n+1}}}$ sind abgeschlossen}
Das Komplement einer Menge $L_{\widetilde{s_{n+1}}}$ ist als Vereinigung offener Mengen offen, also ist $L_{\widetilde{s_{n+1}}}$ abgeschlossen.

\point{Die $L_{\widetilde{s_{n+1}}}$ sind "uberall endlich verzweigt}
Sei $L:=L_{\widetilde{s_{n+1}}}$ f"ur ein $\widetilde{s_{n+1}}\in\omega^{<\omega}$. Es ist zu zeigen, da"s
\begin{equation*}
\FA[s\in T_L]\EX[k<\omega]\FA[t\in\omega^{<\omega}_*](s\kon t\in T_L\Rightarrow t(0)\not > k).
\end{equation*}
\point{Fall 1} 
Sei $s\prec \tilde{s}:=s_0\kon\ldots\kon s_n\kon s_{n+1}$ und $s\neq \tilde{s}$ (d.h. $s$ ist ein echtes Anfangsst"uck von $\tilde{s}$)\gs etwa $s=(n_0,\ldots,n_i)\prec\tilde{s}=(n_0,\ldots,n_m)$ mit $i\lneq m$. Sei $k:=n_{i+1}$, dann gilt $\FA[t\in \omega^{<\omega}_*] (s\kon t\in T_L\Rightarrow t(0)\not> k)$.

\point{Fall 2}
Sei $s\succ s_0\kon\ldots\kon s_n\kon s_{n+1}$. Etwa $s= s_0\kon\ldots\kon s_n\kon s_{n+1}'$ mit $s_{n+1}\prec s_{n+1}'$ (also eventuell auch $s_{n+1}= s_{n+1}'$). Dann gilt wegen $s_{n+1}(0)> k_{n+1}$ und $s_{n+1}\prec s_{n+1}'$: $s_{n+1}'(0)> k_{n+1}$. Spielt nun Spieler II $k_{n+2}':=\tau(u\kon(s_{n+1}'))$, so mu"s wegen (\ref{Ku_02}) f"ur jedes $t\in \omega_*^{<\omega}$ f"ur alle $f'\in K_u$ gelten:
\begin{alignat*}{1}
&s_0\kon\ldots\kon s_n\kon s_{n+1}'\kon t\not\prec f' \ODER t(0)\not> k_{n+2}'.
\end{alignat*}
Also gilt falls $t(0)> k_{n+2}'$, da"s $s\kon t\not\in T_L$. Zu $s$ existiert also ein $k:=k_{n+2}'$ mit $\FA[t\in \omega^{<\omega}_*] (s\kon t\in T_L\Rightarrow t(0)\not> k)$.

\medskip Somit gilt nach Fall 1 und 2, da"s die $L_{\widetilde{s_{n+1}}}$ "uberall endlich verzweigt sind.

\medskip Damit ist die Hinrichtung von $(ii)$ bewiesen.
\end{beweis}

Im Beweis des Theorems \ref{SatzSigmaBeschrCharSpiele} wurde $\AC$ nicht benutzt\gs unter der Annahme von $\AD$ gilt somit:

\begin{korollar}\label{Korollar_Komplement_Sigmabeschr}
Sei $A$ das Komplement einer $\sigma$-beschr"ankten Menge $B\subset\BR$ und es gelte $\AD$, dann enth"alt $A$ eine nicht-leere superperfekte Teilmenge. 
\end{korollar}
\begin{beweis}
Das Komplement $A$ einer $\sigma$-beschr"ankten Menge $B$ kann nicht $\sigma$-besch"ankt sein\gs sonst w"are $\BR=A\cup B$ als endliche Vereinigung $\sigma$-beschr"ankter Mengen $\sigma$-beschr"ankt. Nach Theorem \ref{SatzSigmaBeschrCharSpiele} $(ii)$ hat dann Spieler II in dem Spiel $\BMSb{A}$ keine Gewinnstrategie. Unter der Annahme von $\AD$ hat somit Spieler I eine Gewinnstrategie. Und nach Theorem \ref{SatzSigmaBeschrCharSpiele} $(i)$ folgt, da"s $A$ eine nicht-leere superperfekte Teilmenge enthalten mu"s.
\end{beweis}

\subsection*{Banach-Mazur-Spiel mit Zeugen $\BMSc{B}$}
Wie schon im Falle der $\BMSo{}$-Spiele kann es auch f"ur die $\BMSb{}$-Spiele von Nutzen sein eine gegen"uber der zu charakterisierenden Menge $A$ vereinfachte Gewinnmenge $B$ betrachten zu k"onnen [\siehe etwa S"atze \ref{SatzDefbarkeit_1} und \ref{SatzDefbarkeit_2}].
Analog zu Theorem \ref{SatzMagerCharSpiele_2} soll daher nun eine spieltheoretische Charakterisierung der $\sigma$-Kategorie f"ur Teilmengen $A$ des Baireraumes mit $A=\p(B)$ f"ur eine Teilmenge $B\subset\BR\times\lambda$ (f"ur eine unendliche Ordinalzahl $\lambda$) gegeben werden. Dazu wird f"ur die Menge $B$ ein Spiel $\BMSc{B}$ definiert, in dem Spieler I genau dann eine Gewinnstrategie hat, wenn $A$ eine nicht-leere superperfekte Teilmenge enth"alt und II eine Gewinnstrategie hat genau dann, wenn $A$ $\sigma$-beschr"ankt ist.
 
\begin{definition}\header{Banach-Mazur-Spiel mit Zeugen $\BMSc{B}$}\label{DefBMSMitZeugen}\\
Sei $A\subset\BR$ und $B\subset\BR\times\lambda^\omega$ f"ur eine Ordinalzahl $\lambda$ und es gelte:
\begin{alignat*}{1}
&A=\p(B):=\MNG{f\in{\BR}}{\EX[{\xi\in\lambda^\omega}]((f,\xi)\in B)}.
\end{alignat*}

Dann definiert man das \emph{Banach-Mazur-Spiel $\BMSc{B}$ mit Zeugen} wie folgt:
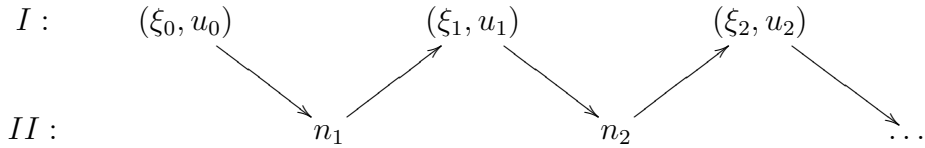
\[
\xymatrix{
I: &(\xi_0,u_0) \ar[dr] & 					 	 &(\xi_1,u_1) \ar[dr]  &						&(\xi_2,u_2) \ar[dr]  &\\
II:&		 				 			&n_1 \ar[ur]	 &						       &n_2 \ar[ur]	&										&\ldots
}
\]
Spieler I und Spieler II spielen abwechselnd\gs Spieler I spielt ein Element $\xi_0\in\lambda$ und eine nicht-leere
endliche Sequenz nat"urlicher Zahlen $u_0\in\omega_*^{<\omega}$, Spieler II spielt eine nat"urliche Zahl $n_1<\omega$, Spieler I spielt ein Element $\xi_1\in\lambda$ und eine nicht-leere endliche Sequenz $u_1\in\omega_*^{<\omega}$, Spieler II spielt eine nat"urliche Zahl $n_2<\omega$ u.s.w.. Sei nun $\xi=(\xi_0,\xi_1,\ldots)\in\lambda^\omega$ und $f=u_0\kon u_1\kon u_2\kon\ldots\in\BR$. 
Man sagt, da"s \emph{Spieler I gewinnt}, falls gilt:
\[
(f,\xi)\in B.
\]
Ansonsten sagt man, da"s \emph{Spieler II gewinnt}. $B$ nennt man dabei auch die \emph{Gewinnmenge}.

Da Spieler I nicht nur endliche Sequenzen sondern auch Elemente aus $\lambda$ spielt, die zu einem \eemph{Zeugen} $\xi\in\lambda^\omega$ zusammengesetzt werden\gs wobei Spieler I $(f,\xi)\in B$ anstrebt\gs nennt man das Spiel $\BMSc{B}$ auch ein Banach-Mazur-Spiel \emph{mit Zeugen}.
\end{definition}

Die Begriffe \emph{(Gewinn-)Strategie f"ur I (bzw. f"ur II)} und \emph{determiniert} sind f"ur die Spiele $\BMSc{B}$ analog zu den Definitionen \ref{Def_Gewinnstrategie_I} (bzw. \ref{Def_Gewinnstrategie_II}) und \ref{Def_determiniert_G_A}  definiert.

\bigskip Ist $\lambda$ abz"ahlbar, l"a"st sich dieses Spiel mit einer Kodierung der in $\BMSc{B}$ gespielten Sequenzen und Zeugen in die nat"urlichen Zahlen auch als ein Spezialfall der Grundvariante [\siehe Definition \ref{Def_SpielGrundvariante}] auffassen und analog zu Lemma \ref{Lemma_Determiniertheit_G**p} und  Lemma \ref{Lemma_Determiniertheit_G_tilde} gilt:

\begin{lemma}\header{Determiniertheit}\label{Lemma_Determiniertheit_G_tilde_p}\\
Falls $\AD$ gilt und $\lambda$ abz"ahlbar ist, so sind auch die obigen Spiele $\BMSc{B}$ determiniert.
\end{lemma}
\begin{beweis}\KOM[KURZ]\KOM[wirklich\\analog?\\$\lambda$]
Analog zu Lemma \ref{Lemma_Determiniertheit_G**p}.
\end{beweis}

Die Begriffe\label{Def_LambdaBeschr} \eemph{$\sigma$-beschr"ankt} und \eemph{$\sigma$-kompakt} lassen sich f"ur beliebige unendliche Ordinalzahlen $\lambda$ verallgemeinern: Eine Teilmenge $A\subset \BR$ hei"se \emph{$\lambda$-beschr"ankt} genau dann, wenn es eine Menge $\MNG{f_i\in\BR}{i<\lambda}$ gibt so, da"s es f"ur jedes $f\in A$ ein $i<\lambda$ gibt mit $f\leq f_i$. Die Menge $\MNG{f_i\in\BR}{i<\lambda}$ hei"st dann eine \emph{$\lambda$-Schranke} von $A$. Eine \emph{$\lambda$-kompakte} Teilmenge eines Raumes sei eine Vereinigung $\card[\lambda]$-vieler kompakter Teilmengen.

\bigskip Analog zu Satz \ref{SatzSigmaBeschr"ankt} $(ii)$ gilt:
\begin{satz}\header{$\lambda$-beschr"ankt}\label{SatzLambdaBeschr}\\
F"ur eine unendliche Ordinalzahl $\lambda$ und eine Teilmenge $A$ des Baireraumes $\BR$ ist "aquivalent:
\begin{ITEMS}[equivalences]
\item $A$ ist $\lambda$-beschr"ankt.
\item $A\subset B$ f"ur eine $\lambda$-kompakte Menge $B\subset\BR$.
\end{ITEMS}
\end{satz}
\begin{beweis}
$(i)\Rightarrow(ii)$: Sei $A$ etwa durch $\MNG{f_i\in {\BR}}{i<\lambda}$ $\lambda$-beschr"ankt. Sei $A_i:=\MNG{g\in A}{g\leq f_i}$. Dann sind die $A_i$ $\leq$-beschr"ankt und wegen
\[
\C[A_i]=\bigcup_{\EX[n<\lng(u)](u(n)>f_i(n))}O_u
\]
auch abgeschlossen\gs nach Satz \ref{Satz_Heine_Borel_1} also kompakt. Somit ist $A\subset \bigcup_{i<\lambda} A_i$ $\lambda$-kompakt.

\bigskip $(ii)\Rightarrow(i)$: Sei $A\subset \bigcup_{i<\lambda} A_i$ mit kompakten $A_i$. Dann ist jedes $A_i$ laut Satz \ref{Satz_Heine_Borel_1} $\leq$-beschr"ankt\gs etwa durch $f_i$. Dann ist $\MNG{f_i\in\BR}{i<\lambda}$ aber eine $\lambda$-Schranke f"ur $A$. Das hei"st $A$ ist $\lambda$-beschr"ankt.
\end{beweis}

\bigskip Aus der Definition f"ur \eemph{$\lambda$-beschr"ankt} folgt sofort, da"s die $\lambda$-beschr"ankten Teilmengen des Baireraumes abgeschlossen sind gegen"uber Teilmengen. Au"serdem ist die Vereinigung $\card[\lambda]$-vieler $\sigma$-beschr"ankter Mengen $\card[\lambda]$-beschr"ankt\footnote{{F"ur beliebiges $\lambda\in\ORD$ gilt $\card[{\omega\times\lambda}]=\card[\lambda]$} (\siehe etwa \cite[3.5]{Jech:2003}).\label{Fussnote_SigmaKat_Def}}.\KOM

\bigskip Damit lassen sich nun \eemph{$\lambda$-Beschr"anktheit} und \eemph{Enthalten einer nicht-leeren superperfekten Teilmenge} f"ur Mengen $A=\textup{p}(B)$ (f"ur ein $B\subset\BR\times\lambda^\omega$) spieltheoretisch wie folgt charakterisieren:

\begin{theorem}\header{Charakterisierung durch Spiele mit Zeugen}\label{SatzSigmaBeschrCharSpiele_2}\\
Sei $A\subset\BR$ und $A=\p(B)$ f"ur ein $B\subset\BR\times\lambda^\omega$ und eine unendliche Ordinalzahl $\lambda$. Dann gilt:
\begin{ITEMS}
\item Spieler I hat eine Gewinnstrategie in $\BMSc{B}$ $\Rightarrow$ $A$ enth"ahlt eine nicht-leere superperfekte Teilmenge $P$.
\item Spieler II hat eine Gewinnstrategie in $\BMSc{B}$ $\Rightarrow$ $A$ ist $\lambda$-beschr"ankt.
\end{ITEMS}
\end{theorem}
\begin{beweis}
$(i)$: Habe I eine Gewinnstrategie $\sigma$ f"ur das Spiel $\BMSc{B}$. Wegen $A=\p(B)$ hat I dann auch eine Gewinnstrategie $\tilde{\sigma}$ f"ur das Spiel $\BMSb{A}$. Nach Theorem \ref{SatzSigmaBeschrCharSpiele} enth"alt $A$ somit eine nicht-leere superperfekte Teilmenge $P$.

\bigskip$(ii)$: Habe Spieler II eine Gewinnstrategie $\tau$.

\point{Gute Sequenz}Eine Sequenz $u=(l_0,s_0,k_1,l_1,s_1,k_2,\ldots,l_n,s_n,k_{n+1})$ mit $s_i\in\omega_*^{<\omega}$ f"ur $i=0,\ldots,n$ sowie $k_i\in\omega$ f"ur $i=1,\ldots,n+1$ und $l_j\in\lambda$ f"ur $j=0,\ldots,n$ hei"se \emph{gut}, falls gilt:
\begin{alignat*}{1}
&s_i(0)>k_i\TT{ f"ur }i=1,\ldots,n\quad\UND\\
&k_j\TT{ mit $\tau$ gespielt f"ur }j=1,\ldots,n+1.
\end{alignat*}

\point{Sequenz gut f"ur $f\in A$}$u$ hei"se \emph{gut f"ur} $f\in A$, falls gilt:
\begin{alignat*}{1}
&\EX[s_{n+1}](s_0\kon\ldots\kon s_n\kon s_{n+1}\prec f\UND s_{n+1}(0)>k_{n+1})
\end{alignat*}
Insbesondere ist die leere Sequenz gut f"ur jedes $f\in A$.

\medskip Sei nun $(f,g)\in B$ beliebig. Dann mu"s es eine Sequenz $u=(\underbrace{l_0,s_0}_{I},\underbrace{k_1}_{II},\underbrace{l_1,s_1}_{I},\underbrace{k_2}_{II},\ldots,\underbrace{l_n,s_n}_{I},\underbrace{k_{n+1}}_{II})$ (eventuell $u=()$) geben mit:
\begin{multline*}
u\TT{ gut }\UND\TT{ gut f"ur $f$}\UND\\ (l_0,\ldots,l_n)=(g(0),\ldots,g(n))\prec g \UND\\
\neg\EX[(l_{n+1},s_{n+1},k_{n+2})](u\kon(l_{n+1},s_{n+1},k_{n+2})\TT{ gut}\UND\TT{ gut f"ur $f$}\UND\\
(l_0,\ldots,l_n, l_{n+1})=(g(0),\ldots,g(n),g(n+1))\prec g)
\end{multline*}
Ansonsten k"onnte Spieler I gewinnen, obwohl Spieler II mit seiner Gewinnstrategie $\tau$ spielt (Widerspruch). 

\medskip Sei nun $r:=l_{n+1}:=g(n+1)$, dann gilt $f\in M_{(u,r)}$ mit 
\begin{alignat}{2}
&M_{(u,r)}&:=\{f'\in \omega^\omega\mid &u\TT{ gut f"ur }f'\UND\notag\\
		&&&\neg\EX[(s_{n+1},k_{n+2})](u\kon (r,s_{n+1},k_{n+2}) \TT{ gut }\UND\TT{ gut f"ur }f'\}\notag\\
&				&=\{f'\in \omega^\omega\mid &\EX[s_{n+1}](s_0\kon\ldots\kon s_n\kon s_{n+1}\prec f'\UND s_{n+1}(0)>k_{n+1})\UND\notag\\
&				 &&\neg\EX[(s_{n+1},k_{n+2})]((s_{n+1}(0)>k_{n+1}\UND k_{n+2}\TT{ $\tau$-gesp.})\UND\notag\\
		&& &\EX[s_{n+2}](s_0\kon\ldots\kon s_n\kon s_{n+1}\kon s_{n+2}\prec f'\UND s_{n+2}(0)> k_{n+2}))\}\label{Mur_02}
\end{alignat}
wobei \glqq$k_{n+2}\TT{ $\tau$-gesp.}$\grqq\ hei"st: $k_{n+2}=\tau(u\kon (r, s_{n+1}))$.
Damit gilt dann
\begin{alignat*}{1}
&A\subset\bigcup_{(u,r)}M_{(u,r)}
\end{alignat*}
wobei "uber die Paare $(u,r)$ vereinigt wird mit $u$ gut und $r\in\lambda$\gs da $(u,r)\in\omega^{<\omega}\times\lambda$ und $\omega^{<\omega}$ abz"ahlbar ist, werden $\card[\lambda]$-viele Mengen vereinigt\footnote{{F"ur beliebiges $\lambda\in\ORD$ gilt $\card[{\omega\times\lambda}]=\card[\lambda]$} (\siehe etwa \cite[3.5]{Jech:2003}, vgl. Fu"snote \ref{Fussnote_SigmaKat_Def}).\label{Fussnote_SBeschr_Spiel}}. 

\medskip Da wir die Mengen $M_{(u,r)}$ [\siehe (\ref{Mur_02})] ganz "ahnlich wie die Mengen $K_u$ [\siehe (\ref{Ku_02})] im Beweis der Hinrichtung von $(ii)$ von Theorem \ref{SatzSigmaBeschrCharSpiele} definiert haben (lediglich die Bedeutung von \glqq$k_{n+2}\TT{ $\tau$-gesp.}$\grqq\ ist eine andere), k"onnen wir von nun an analog dazu fortfahren. 

\medskip Um zu zeigen, da"s $A$ $\lambda$-beschr"ankt ist, gen"ugt es zu zeigen, da"s die $M_{(u,r)}$ $\lambda$-kompakt sind, denn dann ist die abz"ahlbare Vereinigung $\bigcup_{(u,r)}M_{(u,r)}$ $\lambda$-kompakter Mengen $\lambda$-kompakt und $A$ als Teilmenge einer $\lambda$-kompakten Menge $\lambda$-beschr"ankt.

\point{Die $M_{(u,r)}$ sind $\lambda$-kompakt}
Sei $M_{(u,r)}$ wie oben konstruiert mit $u=(l_0,s_0,k_1,l_1,s_1,k_2,\ldots,l_n,s_n,k_{n+1})$. Dann l"a"st sich $M_{(u,r)}$ schreiben als
\[
M_{(u,r)}=\bigcup_{\widetilde{s_{n+1}}(0)>k_{n+1}}\underbrace{O_{s_0\kon\ldots\kon s_n\kon \widetilde{s_{n+1}}}\cap \C[{(\bigcup_{\substack{s_{n+1}(0)>k_{n+1}\\\EX[k_{n+2}](k_{n+2}\TT{\scriptsize{ $\tau$-gesp.}}\UND s_{n+2}(0)>k_{n+2})}} O_{s_0\kon\ldots\kon s_n\kon s_{n+1}\kon s_{n+2}})}]}_{=:L_{\widetilde{s_{n+1}}}}.
\]
Nun gen"ugt es zu zeigen, da"s die $L_{\widetilde{s_{n+1}}}$ kompakt sind. Dies ist nach Satz \ref{Satz_Heine_Borel_2} genau dann der Fall, wenn die $L_{\widetilde{s_{n+1}}}$ abgeschlossen und "uberall endlich verzweigt sind. 

\point{Die $L_{\widetilde{s_{n+1}}}$ sind abgeschlossen}
Das Komplement einer Menge $L_{\widetilde{s_{n+1}}}$ ist als Vereinigung offener Mengen offen, also ist $L_{\widetilde{s_{n+1}}}$ abgeschlossen.

\point{Die $L_{\widetilde{s_{n+1}}}$ sind "uberall endlich verzweigt}
Von hier an besteht der einzige Unterschied zum Beweis der Hinrichtung von $(ii)$ von Theorem \ref{SatzSigmaBeschrCharSpiele} in der Definition von \framebox{$k_{n+2}':=\tau(u\kon(r,s_{n+1}'))$} im zweiten Fall:
Sei $L:=L_{\widetilde{s_{n+1}}}$ f"ur ein $\widetilde{s_{n+1}}\in\omega^{<\omega}$. Es ist zu zeigen, da"s
\begin{equation*}
\FA[s\in T_L]\EX[k<\omega]\FA[t\in\omega^{<\omega}_*](s\kon t\in T_L\Rightarrow t(0)\not > k).
\end{equation*}
\point{Fall 1} 
Sei $s\prec \tilde{s}:=s_0\kon\ldots\kon s_n\kon s_{n+1}$ und $s\neq \tilde{s}$ (d.h. $s$ ist ein echtes Anfangsst"uck von $\tilde{s}$)\gs etwa $s=(n_0,\ldots,n_i)\prec\tilde{s}=(n_0,\ldots,n_m)$ mit $i\lneq m$. Sei $k:=n_{i+1}$, dann gilt $\FA[t\in \omega^{<\omega}_*] (s\kon t\in T_L\Rightarrow t(0)\not> k)$.

\point{Fall 2}
Sei $s\succ s_0\kon\ldots\kon s_n\kon s_{n+1}$. Etwa $s= s_0\kon\ldots\kon s_n\kon s_{n+1}'$ mit $s_{n+1}\prec s_{n+1}'$ (also eventuell auch $s_{n+1}= s_{n+1}'$). Dann gilt wegen $s_{n+1}(0)> k_{n+1}$ und $s_{n+1}\prec s_{n+1}'$: $s_{n+1}'(0)> k_{n+1}$. Spielt nun Spieler II \framebox{$k_{n+2}':=\tau(u\kon(r,s_{n+1}'))$}, so mu"s wegen (\ref{Mur_02}) f"ur jedes $t\in \omega_*^{<\omega}$ f"ur alle $f'\in M_{(u,r)}$ gelten:
\begin{alignat*}{1}
&s_0\kon\ldots\kon s_n\kon s_{n+1}'\kon t\not\prec f' \ODER t(0)\not> k_{n+2}'.
\end{alignat*}
Also gilt falls $t(0)> k_{n+2}'$, da"s $s\kon t\not\in T_L$. Zu $s$ existiert also ein $k:=k_{n+2}'$ mit $\FA[t\in \omega^{<\omega}_*] (s\kon t\in T_L\Rightarrow t(0)\not> k)$.

\medskip Somit gilt nach Fall 1 und 2, da"s die $L_{\widetilde{s_{n+1}}}$ "uberall endlich verzweigt sind.

\medskip Damit ist die Hinrichtung von $(ii)$ bewiesen.
\end{beweis}

\bigskip F"ur $\lambda=\omega$ bedeutet \eemph{$\lambda$-beschr"ankt} das selbe wie \eemph{$\sigma$-beschr"ankt}. In diesem Fall lautet Theorem \ref{SatzSigmaBeschrCharSpiele_2}:

\begin{korollar}\header{Charakterisierung durch Spiele mit Zeugen}\label{SatzSigmaBeschrCharSpiele_3}\\
Sei $A\subset\BR$ und $A=\p(B)$ f"ur eine Menge $B\subset\BR\times\BR$. Dann gilt:
\begin{ITEMS}
\item Spieler I hat eine Gewinnstrategie in $\BMSc{B}$ $\Rightarrow$ $A$ enth"ahlt eine nicht-leere superperfekte Teilmenge $P$.
\item Spieler II hat eine Gewinnstrategie in $\BMSc{B}$ $\Rightarrow$ $A$ ist $\sigma$-beschr"ankt.
\end{ITEMS}
\end{korollar}

Nach Theorem \ref{SatzSigmaBeschrCharSpiele} folgt aus $A\subset\BR$ $\sigma$-beschr"ankt, da"s $A$ keine nicht-leere superperfekte Teilmenge enth"alt und aus $A\subset\BR$ enth"alt eine nicht-leere superperfekte Teilmenge, da"s $A$ nicht $\sigma$-beschr"ankt ist (schlie"slich k"onnen ja nicht gleichzeitig I und II Gewinnstrategien f"ur $\BMSb{A}$ haben). Unter der Annahme von $\AD$ gelten daher in Korollar \ref{SatzSigmaBeschrCharSpiele_3} auch die beiden R"uckrichtungen.

\bigskip Mit Hilfe von Korollar \ref{SatzSigmaBeschrCharSpiele_3} lassen sich nun im n"achsten Kapitel einige Definierbarkeits-Resultate ableiten.

\subsection*{Definierbarkeit}
Mittels Korollar \ref{SatzSigmaBeschrCharSpiele_3} sowie eines Resultates von D.A. Martin und eines weiteren Resultates von Y.N. Moschovakis lassen sich nun Aussagen "uber die Komplexit"at von $\sigma$-Schranken ableiten.

\begin{comment}
\begin{satz}\header{Moschovakis}\label{Moschovakis}\\
Unter der Vorraussetzung, da"s alle projektiven Spiele determinieren, gilt: Ist $\textup{G}(A)$ ein $\DELTA{1}{2n}{}$-Spiel, so hat Spieler I oder Spieler II eine Gewinnstrategie in $\DELTAlf{1}{2n+1}{}$.\KOM[UMFORMEN]
\end{satz}

Mit $\Gamma=\GQ\SIGMAlf{1}{2n}{}(f)$ f"ur ein $f\in\BR$ gilt $\GQ\Gamma=\GQ\SIGMAlf{1}{2n}{}(f)\EQ[(\ref{Frml_G_Versch_Lightface})]\PIlf{1}{2n+1}{}(f)$.\KOM[ERKL"AREN\\$\Gamma(f)$] Hat Spieler I nun eine Gewinnstrategie $\sigma$ f"ur das Spiel $\textup{G}(A)$ mit $A\in\Gamma$, so hat nach dem Dritten Periodizit"atstheorem Spieler I eine $\GQ\Gamma$-rekursive also $\PIlf{1}{2n+1}{}(f)$-rekursive Gewinnstrategie $\sigma$. Damit ist $\sigma$ aber in $\DELTAlf{1}{2n+1}{}(f)$:\ldots\KOM[FEHLT] 
\end{comment}

\bigskip Relativ einfach einzusehen ist, da"s jedes offene (und damit auch jedes abgeschlossene) Spiel $\textup{G}(A)$ determiniert ist\footnote{\emph{Gale-Stewart:} Von David Gale und Frank M. Stewart (\siehe \cite{Gale:1953}).}:

\begin{satz}\header{Gale-Stewart}\label{Gale_Stewart}\\
Sei $A$ eine offene Teilmenge des Baireraumes. Dann ist das Spiel $\textup{G}(A)$ determiniert.
\end{satz}
\begin{beweis}
Gewinnt Spieler I das Spiel $\textup{G}(A)$ (ist also das gespielte $f\in A$), so ist Spieler I schon nach endlich vielen Z"ugen in eine Position gelangt, von der aus er nicht mehr verlieren konnte: Da $A$ offen ist, mu"s Spieler I die Menge $A$ in irgendeiner offenen Basismenge $O_s\subset A$ treffen. Wenn Spieler I dies geschafft hat, so hat er daf"ur aber nur endlich viele Z"uge $\leq\lng(s)$ ben"otigt. In diesem Fall ist das Spiel also schon nach endlich vielen Z"ugen entschieden. Ist das Spiel hingegen erst nach unendlich vielen Z"ugen entschieden, so gewinnt Spieler II, da dies hei"st, da"s Spieler I die offene Menge $A$ nicht getroffen haben kann.

\point{Angenommen Spieler I hat keine Gewinnstrategie} Eine Gewinnstrategie f"ur Spieler II besteht dann darin, N. Lusins Tip: \flqq Avoid losing positions\frqq\ zu befolgen und nur solche Z"uge abzugeben, die garantieren, da"s Spieler I durch seinen n"achsten Zug nicht in eine Position gelangt, von der aus er nicht mehr verlieren kann. Spieler II kann auf diese Weise spielen, da sonst Spieler I von Anfang an eine Gewinnstrategie h"atte (im Widerspruch zur Annahme, da"s Spieler I keine Gewinnstrategie hat). Somit hat Spieler II eine Gewinnstrategie.

\medskip Dies zeigt, da"s $\textup{G}(A)$ f"ur offenes $A$ determiniert ist.
\end{beweis}

\begin{satz}\header{Definierbarkeit einer $\sigma$-Schranke f"ur $\SIGMA{1}{1}{}$-Mengen}\label{SatzDefbarkeit_1}\\
Jede $\SIGMA{1}{1}{}$-Menge in $\BR$ enth"alt eine nicht-leere superperfekte Teilmenge oder ist $\sigma$-beschr"ankt mit einer $\sigma$-Schranke in $\DELTAlf{1}{1}{}$.
\end{satz}
\begin{beweis}
Sei $A\in\SIGMA{1}{1}{}(\BR)$. Dann gilt:
\[
f\in A\Leftrightarrow\EX[g\in{\BR}]((f,g)\in B)
\]
f"ur eine abgeschlossene Teilmenge $B$ des Baireraumes. Das Spiel $\BMSc{B}$ ist also determiniert. 

\medskip Enth"alt nun $A$ keine nichtleere superperfekte Teilmenge, dann hat nach Korollar \ref{SatzSigmaBeschrCharSpiele_3} Spieler II eine Gewinnstrategie in $\BMSc{B}$. Nach dem dritten Periodizit"atstheorem von Y. N. Moschovakis\footnote{\emph{Drittes Periodizit"atstheorem}:\label{FN_Moschovakis} Moschovakis definiert in \cite[6D.]{Moschovakis:1980} einen \emph{Spielequantoren} $\GQ$ f"ur Punktemengen $P\subset \omega\times\BR$:
\begin{alignat*}{1}
n\in\GQ P 	&:\Leftrightarrow \GQ g ((n,g)\in P)\\
						&:\Leftrightarrow\TT{Spieler I gewinnt das Spiel }\textup{G}(A)\TT{ mit } A:=\MNG{g\in{\BR}}{(n,g)\in P} 
\end{alignat*}
und f"ur Punktklassen $\Gamma$: $\GQ\Gamma:=\MNG{\GQ P}{P\in\Gamma({\omega\times\BR})}$. Der so definierte Spielquantor $\GQ$ \glqq verschiebt\grqq\ die Lightface-Hierarchie (\siehe \cite[6D.]{Moschovakis:1980}):
\begin{align}
\GQ\SIGMAlf{1}{1}{}=\PIlf{1}{1}{},\quad\GQ\PIlf{1}{1}{}=\SIGMAlf{1}{2}{},\quad\GQ\SIGMAlf{1}{2}{}=\PIlf{1}{3}{},\quad\GQ\PIlf{1}{3}{}=\SIGMAlf{1}{4}{},\ldots .\label{Frml_G_Versch_Lightface}
\end{align}
Unter der Annahme von projektiver Determiniertheit erf"ullt die Boldface-Hierarchie die Voraussetzungen des 3. Periodizit"atstheorems. Mit $\Gamma=\SIGMAlf{1}{2n}{}$ gilt $\GQ\Gamma=\GQ\SIGMAlf{1}{2n}{}\EQ[(\ref{Frml_G_Versch_Lightface})]\PIlf{1}{2n+1}{}$. Hat Spieler I nun eine Gewinnstrategie $\sigma$ f"ur das Spiel $\textup{G}(A)$ mit $A\in\Gamma$, so hat nach dem Dritten Periodizit"atstheorem Spieler I eine $\GQ\Gamma$-rekursive also $\PIlf{1}{2n+1}{}$-rekursive Gewinnstrategie $\sigma$ (d.h. der Graph $\textup{G}_\sigma$ von $\sigma$ ist $\PIlf{1}{2n+1}{}$). $\PIlf{1}{2n+1}{}$-rekursive Funktionen liegen aber in $\DELTAlf{1}{2n+1}{}$ (Sei $\sigma:\PFEIL{\omega}{}{\omega}$ eine Strategie und $\textup{G}_\sigma\in\PIlf{1}{2n+1}{}$, dann sieht man  $\C[{\textup{G}_\sigma}]=\MNG{(m,n)}{\EX[l<\omega](\sigma(m)=l\UND l\neq n)}\in\PIlf{1}{2n+1}{}$ \siehe etwa \cite[3E.2]{Moschovakis:1980}). Aus Symmetriegr"unden erh"alt man dieses Resultat auch f"ur Spieler II.} (\siehe \cite{Moschovakis:1973} oder \cite[6E.1.]{Moschovakis:1980}) hat Spieler II dann eine $\DELTAlf{1}{1}$-Gewinnstrategie. 
Da diese Gewinnstrategie eine $\sigma$-Schranke f"ur $A$ definiert folgt, da"s $A$ eine $\sigma$-Schranke in $\DELTAlf{1}{1}$ hat. 
\end{beweis}

Nimmt man ein Resultat von D.A. Martin hinzu, l"a"st sich die Komplexit"at von $\sigma$-Schranken f"ur $\SIGMA{1}{2n+1}{}$-Mengen absch"atzen: 

\begin{satz}\header{Definierbarkeit einer $\sigma$-Schranke f"ur $\SIGMA{1}{2n+1}{}$-Mengen}\label{SatzDefbarkeit_2}\\
Ist jede $\DELTA{1}{2n}{}$-Menge in $\BR$ determiniert, so enth"alt jede $\SIGMA{1}{2n+1}{}$-Menge in $\BR$ eine nicht-leere superperfekte Teilmenge oder ist $\sigma$-beschr"ankt mit einer $\sigma$-Schranke in $\DELTAlf{1}{2n+1}{}$.
\end{satz}
\begin{beweis}
Sei $A\in\SIGMA{1}{2n+1}{}$. Dann gilt:
\[
f\in A\Leftrightarrow\EX[g\in{\BR}]((f,g)\in B)
\]
f"ur eine $\PI{1}{2n}{}$-Menge $B$. Das Spiel $\BMSc{B}$ ist also ein $\PI{1}{2n}{}$-Spiel. 

\medskip Sei nun jede $\DELTA{1}{2n}{}$-Menge in $\BR$ determiniert. Aus $\DELTA{1}{2n}{}$-Determiniertheit folgt nach einem Resultat von D. A. Martin\footnote{\emph{Satz(Martin):} Unter der Voraussetzung von $\DC$ gilt f"ur alle $n<\omega$: Ist jede $\DELTA{1}{2n}{}$-Menge determiniert, so ist auch jede $\SIGMA{1}{2n}{}$-Menge determiniert. (\siehe \cite[196ff]{KechrisSolovay:1985} oder \cite[30.10]{Kanamori:2003})

Das Resultat von D.A. Martin setzt das Axiom der abh"angigen Auswahl $\DC$ voraus. $\DC$ besagt, da"s jeder nicht-leere beschnittene Baum auf einer Menge $A$ einen unendlichen Pfad hat (\siehe etwa \cite[{20.3f}, S. 139]{Kechris:1995}). Eine gleichwertige Formulierung lautet:
\begin{alignat*}{1}
\FA[X]\FA[R](	&X\neq\LM\UND R\subset X\times X\UND \FA[x\in X]\EX[y\in X](\left\langle x,y\right\rangle\in R)\\
							&\Rightarrow\EX[f\in X^\omega]\FA[n<\omega](\left\langle f(n),f(n+1)\right\rangle\in R) )
\end{alignat*}
(\siehe etwa \cite[5.4f]{Jech:2003} oder \cite[11.1f]{Kanamori:2003}).

$\DC$ ist st"arker als $\AC[\omega]$, d.h. $\DC\Rightarrow\AC[\omega]$ (\siehe etwa \cite[5.6, S. 60]{Jech:2003}). Wie $\AC[\omega]$ ist jedoch auch $\DC$ vertr"aglich mit $\AD$ (\siehe etwa \cite[30.29]{Kanamori:2003} vgl. Satz \ref{Lemma_Konsistenz_AD_AC_omega}).\label{FN_GaleMartin}}
$\SIGMA{1}{2n}{}$-Determiniertheit. Somit ist dann auch $\BMSc{B}$ als $\PI{1}{2n}$-Spiel determiniert. Enth"alt nun $A$ keine nichtleere superperfekte Teilmenge, dann hat nach Korollar \ref{SatzSigmaBeschrCharSpiele_3} Spieler II eine Gewinnstrategie in $\BMSc{B}$. Nach einem Resultat von Y. N. Moschovakis\footnote{\siehe Fu"snote \ref{FN_Moschovakis}.} hat Spieler II dann eine Gewinnstrategie in $\DELTAlf{1}{2n+1}$ (\siehe \cite{Moschovakis:1973} oder \cite[6E.1.]{Moschovakis:1980}). Da diese Gewinnstrategie eine $\sigma$-Schranke f"ur $A$ definiert, hei"st dies, da"s $A$ eine $\sigma$-Schranke in $\DELTAlf{1}{2n+1}$ hat.
\end{beweis}

\clearemptydoublepage

\chapter{Verallgemeinerte Kategorie}\label{Kap_VerallgKategorie}
\setcounter{lemma}{0}
Die Baire-Kategorie und die $\sigma$-Kategorie haben im Baireraum eine ganz "ahnliche spieltheoretische Charakterisierung [\siehe Kapitel \ref{KapSpielCharMager} und \ref{KapSpielCharSigmaBeschr"ankt}]. Diese Gemeinsamkeit l"a"st sich nun dazu nutzen, um die verschiedenen Kategorien-Konzepte des Baire-Raumes in einer spieltheoretischen Definition mittels sogenannter \eemph{Bedingungsmengen} zu verallgemeinern [\siehe Kapitel \ref{Kapitel_BMager}].

\bigskip Zun"achst definieren wir in Kapitel \ref{Kap_VerallgSpiele} f"ur die verallgemeinerte Kategorie ein verallgemeinertes Spiel $\BMSd{\BM{B}}(A)$ f"ur Teilmengen $A$ des Raumes $X^\omega$. In dem verallgemeinerten Spiel werden die vorherigen Spiele $\BMSo{A}$ und $\BMSb{A}$  dahingehend verallgemeinert, da"s Spieler II nun keine endlichen Sequenzen oder nat"urliche Zahlen mehr spielt, sondern \eemph{Bedingungen} f"ur den weiteren Spielverlauf.

\bigskip In Kapitel \ref{Kapitel_BMager} dienen anschlie"send die Charakterisierungen aus Satz \ref{Satz_Heine_Borel_2} f"ur kompakte und aus Satz \ref{SatzSigmaBeschr"ankt} $(ii)$ f"ur $\sigma$-beschr"ankte Teilmengen des Baireraumes als Vorlagen f"ur die Definition der allgemeineren Begriffe \eemph{$\BM{B}$-dirgends dicht} und \eemph{$\BM{B}$-mager} [\siehe Definitionen \ref{DefBNirgendsDicht} und \ref{DefBMager}]. 

\bigskip In Kapitel \ref{Kap_BPerf} wird der Begriff der \eemph{$\BM{B}$-perfekten} Teilmenge des Bairerauems eingef"uhrt. Beispiel \ref{Bsp_BPerfEnthSuperperfTM} zeigt, da"s die nicht-leeren $\BM{B}$-perfekten Teilmengen eine Verallgemeinerung von Teilmengen des Baireraumes sind, die eine nicht-leere superperfekte Teilmenge enthalten.

\bigskip Abschlie"send wird in Kapitel \ref{KapSpielCharBMager} gezeigt, da"s sich die Theoreme aus den Kapiteln \ref{KapSpielCharMager} und \ref{KapSpielCharSigmaBeschr"ankt} analog auf die verallgemeinerte Kategorie "ubertragen lassen [\siehe Theoreme \ref{Theorem_Char_AllgSpiel_OhneZeugen} und \ref{Theorem_Char_AllgSpiel_MitZeugen}].

\bigskip Kapitel \ref{Kap_Resume} gibt eine "Ubersicht "uber die wichtigsten im Laufe der Arbeit entwickelten Begriffe und wie sie sich zueinander verhalten.

\bigskip Die Definitionen sowie teilweise auch die Aussagen der Lemmata und S"atze (ohne Beweise) finden sich in \eemph{\glqq On a notion of smallness for subsets of the Baire space\grqq\ } von A. Kechris \cite[6]{Kechris:1977}. Die Beispiele (ohne Beweise) stammen alle aus \cite[6]{Kechris:1977}. Theorem \ref{Theorem_Char_AllgSpiel_OhneZeugen} sowie die Aussage von Theorem \ref{Theorem_Char_AllgSpiel_MitZeugen} finden sich in \cite[3.1, 3.3]{Kechris:1977}.

\bigskip Sei $X$ in diesem Kapitel eine Menge mit mehr als einem Element, $(X,\POW[X])$ der topologische Raum $X$ mit der diskreten Topologie und $X^{\omega}$ der dazugeh"orige Produktraum der Folgen in $X$. F"ur $X=\omega$ erh"alt man so etwa den Baire-Raum und f"ur $X=\MN{0,1}$ den Cantorraum. Die Menge $X^{<\omega}$ sei analog zum Fall $X=\omega$ definiert als Menge der endlichen Sequenzen von Elementen aus $X$. Die Menge $X_*^{<\omega}$ bezeichne $X^{<\omega}$ ohne die leere Sequenz $()$. Die offenen Basismengen seien in $X^\omega$ analog zu den offenen Basismengen des Baireraumes definiert:
\[
O_u:=\{f\in X^\omega\mid u\prec f\}
\]
mit $u\in X^{<\omega}$.

\bigskip Mit der Metrik
\begin{equation*}
\dd(f,g):=\sum_{n=0}^{\infty}\frac{1}{2^{n+1}}\Kronecker(f(n),g(n))\\
\Kronecker(i,j):=
\begin{cases}
0 &\TT{falls }i=j\\
1 &\TT{falls }i\neq j
\end{cases}
\end{equation*}
f"ur $f,g\in X^\omega$ l"ast sich $X^\omega$ in gleicher Weise wie der Baireraum vollst"andig metrisieren. Nach dem Baire'schen Kategoriensatz \ref{BKS} ist $X^\omega$ daher ein Baire'scher Raum. 

\bigskip Ist $X$ abz"ahlbar, so hat $X^\omega$ eine abz"ahlbare Basis (\siehe vgl. Lemma \ref{BR->2.Axiom}) und ist somit nach Lemma \ref{2.Axiom->separabel} separabel und somit polnisch (\siehe vgl. Satz \ref{BR->polnisch}).

\bigskip In $X^\omega$ gilt analog zu Lemma \ref{Lemma_Baireraum_OffeneBasismengen}, da"s zwei offene Basismengen in $X^\omega$ entweder ineinander liegen oder einen leeren Durschnitt haben:
Seien $O_s$ und $O_t$ offene Basismengen in $X^\omega$ mit $O_s\cap O_t \neq \LM$, dann gilt
\begin{alignat*}{1}
&O_s\subset O_t \TT{ oder } O_t\subset O_s.
\end{alignat*}

\bigskip Die offenen Basismengen $O_u$ von $X^\omega$ sind auch abgeschlossen (analog zu Lemma \ref{BR:offeneBasismengen->abgeschlossen}).

\bigskip Ist die Menge $X$ nicht endlich, so sind die offenen Teilmengen des Raumes $X^\omega$ nicht kompakt (\siehe vgl. Lemma \ref{BR:offeneBasismengen->nichtkompakt}).

\section{Verallgemeinerte Spiele}\label{Kap_VerallgSpiele}
In den Spielen $\BMSo{A}$, $\BMSa{B}$, $\BMSb{A}$ und $\BMSc{B}$ l"a"st sich die Spielweise von II als das Stellen von gewissen Bedingungen an den weiteren Spielverlauf auffassen. Beispielsweise spielt II in der Situation $(s_0,k_1,s_1,\ldots,k_n,s_n)$ des Spieles $\BMSb{A}$ ein $k_{n+1}<\omega$. Nach den Regeln von $\BMSb{A}$ entspricht dies der Bedingung $s_{n+1}(0)>k_{n+1}$ f"ur den n"achsten Zug $s_{n+1}$ von I. In dem verallgemeinerten Spiel $\BMSd{\BM{B}}(A)$ wird dies nun dahingehend verallgemeinert, da"s Spieler II nicht mehr konkrete Bedingungen, sondern nur noch abstrakte Elemente $b$ (genannt \eemph{Bedingungen}) aus einer Menge $\BED{B}$ (genannt \eemph{Bedingungsmenge}) spielt. Die Bedingungen gen"ugen dabei drei nat"urlichen Eigenschaften: Monotonie, Unterscheidbarkeit und Reduzierbarkeit (diese Eigenschaften sind \eemph{nat"urlich} in dem Sinne, da"s einige bereits aus den vorherigen Kapiteln bekannte Konzepte diese Eigenschaften erf"ullen \siehe Beispiele \ref{bsp_6_2} und \ref{bsp_6_3}).  

\bigskip Zun"achst mu"s an dieser Stelle also der Begriff der \glqq Bedingungsmenge\grqq\ pr"azisiert werden.

\begin{definition}\header{Bedingungsmenge}\label{DefBedingungsmenge}\\
Sei $\BED{B}$ eine nicht-leere Menge\gs die Elemente von $\BED{B}$ hei"sen \emph{Bedingungen}\gs 
und $\BF{E}$ eine Funktion, die jedem Element $b\in\BED{B}$ eine nicht-leere Menge $\BF{\TT{\textsl{Erf}}}(b)$ oder $\BF{E}(b)$ von nicht-leeren endlichen Sequenzen von Elementen aus $X$ (mit $X$ wenigstens zweielementig) zuordnet:
\[%\begin{equation*}
\BF{E}:\PFEIL{\BED{B}}{}{ X^{<\omega}_*}.
\]%\end{equation*}
Man schreibt auch 
\[
u\erf_{\BF{E}} b \TT{ statt } u\in\BF{E}(b) \\(u\in X_*^{<\omega}, b\in\BED{B})
\]
und sagt \emph{$u$ erf"ullt $b$}. Der Einfachheit halber schreibt man auch $u\erf b$. Dann hei"st das Paar $\BM{B}:=(\BED{B},\BF{E})$ eine \emph{Bedingungsmenge f"ur $X^\omega$} oder einfach \emph{Bedingungsmenge} wenn gilt:
\begin{ITEMS}[arabic)]
\item (Monotonie:)\\
$\FA[u,v\in X_*^{<\omega}]\FA[b\in\BED{B}](u\prec v\UND u\erf_{\BF{E}} b\Rightarrow v\erf_{\BF{E}} b)$,
\item (Unterscheidbarkeit:)\\
$\FA[x\in X]\EX[b\in \BED{B}]\FA[u\in X_*^{<\omega}](u\erf_{\BF{E}}b\Rightarrow u(0)\neq x)$.

\item (Reduzierbarkeit:)\\%(Charakterisierbarkeit:)\\%(Bestimmtheit:)\\
Es gibt eine Abbildung $\l:\PFEIL{\BED{B}}{}{\omega}$
so, da"s gilt
\begin{multline*}
\FA[b\in\BED{B}]\FA[u,v\in X_*^{<\omega}][(u\prec v\UND u\not\erf_{\BF{E}} b\UND v\erf_{\BF{E}} b)\Rightarrow \\
\EX[b'\in\BED{B}](\underbrace{l(b')\lneq l(b)}_!\UND \FA[w\in X_*^{<\omega}](w\erf_{\BF{E}} b'\Leftrightarrow u\kon w\erf_{\BF{E}} b))].
\end{multline*}
\end{ITEMS} 
\end{definition}
Die Intuition hierzu ist folgende:
Die Bedingungsmenge $\BM{B}=(\BED{B},\BF{E})$ legt f"ur das noch zu definierende Spiel $\BMSd{\BM{B}}(A)$ 
(mit $A\subset X^{\omega}$) fest, welche die weiteren Zugm"oglichkeiten von Spieler I einschr"ankenden Bedingungen
Spieler II spielen kann. F"ur $\BM{B}$ fordert die obige Definition dabei im einzelnen:
\begin{ITEMS}[arabic)]
\item (Monotonie:)\\
Die Mengen $\BF{E}(b)\subset X_*^{<\omega}$ sind abgeschlossen gegen Erweiterungen $u\prec v$ von Elementen
$u\in\BF{E}(b)$.
\item (Unterscheidbarkeit:)\\
F"ur jedes $x\in X$ gibt es eine Bedingung $b_x\in\BED{B}$, die aussagt, da"s das erste Glied $u(0)$ einer
endlichen Sequenz $u\neq()$ ungleich $x$ ist. Daraus folgt sp"ater unter anderem, da"s jede einelementige Menge \eemph{klein} ist [\siehe Satz \ref{SatzNDichtUndBMager}].
\item (Reduzierbarkeit:)\\%Bestimmtheit:)\\
Im Falle einer Erweiterung $v$ von $u\neq()$, die eine Bedingung $b$ erf"ullt, die $u$ aber nicht erf"ullt, gibt es
eine Bedingung $b'$, die "uber endliche Sequenzen $w\neq()$ aussagt, ob $u\kon w$ die Bedingung $b$ erf"ullt oder nicht. Diese Bedingung $b'$ ist \eemph{reduziert} gegen"uber $b$ in dem Sinne, da"s $\l(b')\lneq \l(b)$ gilt.  
\end{ITEMS}

Von nun an stehen\gs falls nicht ausdr"ucklich anders erw"ahnt\gs $s, t, u, v, w$ und $s_i,t_i,\ldots$ (mit $i<\omega$) f"ur endliche sowie $f, g, h$ und $f_i,g_i,\ldots$ (mit $i<\omega$) f"ur unendliche Sequenzen von Elementen in einer vorgegebenen Menge $X$. Au"serdem stehen $b$ und $b_i$ (mit $i<\omega$) f"ur Bedingungen (Elemente aus einer vorgegebenen Bedingungsmenge $\BED{B}$). 

\begin{beispiel}\header{$X=\BED{B}=2$}\label{bsp_6_1}\\
Das Tupel $\BM{B}:=(\BED{B},\BF{E})$ mit
\begin{alignat*}{1}
X&:=2:=\MN{0,1},\\ 
\BED{B}&:=2,\\
u\erf_{\BF{E}}b&:\AQ u(0)=b \TT{ (f"ur }b\in\BED{B}),\\
\l(b)&:=0 \TT{ f"ur alle }b\in\BED{B}
\end{alignat*}
ist eine Bedingungsmenge f"ur $X^\omega$.
\end{beispiel}

\begin{beweis}
\medskip{\eemph{Monotonie:}} $u\prec v\UND u(0)=b\Rightarrow v(0)=b$.

\medskip{\eemph{Unterscheidbarkeit:}} 
Sei f"ur beliebiges $x\in 2$
\begin{equation*}
b_x:=
\begin{cases}
0 &\TT{falls }x=1,\\
1 &\TT{falls }x=0.\\
\end{cases}
\end{equation*}
Dann gilt f"ur beliebiges $u\in 2^{<\omega}$: $u\erf b_x\Rightarrow u(0)\neq x$.

\medskip{\eemph{Reduzierbarkeit:}}
Die Reduzierbarkeit ist trivialerweise erf"ullt, da der Fall $u\prec v$ mit $u(0)\neq b\UND v(0)=b$ nicht eintreten kann. 
\end{beweis}

\begin{beispiel}\header{$X$ $\TT{beliebig}^*$, $\BED{B}= X^{<\omega}$}\label{bsp_6_2}\\
Das Tupel $\BM{B}:=(\BED{B},\BF{E})$ mit
\begin{alignat*}{1}
X&\TT{ sei beliebig ($*:$ mit wenigstens zwei Elementen)},\\
\BED{B}&:=X_*^{<\omega},\\
u\erf_{\BF{E}}b&:\AQ u\succ b \TT{ (f"ur }b\in\BED{B}),\\
\l(b)&:=\lng(b) \TT{ f"ur alle }b\in\BED{B}
\end{alignat*}
ist eine Bedingungsmenge f"ur $X^\omega$.
\end{beispiel}
\begin{beweis}
\medskip{\eemph{Monotonie:}} $u\prec v\UND b\prec u\Rightarrow b\prec v$.

\medskip{\eemph{Unterscheidbarkeit:}} 
Sei f"ur beliebiges $x\in X$
\begin{equation*}
b_x:= (y) \TT{ mit }y\neq x.
\end{equation*}
Dann gilt f"ur beliebiges $u\in X^{<\omega}$: $u\erf b_x\Rightarrow ((y)\prec u\UND y\neq x) \Rightarrow u(0)\neq x$.

\medskip{\eemph{Reduzierbarkeit:}}
Seien $u,v,b\in X_*^{<\omega}=\BED{B}$ mit $u\prec v$ und $b\not\prec u\UND b\prec v$. Sei dann $b'\in\BED{B}$ das Erg"anzugsst"uck von $u$ zu $b$:
\[
u\kon b'=b.
\]
Dann gilt $\lng(b')\lneq\lng(b)$ ( $\lng(b')\neq\lng(b)$ da $u$ nicht leer ist) und f"ur $w\in X_*^{<\omega}$ erh"alt man:
\[
b'\prec w\Leftrightarrow b=u\kon b'\prec u\kon w.
\]
\end{beweis}

\begin{beispiel}\header{$X=\BED{B}=\omega$}\label{bsp_6_3}\\
Das Tupel $\BM{B}:=(\BED{B},\BF{E})$ mit
\begin{alignat*}{1}
X&:=\omega,\\
\BED{B}&:=\omega,\\
u\erf_{\BF{E}}b&:\AQ u(0)>b \TT{ (f"ur }b\in\BED{B}),\\
\l(b)&:=0 \TT{ f"ur alle }b\in\BED{B}
\end{alignat*}
ist eine Bedingungsmenge f"ur $X^\omega$.
\end{beispiel}
\begin{beweis}
\medskip{\eemph{Monotonie:}} $u\prec v\UND u(0)>b\Rightarrow v(0)>b$.

\medskip{\eemph{Unterscheidbarkeit:}} 
Sei f"ur beliebiges $x\in X$
\begin{equation*}
b_x:= x.
\end{equation*}
Dann gilt f"ur beliebiges $u\in X^{<\omega}$: $u\erf b_x\Rightarrow u(0)> x \Rightarrow u(0)\neq x$.

\medskip{\eemph{Reduzierbarkeit:}}
Die Reduzierbarkeit ist trivialerweise erf"ullt, da der Fall $u\prec v$ mit $u(0)\not >b\UND v(0)>b$ nicht eintreten kann. 
\end{beweis}

Nun wird es Zeit f"ur die Definition der verallgemeinerten Spiele:

\begin{definition}\header{verallgemeinertes Spiel ohne Zeugen $\BMSd{\BM{B}}(A)$}\\
Sei $\BM{B}:=(\BED{B},\BF{E})$ eine Bedingungsmenge. F"ur eine beliebige Menge $A\subset X^\omega$ definiert man dann
das \emph{verallgemeinerte Banach-Mazur-Spiel $\BMSd{\BM{B}}(A)$} wie folgt:
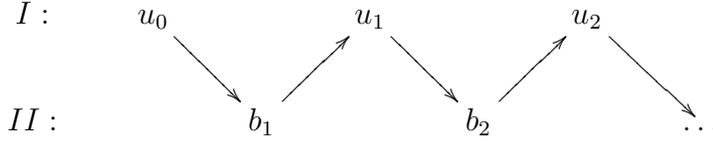
\[
\xymatrix{
I: &u_0 \ar[dr]  & 					 	 &u_1 \ar[dr]  &						&u_2 \ar[dr]\\
II:&		 				 &b_1 \ar[ur]	 &						 &b_2 \ar[ur]	&						&\ldots
}
\]
Spieler I und Spieler II spielen abwechselnd\gs Spieler I spielt eine nicht-leere endliche Sequenz $u_0\in X_*^{<\omega}$, Spieler II spielt eine Bedingung $b_1\in\BED{B}$, Spieler I spielt eine nicht-leere endliche Sequenz $u_1\in X^{<\omega}_*$, Spieler II spielt eine Bedingung $b_2\in\BED{B}$ u.s.w.. Sei nun $f:=u_0\kon u_1\kon u_2\kon\ldots$. Man sagt, da"s \emph{Spieler I gewinnt}, falls gilt:
\begin{ITEMS}[arabic)]
\item $f\in A$,
\item $\FA[i\geq 1](u_i\erf_{\BF{E}}b_i)$.
\end{ITEMS}
Ansonsten sagt man, da"s \emph{Spieler II gewinnt}. $A$ nennt man dabei auch die \emph{Gewinnmenge}.
\end{definition}

Die Begriffe \emph{(Gewinn-)Strategie f"ur I (bzw. f"ur II)} und \emph{determiniert} sind f"ur die Spiele $\BMSd{\BM{B}}(A)$ analog zu den Definitionen \ref{Def_Gewinnstrategie_I} (bzw. \ref{Def_Gewinnstrategie_II}) und \ref{Def_determiniert_G_A}  definiert.

\bigskip Sind $X$ und \BM{B} abz"ahlbar, l"a"st sich dieses Spiel mit einer Kodierung der in $\BMSd{\BM{B}}(A)$ gespielten Sequenzen und Bedingungen in die nat"urlichen Zahlen auch als ein Spezialfall der Grundvariante [\siehe Definition \ref{Def_SpielGrundvariante})] auffassen. Analog zu Lemma \ref{Lemma_Determiniertheit_G**p}, Lemma \ref{Lemma_Determiniertheit_G_tilde} und Lemma \ref{Lemma_Determiniertheit_G_tilde_p} gilt:

\begin{lemma}\header{Determiniertheit}\label{Lemma_Determiniertheit_G_B}\\
Falls $\AD$ gilt und die Mengen $X$ und $\BED{B}$ abz"ahlbar sind, so sind auch die obigen Spiele $\BMSd{\BM{B}}(A)$ determiniert.
\end{lemma}
\begin{beweis}\KOM[KURZ]
Analog zu Lemma \ref{Lemma_Determiniertheit_G**p}.
\end{beweis}

F"ur die Bedingungsmenge $\BM{B}$ aus Beispiel \ref{bsp_6_2}, $X:=\omega$ und $A\subset\BR$ erh"alt man mit $\BMSd{\BM{B}}(A)$ das Spiel $\BMSo{A}$ (ob nun wie in $\BMSd{\BM{B}}(A)$ inneinanderliegende endliche Sequenzen gespielt oder die Fortsetzungen aneinandergereiht werden wie in $\BMSo{A}$ ist unerheblich). Dieses Spiel diente in Kapitel \ref{KapSpielCharMager} zur Beschreibung der Baire-Kategorie. 

\bigskip F"ur die Bedingungsmenge $\BM{B}$ aus Beispiel \ref{bsp_6_3} und $A\subset\BR$ erh"alt man mit $\BMSd{\BM{B}}(A)$ das Spiel $\BMSb{A}$. Dieses Spiel diente in Kapitel \ref{KapSpielCharSigmaBeschr"ankt} zur Beschreibung der $\sigma$-Kategorie.

\section{$\BM{B}$-magere Mengen}\label{Kapitel_BMager}
\begin{comment}
Nachfolgend bezeichne $\BM{B}=(\BED{B},\BF{E})$ stets eine Bedingungsmenge f"ur $X^\omega$ (wobei $X$ wenigstens zwei Elemente habe).
\end{comment}

Die beiden Begriffe \eemph{mager} (aus Kapitel \ref{Kap_Magere_Mengen}) und \eemph{$\sigma$-beschr"ankt} (aus Kapitel \ref{Kapitel_SigmaBeschr}), die in der Baire- bzw. $\sigma$-Kategorie die kleinen Mengen (im Sinne von Definition \ref{Def_SigmaIdeal}) kennzeichnen, sollen nun in dem allgemeineren spieltheoretischen Begriff \eemph{$\BM{B}$-mager} vereinheitlicht werden. Dabei gehen wir so vor, da"s wir zun"achst den Begriff der \eemph{kompakten} Teilmenge mittels der Charakterisierung in Satz \ref{Satz_Heine_Borel_2} zum Begriff der \eemph{$\BM{B}$-nirgends dichten} Teilmenge verallgemeinern. Darunter fallen dann auch die \eemph{nirgends dichten und abgeschlossen} Teilmengen. Die $\BM{B}$-nirgends dichten Teilmengen haben die Idealeigenschaft. Um ein $\sigma$-Ideal zu erhalten, verallgemeinern wir in einem n"achsten Schritt die Charakterisierung der $\sigma$-beschr"ankten Teilmengen aus Satz \ref{SatzSigmaBeschr"ankt} $(ii)$ und erhalten den Begriff der \eemph{$\BM{B}$-mageren} Teilmenge, der sowohl die  mageren als auch die $\sigma$-beschr"ankten Teilmengen umfa"st.

\bigskip Die kompakten Teilmengen des Baireraumes lassen sich nach Satz \ref{Satz_Heine_Borel_2} vollst"andig charakterisieren als diejenigen abgeschlossenen Teilmengen $A$, deren B"aume $T_A$ "uberall endlich verzweigt sind:
\begin{alignat}{1}
&\FA[s\in T_A]\EX[k<\omega]\FA[t\in\omega^{<\omega}_*](s\kon t\in T_A\Rightarrow t(0)\not > k).
\end{alignat}
Als eine Verallgemeinerung dieser Charakterisierung erh"alt man den Begriff der $\BM{B}$-nirgends dichten Teilmengen:\KOM[WICHTIG]

\begin{definition}\header{$\BM{B}$-nirgends dicht}\label{DefBNirgendsDicht}\\
Sei $A\subset X^\omega$ abgeschlossen und $\BM{B}=(\BED{B},\BF{E})$ eine Bedingungsmenge. Dann hei"st $A$
\emph{abgeschlossen $\BM{B}$-nirgends dicht} oder einfach \emph{$\BM{B}$-nirgends dicht}, falls gilt:
\begin{equation*}
\FA[u\in T_A]\EX[b\in\BED{B}]\FA[v\in X^{<\omega}_*](u\kon v\in T_A\Rightarrow v\not\erf b).
\end{equation*}
\end{definition}

Die $\sigma$-beschr"ankten Teilmengen des Baireraumes sind nun laut Satz \ref{SatzSigmaBeschr"ankt} gerade diejenigen Teilmengen $A$, die sich als Teilmengen $\sigma$-kompakter Mengen darstellen lassen:
\begin{equation*}
A\subset\bigcup_{n<\omega}A_n\\\TT{wobei }\FA[n<\omega](A_n\subset \omega^\omega \TT{ kompakt}).
\end{equation*}
Analog dazu definiert man den allgemeineren Begriff der $\BM{B}$-mageren Teilmenge:

\begin{definition}\header{$\BM{B}$-mager}\label{DefBMager}\\
Sei $A\subset X^\omega$ beliebig und $\BM{B}=(\BED{B},\BF{E})$ eine Bedingungsmenge. Dann hei"st $A$
\emph{$\BM{B}$-mager}, falls $A$ Teilmenge einer abz"ahlbaren Vereinigung abgeschlossen $\BM{B}$-nirgends dichter 
Teilmengen von $X^\omega$ ist:
\begin{equation*}
A\subset\bigcup_{n<\omega}A_n\\\TT{ wobei }\FA[n<\omega](A_n\subset X^\omega \TT{ $\BM{B}$-nirgends dicht}).
\end{equation*}
\end{definition}

Die abz"ahlbare Vereinigung $\BM{B}$-nirgends dichter Teilmengen von $X^\omega$ ist i.a. nicht abgeschlossen und somit nicht wieder $\BM{B}$-nirgends dicht.\KOM[liegt nicht nur an abgeschl.] Offene Basismengen $O_u\subset X^\omega$ sind nicht $\BM{B}$-nirgends dicht, da es sonst f"ur $u\in T_{O_u}$ eine Bedingung $b\in\BED{B}$ gibt, mit $\FA[v\in X^{<\omega}_*](u\kon v\in T_{O_u}\Rightarrow v\not\erf b)$, da nun f"ur alle $v\in X^{<\omega}_*$ gilt, da"s $u\kon v\in T_{O_u}$, ist dies gleichbedeutend mit $\FA[v\in X^{<\omega}_*](v\not\erf b)$\gs laut Definition \ref{DefBedingungsmenge} ist die Menge der $v\in X^{<\omega}_*$, die eine Bedingung $b\in\BED{B}$ erf"ullen, aber nie leer. Da Teilmengen $\BM{B}$-nirgends dichter Mengen wieder $\BM{B}$-nirgends dicht sind (dies folgt direkt aus der Definition \ref{DefBNirgendsDicht}), k"onnen keine offenen Mengen $\BM{B}$-nirgends dicht sein (eine offene Menge enth"alt ja immer eine offene Basismenge).

\bigskip Eine Bedingungsmenge $\BM{B}$ hei"se \emph{gegen $\UND$ abgeschlossen}, falls:
\[
\FA[b_0,b_1\in\BED{B}]\EX[b\in\BED{B}]\FA[f\in X_*^{<\omega}][(f\erf b_0\UND f\erf b_1)\Leftrightarrow f\erf b]
\]

\begin{satz}\header{Ideal der $\BM{B}$-nirgends dichten Teilmengen}\\
Sei die Bedingungsmenge $\BM{B}$ gegen $\UND$ abgeschlossen. Dann bilden die $\BM{B}$-nirgends dichten Teilmengen in $X^\omega$ ein Ideal.
\end{satz}
\begin{beweis}Seien $A,B\subset X^\omega$.

\medskip Mit $B\subset A$ ist auch $T_B\subset T_A$. Aus der Definition \ref{DefBNirgendsDicht} von $\BM{B}$-nirgends dicht folgt daher direkt, da"s mit $A$ auch $B$ $\BM{B}$-nirgends dicht ist.

\medskip Wegen $T_{A\cup B}=T_A\cup T_B$ folgt direkt aus der Definition \ref{DefBNirgendsDicht} von $\BM{B}$-nirgends dicht, da"s mit $A,B$ $\BM{B}$-nirgends dicht auch $A\cup B$ $\BM{B}$-nirgends dicht sein mu"s.

\medskip Demnach bilden die $\BM{B}$-nirgends dichten Teilmengen in $X^\omega$ ein Ideal.
\end{beweis}

Die $\BM{B}$-mageren Teilmengen in $X^\omega$ verhalten sich nun so, wie wir es im einf"uhrenden Kapitel \ref{Kap_Kleine_Gro"se_Mengen} von \glqq kleinen\grqq\ Mengen gefordert haben: 

\begin{satz}\header{$\sigma$-Ideal der $\BM{B}$-mageren Teilmengen}\\
Die $\BM{B}$-mageren Teilmengen in $X^\omega$ bilden ein $\sigma$-Ideal.
\end{satz}
\begin{beweis}
Aus der Definition \ref{DefBMager} von $\BM{B}$-mager folgt direkt, da"s das System $\BM{B}$-magerer Teilmengen von $X^\omega$ abgeschlossen ist gegen"uber Teilmengen.

\medskip Da eine abz"ahlbare Vereinigung abz"ahlbarer Vereinigungen wieder eine abz"ahlbare Vereinigung ergibt, folgt aus der Definition \ref{DefBMager} von $\BM{B}$-mager, da"s das System $\BM{B}$-magerer Teilmengen von $X^\omega$ abgeschlossen gegen"uber abz"ahlbaren Vereinigungen ist.

\medskip Demnach bilden die $\BM{B}$-mageren Teilmengen in $X^\omega$ ein $\sigma$-Ideal.
\end{beweis}

\begin{beispiel}\header{In Bsp. \ref{bsp_6_1}: $\BM{B}$-mager $\AQ$ abz"ahlbar}\\
In Beispiel \ref{bsp_6_1} gilt mit $X:=\BED{B}:=2:=\{0,1\}$ und $u\erf b:\AQ u(0)=b$ f"ur $A\subset 2^\omega$:
\begin{ITEMS}
\item $A$ ist $\BM{B}$-nirgends dicht genau dann, wenn $A=\MN{f}$ f"ur ein $f\in X^\omega$.
\item $A$ ist $\BM{B}$-mager genau dann, wenn $A$ abz"ahlbar ist.
\end{ITEMS}
\end{beispiel}
\begin{beweis}
$(i)$ In Beispiel \ref{bsp_6_1} gilt:
\begin{alignat*}{2}
&&&A\TT{ $\BM{B}$-nirgends dicht}\\
&&\AQ[Def. \ref{DefBNirgendsDicht}]&\FA[u\in T_A]\EX[k<2]\FA[v\in 2^{<\omega}_*](u\kon v\in T_A\Rightarrow v(0)\neq k)\TT{ und }A \TT{ abg.}\\
&&\Leftrightarrow&A=\{f\}\TT{ f"ur ein $f\in 2^\omega$}
\end{alignat*}
Dabei sind Singleton-Mengen in $2^\omega$ abgeschlossen, da ihr Komplement $\C[\{f\}]=\bigcup_{u\not\prec f}O_u$ offen ist.

\medskip $(ii)$ In Beispiel \ref{bsp_6_1} gilt:
\begin{alignat*}{2}
&&&A\TT{ $\BM{B}$-mager}\\
&&\AQ[$(i)$]&A\subset\bigcup_{i<\omega}\{f_i\}\TT{ f"ur $(f_i\in 2^\omega)_i$}\\
&&\Leftrightarrow&A\TT{ ist abz"ahlbar}
\end{alignat*}
\end{beweis}

\begin{beispiel}\header{In Bsp. \ref{bsp_6_2}: $\BM{B}$-mager $\AQ$ mager}\label{BspBMagerMager}\\
In Beispiel \ref{bsp_6_2} gilt mit einer beliebigen Menge $X$ mit wenigstens zwei Elementen, $\BED{B}:=X^{<\omega}$ und $u\erf b:\AQ u\succ b$ f"ur $A\subset X^\omega$:
\begin{ITEMS}
\item $A$ ist $\BM{B}$-nirgends dicht genau dann, wenn $A$ abgeschlossen und nirgends dicht ist.
\item $A$ ist $\BM{B}$-mager genau dann, wenn $A$ mager ist.
\end{ITEMS}
\end{beispiel}
\begin{beweis}
$(i)\Leftrightarrow:$ Sei in Beispiel \ref{bsp_6_2} $A$ $\BM{B}$-nirgends dicht\gs d.h. (laut Definition \ref{DefBNirgendsDicht}):
\begin{alignat*}{1}
&\FA[u\in T_A]\EX[s\in X_*^{<\omega}]\FA[v\in X^{<\omega}_*](u\kon v\in T_A\Rightarrow s\not\prec v)\TT{ und }A \TT{ abg.}\\
\Leftrightarrow	&\FA[u\in T_A]\EX[s\in X^{<\omega}]\FA[v\in X^{<\omega}_*](u\kon v\in T_A\Rightarrow u\kon s\not\prec u\kon v)\TT{ und }A \TT{ abg.}\\
\Leftrightarrow	&\FA[O_u]\EX[O_{u\kon s}\subset O_u](O_{u\kon s}\not\subset A)\TT{ und }A \TT{ abg.}\\
\Leftrightarrow	&A \TT{ nirgends dicht}\TT{ und }A \TT{ abg.}
\end{alignat*}

\medskip $(ii)\Rightarrow:$ Nach $(i)$ sind in Beispiel \ref{bsp_6_2} $\BM{B}$-nirgends dichte abgeschlossene Mengen auch nirgends dicht. Also sind in Beispiel \ref{bsp_6_2} $\BM{B}$-magere Mengen (also Teilmengen abz"ahlbarer Vereinigungen $\BM{B}$-nirgends dichter Mengen) insbesondere mager (also Teilmengen abz"ahlbarer Vereinigungen nirgends dichter Mengen).

\medskip $(ii)\Leftarrow:$ In Beispiel \ref{bsp_6_2} sind magere Mengen auch $\BM{B}$-mager, da gilt:
\[
A\TT{ nirgends dicht}\Rightarrow \ABS[A]\TT{ nirgends dicht}.
\]

\medskip Es gilt $A\TT{ nirgends dicht}\Rightarrow \ABS[A]\TT{ nirgends dicht}$: 
Sei $A\subset X^\omega$ nirgends dicht. Sei etwa $G\subset \C[A]$ dicht und offen. Dann ist $G\cap\C[{\ABS[A]}]\subset\C[{\ABS[A]}]$ dicht und offen. Also ist $\ABS[A]$ nirgends dicht.
\begin{comment}%geht einfacher
Angenommen eine Teilmenge $A\subset X^\omega$ ist nirgends dicht (sei etwa $G\subset \C[A]$ dicht und offen) und $\ABS[A]$ ist nicht nirgends dicht. Dann folgt:
\[
\EX[{O_u\subset \ABS[A]}](O_u\cap G \TT{ offen}\subset\ABS[A]\o A).
\]

\medskip In $\ABS[A]\o A$ pa"st aber keine offene Menge $O_u$: Angenommen $O_u\subset \ABS[A]\cap \C[A]$. Dann ist insbesondere $A\subset \C[O_u]$ $(*)$ und somit folgt:
\[
O_u\subset\ABS[A]\o A\subset\ABS[A]\underset{\TT{\tiny{Def.}}}{=}\bigcap_{\substack{B\supset A\\B\TT{ \scriptsize{abg.}}}}B\underset{(*)}{\subset} \C[O_u].
\]
Das ist ein Widerspruch. Also pa"st in $\ABS[A]\o A$ keine offene Menge $O_u$ und es ist damit gezeigt da"s falls eine Menge $A$ nirgends dicht ist, dann auch $\ABS[A]$ nirgends dicht sein mu"s.
\end{comment}

\medskip Insgesamt gilt also, da"s in Beispiel \ref{bsp_6_2} magere Mengen auch $\BM{B}$-mager sind.
\end{beweis}

\begin{beispiel}\header{In Bsp. \ref{bsp_6_3}: $\BM{B}$-mager $\AQ$ $\sigma$-beschr"ankt}\label{Bsp_BNirgDicht_Kompakt}\\
In Beispiel \ref{bsp_6_3} gilt mit $X:=\BED{B}:=\omega$ und $u\erf b:\AQ u(0)>b$ f"ur $A\subset \omega^\omega$:
\begin{ITEMS}
\item $A$ ist $\BM{B}$-nirgends dicht genau dann, wenn $A$ kompakt ist.
\item $A$ ist $\BM{B}$-mager genau dann, wenn $A$ $\sigma$-beschr"ankt ist.
\end{ITEMS}
\end{beispiel}
\begin{beweis}
$(i)$: Da der Kompaktheitsbegriff gerade die Vorlage f"ur unsere Definition von \eemph{$\BM{B}$-nirgends dicht} war ist klar, da"s \eemph{$\BM{B}$-nirgends dicht} und \eemph{kompakt} in Beispiel \ref{bsp_6_3} "aquivalent sind [\siehe
Satz \ref{Satz_Heine_Borel_2} und Definition \ref{DefBNirgendsDicht}]. 

\medskip $(ii)$: Ebenso war der Begriff \eemph{$\sigma$-beschr"ankt} unser Vorbild f"ur die Definition des Begriffes \eemph{$\BM{B}$-mager}:
\begin{alignat*}{1}
&A \TT{ $\BM{B}$-mager}\\
&\AQ[Def. \ref{DefBMager}] A\subset\bigcup_{n<\omega}A_i \TT{ mit $A_i$ $\BM{B}$-nirgends dicht}\\
\intertext{und}
&A \TT{ $\sigma$-beschr"ankt}\\
&\AQ[Satz \ref{SatzSigmaBeschr"ankt} $(ii)$] A\subset\bigcup_{n<\omega}A_i \TT{ mit $A_i$ kompakt}
\end{alignat*}
Somit folgt $(ii)$ dann aus $(i)$.
\end{beweis}

\begin{satz}\header{$\BM{B}$-nirgends dicht, $\BM{B}$-mager}\label{SatzNDichtUndBMager}\\
Sei $A\subset X^\omega$ beliebig und $\BM{B}=(\BED{B},\BF{E})$ eine Bedingungsmenge. Dann gilt:
\begin{ITEMS}
\item $A$ einelementig $\Rightarrow$ $A$ $\BM{B}$-nirgends dicht.
\item $A$ abz"ahlbar $\Rightarrow$ $A$ $\BM{B}$-mager.
\item $A$ $\BM{B}$-nirgends dicht $\Rightarrow$ $A$ nirgends dicht.
\item $A$ $\BM{B}$-mager $\Rightarrow$ $A$ mager.
\end{ITEMS}
\end{satz}
\begin{beweis}
\point{(i) einelementig $\Rightarrow$ $\BM{B}$-nirgends dicht}

Sei $A\subset X^\omega$ einelementig\gs etwa $A=\MN{f}$ mit $f\in X^\omega$. Einelementige Teilmengen 
von $X^\omega$ sind abgeschlossen (s.o.). 

Es bleibt zu zeigen, da"s f"ur $A$ gilt: 
\[
\FA[u\in T_A]\EX[b\in\BED{B}]\FA[v\in X^{<\omega}_*](u\kon v\in T_A\Rightarrow v\not\erf b).
\] 
Es gilt $T_A=\MNG{u\in X^{<\omega}}{u\prec f}$. Dann gibt es nach Definition \ref{DefBedingungsmenge} 2) 
zu $f(\lng(u))\in X$\gs und somit zu jedem $u\in T_A$\gs eine Bedingung $b_u$ mit:  
\[
\FA[v\in X^{<\omega}_*](v\erf_{\BF{E}}b_u\Rightarrow v(0)\neq f(\lng(u))).
\]
Die Bedingung $b_u$ stellt also sicher, da"s ein $v$, das $b_u$ erf"ullt, $u$ nicht innerhalb von $f$ fortsetzt.
Es gilt:
\begin{alignat*}{1}
\FA[u\in T_A]\EX[b_u\in\BED{B}]\FA[v\in X^{<\omega}_*]&(v\erf_{\BF{E}}b_u\Rightarrow v(0)\neq f(\lng(u)))\\
\ldots																					&(v\erf_{\BF{E}}b_u\Rightarrow u\kon v\not\prec f)\\
\ldots																					&(v\erf_{\BF{E}}b_u\Rightarrow u\kon v\not\in T_A)\\
\ldots																					&(u\kon v\in T_A\Rightarrow v\not\erf_{\BF{E}}b_u)
\end{alignat*}
wobei ja f"ur $u\in X^{<\omega}_*$ gilt: $u=(u(0),\ldots u(\lng(u)-1))$.

\medskip Damit ist laut Definition \ref{DefBNirgendsDicht} gezeigt, da"s $A$ $\BM{B}$-nirgends dicht.

\point{(ii) abz"ahlbar $\Rightarrow$ $\BM{B}$-mager}

Sei $A\subset X^\omega$ abz"ahlbar. Dann ist $A$ abz"ahlbare Vereinigung einelementiger Teilmengen von $X^\omega$ und
somit nach (i) abz"ahlbare Vereinigung abgeschlossen $\BM{B}$-nirgends dichter Teilmengen von $X^\omega$. Nach 
Definition \ref{DefBMager} ist $A$ also $\BM{B}$-mager.

\point{(iii) $\BM{B}$-nirgends dicht $\Rightarrow$ nirgends dicht}
Sei $A$ $\BM{B}$-nirgends dicht. Also gilt:
\begin{equation}
\FA[u\in T_A]\EX[b\in\BED{B}]\FA[v\in X^{<\omega}_*](u\kon v\in T_A\Rightarrow v\not\erf b)\label{Frml_BNirgDicht}.
\end{equation}
Angenommen es gibt ein $O_u\subset A$. F"ur $u$ sei dann $b_u$ die Bedingung, f"ur die gilt:
\begin{equation}
\FA[v\in X^{<\omega}_*](u\kon v\in T_A\Rightarrow v\not\erf b).\label{Frml_BNirgDicht_2}
\end{equation}
Solch ein $b_u$ exitiert nach (\ref{Frml_BNirgDicht}). Laut Definition \ref{DefBedingungsmenge} enth"alt $\BF{E}(b_u)$ dann eine nicht-leere Sequenz $v\in X^{<\omega}_*$, d.h. es gibt eine nicht-leere Sequenz $v$, die $b_u$ erf"ullt. 
Aus $v\erf b_u$ folgt dann aber $u\kon v\not\in T_A$ (nach (\ref{Frml_BNirgDicht_2})). Da $A$ abgeschlossen ist (da $\BED{B}$-nirgends dicht) ist $T_A$ blattlos und es gibt ein $f\succ u\kon v$ in $O_u\o A$. Widerspruch zu $O_u\subset A$. Es gibt also kein $O_u\subset A$.

\medskip Damit ist $\INT[A]$ leer und $A$ nach Satz \ref{satz_nirgdicht} $(iv)$ nirgends dicht.

\point{(iv) $\BM{B}$-mager $\Rightarrow$ mager}

Sei $A$ $\BM{B}$-mager. Dann ist $A$ laut Definition \ref{DefBMager} Teilmenge einer abz"ahlbaren Vereinigung 
abgeschlossen $\BM{B}$-nirgends dichter Teilmengen von $X^\omega$, somit wegen (iii) Teilmenge einer abz"ahlbaren
Vereinigung nirgends dichter Teilmengen von $X^\omega$ und somit als Teilmenge einer mageren Teilmenge von $X^\omega$
selber mager (s.o.).
\end{beweis}

\begin{korollar}\header{kompakt $\Rightarrow$ nirgends dicht}\\
Kompakte Teilmengen des Baireraumes sind nirgends dicht.
\end{korollar}
\begin{beweis} 
\siehe Satz \ref{SatzNDichtUndBMager} $(iii)$ und Beispiel \ref{Bsp_BNirgDicht_Kompakt} $(i)$.
\end{beweis}

\section{$\BM{B}$-perfekte Mengen}\label{Kap_BPerf}
Die \eemph{gro"sen} Mengen der verallgemeinerten Kategorie sind die Komplemente $\BM{B}$-magerer Teilmengen des Baireraumes. Der Begriff der \eemph{$\BM{B}$-perfekten} Teilmenge des Bairerauems kennzeichnet analog zum Begriff der \eemph{superperfekten} (und nicht-leeren) Teilmenge \eemph{relativ gro"se} Teilmengen des Baireraumes. Beispiel \ref{Bsp_BPerfEnthSuperperfTM} zeichnet die nicht-leeren $\BM{B}$-perfekten Teilmengen aus als Verallgemeinerung von Teilmengen des Baireraumes, die eine nicht-leere superperfekte Teilmenge enthalten.

\begin{definition}\header{$\BM{B}$-dichte Teilmenge von $X^{<\omega}$}\label{Def_BDichteMenge}\\
Eine Menge $Q\subset X_*^{<\omega}$ hei"st \emph{$\BM{B}$-dicht} genau dann, wenn gilt:
\[
\FA[b\in\BED{B}]\EX[u\in X_*^{<\omega}]\EX[u'\succ u](u\in Q \UND u'\erf b).
\]
%f"ur jedes $b\in \BED{B}$ ein $u\in Q$ existiert so, da"s $u\erf b$.
\end{definition}

F"ur einen Baum $J$ auf $X^{<\omega}$ (also $J\subset(X^{<\omega})^{<\omega}$) bezeichne $J_*$ die Menge $J$ ohne die leere Sequenz $()\in (X^{<\omega})^{<\omega}$.

\begin{definition}\header{$\BM{B}$-perfekter Baum auf $X^{<\omega}$}\label{Def_BPerfekterBaum}\\
Sei $J$ ein Baum auf $X^{<\omega}$. $J$ hei"se \emph{$\BM{B}$-perfekt}, falls gilt:
\begin{ITEMS}[arabic)]
\item F"ur jedes $p\in J_*$ ist die Menge $\MNG{u\in X^{<\omega}_*}{p\kon(u)\in J}$ $\BM{B}$-dicht. Das hei"st laut Definition \ref{Def_BDichteMenge}: 
\[
\FA[p\in J_*]\FA[b\in\BED{B}]\EX[u\in X_*^{<\omega}]\EX[u'\succ u](p\kon (u)\in J \UND u'\erf b).
\]
\item F"ur jedes $p\in J$ gilt: falls $p\kon (u), p\kon(v) \in J$ und $u\neq v$, dann sind $u$ und $v$ inkompatibel ($u,v\in X^{<\omega}$ beliebig).
\end{ITEMS}
\end{definition}

Punkt $2)$ in Definition \ref{Def_BPerfekterBaum} garantiert, da"s f"ur ein $f\in X^\omega$ die Sequenz $(s_0,\ldots,s_n)\in J$ der L"ange $n+1$ mit $s_0\kon \ldots\kon s_n\prec f$ eindeutig bestimmt sind [\siehe Definition \ref{Def_[J]} und Lemma \ref{Lemma_Eindeutigkeit_Def_[J]}].

\begin{definition}\header{$[J]$ f"ur einen $\BM{B}$-perfekten Baum $J$ auf $X^{<\omega}$}\label{Def_[J]}\\
F"ur einen $\BM{B}$-perfekten Baum $J$ definiert man 
\begin{alignat*}{1}
[J]	&= \MNG{f\in X^\omega}{\FA[n<\omega]\EX[(s_0,\ldots,s_n)\in J](s_0\kon\ldots\kon s_n\prec f)}
\end{alignat*} 
\end{definition}

Man beachte den Unterschied zu den herk"ommlichen B"aumen $T$ auf einer Menge $X$ [\siehe Kapitel \ref{Kap_Baire_Raum}]. Ein Baum $T$ auf $X$ ist eine Teilmenge von $X^{<\omega}$, wohingegen in diesem Kapitel B"aume $J$ auf $X^{<\omega}$ betrachtet werden\gs also Teilmengen $J\subset (X^{<\omega})^{<\omega}$.

\begin{lemma}\header{Baum $\BM{B}$-perfekt $\Rightarrow$ Baum blattlos}\label{Blattlosigkeit}\\
Sei $J$ ein $\BM{B}$-perfekter Baum auf $X^{<\omega}$. Dann enth"alt $J$ keine Bl"atter.
\end{lemma}
\begin{beweis}
Sei $J$ ein $\BM{B}$-perfekter Baum auf $X^{<\omega}$. Laut Definition \ref{Def_BPerfekterBaum} ist f"ur jedes $p\in J_*$ die Menge $\MNG{u\in X^{<\omega}_*}{p\kon(u)\in J}$ $\BM{B}$-dicht. Das hei"st: F"ur jedes $p\in J_*$ und f"ur jedes $b\in\BED{B}$ existiert ein $u\in X^{<\omega}_*$ mit $p\kon(u)\in J$. Da f"ur eine Bedingungsmenge $\BM{B}:=(\BED{B},\BF{E})$ laut Definition \ref{DefBedingungsmenge} die Menge $\BED{B}$ als nicht-leer vorrausgesetzt wird, enth"alt $J$ keine Bl"atter. 
\end{beweis}

\begin{lemma}\header{Eindeutigkeit der $(s_0,\ldots,s_n)$ in Definition \ref{Def_[J]}}\label{Lemma_Eindeutigkeit_Def_[J]}\\
Sei $J$ ein $\BM{B}$-perfekter Baum auf $X^{<\omega}$ und $f\in [J]$. Dann gilt f"ur
$n < m<\omega$ und $(s_0,\ldots,s_n), (t_0,\ldots,t_m) \in J$ mit $s_0\kon\ldots\kon s_n\prec f$ sowie $t_0\kon\ldots\kon t_m\prec f$, da"s $(s_0,\ldots,s_n) = (t_0,\ldots,t_n)$.
\end{lemma}
\begin{beweis}
Seien $n < m<\omega$ und $(s_0,\ldots,s_n), (t_0,\ldots,t_m) \in J$ mit $s_0\kon\ldots\kon s_n\prec f$ sowie $t_0\kon\ldots\kon t_m\prec f$. Dann sind $s_0$ und $t_0$ nicht inkompatibel und somit $s_0=()\kon s_0 \EQ[\ref{Def_BPerfekterBaum} 2)] ()\kon t_0= t_0$. Dann sind $s_1$ und $t_1$ ebenfalls nicht inkompatibel und somit $s_0\kon s_1 \EQ[\ref{Def_BPerfekterBaum} 2)] t_0\kon t_1$ und wegen $s_0=t_0$ gilt $s_1=t_1$, u.s.w.\gs insgesamt gilt also $(s_0,\ldots,s_n) = (t_0,\ldots,t_n)$.
\end{beweis}

Sei $J$ ein $\BM{B}$-perfekter Baum auf $X^{<\omega}$ und $g\in [J]$. Dann gibt es laut Definition \ref{Def_[J]} f"ur jedes $n<\omega$ ein $(s_0,\ldots,s_n)\in J$ mit $s_0\kon\ldots\kon s_n\prec g$ und die $s_i\in X^{<\omega}$ sind gem"a"s Lemma \ref{Lemma_Eindeutigkeit_Def_[J]} unabh"angig von $n$ eindeutig bestimmt. Es bezeichne \[ g[i]:=s_i \] die $i$-te endliche Sequenz in $(s_0,\ldots,s_n)\in J$ mit $s_0\kon\ldots\kon s_n\prec g$.

\begin{definition}\header{$\BM{B}$-perfekte Teilmenge von $X^\omega$}\label{Def_BPerfTM}\\
Eine Teilmenge $A$ von $X^\omega$ hei"st \emph{$\BM{B}$-perfekt}, falls $A=[J]$ f"ur einen $\BM{B}$-perfekten
Baum $J$ auf $X^{<\omega}$.
\end{definition}

Man beachte, da"s eine $\BM{B}$-perfekte Menge $A$ nicht abgeschlossen sein mu"s, denn aus $A=[J]$ f"ur einen $\BM{B}$-perfekten Baum $J$ folgt (im Gegensatz zu den herk"ommlichen B"aumen) i.a. nicht, da"s $A$ abgeschlossen ist.

\bigskip Nachdem wir die folgenden beiden Lemmata gezeigt haben, k"onnen wir damit anschlie"send den Begriff der $\BM{B}$-perfekten Teilmenge in einer Reihe von Beispielen zu einigen anderen Begriffen von \eemph{relativ gro"sen} Mengen in Beziehung setzen. Insbesondere zeigt Beispiel \ref{Bsp_BPerfEnthSuperperfTM}, da"s nicht-leere $\BM{B}$-perfekte Teilmengen eine Verallgemeinerung von Teilmengen des Baireraumes sind, die eine nicht-leere superperfekte Teilmenge enthalten.

\begin{lemma}\header{$\BM{B}$-perfekt $\Rightarrow$ $\G_\delta$}\label{lemma_Bperfekt_G_delta}\\
Jede $\BM{B}$-perfekte Menge $A$ ist eine $\G_\delta$-Menge\footnote{\emph{$\G_\delta$-Menge:} Anstatt $\PI{0}{2}{}$ schreibt man auch $\G_\delta$. Eine $\G_\delta$-Menge ist also ein abz"ahlbarer Schnitt offener Mengen [\siehe vgl. Kapitel \ref{KapSpielCharMager} Fu"snote \ref{Def_GDeltaFSigma}].}.
\end{lemma}

\begin{beweis}
Sei $A=[J]$ f"ur einen $\BM{B}$-perfekten Baum J auf $X^{<\omega}$\gs dann gilt:
\begin{alignat*}{1}
A	&=[J]\\
	&\EQ[\ref{Def_[J]}] \MNG{f\in X^\omega}{\FA[n<\omega]\EX[(s_0,\ldots,s_n)\in J](s_0\kon\ldots\kon s_n\prec f)}\\
	&=\bigcap_{n<\omega}\bigcup_{(s_0,\ldots,s_n)\in J} O_{s_0\kon\ldots\kon s_n}.
\end{alignat*} 
\end{beweis}

\newpage%%%%%%%%%%%%%%%%%%%%%%

\begin{lemma}\header{Abgeschlossene $\BM{B}$-perfekte Teilmenge}\label{lemma_abg_Bperf_TM}\\
Falls f"ur $\BM{B}$ gilt, da"s die Erf"ullung einer Bedingung $b\in \BED{B}$ durch eine endliche Sequenz $u\in X^{<\omega}$ nur von der ersten Stelle von $u$ abh"angt\gs genauer:
\begin{alignat*}{1}
&\FA[u\in X^{<\omega}](u\erf b \AQ (u(0))\erf b),
\end{alignat*}
dann enth"alt jede $\BM{B}$-perfekte Menge $A$ eine abgeschlossene $\BM{B}$-perfekte Teilmenge $B\subset A$ und falls $A\neq\LM$ so ist auch $B\neq\LM$.
\end{lemma}

\begin{beweis}
Sei $A=[J]$ f"ur einen $\BM{B}$-perfekten Baum $J$ auf $X^{<\omega}$. Ist $A$ abgeschlossen, so folgt die Behauptung. Sei nun $A$ nicht abgeschlossen. Das ist gleichbedeutend damit, da"s es eine Folge $g_0, g_1, \ldots$ in $A$ gibt, die gegen ein $f\in X^{\omega}$ konvergiert mit $f\not\in A$. Um eine abgeschlossene $\BM{B}$-perfekte Teilmenge $B\subset A$ zu erhalten definiert man einen ebenfalls $\BM{B}$-perfekten Baum $\widetilde{J}$ so, da"s keine Folge mehr in $[\widetilde{J}]$ auftreten kann, die gegen ein Element au"serhalb von $[\widetilde{J}]$ konvergiert\gs dabei benutzen wir das Auswahlaxiom:

Falls eine Folge $(g_i)_i<\omega$ in $[J]$ gegen ein $f\not\in [J]$ konvergiert, so mu"s wegen $f\not\in [J]$ (laut Definition von $[J]$) gelten
\begin{alignat*}{1}
&\EX[n<\omega]\FA[(s_0,\ldots,s_n)\in J](s_0\kon\ldots s_n\not\prec f),
\end{alignat*}
und wegen $\FA[g\in\MNG{g_i}{i<\omega}]\FA[n<\omega]((g[0],\ldots, g[n])\in J)$ somit insbesondere
\begin{alignat}{1}
&\EX[n<\omega]\FA[g\in\MNG{g_i}{i<\omega}](g[0]\kon\ldots\kon g[n]\not\prec f)\label{Frml_n_Verzw}.
\end{alignat}
W"ahlt man $n$ in (\ref{Frml_n_Verzw}) minimal, so folgt in $J$ auf $(s_0\ldots s_{n-1})\in J$ mit $(s_0\ldots s_{n-1})\prec f$ eine unendliche Verzweigung $\MNG{g_i[n]}{i<\omega}$ (unendlich viele $g_i[n]$ m"ussen paarweise verschieden sein, da die $g_i$ gegen $f$ konvergieren). 

%\nopagebreak%%%%%%%%%%%%%%%%%%%%%%%%

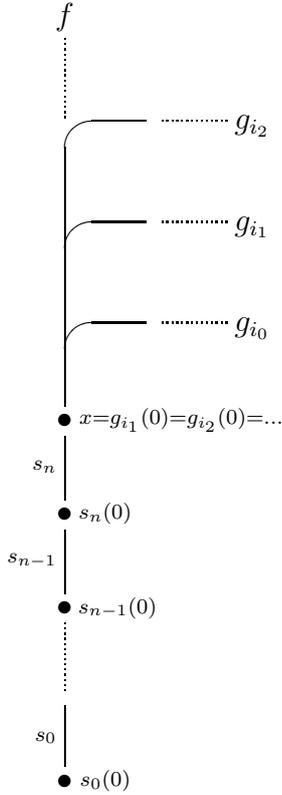
\begin{equation*}\xymatrix{ 
%{}&{}&{}&{}\\
{f}
&{}&{}&{}\\
{}\ar@{.} [u] 
&{}\ar@{.}[r]&{g_{i_2}}&{}\\
{}&{}\ar@{.}[r]&{g_{i_1}}&{}\\
{}&{}\ar@{.}[r]&{g_{i_0}}&{}\\
{\bullet}\ar@{}[r]|>(1.43){x=g_{i_1}(0)=g_{i_2}(0)=\ldots}
\ar@{-} `u[u] [ur] 
\ar@{-} `u[uur] [uur]
\ar@{-} `u[uuur] [uuur]
\ar@{-}[d] _{s_{n}}
&{}&{}&{}\\
{\bullet}\ar@{}[r]|>(0.5){s_{n}(0)}
\ar@{-}[d] _{s_{n-1}}
&{}&{}&{}\\
{\bullet}\ar@{}[r]|>(0.65){s_{n-1}(0)}
\ar@{.}[d] 
&{}&{}&{}\\
{}
\ar@{-}[d] _{s_{0}}
&{}&{}&{}\\
{\bullet}\ar@{}[r]|>(0.5){s_0(0)}
&{}&{}&{}\\
}
\end{equation*}

Um dies f"ur einen Teilbaum $\widetilde{J}$ von $J$ auszuschlie"sen, definiert man f"ur $(s_0,\ldots,s_n)\in J$ und $x\in X$ eine Teilmenge $Y_{(s_0,\ldots, s_n),x}\subset X_*^{<\omega}$:
\begin{alignat*}{1}
&Y_{(s_0,\ldots ,s_n),x}:=\MNG{s\in X_*^{<\omega}}{(s_0,\ldots ,s_n,s)\in J \UND s(0)=x}.
\end{alignat*}
Zu jedem $(s_0,\ldots, s_n)\in J$ und der Familie der nicht-leeren $Y_{(s_0,\ldots, s_n),x}$ erh"alt man nun mittels Auswahlaxiom eine zugeh"orige Auswahlfunktion
\begin{alignat*}{1}
&\varphi_{(s_0,\ldots, s_n)}: \PFEIL{ \MNG{Y_{ (s_0,\ldots ,s_n) ,x}\neq\LM}{x\in X} }{}{ X_*^{<\omega} }.
\end{alignat*}
Sei nun 
\begin{alignat*}{1}
&\widetilde{J}:=\MNG{(s_0,\ldots ,s_n)\in J}{\FA[i\leq n](s_i\in \IMG( \varphi_{(s_0,\ldots, s_{i-1})})}.
\end{alignat*}
Dann gilt:

\point{$[\widetilde{J}]$ ist $\BM{B}$-perfekt}
Zu zeigen ist, da"s $\widetilde{J}$ ein $\BM{B}$-perfekter Baum ist. Da $J$ ein $\BM{B}$-perfekter Baum ist, ist f"ur jeden Verzweigungspunkt $p\in J_*$ die Menge $\MNG{u\in X_*^{<\omega}}{p\kon (u)\in J}$ der Verzweigungen in $J$ $\BM{B}$-dicht. Das hei"st:
\[
\FA[p\in J_*]\FA[b\in\BED{B}]\EX[u\in X_*^{<\omega}]\EX[u'\succ u](p\kon (u)\in J\UND u'\erf b).
\]
Da laut Vorraussetzung $u\erf b$ nur von $u(0)=u'(0)$ abh"angt, ist dies "aquivalent zu:
\[
\FA[p\in J_*]\FA[b\in\BED{B}]\EX[u\in X_*^{<\omega}](p\kon (u)\in J\UND u\erf b).
\]
Da $u\erf b$ nur von $x:=u(0)$ abh"angt, und wir f"ur die Konstruktion von $\widetilde{J}$ mittels der Auswahlfunktionen an jedem Verzweigungspunkt $p$ in $J$ f"ur jedes $x\in X$, f"ur das in $J$ eine Verzweigung $u_x$ mit $u_x(0)=x$ exisiert, genau ein solches $u_x$ ausgew"ahlt haben, gibt es auch in $\widetilde{J}$ an jedem Verzweigungspunkt $p\in\widetilde{J}_*$ f"ur jede Bedingung $b\in \BED{B}$ eine Verzweigung $u\in X_*^{<\omega}$  mit $p\kon(u)\in \widetilde{J}$ und $u\erf b$. Also ist $\widetilde{J}$ ebenfalls $\BM{B}$-perfekt.

\point{$[\widetilde{J}]$ ist abgeschlossen}
Falls eine Folge $(g_i)_i<\omega$ in $[\widetilde{J}]$ gegen ein $f\not\in [\widetilde{J}]$ konvergiert, so mu"s wegen $f\not\in [\widetilde{J}]$ (laut Definition von $[\widetilde{J}]$) gelten
\begin{alignat*}{1}
&\EX[n<\omega]\FA[(s_0,\ldots,s_n)\in \widetilde{J}](s_0\kon\ldots s_n\not\prec f),
\end{alignat*}
und wegen $\FA[g\in\MNG{g_i}{i<\omega}]\FA[n<\omega]((g[0],\ldots, g[n])\in \widetilde{J})$ somit insbesondere
\begin{alignat}{1}
&\EX[n<\omega]\FA[g\in\MNG{g_i}{i<\omega}](g[0]\kon\ldots\kon g[n]\not\prec f)\label{Frml_n_Verzw_1}.
\end{alignat}
W"ahlt man $n$ in (\ref{Frml_n_Verzw_1}) minimal, so folgt in $\widetilde{J}$ auf $(s_0\ldots s_{n-1})\in \widetilde{J}$ mit $(s_0\ldots s_{n-1})\prec f$ eine unendliche Verzweigung $\MNG{g_i[n]}{i<\omega}$ (unendlich viele $g_i[n]$ m"ussen paarweise verschieden sein, da die $g_i$ gegen $f$ konvergieren).  Da die $g_i$ gegen $f$ konvergieren, m"ussen auch unendlich viele der $g_i[n]$ die selbe erste Stelle $g_i[n](0)$ haben\gs im Widerspruch zur Konstruktion von $\widetilde{J}$.

\point{$B\neq\LM$ falls $A\neq\LM$} Da wir bei der Konstruktion von $\widetilde{J}$ Elemente aus $J$ ausgew"ahlt haben, falls $J$ nicht leer war, gilt $B\EQ[\scriptsize{def.}][\widetilde{J}]\neq\LM$ falls $A\EQ[\scriptsize{def.}][J]\neq\LM$.
\end{beweis}

\begin{beispiel}\header{In Bsp. \ref{bsp_6_1}: $\BM{B}$-perfekt $\neq\LM$ $\Rightarrow$ enth. perf. Teilm. $\neq\LM$}\\
In Beispiel \ref{bsp_6_1} gilt mit $X:=\BED{B}:=2:=\{0,1\}$ und $u\erf b:\AQ u(0)=b$ f"ur $A\subset 2^\omega$:
\begin{alignat*}{1}
\TT{$A\neq\LM$ ist $\BM{B}$-perfekt $\Rightarrow$} &\TT{ $A$ enth"alt eine perfekte Teilmenge $B\neq\LM$}
\end{alignat*}
wobei $B$ perfekt hei"se, da"s $B$ perfekt ist bzgl. der Relativtopologie und abgeschlossen in $2^\omega$ [\siehe S. \pageref{Def_perfekteTM}].
\end{beispiel}
\begin{beweis}
Sei $A\neq\LM$ $\BM{B}$-perfekt\gs etwa $A=[J]$ mit $J$ $\BM{B}$-perfekter Baum auf $X^{<\omega}$. In Beispiel \ref{bsp_6_1} gilt $u\erf b\AQ (u(0))\erf b$\gs nach Lemma \ref{lemma_abg_Bperf_TM} enth"alt $A$ also eine abgeschlossene $\BM{B}$-perfekte Teilmenge $B\subset A$ und $B\neq\LM$, da $A\neq\LM$. Dann ist $B$ perfekt (d.h. abgeschlossen und perfekt bzgl. der Relativtopologie):

\point{$B$ abgeschlossen}
Dies gilt nach Lemma \ref{lemma_abg_Bperf_TM}.

\point{$B$ perfekt bzgl. der Relativtopologie}
Sei $B$ $\BM{B}$-perfekt\gs etwa $B=[\widetilde{J}]$. Dann enth"alt $B$ bez"uglich der Relativtopologie in $2^\omega$ keine isolierten Punkte: 

Sei ein $f\in O_t\cap B$ f"ur eine offene Basismenge $O_t\subset 2^\omega$. Da $t\in T_B$ eine Anfangssequenz einer Sequenz $p[0]\kon\ldots\kon p[n]$ mit $p\in \widetilde{J}$ ist und $\widetilde{J}$ $\BM{B}$-perfekt ist gilt:
\[
\FA[b\in \BED{B}]\EX[u\in X_*^{<\omega}]\EX[u'\succ u](p\kon(u)\in J\UND u'\erf b).
\]
Da $u\erf b\Leftrightarrow (u(0))=(u'(0))\erf b$ gilt, gibt es zu jedem $b\in \BED{B}=2$ ein $u\in X_*^{<\omega}$ so, da"s $u\erf b$ (d.h. $u(0)=b$) und $p\kon (u)\in J$. W"ahlt man $b\neq f[n+1](0)$, so ist $p[0]\kon\ldots\kon p[n]\kon u\not\prec f$. Da es in $\widetilde{J}$ nach Lemma \ref{Blattlosigkeit} keine Bl"atter gibt ist dann $p[0]\kon\ldots\kon p[n]\kon u\prec g$ f"ur ein $g\neq f$ aus $O_t\cap [\widetilde{J}]$. Damit ist gezeigt, da"s $B$ bez"uglich der Relativtopologie in $2^\omega$ keine isolierten Punkte enth"alt.
\end{beweis}

\begin{comment}
\begin{beispiel}\header{In Bsp. \ref{bsp_6_2}: $\BM{B}$-perfekt $\neq\LM$ $\Rightarrow$ dicht in off. Umgeb. $\neq\LM$ und $\G_\delta$}\\
In Beispiel \ref{bsp_6_2} gilt mit einer beliebigen Menge $X$ mit wenigstens zwei Elementen, $\BED{B}:=X^{<\omega}$ und $u\erf b:\AQ u\succ b$ f"ur $A\subset X^\omega$:
\[
\TT{$A\neq\LM$ ist $\BM{B}$-perfekt $\Rightarrow$ $\EX[U\subset X^\omega \TT{ offen}] (A\subset U \UND A \TT{ dicht in }U)$ und $A$ ist $\G_\delta$.}
\]
\end{beispiel}
\begin{beweis}
Sei $A\neq\LM$ $\BM{B}$-perfekt\gs etwa $A=[J]$ mit $J$ $\BM{B}$-perfekter Baum auf $X^{<\omega}$.

\point{$\EX[U\subset X^\omega \TT{ offen}] (A\subset U \UND A \TT{ dicht in }U)$} 
Da $A=[J]$ nicht leer ist, ist auch die offene Menge $U:=\bigcup_{(s)\in J_*}O_s$ nicht leer. Es gilt $[J]\subset U$. F"ur eine beliebige offene Basismenge $O_u\subset U$ mu"s $u=s\kon u'$ gelten f"ur ein $(s)\in J_*$ und $u`\in X^{<\omega}$ (laut Definition von $U$). Da $J$ $\BM{B}$-perfekt ist, gibt es insbesondere zu $p:=(s)\in J_*$ und $b:=u'\in \BED{B}=X^{<\omega}$ ein $s'\in X^{<\omega}$ mit $s'\erf b$ (also $s\kon u'\prec s\kon s'$) und $(s)\kon (s')\in J$. Da $J$ blattlos ist, mu"s es ein $g\in [J]$ geben mit $s\kon s'\prec g$\gs f"ur dieses $g$ gilt: $g\in [J]\cap O_u$ ($g\in O_u$ wegen $u=s\kon u'\prec s\kon s'\prec g$). Da $O_u\subset U$ als offene Basismenge in $U$ beliebig gew"ahlt war ist somit $A=[J]$ dicht in $U$.

\point{$A$ ist $\G_\delta$-Menge} 
Nach Lemma \ref{lemma_Bperfekt_G_delta} ist $A$ eine $\G_\delta$-Menge.
\end{beweis}
\end{comment}

\begin{beispiel}\header{In Bsp. \ref{bsp_6_3}: $\BM{B}$-perfekt $\neq\LM$ $\Rightarrow$ enth. superperf. Teilmenge $\neq\LM$}\label{Bsp_BPerfEnthSuperperfTM}
In Beispiel \ref{bsp_6_3} gilt mit $X:=\BED{B}:=\omega$ und $u\erf b:\AQ u(0)>b$ f"ur $A\subset \omega^\omega$:
\[
\TT{$A\neq\LM$ ist $\BM{B}$-perfekt $\Rightarrow$ $A$ enth"alt eine superperfekte Teilmenge $B\neq\LM$.}
\]
\end{beispiel}
\begin{beweis}
Sei $A\subset \omega^\omega$ nicht-leer und $\BM{B}$-perfekt. In Beispiel \ref{bsp_6_3} gilt $u\erf b\AQ (u(0))\erf b$\gs nach Lemma \ref{lemma_abg_Bperf_TM} enth"alt $A$ also eine abgeschlossene $\BM{B}$-perfekte Teilmenge $B\subset A$ und $B\neq\LM$, da $A\neq\LM$. Dann ist $B$ superperfekt (d.h. abgeschlossen und der Baum $T_B$ ist superperfekt):

\point{$T_B$ superperperfekt}
$B$ ist $\BM{B}$-perfekt\gs etwa $B=[J]$. Somit gilt
\[
\FA[p\in J_*]\FA[b\in \BED{B}]\EX[u\in X_*^{<\omega}]\EX[u'\succ u](p\kon(u)\in J\UND u'\erf b).
\]
Da in Bsp. \ref{bsp_6_3} $u\erf b\Leftrightarrow (u(0))\erf b$ gilt, gibt es dann zu jedem $p\in J_*$ und jedem $b\in \BED{B}=\omega$ ein $u\in \omega_*^{<\omega}$ mit $p\kon(u)\in J$ und $u\erf b$ (d.h. $u(0)>b$)\gs es existieren also unendlich viele Verzweigungen $u\in \omega_*^{<\omega}$ mit $p[0]\kon\ldots\kon p[n]\kon u\in T_B$. Da jedes $t\in T_B$ ein Anfangsst"uck eines $p\in J$ ist, ist $T_B$ superperperfekt.
\end{beweis}

\section{Charakterisierung durch verallgemeinerte Spiele}\label{KapSpielCharBMager}
Abschlie"send soll nun gezeigt werden, da"s sich die Theoreme aus den Kapiteln \ref{KapSpielCharMager} und \ref{KapSpielCharSigmaBeschr"ankt} analog auf die verallgemeinerte Kategorie "ubertragen lassen.
Die verallgemeinerten Spiele $\BMSd{\BM{B}}(A)$ ohne Zeugen sind bereits in Kapitel \ref{Kap_VerallgSpiele} definiert worden. Dort wurden die Begriffe \eemph{kompakt} und \eemph{$\sigma$-beschr"ankt} ja gerade im Hinblick auf dieses Spiel verallgemeinert. Das entsprechende projektive Spiel ist in Kapitel \ref{Kap_VerallgSpiel_2} noch zu definieren. 

\subsection*{Verallgemeinerte Spiele ohne Zeugen $\BMSd{\BM{B}}(A)$}
Das folgende Theorem bietet nun eine spieltheoretische Charakterisierung der verallgemeinerten Kategorie analog zu den spieltheoretischen Charakterisierungen der Baire- und der $\sigma$-Kategorie durch die Theoreme \ref{SatzMagerCharSpiele} bzw. \ref{SatzSigmaBeschrCharSpiele}:

\begin{theorem}\header{Charakterisierung durch Spiele ohne Zeugen}\label{Theorem_Char_AllgSpiel_OhneZeugen}\\
Sei $\BM{B}$ eine abz"ahlbare Bedingungsmenge auf $X$ (das hei"se: $X$ und $\BED{B}$ abz"ahlbar) und $A\subset X^\omega$. Dann gilt:
\begin{ITEMS}
\item Spieler I hat eine Gewinnstrategie in $\BMSd{\BM{B}}(A)$ $\Leftrightarrow$ $A$ enth"alt eine nicht-leere $\BM{B}$-perfekte Teilmenge $B$.
\item Spieler II hat eine Gewinnstrategie in $\BMSd{\BM{B}}(A)$ $\Leftrightarrow$ $A$ ist $\BM{B}$-mager.
\end{ITEMS}
\end{theorem}
\begin{beweis}\KOM[FERTIG]
(i)$\Rightarrow$:
Habe I eine Gewinnstrategie $\sigma$. Sei
\begin{alignat*}{2}
s_0	&:=s_{()}						&&:=\sigma()\\
s_1	&:=s_{(b_0)}				&&:=\sigma(s_0,b_0)\\
s_2	&:=s_{(b_0,b_1)}		&&:=\sigma(s_0,b_0,s_1,b_1)\\
&\ldots
\end{alignat*}
$s_i$ sei also eine abk"urzende Schreibweise f"ur $s_{(b_0,\ldots,b_{i-1})}$ und  h"ange von $(b_0,\ldots,b_{i-1})$ ab. Dann erf"ullt der Baum der $\sigma$-gespielten endlichen Sequenzen $J_\sigma:=\MNG{(s_0,\ldots,s_n)}{b_0,\ldots,b_{n-1}\in\BED{B}}$ die erste Bedingung f"ur einen $\BM{B}$-perfekten Baum auf $X^{<\omega}$ [\siehe Definition \ref{Def_BPerfekterBaum}], weil $\sigma$ eine Gewinnstrategie von I ist\gs das hei"st ja: an jeder Verzweigungsstelle $p:=(s_0,\ldots,s_n)\in {J_\sigma}_*$ und f"ur jede Bedingung $b_n\in\BED{B}$ gibt es ein $s\in X_*^{<\omega}$ mit $s\erf b_n$ und $p\kon (s)\in {J_\sigma}$, insbesondere gilt mit $s':=s$ also Bedingung $1)$ aus Definition \ref{Def_BPerfekterBaum}:
\begin{equation}
\FA[p\in {J_\sigma}_*]\FA[b\in\BED{B}]\EX[s\in X_*^{<\omega}]\EX[s'\succ s](p\kon (s)\in J \UND s'\erf b).\label{B_Perf_1}
\end{equation}

Der Baum $J_\sigma$ l"a"st sich nun so beschneiden, da"s der daraus entstehende Teilbaum $\widetilde{J_\sigma}$ zus"atzlich zur ersten auch die zweite Bedingung aus der Definition eines $\BM{B}$-perfekten Baumes erf"ullt. F"ur ein $p\in J_\sigma$, $s\in X_*^{<\omega}$ und $b\in\BED{B}$ sei
\[
\varphi(s,b,p):\equiv\EX[s'\succ s](p\kon(s)\in J_\sigma\UND s'\erf b).
\]
Da $J_\sigma$ die erste Bedingung $(\ref{B_Perf_1})$ aus der Definition eines $\BM{B}$-perfekten Baumes erf"ullt, gilt
f"ur alle $p\in J_\sigma$:
\[
\FA[b\in\BED{B}]\EX[s\in X_*^{<\omega}]\varphi(s,b,p).
\]

\medskip Sei nun $J_\sigma^0:=\{()\}\subset J_\sigma$. Sei $J_\sigma^n\subset J_\sigma$ bereits definiert, dann erh"alt man den Baum $J_\sigma^{n+1}\supset J_\sigma^n$ wie folgt: Wir definieren f"ur jedes Blatt $p\in J_\sigma^n$ eine Menge $A_p\subset J_\sigma$ und setzen dann 
\[
J_\sigma^{n+1}:=(\bigcup_{p\TT{ \scriptsize{Blatt in }} J_\sigma^n} A_p) \cup J_\sigma^n.
\]
Dabei erh"alt man die Menge $A_p$ f"ur ein Blatt $p$ in $J_\sigma^n$ wie folgt:

\medskip Sei $p$ ein Blatt in $J_\sigma^n$ und $b_0,b_1,\ldots$ eine Aufz"ahlung der Bedingungen aus $\BED{B}$. Sei zun"achst $A_p=\LM$. Wir f"ugen iterativ Knoten aus $J_\sigma$ zu $A_p$ hinzu:

\point{Schritt $0$} F"ur $b_0$ gilt $\EX[s\in X_*^{<\omega}]\varphi(s,b_0,p)$: w"ahle ein $s_0\in X_*^{<\omega}$ mit $\varphi(s_0,b_0,p)$ und f"uge $p\kon(s_0)$ zu $A_p$ hinzu.

\point{Schritt $i$ f"ur $i\geq 1$} F"ur $b_i$ gilt $\EX[s\in X_*^{<\omega}]\varphi(s,b_i,p)$\gs w"ahle ein $s_i\in X_*^{<\omega}$ mit $\varphi(s_i,b_i,p)$ und unterscheide:

\point{Fall 1} Gilt $\FA[p\kon (s)\in A_p](s\perp s_i)$, so f"uge auch $p\kon(s_i)$ zu $A_p$ hinzu.

\point{Fall 2} Gilt $\EX[p\kon (s)\in A_p](s\not\perp s_i)$\gs unterscheide:

\point{Fall 2.1} Gilt $\EX[p\kon (s)\in A_p](s\prec s_i)$: f"uge $p\kon(s_i)$ \eemph{nicht} zu $A_p$ hinzu.

\point{Fall 2.2} Gilt $\EX[p\kon (s)\in A_p](s\succ s_i\UND s\neq s_i)$: entferne alle $p\kon(s)$ aus $A_p$, f"ur die gilt, da"s $s\succ s_i\UND s\neq s_i$ und f"uge statt ihrer $p\kon(s_i)$ zu $A_p$ hinzu.

\medskip Sei nun 
\[
\widetilde{J_\sigma}:=\bigcup_{n<\omega} J_\sigma^n.
\]
Mit dieser Konstruktion ist dann sichergestellt, da"s $\widetilde{J_\sigma}$ zus"atzlich zur ersten auch die zweite Bedingung aus Definition \ref{Def_BPerfekterBaum} erf"ullt und somit $\BM{B}$-perfekt ist.

\begin{comment}%Beschneidungs-Vorschrift angeben
F"ur jedes $n\geq 1$ und jeden Knoten $(s_0,\ldots,s_{n-1})\in J_\sigma$ und jede aufsteigende Kette $s^0_n\prec s^1_n\prec s^2_n\prec\ldots$ mit $\FA[i<\omega]((s_0,\ldots,s_{n-1},s^i_{n})\in J_\sigma)$ entferne man aus der Menge $\MNG{(s_0,\ldots,s_{n-1},s^i_{n})\in J_\sigma}{i<\omega}$ alle bis auf $(s_0,\ldots,s_{n-1},s^0_{n})\in J_\sigma$. Geht man dabei iterativ bei $n=1$ beginnend vor, so erh"alt man auf diese Weise einen Teilbaum $\widetilde{J_\sigma}$ von $J_\sigma$. Dieser erf"ullt trivialerweise immer noch die erste Bedingung der Definition eines $\BM{B}$-perfekten Baumes: Sei etwa $b_i\in\BED{B}$ eine Bedingung, die von der durch die Beschneidung an der Stelle $(s_0,\ldots,s_{n-1})$ weggefallenen Sequenz $s^i_n\succ s^0_n$ erf"ullt wird\gs dann gilt $\EX[s'\succ s^i_n](s^0_n\prec s^i_n\prec s'\erf b_i)$ mit $(s_0,\ldots,s_{n-1},s^0_{n})\in \widetilde{J_\sigma}$. Zus"atzlich ist durch die Beschneidung auch Bedingung $(ii)$ aus der Definition eines $\BM{B}$-perfekten Baumes erf"ullt. Also ist $[\widetilde{J_\sigma}]$ eine $\BM{B}$-perfekte Menge. 
\end{comment}

\medskip Da $J_\sigma$ mit einer Gewinnstrategie von I definiert ist, gilt $J_\sigma\neq\LM$ und $[J_\sigma]\subset A$ und somit auch $\widetilde{J_\sigma}\neq\LM$ und $[\widetilde{J_\sigma}]\subset A$.

\bigskip(i)$\Leftarrow$:
Enthalte $A$ eine nicht-leere $\BM{B}$-perfekte Menge $[J]$ ($J$ $\BM{B}$-perfekter Baum auf $X^{<\omega}$). Es ist zu zeigen, da"s I eine Gewinnstrategie besitzt:

\medskip Spiele I zun"achst ein $(u_0)\in J$. Habe I nun bereits $p=(u_0,\ldots,u_{n-1})\in J$ gespielt und II die Bedingungen $(b_1,\ldots,b_{n-1})$ und gelte $u_i\erf b_i$ f"ur $i=1,\ldots,n-1$. Spielt nun II $b_n\in\BED{B}$, so gen"ugt zu zeigen:
\begin{equation}\label{Frml00}
\EX[u_n\in X_*^{<\omega}](p\kon(u_n)\in J\UND u_n\erf b_n),
\end{equation}
denn dann bewegt sich I wegen $[J]\subset A$ weiterhin in der Gewinnmenge $A$ und erf"ullt die n"achste Bedingung $b_n$. 

\medskip Da $J$ $\BM{B}$-perfekt ist, gilt f"ur $b_n$:
\begin{equation*}\label{Frml01}
\EX[u_n^0\in X_*^{<\omega}]\EX[u'\succ u_n^0](p\kon(u_n^0)\in J\UND u'\erf b_n).
\end{equation*}
Falls $u_n^0\erf b_n$, so gilt (\ref{Frml00}) mit $u_n:=u_n^0$. 

Falls hingegen $u_n^0\not\erf b_n$ gilt, folgt wegen der \eemph{Reduzierbarkeit} von $b_n$ [\siehe Definition \ref{DefBedingungsmenge} Punkt $3)$]:
\begin{equation}\label{Frml02}
\EX[b_n^1][\l(b_n^1)\lneq\l(b_n)\UND\FA[v\in X_*^{<\omega}](v\erf b_n^1\Leftrightarrow u_n^0\kon v\erf b_n)].
\end{equation}

\medskip Da $J$ $\BM{B}$-perfekt ist, gilt nun wiederum f"ur $b_n^1$:
\begin{equation*}\label{Frml03}
\EX[u_n^1\in X_*^{<\omega}]\EX[u'\succ u_n^1](p\kon(u_n^0)\kon(u_n^1)\in J\UND u'\erf b_n^1).
\end{equation*}
Falls $u_n^1\erf b_n^1$, so gilt (\ref{Frml00}) mit $u_n:=u_n^0\kon u_n^1$ wegen 
\[
u_n^1\erf b_n^1\AQ[(\ref{Frml02})]u_n^0\kon u_n^1\erf b_n.
\]
Falls hingegen $u_n^1\not\erf b_n^1$ gilt, folgt wie oben wegen der \eemph{Reduzierbarkeit} von $b_n^1$:
\begin{equation}\label{Frml04}
\EX[b_n^2][\l(b_n^2)\lneq\l(b_n^1)\UND\FA[v\in X_*^{<\omega}](v\erf b_n^2\Leftrightarrow u_n^1\kon v\erf b_n^1)]
\end{equation}
u.s.w.. 

Da bei jedem Reduzierungsschritt $\l(b_n^{i+1})\lneq \l(b_n^{i})$ gilt, l"a"st sich auf diese Weise nur endlich oft fortfahren. Falls eine derartige Kette von Reduzierungen nicht schon vorher dadurch abbricht, da"s eines der $u_n^i$ die Bedingung $b_n^i$ erf"ullt (wodurch dann (\ref{Frml00}) folgt), steht an ihrem Ende:
\begin{equation}\label{Frml05}
\EX[b_n^m][\l(b_n^m)=0\lneq\l(b_n^{m-1})\UND\FA[v\in X_*^{<\omega}](v\erf b_n^m\Leftrightarrow u_n^{m-1}\kon v\erf b_n^{m-1})].
\end{equation}
Da $J$ $\BM{B}$-perfekt ist, gilt f"ur $b_n^{m}$:
\begin{equation*}\label{Frml06}
\EX[u_n^{m}\in X_*^{<\omega}]\EX[u'\succ u_n^{m}](p\kon(u_n^0)\kon\ldots\kon(u_n^{m})\in J\UND u'\erf b_n^{m}).
\end{equation*}
Falls $u_n^{m}\erf b_n^{m}$, so gilt (\ref{Frml00}) mit $u_n:=u_n^0\kon\ldots\kon u_n^{m}$, wegen
\begin{alignat*}{1}
u_n^{m}\erf b_n^{m}&\AQ[(\ref{Frml05})] u_n^{m-1}\kon u_n^{m}\erf b_n^{m-1}\\
									 &\Leftrightarrow u_n^{m-2}\kon u_n^{m-1}\kon u_n^{m}\erf b_n^{m-2}\\
									 &\ldots\\
									 &\AQ[(\ref{Frml04})] u_n^1\kon\ldots\kon u_n^{m}\erf b_n^1\\
									 &\AQ[(\ref{Frml02})] u_n^0\kon\ldots\kon u_n^{m}\erf b_n.
\end{alignat*}

Der Fall $u_n^{m}\not\erf b_n^{m}$ kann nun nicht mehr eintreten, da dies die Voraussetzung f"ur die Reduzierung von $b_n^{m}$ w"are [\siehe Definition \ref{DefBedingungsmenge} Punkt $3)$] und somit gelten m"u"ste: $\EX[b_n^{m+1}](\l(b_n^{m+1})\lneq 0=\l(b_n^{m}))$ (im Widerspruch zu $\l:\PFEIL{\BED{B}}{}{\omega}$).

Damit ist (\ref{Frml00}) gezeigt.

\bigskip(ii)$\Leftarrow$:
Sei $A$ $\BM{B}$-mager\gs etwa $A\subset\bigcup_{i<\omega}A_i$ mit $A_i$ abgeschlossen $\BM{B}$-nirgends dicht (d.h. $\FA[u\in T_{A_i}]\EX[b\in\BED{B}]\FA[v\in X^{<\omega}_*](v\erf b\Rightarrow u\kon v\not\in T_{A_i})$). Spieler I spiele als erste endliche Sequenz $u_0\in X^{<\omega}$

\point{Liege $u_0$ in $\bigcap_{j\geq 1}T_{A_{i_j}}$ und in keinem anderen $T_{A_i}$}
Dann gibt es nach Definition \ref{DefBNirgendsDicht} (f"ur \eemph{$\BM{B}$-nirgends dicht}) Bedingungen $b_1,b_2,\ldots\in\BED{B}$ mit
\begin{alignat*}{1}
b_1\in\BED{B}:&\FA[v\in X^{<\omega}_*](v\erf b_1\Rightarrow u_0\kon v\not\in T_{A_{i_1}})\\
b_2\in\BED{B}:&\FA[v\in X^{<\omega}_*](v\erf b_2\Rightarrow u_0\kon u_1\kon v\not\in T_{A_{i_2}})\\
\ldots
\end{alignat*}
Wegen der \eemph{Monotonie} von $\erf$ in Definition \ref{DefBedingungsmenge} ($u\prec u'\UND u\erf b\Rightarrow u'\erf b$) gilt dann f"ur $f:=u_0\kon u_1\kon\ldots$ und beliebiges $n<\omega$:
\begin{alignat*}{1}
&u_0\kon\ldots\kon u_{n} \not\in T_{A_{i_1}}\cup\ldots\cup T_{A_{i_n}}
\end{alignat*}
und somit (da $A_{i_j}=[T_{A_{i_j}}]$) f"ur beliebiges $n<\omega$:
\begin{alignat*}{1}
&f \not\in A_{i_1}\cup\ldots\cup A_{i_n}
\end{alignat*}
also
\begin{alignat*}{1}
&f \not\in \bigcup_{j\geq 1}A_{i_j}
\end{alignat*}
und somit auch $f \not\in A$ (nach der Wahl der $T_{A_{i_j}}$).
Spielt Spieler II die $b_1,b_2,\ldots$ wie oben angegeben, stellt dies also eine Gewinnstrategie f"ur ihn dar.

\bigskip(ii)$\Rightarrow$:
Habe II eine Gewinnstrategie $\tau$. 

\point{Gute Sequenz} 
Eine Sequenz $p:=(u_0,b_1,u_1,b_2,\ldots,u_n,b_{n+1})$ mit $u_i\in X^{<\omega}_*$ f"ur $i=0,\ldots, n$ und $b_i\in\BED{B}$ f"ur $i=1,\ldots,n+1$ hei"se \emph{gut}, falls gilt:
\begin{alignat*}{1}
&u_i\erf b_i \TT{ f"ur }i=1,\ldots,n\UND\\
&b_j\TT{ mit $\tau$ gespielt f"ur }j=1,\ldots,n+1.
\end{alignat*}
Die leere Sequenz sei definitionsgem"a"s gut. 

\point{Sequenz gut f"ur $f\in A$}
$p$ hei"se \emph{gut f"ur $f\in A$}, falls gilt:
\begin{alignat*}{1}
&\EX[u_{n+1}](u_0\kon\ldots\kon u_n\kon u_{n+1}\prec f \UND u_{n+1}\erf b_{n+1}).
\end{alignat*}
Insbesondere ist die leere Sequenz gut f"ur jedes $f\in A$.

\medskip Sei nun $f\in A$ beliebig. Dann mu"s es eine Sequenz $p:=(u_0,b_1,u_1,b_2,\ldots,u_n,b_{n+1})$ geben (eventuell $p=()$) mit $p$ ist gut und $p$ ist gut f"ur $f$ und:
\begin{alignat*}{1}
&\neg\EX[(u_{n+1},b_{n+2})\in X^{<\omega}\times\BED{B}](p\kon (u_{n+1},b_{n+2}) \TT{ gut und gut f"ur } f).
\end{alignat*}
Ansonsten k"onnte Spieler I gewinnen, obwohl Spieler II mit seiner Gewinnstrategie $\tau$ spielt (Widerspruch). Dann gilt $f\in K_p$ mit
\begin{alignat}{2}
K_p &&:=\{f'\in X^\omega\mid &p\TT{ gut f"ur }f'\TT{ und }\notag\\
		&&					 &\neg\EX[(u_{n+1},b_{n+2})](p\kon (u_{n+1},b_{n+2}) \TT{ gut und gut f"ur } f') \}\notag\\
		&&=\{f'\in X^\omega\mid &\EX[u_{n+1}](u_0\kon\ldots\kon u_n\kon u_{n+1}\prec f' \UND u_{n+1}\erf b_{n+1})\UND\notag\\
		&&					 &\neg\EX[(u_{n+1},b_{n+2})]((u_{n+1}\erf b_{n+1}\UND b_{n+2}\TT{ $\tau$-gespielt }) \UND\notag\\
		&&					 &\EX[u_{n+2}](u_0\kon\ldots\kon u_n\kon u_{n+1}\kon u_{n+2}\prec f' \UND u_{n+2}\erf b_{n+2})) \}\label{Kp_02}
\end{alignat}
und
\begin{alignat*}{1}
&A\subset\bigcup_{p\TT{ \scriptsize{gut}}}K_p
\end{alignat*}
Da die Bedingungsmenge $\BED{B}$ als abz"ahlbar vorausgesetzt wurde, ist diese Vereinigung abz"ahlbar. 

\medskip Um zu zeigen, da"s $A$ $\BM{B}$-mager ist, gen"ugt es nun zu zeigen, da"s die $K_p$ $\BM{B}$-mager sind, denn dann ist die abz"ahlbare Vereinigung $\bigcup_{p}K_p$ $\BM{B}$-magerer Mengen $\BM{B}$-mager und $A$ als Teilmenge einer $\BM{B}$-mageren Menge $\BM{B}$-mager.

\point{Die $K_p$ sind $\BM{B}$-mager}
Sei $K_p$ wie oben konstruiert mit $p=(u_0,b_1,u_1,b_2,\ldots,u_n,b_{n+1})$. Dann l"a"st sich $K_p$ schreiben als
\[
K_p=\bigcup_{\widetilde{u_{n+1}}\erf b_{n+1}}\underbrace{O_{u_0\kon\ldots\kon u_n\kon \widetilde{u_{n+1}}}\cap \C[{(\bigcup_{\substack{u_{n+1}\erf b_{n+1}\\\EX[b_{n+2}](b_{n+2}\TT{\scriptsize{ $\tau$-gesp.}}\UND u_{n+2}\erf b_{n+2})}} O_{u_0\kon\ldots\kon u_n\kon u_{n+1}\kon u_{n+2}})}]}_{=:L_{\widetilde{u_{n+1}}}}.
\]
Nun gen"ugt es zu zeigen, da"s die $L_{\widetilde{u_{n+1}}}$ $\BM{B}$-nirgends dicht sind: 

\point{Die $L_{\widetilde{u_{n+1}}}$ sind abgeschlossen}
Das Komplement einer Menge $L_{\widetilde{u_{n+1}}}$ ist als Vereinigung offener Mengen offen, also ist $L_{\widetilde{u_{n+1}}}$ abgeschlossen.

\point{Die $L_{\widetilde{s_{n+1}}}$ sind $\BM{B}$-nirgends dicht}
Sei $L:=L_{\widetilde{u_{n+1}}}$ f"ur $\widetilde{u_{n+1}}\in X_*^{<\omega}$. Es ist zu zeigen, da"s
\begin{equation*}
\FA[u\in T_L]\EX[b\in\BED{B}]\FA[v\in X^{<\omega}_*](u\kon v\in T_L\Rightarrow v\not\erf b)\label{Frml_X1}.
\end{equation*}

\point{Fall 1} 
Sei $u\prec \tilde{u}:=u_0\kon\ldots\kon u_n\kon u_{n+1}$ und $u\neq \tilde{u}$ (d.h. $u$ ist ein echtes Anfangsst"uck von $\tilde{u}$)\gs etwa $u=(x_0,\ldots,x_i)\prec\tilde{u}=(x_0,\ldots,x_m)$ mit $i\lneq m$. Dann gibt es nach \eemph{Unterscheidbarkeit} [\siehe Definition \ref{DefBedingungsmenge} $2)$] f"ur $x=x_{i+1}$ ein $b_x$ mit $\FA[v\in X_*^{<\omega}] (v\erf b_x\Rightarrow v(0)\neq x)$. Aus $v(0)\neq x$ folgt aber $u\kon v\not\in T_L$. Somit gibt es zu $u$ eine Bedingung $b:=b_x$ mit $\FA[v\in X^{<\omega}_*](u\kon v\in T_L\Rightarrow v\not\erf b)$.

\point{Fall 2}
Sei $u\succ u_0\kon\ldots\kon u_n\kon u_{n+1}$. Etwa $u= u_0\kon\ldots\kon u_n\kon u_{n+1}'$ mit $u_{n+1}\prec u_{n+1}'$ (also eventuell auch $u_{n+1}= u_{n+1}'$). Dann gilt wegen $u_{n+1}\erf b_{n+1}$ und der \eemph{Monotonie} f"ur die Relation $\erf$ [\siehe Definition \ref{DefBedingungsmenge} $1)$]: $u_{n+1}'\erf b_{n+1}$. Spielt nun Spieler II $b_{n+2}':=\tau(p\kon (u_{n+1}'))$, so mu"s wegen (\ref{Kp_02}) f"ur jedes $v\in X_*^{<\omega}$ f"ur alle $f'\in K_p$ gelten:
\begin{alignat*}{1}
&u_0\kon\ldots\kon u_n\kon u_{n+1}'\kon v\not\prec f' \ODER v\not\erf b_{n+2}'.
\end{alignat*}
Also gilt falls $v\erf b_{n+2}'$, da"s $u\kon v\not\in T_L$. Zu $u$ existiert also eine Bedingung $b:=b_{n+2}'$ mit $\FA[v\in X^{<\omega}_*](u\kon v\in T_L\Rightarrow v\not\erf b)$.

\medskip Somit gilt dann nach Fall 1 und 2, da"s die $L_{\widetilde{u_{n+1}}}$ $\BM{B}$-nirgends dicht sind.

\medskip Damit ist die Hinrichtung von $(ii)$ bewiesen.
\end{beweis}

\subsection*{Verallgemeinerte Spiele mit Zeugen $\BMSe{\BM{B}}(C)$}\label{Kap_VerallgSpiel_2}
Analog zu den spieltheoretischen Charakterisierungen der Baire- und der $\sigma$-Kategorie durch die Theoreme \ref{SatzMagerCharSpiele_2} bzw. \ref{SatzSigmaBeschrCharSpiele_2} soll nun eine spieltheoretische Charakterisierung der verallgemeinerten Kategorie f"ur Teilmengen $A$ des Bairerauems mit $A=\p(C)$ f"ur eine Teilmenge $C\subset X^\omega\times\lambda$ (f"ur eine unendliche Ordinalzahl $\lambda$) gegeben werden. Der Nutzen hiervon ist eine gegen"uber $A$ vereinfachte Gewinnmenge $C$. Dazu wird f"ur eine Bedingungsmenge $\BM{B}:=(\BED{B},\BF{E})$ auf $X$ und die Menge $C$ ein Spiel $\BMSe{\BM{B}}(C)$ definiert, in dem Spieler I genau dann eine Gewinnstrategie hat, wenn $A$ eine nicht-leere $\BM{B}$-perfekte Teilmenge enth"alt und II eine Gewinnstrategie hat genau dann, wenn $A$ $\BM{B}^\lambda$-mager ist.

\begin{definition}\header{verallgemeinertes Spiel mit Zeugen $\BMSe{\BM{B}}(C)$}\\
Sei $\BM{B}:=(\BED{B},\BF{E})$ eine Bedingungsmenge auf $X$, sei $A\subset X^\omega$ und $C\subset X^\omega\times\lambda^\omega$ f"ur eine Ordinalzahl $\lambda$ und es gelte:
\begin{alignat*}{1}
&A=\p(C):=\MNG{f\in{X^\omega}}{\EX[{\xi\in\lambda^\omega}]((f,\xi)\in C)}.
\end{alignat*}

Dann definiert man das \emph{verallgemeinerte Spiel $\BMSe{\BM{B}}(C)$ mit Zeugen} wie folgt:
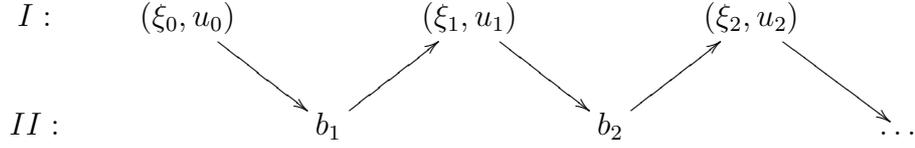
\[
\xymatrix{
I: &(\xi_0,u_0) \ar[dr] & 					 	 &(\xi_1,u_1) \ar[dr]  &						&(\xi_2,u_2) \ar[dr]  &\\
II:&		 				 			&b_1 \ar[ur]	 &						       &b_2 \ar[ur]	&										&\ldots
}
\]
Spieler I und Spieler II spielen abwechselnd\gs Spieler I spielt ein Element $\xi_0\in\lambda$ und eine nicht-leere endliche Sequenz $u_0\in X_*^{<\omega}$, Spieler II spielt eine Bedingung $b_1\in\BED{B}$, Spieler I spielt ein Element $\xi_1\in\lambda$ und eine nicht-leere endliche Sequenz $u_1\in X^{<\omega}_*$, Spieler II spielt eine Bedingung $b_2\in\BED{B}$ u.s.w.. Sei nun $\xi=(\xi_0,\xi_1,\ldots)\in\lambda^\omega$ und $f:=u_0\kon u_1\kon u_2\kon\ldots$. 
Man sagt, da"s \emph{Spieler I gewinnt}, falls gilt:
\begin{ITEMS}[arabic)]
\item $(f,\xi)\in C$,
\item $\FA[i\geq 1](u_i\erf b_i)$.
\end{ITEMS}
Ansonsten sagt man, da"s \emph{Spieler II gewinnt}. $C$ nennt man dabei auch die \emph{Gewinnmenge}.
\end{definition}

Die Begriffe \emph{(Gewinn-)Strategie f"ur I (bzw. f"ur II)} und \emph{determiniert} sind f"ur die Spiele $\BMSe{\BM{B}}(C)$ analog zu den Definitionen \ref{Def_Gewinnstrategie_I} (bzw. \ref{Def_Gewinnstrategie_II}) und \ref{Def_determiniert_G_A}  definiert.

\bigskip Sind $X$, $\lambda$ und \BM{B} abz"ahlbar, l"a"st sich dieses Spiel (wie beim Spiel $\BMSd{\BM{B}}(A)$) mit einer Kodierung der in $\BMSe{\BM{B}}(A)$ gespielten Sequenzen, Zeugen und Bedingungen in die nat"urlichen Zahlen auch als ein Spezialfall der Grundvariante [\siehe Definition \ref{Def_SpielGrundvariante}] auffassen. Analog zu Lemma \ref{Lemma_Determiniertheit_G**p}, Lemma \ref{Lemma_Determiniertheit_G_tilde} und Lemma \ref{Lemma_Determiniertheit_G_tilde_p} gilt:

\begin{lemma}\header{Determiniertheit}\label{Lemma_Determiniertheit_G_B_p}\\
Falls $\AD$ gilt und die Mengen $X$, $\lambda$ und $\BED{B}$ abz"ahlbar sind, so sind auch die obigen Spiele $\BMSe{\BM{B}}(C)$ determiniert.
\end{lemma}
\begin{beweis}
Analog zu Lemma \ref{Lemma_Determiniertheit_G**p}.
\end{beweis}

\bigskip Analog zur Verallgemeinerung der Begriffe \eemph{$\sigma$-beschr"ankt} und \eemph{$\sigma$-kompakt} f"ur beliebige unendliche Ordinalzahlen $\lambda$ [\siehe Satz \ref{SatzLambdaBeschr}] l"a"st sich nun auch der Begriff \eemph{$\BM{B}$-mager} verallgemeinern: Sei $A\subset X^\omega$ beliebig und $\BM{B}=(\BED{B},\BF{E})$ eine Bedingungsmenge. Dann hei"se $A$
\emph{$\BM{B}^\lambda$-mager}, falls $A$ Teilmenge einer Vereinigung $\card[\lambda]$-vieler abgeschlossener $\BM{B}$-nirgends dichter Teilmengen von $X^\omega$ ist:
\begin{equation*}
A\subset\bigcup_{n<\lambda}A_n\\\TT{ wobei }\FA[n<\lambda](A_n\subset X^\omega \TT{ $\BM{B}$-nirgends dicht}).
\end{equation*}

Aus der Definition f"ur \eemph{$\BM{B}^\lambda$-mager} folgt sofort, da"s die $\BM{B}^\lambda$-mageren Teilmengen des Baireraumes abgeschlossen sind gegen"uber Teilmengen. Au"serdem ist die Vereinigung $\card[\lambda]$-vieler $\BM{B}$-mager Mengen $\BM{B}^\lambda$-mager\footnote{{F"ur beliebiges $\lambda\in\ORD$ gilt $\card[{\omega\times\lambda}]=\card[\lambda]$} (\siehe etwa \cite[3.5]{Jech:2003}, vgl. Kap. \ref{Kap_SigmaKategorie} Fu"snote \ref{Fussnote_SigmaKat_Def}).\label{Fussnote_VerallgemKat_Spiel_1}}.

\bigskip Damit erh"alt man nun:

\begin{theorem}\header{Charakterisierung durch Spiele mit Zeugen}\label{Theorem_Char_AllgSpiel_MitZeugen}\\
Sei $\BM{B}$ eine abz"ahlbare Bedingungsmenge auf $X$ (das hei"se: $X$ und $\BED{B}$ abz"ahlbar) und sei $A\subset X^\omega$ mit $A=\p(C)$ f"ur eine Menge $C\subset X^\omega\times\lambda^\omega$ und eine unendliche Ordinalzahl $\lambda$. Dann gilt:
\begin{ITEMS}
\item Spieler I hat eine Gewinnstrategie in $\BMSe{\BM{B}}(C)$ $\Rightarrow$ $A$ enth"alt eine nicht-leere $\BM{B}$-perfekte Teilmenge $B$.
\item Spieler II hat eine Gewinnstrategie in $\BMSe{\BM{B}}(C)$ $\Rightarrow$ $A$ ist $\BM{B}^\lambda$-mager.
\end{ITEMS}
\end{theorem}
\begin{beweis}
$(i)$
Habe I eine Gewinnstrategie $\sigma$ f"ur das Spiel $\BMSe{\BM{B}}(C)$. Wegen $A=\p(C)$ hat I dann auch eine Gewinnstrategie $\tilde{\sigma}$ f"ur das Spiel $\BMSd{\BM{B}}(A)$. Nach Theorem \ref{Theorem_Char_AllgSpiel_OhneZeugen} $(i)\Rightarrow$ enth"alt $A$ somit eine nicht-leere $\BM{B}$-perfekte Teilmenge $B$.

\bigskip$(ii)$
Habe II eine Gewinnstrategie $\tau$. 

\point{Gute Sequenz} 
Eine Sequenz $p:=(l_0,u_0,b_1,l_1,u_1,b_2,\ldots,l_n,u_n,b_{n+1})$ mit $u_i\in X^{<\omega}_*$ f"ur $i=0,\ldots, n$, sowie und $l_j\in\lambda$ f"ur $j=0,\ldots,n$ und $b_i\in\BED{B}$ f"ur $i=1,\ldots,n+1$ hei"se \emph{gut}, falls gilt:
\begin{alignat*}{1}
&u_i\erf b_i \TT{ f"ur }i=1,\ldots,n\UND\\
&b_j\TT{ mit $\tau$ gespielt f"ur }j=1,\ldots,n+1.
\end{alignat*}
Die leere Sequenz sei definitionsgem"a"s gut. 

\point{Sequenz gut f"ur $f\in A$}
$p$ hei"se \emph{gut f"ur $f\in A$}, falls gilt:
\begin{alignat*}{1}
&\EX[u_{n+1}](u_0\kon\ldots\kon u_n\kon u_{n+1}\prec f \UND u_{n+1}\erf b_{n+1}).
\end{alignat*}
Insbesondere ist die leere Sequenz gut f"ur jedes $f\in A$.

\medskip Sei nun $(f,g)\in C$ beliebig. Dann mu"s es eine Sequenz $p:=(\underbrace{l_0,u_0}_I,\underbrace{b_1}_{II},\underbrace{l_1,u_1}_I,\underbrace{b_2}_{II},\ldots,\underbrace{l_n,u_n}_I,\underbrace{b_{n+1}}_{II})$ (eventuell $p=()$) geben mit:
\begin{multline*}
p\TT{ gut }\UND\TT{ gut f"ur $f$}\UND\\ (l_0,\ldots,l_n)=(g(0),\ldots,g(n))\prec g \UND\\
\neg\EX[(l_{n+1},u_{n+1},b_{n+2})](p\kon (l_{n+1},u_{n+1},b_{n+2}) \TT{ gut}\UND \TT{ gut f"ur $f$}\UND\\
(l_0,\ldots,l_n, l_{n+1})=(g(0),\ldots,g(n),g(n+1))\prec g).
\end{multline*}
Ansonsten k"onnte Spieler I gewinnen, obwohl Spieler II mit seiner Gewinnstrategie $\tau$ spielt (Widerspruch). 

\medskip Sei nun $r:=l_{n+1}:=g(n+1)$, dann gilt $f\in M_{(p,r)}$ mit 
\begin{alignat}{2}
&M_{(p,r)}&:=\{f'\in \omega^\omega\mid &p\TT{ gut f"ur }f'\UND\notag\\
		&&&\neg\EX[(u_{n+1},b_{n+2})](p\kon (r,u_{n+1},b_{n+2}) \TT{ gut }\UND\TT{ gut f"ur }f'\}\notag\\
&				&=\{f'\in \omega^\omega\mid &\EX[u_{n+1}](u_0\kon\ldots\kon u_n\kon u_{n+1}\prec f'\UND u_{n+1}\erf b_{n+1})\UND\label{Mpr_01}\\
&				 &&\neg\EX[(u_{n+1},b_{n+2})]((u_{n+1}\erf b_{n+1}\UND b_{n+2}\TT{ $\tau$-gesp.})\UND\notag\\
		&& &\EX[u_{n+2}](u_0\kon\ldots\kon u_n\kon u_{n+1}\kon u_{n+2}\prec f'\UND u_{n+2}\erf b_{n+2}))\}\label{Mpr_02}
\end{alignat}
wobei \glqq$b_{n+2}\TT{ $\tau$-gesp.}$\grqq\ hei"st: $b_{n+2}=\tau(p\kon (r, u_{n+1}))$.
Damit gilt dann
\begin{alignat*}{1}
&A\subset\bigcup_{(p,r)}M_{(p,r)}
\end{alignat*}
wobei "uber die Paare $(p,r)$ vereinigt wird mit $p$ wie oben konstruiert und $r\in\lambda$\gs da $(p,r)\in\omega^{<\omega}\times\lambda$ und $\omega^{<\omega}$ abz"ahlbar ist, werden $\card[\lambda]$-viele Mengen vereinigt\footnote{{F"ur beliebiges $\lambda\in\ORD$ gilt $\card[{\omega\times\lambda}]=\card[\lambda]$} (\siehe etwa \cite[3.5]{Jech:2003}, vgl. Fu"snote \ref{Fussnote_VerallgemKat_Spiel_1}).\label{Fussnote_VerallgemKat_Spiel}}. 

\medskip Da wir die Mengen $M_{(p,r)}$ ganz "ahnlich wie die Mengen $K_p$ im Beweis der Hinrichtung von $(ii)$ von Theorem \ref{Theorem_Char_AllgSpiel_OhneZeugen} definiert haben (lediglich die Bedeutung von \glqq$b_{n+2}\TT{ $\tau$-gesp.}$\grqq\ ist eine andere), k"onnen wir von nun an analog dazu fortfahren. 

\medskip Um zu zeigen, da"s $A$ $\BM{B}^\lambda$-mager ist, gen"ugt es zu zeigen, da"s die $M_{(p,r)}$ $\BM{B}$-mager sind, denn eine Vereinigung  $\card[\lambda]$-vieler $\BM{B}$-magerer Mengen sowie Teilmengen $\BM{B}^\lambda$-magerer Mengen sind $\BM{B}^\lambda$-mager.

\point{Die $M_{(p,r)}$ sind $\BM{B}$-mager}
Sei $M_{(p,r)}$ wie oben konstruiert mit $p:=(l_0,u_0,b_1,l_1,u_1,b_2,\ldots,l_n,u_n,b_{n+1})$. Dann l"a"st sich $M_{(p,r)}$ schreiben als
\[
M_{(p,r)}=\bigcup_{\widetilde{u_{n+1}}\erf b_{n+1}}\underbrace{O_{u_0\kon\ldots\kon u_n\kon \widetilde{u_{n+1}}}\cap \C[{(\bigcup_{\substack{u_{n+1}\erf b_{n+1}\\\EX[b_{n+2}](b_{n+2}\TT{\scriptsize{ $\tau$-gesp.}}\UND u_{n+2}\erf b_{n+2})}} O_{u_0\kon\ldots\kon u_n\kon u_{n+1}\kon u_{n+2}})}]}_{=:L_{\widetilde{u_{n+1}}}}.
\]
Nun gen"ugt es zu zeigen, da"s die $L_{\widetilde{u_{n+1}}}$ $\BM{B}$-nirgends dicht sind:

\point{Die $L_{\widetilde{u_{n+1}}}$ sind abgeschlossen}
Das Komplement einer Menge $L_{\widetilde{u_{n+1}}}$ ist als Vereinigung offener Mengen offen, also ist $L_{\widetilde{u_{n+1}}}$ abgeschlossen.

\point{Die $L_{\widetilde{u_{n+1}}}$ sind $\BM{B}$-nirgends dicht}

Sei $L:=L_{\widetilde{u_{n+1}}}$ f"ur $\widetilde{u_{n+1}}\in X_*^{<\omega}$. F"ur die beiden F"alle 

\point{Fall 1} 
$u\prec \tilde{u}:=u_0\kon\ldots\kon u_n\kon u_{n+1}$ und $u\neq \tilde{u}$ und

\point{Fall 2}
$u\succ u_0\kon\ldots\kon u_n\kon u_{n+1}$

zeigt man nun analog zum Beweis der Hinrichtung von $(ii)$ von Theorem \ref{Theorem_Char_AllgSpiel_OhneZeugen}, da"s jeweils 
\[
\EX[b\in\BED{B}]\FA[v\in X^{<\omega}_*](u\kon v\in T_L\Rightarrow v\not\erf b)
\]
gilt (der einzige Unterschied zum Beweis der Hinrichtung von $(ii)$ von Theorem \ref{Theorem_Char_AllgSpiel_OhneZeugen} besteht in der Definition von $b_{n+2}':=\tau(p\kon(r,b_{n+1}'))$ im zweiten Fall).

\medskip Somit gilt nach Fall 1 und 2, da"s die $L_{\widetilde{u_{n+1}}}$ $\BM{B}$-nirgends dicht sind.

\medskip Damit ist die Hinrichtung von $(ii)$ bewiesen.
\end{beweis}

F"ur $\lambda=\omega$ bedeutet \eemph{$\BM{B}^\lambda$-mager} das selbe wie \eemph{$\BM{B}$-mager}. In diesem Fall lautet Theorem \ref{Theorem_Char_AllgSpiel_MitZeugen}:

\begin{korollar}\header{Charakterisierung durch Spiele mit Zeugen}\label{Theorem_Char_AllgSpiel_MitZeugen_2}\\
Sei $\BM{B}$ eine abz"ahlbare Bedingungsmenge auf $X$ (das hei"se: $X$ und $\BED{B}$ abz"ahlbar) und sei $A\subset X^\omega$ mit $A=\p(C)$ f"ur eine Menge $C\subset X^\omega\times X^\omega$. Dann gilt:
\begin{ITEMS}
\item Spieler I hat eine Gewinnstrategie in $\BMSe{\BM{B}}(C)$ $\Rightarrow$ $A$ enth"alt eine nicht-leere $\BM{B}$-perfekte Teilmenge $B$.
\item Spieler II hat eine Gewinnstrategie in $\BMSe{\BM{B}}(C)$ $\Rightarrow$ $A$ ist $\BM{B}$-mager.
\end{ITEMS}
\end{korollar}

\section{Res"umee}\label{Kap_Resume}
Die eingangs unter \ref{Kap_VereinhProblem} gestellte Frage nach einer Vereinheitlichung der Baire- und der $\sigma$-Kategorie ist in Kapitel \ref{Kap_VerallgKategorie} durch ein verallgemeinertes Kategorien-Konzept f"ur den Baireraum beantwortet worden:

\bigskip Dazu wurden in den Kapiteln \ref{Kap_BaireKategorie} und \ref{Kap_SigmaKategorie} die Baire- und die $\sigma$-Kategorie vorgestellt und spieltheoretisch charakterisiert. In Kapitel \ref{Kap_VerallgKategorie} konnte dann auf Basis dieser Charakterisierungen ein verallgemeinertes spieltheoretisches Kategorienkonzept definiert werden, das Baire- und $\sigma$-Kategorie vereinheitlicht und zudem noch weitere Spezialf"alle beinhaltet [\siehe die Beispiele in den Kapiteln \ref{Kapitel_BMager} und \ref{Kap_BPerf}]. Dabei verhalten sich die Begriffe des verallgemeinerten Konzeptes zu denen der Baire- und der $\sigma$-Kategorie wie in folgender Tabelle dargestellt:\\

\begin{tabular}{ccc}
\textbf{Verallg. Kategorie}	
&\textbf{Baire-Kategorie}	
&\textbf{$\sigma$-Kategorie}\\
$\BM{B}$-nirgends-dicht	
&$\underset{\TT{\scriptsize{\ref{BspBMagerMager} $(i)$}}}{\TT{abg. $+$ nirg. dicht}}$ 	&$\underset{\TT{\scriptsize{\ref{Bsp_BNirgDicht_Kompakt} $(i)$}}}{\TT{kompakt}}$\\
$\BM{B}$-mager				
&$\underset{\TT{\scriptsize{\ref{BspBMagerMager} $(ii)$}}}{\TT{mager}}$									&$\underset{\TT{\scriptsize{\ref{Bsp_BNirgDicht_Kompakt} $(ii)$}}}{\TT{$\sigma$-beschr"ankt}}$\\
$\BM{B}$-perfekt$\neq\LM$			
&							
&$\underset{\TT{\scriptsize{\ref{Bsp_BPerfEnthSuperperfTM}}}}{\TT{enth. superperf. Teilm.$\neq\LM$}}$
\end{tabular}

\bigskip wobei die Begriffe der linken Spalte die spieltheoretischen Verallgemeinerungen der entsprechenden Begriffe der mittleren und der rechten Spalte darstellen.

\bigskip In Kapitel \ref{KapSpielCharBMager} wurde dann gezeigt, da"s die Theoreme in den Kapiteln \ref{KapSpielCharMager} und \ref{KapSpielCharSigmaBeschr"ankt}, die die Baire- bzw. die $\sigma$-Kategorie spieltheoretisch charakterisieren, sich analog auch f"ur die verallgemeinerte Kategorie gewinnen lassen. Dar"uber hinaus lassen sich auch andere Resultate wie etwa diejenigen zur Definierbarkeit in Kapitel \ref{KapSpielCharSigmaBeschr"ankt} auf die allgemeine Kategorie "ubertragen \siehe \cite{Kechris:1977}.

\nocite{Alexandroff:1994}
\nocite{Alsmeyer:1998}
\nocite{Baire:1899}
\nocite{Bauer:1974}
\nocite{Bourbaki:1998:1}
\nocite{Bourbaki:1998:2}
\nocite{Cantor:1883}
\nocite{Elstrodt:1996}
\nocite{Frechet:1906}
\nocite{Gale:1953}
\nocite{Hausdorff:1914}
\nocite{Haworth:1977}
\nocite{Jech:2003}
\nocite{Kanamori:2003}
\nocite{Kechris:1977}
\nocite{Kechris:1995}
\nocite{KechrisSolovay:1985}
\nocite{Koehnen:1988}
\nocite{Kuratowski:1958}
\nocite{Kuratowski:1976}
\nocite{Lebesgue:1904}
\nocite{Lebesgue:1972:2}
%\nocite{Manna:1992}
%\nocite{Manna:1995}
%\nocite{Martin:1973}
\nocite{Moschovakis:1973}
\nocite{Moschovakis:1980}
\nocite{Mycielski:1962}
\nocite{Mycielski:1964:1}
%\nocite{Mycielski:1964:2}
%\nocite{Myerson:1991}
\nocite{Oxtoby:1980}
\nocite{Preuss:1970}
%\nocite{Reimann:2005}
\nocite{Steinhaus:1965}
\nocite{Telarsky:1987}
\nocite{vNeumann:1944}
\nocite{Zermelo:1904}
%==============================backmatter===========================
\backmatter

\chapter{Konventionen}\label{Kap_Konventionen}
%\chapter*{Konventionen}\label{Kap_Konventionen}
%\addcontentsline{toc}{chapter}{\numberline{}Konventionen}

Im Allgemeinen werden (falls nicht ausdr"ucklich anders vermerkt) die folgenden Konventionen befolgt:

%BRAUCHT MAN NICHT
\begin{comment}\KOM[FEHLT:\\\eemph{Produktraum\\WICHTIG F"UR\\BEGRIFFE "UBER\\PROJ. MENGEN\\UND SPIELE}]\end{comment}

\Point{Ableitungen}
Bei l"angeren logischen Ableitungen werden der "Ubersicht halber die Ausdr"ucke nicht immer mit \glqq$\Rightarrow$\grqq\ getrennt. Wenn es der "Ubersichtlichkeit dient, werden "ahnliche Ausdr"ucke nicht immer wieder von neuem aufgef"uhrt, sondern nur diejenigen Teile, die sich ver"andern - etwa:
\begin{alignat*}{4}
\EX[n_0\in\omega]\FA[n>n_0]	&(A(n))	&																&						&\EX[n_0\in\omega]\FA[n>n_0]	&(A(n))\\
\ldots											&(B(n))	&\TT{\quad anstelle von\quad }	&\Rightarrow&\EX[n_0\in\omega]\FA[n>n_0]	&(B(n))\\
\ldots											&(C(n)) &																&\Rightarrow&\EX[n_0\in\omega]\FA[n>n_0]	&(C(n)).
\end{alignat*}

\Point{Abschluss} \siehe \eemph{Inneres}.

\Point{Absolutbetrag}
F"ur eine reelle Zahl $x$ sei $|x|\in\RZ$ ihr \eemph{Absolutbetrag}:
\begin{equation*}
|x|:=
\begin{cases}
x &\TT{falls }x\geq 0,\\
-x &\TT{falls }x< 0.
\end{cases}
\end{equation*}

\Point{Aufz"ahlungen}
Bei Aufz"ahlungen wird oft die verk"urzende Schreibweise
\begin{alignat*}{4}
(x_i>y)_{i\geq 1}			&\TT{\quad anstelle von\quad }	&\FA[i\geq 1](x_i>y)
\end{alignat*}
verwendet.

\Point{Differenz} Sei $X$ eine Menge und seien $A, B\subset X$. Dann bezeichnet $A\o B$ die \eemph{Differenz} $A$ ohne $B$, d.h. die Menge aller $x\in A$ mit $x\not\in B$. Das ist gleichbedeutend mit
\[
A-B:=A\cap\C[B].
\]

\Point{Disjunkte Vereinigung} 
Sei $X$ eine Menge, f"ur $A,B\subset X$, $I$ eine beliebige Indexmenge und $A_i\subset X$ f"ur alle $i\in I$. F"ur paarweise disjunkte Mengen $A_1,A_2,\ldots$ schreiben wir auch $\sum_{i\geq 1}A_i$ statt $\bigcup_{i\geq 1}A_i$.

\Point{Eigenschaft} Eine Teilmenge $E\subset X$ einer Menge $X$ nennen wir auch eine \eemph{Eigenschaft} von Elementen von $X$. Entsprechend hei"st eine Teilmenge $\SYS{E}\subset\POW[X]$ eine \eemph{Eigenschaft} von Teilmengen von $X$.

\Point{Elemente}
\noindent $m, n, i, j, k$ bezeichnen nat"urliche Zahlen,

\noindent $s, t, u, v, w$ bezeichnen endliche und $f, g, h$ bezeichnen unendliche Sequenzen von Elementen in $\omega$ oder in einer vorgegebenen Menge $X$.

\Point{Erzeugte $\sigma$-Algebra} F"ur ein System $\SYS{E}$ von Teilmengen einer nicht-leeren Menge $\Omega$ sei $\sigma(\SYS{E})$ die kleinste $\sigma$-Algebra, die $\SYS{E}$ enth"alt, d.h.:
\[
\sigma(\SYS{E}):=\bigcap_{\substack{\SYS{E}\subset\SYS{A}\\\SYS{A}\TT{\scriptsize{ $\sigma$-Alg.}}}} \SYS{A}.
\]
$\sigma(\SYS{E})$ hei"st die von $\SYS{E}$ \eemph{erzeugte $\sigma$-Algebra}.

\Point{Euklidischer Raum}\label{Def_Eukl_Raum} Das Tupel $(\RZ^n,\textup{d}^n)$ mit 
\[
\textup{d}^n(x,y):=\sqrt{(x_1-y_1)^2+\ldots+(x_n-y_n)^2}
\]
f"ur $x=(x_1,\ldots,x_n),y=(y_1,\ldots,y_n)\in\RZ$ ist ein metrischer Raum\gs genannt \eemph{$n$-dimensionaler Euklidischer Raum}.

\Point{$\G_\delta$- und $\F_\sigma$-Mengen} Eine Teilmenge $A$ eines topologischen Raumes bezeichnet man als \eemph{$\G_\delta$-} bzw. \eemph{$\F_\sigma$-Menge}, falls $A$ Durchschnitt abz"ahlbar vieler offener Mengen bzw. Vereinigung abz"ahlbar vieler abgeschlossener Mengen ist:
\begin{alignat*}{2}
&A=\bigcap_{i<\omega} O_i\TT{ mit $O_i$ offen } &\Leftrightarrow : A\TT{ ist } G_\delta,\\
&A=\bigcup_{i<\omega} A_i\TT{ mit $A_i$ abg. } &\Leftrightarrow : A\TT{ ist } F_\sigma.
\end{alignat*}

\Point{Infimum} \siehe \eemph{Supremum}.

\Point{Inneres, Abschluss} F"ur einen Teilmenge $A$ eines topologischen Raumes $(X,\TOP{X})$ bezeichnet $\INT[A]$ das \eemph{Innere} von $A$, d.h. die maximale offene Teilmenge von $A$:
\begin{alignat*}{1}
&\INT[A]:=\bigcup_{\substack{U\subset A\\U\TT{ \tiny{offen}}}}U.
\end{alignat*}
Mit $\ABS[A]$ bezeichnet man den \eemph{Abschlu"s} von $A$, d.h. die minimale abgeschlossene Obermenge von $A$:
\begin{alignat*}{1}
&\ABS[A]:=\bigcap_{\substack{B\supset A\\B\TT{ \tiny{abgeschl.}}}}B.
\end{alignat*}

\begin{comment}
\Point{Intervalle}
F"ur Beispiele brauchen wir folgende Begriffe: Seien $a,b\in\RZ$ und $a\leq b$. Dann seien \eemph{abgeschlossene}, \eemph{offene}, und \eemph{halboffene} Intervalle (wie "ublich) definiert:
\begin{alignat*}{1}
&[a,b]:=\MNG{x\in\RZ}{a\leq x \leq b}\\
&(a,b):=\MNG{x\in\RZ}{a< x < b}\\
&[a,b):=\MNG{x\in\RZ}{a\leq x < b}\\
&(a,b]:=\MNG{x\in\RZ}{a< x \leq b}.
\end{alignat*}
Falls $a\neq b$ gilt, bezeichnet man diese Intervalle auch als \eemph{eigentliche Intervalle}, und falls $a=b$ gilt als \eemph{uneigentliche Intervalle}.
\end{comment}

\Point{Kardinalit"at} Die \eemph{Kardinalit"at} einer Menge $X$ werde mit $\card[X]$ bezeichnet. F"ur zwei Mengen $X$, $Y$ gilt 
\begin{alignat*}{2}
&\card[X]=\card[Y]	&&:\Leftrightarrow \TT{ Es gibt eine Bijektion zwischen $X$ und $Y$},\\
&\card[X]\leq\card[Y]	&&:\Leftrightarrow \TT{ Es gibt eine Injektion von $X$ nach $Y$},\\
&\card[X]\geq\card[Y]	&&:\Leftrightarrow \TT{ Es gibt eine Surjektion von $X$ nach $Y$},\\
&\card[X]=\card[Y]	&&\AQ[$\substack{\TT{Cantor}\\\TT{Bernstein}}$] \card[X]\leq\card[Y]\UND\card[X]\geq\card[Y],\\
&\card[X]\underset{\TT{\tiny{Cantor}}}{\lneq}\card[{\POW[X]}]&&.
\end{alignat*}
(\siehe etwa \cite[3]{Jech:2003})

\Point{Komplement} Sei $X$ eine Menge und $A\subset X$. Dann bezeichnet $\C[A]$ das \eemph{Komplement} von $A$ in $X$, d.h. die Menge aller $x\in X$ mit $x\not\in A$.

\Point{Leere Menge} Die \eemph{leere Menge} werde mit $\LM$ bezeichnet.

\begin{comment}
\Point{Lineare Ordnung}\label{Def_LinOrd} Eine \eemph{lineare Ordnung} ist eine zweistellige Relation $<$ auf einer Menge $P$ mit:
\begin{ITEMS}[arabic)]
\item $p\not< p$,
\item $p< q\UND q< r\Rightarrow p< r$,
\item $p< q \ODER p=q \ODER q<p$.
\end{ITEMS}
f"ur alle $p,q,r\in P$. Die Menge $P$ hei"st dann \eemph{linear geordnet}.
\end{comment}

\Point{Ma"s} \siehe \eemph{Ma"sraum}.

\Point{Ma"sraum} Ein Tripel $(\Omega,\SYS{A},\mu)$ bestehend aus einer nicht-leeren Menge $\Omega$, einer \eemph{$\sigma$-Algebra}  $\SYS{A}$ auf $\Omega$, d.h.:
\begin{ITEMS}[arabic)]
\item $\Omega\in\SYS{A}$,
\item $A\in\SYS{A}\Rightarrow\C[A]\in\SYS{A}$,
\item $A_1,A_2,\ldots\in\SYS{A}\Rightarrow\bigcup_{i\geq 1}A_i\in\SYS{A}$.
\end{ITEMS}
und einem \eemph{Ma"s} $\mu:\PFEIL{\SYS{A}}{}{[0,\infty]}$ auf $(\Omega,\SYS{A})$, d.h.:
\begin{ITEMS}[arabic)]
\item $\mu(\LM)=0$,
\item $A_1, A_2,\ldots\in\SYS{A}$ paarweise disjunkt $\Rightarrow \mu(\sum_{i\geq 1}A_i)=\sum_{i\geq 1}\mu(A_i)$. 
\end{ITEMS}
hei"st ein \eemph{Ma"sraum}. Dabei hei"sen $A_1, A_2,\ldots\in\SYS{A}$ \eemph{paarweise disjunkt}, falls $A_i\cap A_j=\LM$ f"ur alle $i\neq j$ gilt. (\siehe etwa \cite[1--2]{Alsmeyer:1998})

\Point{Mengen und R"aume}\label{DefMengenR"aumeEtc}
\noindent $M,N,X,Y$ bezeichnen Mengen,

\noindent $X^\omega$ bezeichnet f"ur eine Menge $X$ die Menge aller Folgen in $X$,

\noindent $X^{<\omega}$ die Menge aller endlichen Sequenzen in $X$ und

\noindent $X_*^{<\omega}$ die Menge aller nicht-leeren endlichen Sequenzen in $X$,

\noindent $\Omega$ steht in ma"stheoretischen Beispielen f"ur eine nicht-leere Menge mit der zus"atzlichen Struktur eines Ringes oder einer $\sigma$-Algebra,

\noindent $\omega$ bezeichnet die Menge der nat"urlichen Zahlen\gs dabei wird die Schreibweise
\begin{alignat*}{2}
n<\omega	&\TT{\quad analog zu\quad}	&n\in\NZ
\end{alignat*} 
verwendet,

$\QZ$ und $\RZ$ bezeichnen die Mengen der rationalen bzw. der reellen Zahlen.

\Point{(Mengen-)System, Familie} Eine Teilmenge $\SYS{Y}\subset \POW[X]$ hei"st auch ein \eemph{System von Teilmengen von $X$} oder (falls klar ist, in welcher Menge $X$ sich die Mengen befinden) einfach ein \eemph{(Mengen-)System}. Ist $\SYS{Y}=\MNG{A_i}{i\in I}$ ein indiziertes Teilmengensystem, so schreibt man stattdessen auch $(A_i)_{i\in I}$ und nennt es eine \eemph{Familie} von Teilmengen von $X$. Ist $I$ abz"ahlbar, nennt man die Familie auch eine \eemph{Folge} von Teilmengen von $X$.

\Point{Potenzmenge} Sei $X$ eine Menge. Dann bezeichnet $\POW[X]$ ihre \eemph{Potenzmenge}, d.h. die Menge aller Teilmengen von $X$.

\Point{Produktmenge} Seien $X$ und $Y$ Mengen. Dann bezeichnet $X\times Y$ die \eemph{Produktmenge} oder das \eemph{Produkt} von $X$ und $Y$, d.h. die Menge aller Paare $(x,y)$ mit $x\in X$ und $y\in Y$.

\Point{Ring} Ein System von Teilmengen $\SYS{R}\subset X$ einer nicht-leeren Menge $\Omega$ hei"st \eemph{Ring} auf $\Omega$, falls gilt:
\begin{ITEMS}[arabic)]
\item $\LM\in\SYS{R}$,
\item $A,B\in\SYS{R}\Rightarrow A\o B\in\SYS{R}$,
\item $A,B\in\SYS{R}\Rightarrow A\cup B\in\SYS{R}$.
\end{ITEMS}
(\siehe etwa \cite[3.1]{Alsmeyer:1998})

\Point{$\sigma$-Algebra} \siehe \eemph{Ma"sraum}.

\Point{Supremum, Infimum}
Sei $A\subset\RZ$. Eine reelle Zahl $z$ hei"st \eemph{Supremum} von $A$, falls $z$ kleinste obere Schranke von $A$ ist:
\begin{alignat*}{1}
&z\TT{ ist obere Schranke von $A$}\TT{ (d.h.: }\FA[x\in A](x\leq z)\TT{) und}\\
&\TT{falls $z'$ obere Schranke von $A$ ist, so ist $z<z'$}.
\end{alignat*} 
Analog ist das \eemph{Infimum} von $A$ definiert als die gr"o"ste untere Schranke von $A$.

\Point{Symmetrische Differenz} Sei $X$ eine Menge und seien $A, B\subset X$. Dann bezeichnet $A\sd B$ die \eemph{symmetrische Differenz} von $A$ und $B$, d.h. die Menge aller $x\in A\cup B$ mit $x\not\in A\cap B$. Das ist gleichbedeutend mit
\[
A\sd B := (A-B)\cup (B-A).
\]

\Point{Quantoren}
\noindent In Formeln mit Quantoren schreibt man auch:
\begin{alignat*}{2}
\FA[n>n_0](A(n)) 	&\TT{\quad anstelle von\quad }	&\FA[n<\omega](n>n_0\Rightarrow A(n)),\\
\EX[n>n_0](A(n)) 	&\TT{\quad anstelle von\quad }	&\EX[n<\omega](n>n_0\Rightarrow A(n)).
\end{alignat*} 

\begin{comment}
\Point{Vereinigung, Schnitt} 
Sei $X$ eine Menge. Dann seien \eemph{Vereinigungen} und \eemph{Schnitte} auf Teilmengen von $X$ wie "ublich definiert:
\begin{alignat*}{1}
&A\cup B:=\MNG{x\in X}{x\in A \TT{ oder } x\in B}\\
&\bigcup_{i\in I}A_i:=\MNG{x\in X}{\EX[i\in I](x\in A_i)}\\
&A\cap B:=\MNG{x\in X}{x\in A \TT{ und } x\in B}\\
&\bigcap_{i\in I}A_i:=\MNG{x\in X}{\FA[i\in I](x\in A_i)}
\end{alignat*}
f"ur $A,B\subset X$, $I$ eine beliebige Indexmenge und $A_i\subset X$ f"ur alle $i\in I$. F"ur paarweise disjunkte Mengen $A_1,A_2,\ldots$ schreiben wir auch $\sum_{i\geq 1}A_i$ statt $\bigcup_{i\geq 1}A_i$. 
\end{comment}

\newpage
\point{Abk"urzungen} Gelegentlich werden folgende Abk"urzungen verwendet:
\begin{alignat*}{2}
&\TT{abg.} &=	&\TT{\glqq abgeschlossen(e)\grqq}\\
&\TT{abzb.} &=	&\TT{\glqq abz"ahlbar(e)\grqq}\\
&\TT{Bsp.} &=	&\TT{\glqq Beispiel\grqq}\\
&\TT{bzw.} &=	&\TT{\glqq beziehungsweise\grqq}\\
&\TT{Def.} &=	&\TT{\glqq Definition\grqq}\\
&\TT{d.h.} &=	&\TT{\glqq das hei"st\grqq}\\
&\TT{endl.} &=	&\TT{\glqq endlich\grqq}\\
&\TT{gesp.} &=	&\TT{\glqq gespielt\grqq}\\
&\TT{hinr. Krit.} &=	&\TT{\glqq hinreichendes Kriterium\grqq}\\
&\TT{i.a.} &=	&\TT{\glqq im allgemeinen\grqq}\\
&\TT{insb.} &=	&\TT{\glqq insbesondere\grqq}\\
&\TT{nirg.} &=	&\TT{\glqq nirgends\grqq}\\
&\TT{notw. Krit.} &=	&\TT{\glqq notwendiges Kriterium\grqq}\\
&\TT{off.} &=	&\TT{\glqq offen(e)\grqq}\\
&\TT{S.} &=	&\TT{\glqq Seite(n)\grqq}\\
&\TT{Teilm.} &=	&\TT{\glqq Teilmenge\grqq}\\
&\TT{Umgeb.} &=	&\TT{\glqq Umgebung\grqq}\\
&\TT{vgl.} &=	&\TT{\glqq vergleiche\grqq}\\
&\TT{z.B.} &=	&\TT{\glqq zum Beispiel\grqq}\\
\end{alignat*}

\newpage{
\addcontentsline{toc}{chapter}{Literaturverzeichnis}
%\addcontentsline{toc}{chapter}{\numberline{}Literaturverzeichnis}
\bibliographystyle{alpha}
\bibliography{bibtex-math-logik}
%\bibliography{bibtex-math-logik}
}

\end{document}